\documentclass[a4paper,fleqn,leqno]{article}
\usepackage[utf8]{inputenc}
\usepackage[english]{babel}
\usepackage{latexsym,amssymb,amsthm,xcolor}
\usepackage[centertags]{amsmath}
\usepackage{mathtools}
\usepackage{graphicx}
\usepackage{bbm}
\usepackage{framed}
\usepackage{tikz}
\usepackage{pdflscape}
\usepackage{hyperref}

\DeclareSymbolFont{yhmathlsym}{OMX}{yhex}{m}{n}
\DeclareMathAccent{\wideparen}{\mathord}{yhmathlsym}{"F3}
\usepackage[normalem]{ulem}

\newcounter{secnum}
\numberwithin{equation}{section}

\newtheorem{xx}{\bf xxx}

\newtheorem{yy}{\bf yyy}

\newcommand{\Rand}[1]{\marginpar{#1}}
\renewcommand{\Rand}[1]{}

\newcommand{\be}[1]{\Rand{\vspace{0,6cm}\tt #1}\begin{equation}\label{#1}}
\newcommand{\beL}[1]{\Rand{\vspace{0,6cm}\tt #1}\begin{lemma}\label{#1}}
\newcommand{\belC}[2]{\Rand{\vspace{0,6cm}\tt #1}\begin{lemma}[#2]\label{#1}}
\newcommand{\beP}[1]{\Rand{\vspace{0,6cm}\tt #1}\begin{proposition}\label{#1}}
\newcommand{\bePC}[2]{\Rand{\vspace{0,6cm}\tt #1}\begin{proposition}[#2]\label{#1}}
\newcommand{\beD}[1]{\Rand{\vspace{0,6cm}\tt #1}\begin{definition}\label{#1}}
\newcommand{\beT}[1]{\Rand{\vspace{0,6cm}\tt #1}\begin{theorem}\label{#1}}
\newcommand{\beC}[1]{\Rand{\vspace{0,6cm}\tt #1}\begin{corollary}\label{#1}}
\newcommand{\ep}{\end{proposition}}
\newcommand{\bea}[1]{\Rand{\vspace{0,7cm}\tt #1\vspace{-0,7cm}}\begin{eqnarray}\label{#1}}

\newcommand{\wt}{\widetilde}

\def\CB{\mathcal{B}}
\def\CR{\mathcal{R}}
\def\CD{\mathcal{D}}
\def\CF{\mathcal{F}}
\def\CI{\mathcal{I}}
\def\CM{\mathcal{M}}

\def\CL{\mathcal{L}}

\def\CP{\mathcal{P}}

\def\CU{\mathcal{U}}

\def\CW{\mathcal{W}}

\def\mfA{\mathfrak{A}}
\def\mfC{\mathfrak{C}}
\def\mfh{\mathfrak{h}}

\def\mfL{\mathfrak{L}}

\def\mfN{\mathfrak{N}}

\def\mfY{\mathfrak{Y}}
\def\mfZ{\mathfrak{Z}}

\def\mfa{\mathfrak{a}}
\def\mfc{\mathfrak{c}}
\def\mfn{\mathfrak{n}}

\def\K{\mathbb{K}}

\def\N{\mathbb{N}}
\def\R{\mathbb{R}}
\def\U{\mathbb{U}}
\def\Z{\mathbb{Z}}
\def\E{\mathbb{E}}
\def\PP{\mathbb{P}}

\newcommand{\Ltoo}{{_{\D \longrightarrow \atop L \to \infty}}}
\newcommand{\Lto}{{_{\D \Longrightarrow \atop L \to \infty}}}
\newcommand{\kto}{{_{\D \Longrightarrow \atop k \to \infty}}}

\newcommand{\nto}{{_{\D \Longrightarrow \atop n \to \infty}}}
\newcommand{\ntoo}{{_{\D \longrightarrow \atop n \to \infty}}}

\newcommand{\D}{\displaystyle}
\newcommand{\wh}{\widehat}

\DeclareMathAlphabet{\mathpzc}{OT1}{pzc}{m}{it}

\newcommand{\dd}{\mathrm{d}}
\newcommand{\eee}{\mathrm{e}}
\newcommand{\ee}{\end{equation}}
\newcommand{\eea}{\end{eqnarray}}
\newcommand{\bean}{\begin{eqnarray*}}
\newcommand{\eean}{\end{eqnarray*}}

\newcommand{\prodl}{\prod\limits}
\newcommand{\liml}{\lim\limits}
\newcommand{\suml}{\sum\limits}

\newcommand{\cadlag}{c\`adl\`ag}

\makeatletter \@addtoreset{equation}{section} \makeatother
\makeatletter \@addtoreset{enunciato}{section} \makeatother
\newcounter{enunciato}[section]

\newtheorem{ittheorem}{Theorem}
\newtheorem{itlemma}{Lemma}
\newtheorem{itproposition}{Proposition}
\newtheorem{itdefinition}{Definition}
\newtheorem{itremark}{Remark}
\newtheorem{itexample}{Example}
\newtheorem{itclaim}{Claim}
\newtheorem{itcorollary}{Corollary}

\newenvironment{theorem}{\addtocounter{enunciato}{1}
\begin{ittheorem}}{\end{ittheorem}}
\newenvironment{lemma}{\addtocounter{enunciato}{1}
\begin{itlemma}}{\end{itlemma}}
\newenvironment{proposition}{\addtocounter{enunciato}{1}
\begin{itproposition}}{\end{itproposition}}
\newenvironment{definition}{\addtocounter{enunciato}{1}
\begin{itdefinition}}{\end{itdefinition}}
\newenvironment{remark}{\addtocounter{enunciato}{1}
\begin{itremark}}{\end{itremark}}
\newenvironment{example}{\addtocounter{enunciato}{1}
\begin{itexample}}{\end{itexample}}

\newenvironment{corollary}{\addtocounter{enunciato}{1}
\begin{itcorollary}}{\end{itcorollary}}

\usepackage{stackengine}

\begin{document}

\title{Continuum hierarchical Cannings process:\\
I.\ Construction, duality and dichotomy}

\author{A. Greven$^1$, F. den Hollander$^2$}

\date{\today}

\maketitle

\begin{abstract}
We construct spatial measure-valued Cannings processes on the \emph{continuum} hierarchical group $\Omega_N^\infty$ of order $N\in\N$ (an uncountable topological group with an ultrametric distance). These processes arise as the continuum-space limit of spatial Cannings processes on the \emph{discrete} hierarchical group $\Omega_N$ of order $N\in\N$ (a countable topological group with an ultrametric distance) analysed in \cite{GHKK14}. On $\Omega_N$, each site represents the \emph{geographic location} of a colony of individuals, each carrying a \emph{type} drawn from a compact type space $\K$, and the state in each colony is a probability measure on $\K$ describing the type distribution in that colony in the large-colony-size limit of individual-based models. The continuum-space processes on $\Omega_N^\infty$ take values in the set of locally finite measures on $\Omega_N^\infty \times \K$ and are to be thought of as describing a population of individuals, each equipped with a geographic location in $\Omega_N^\infty$ and a type in $\K$, and each being subjected to \emph{migration} and three mechanisms of \emph{resampling}: (1) Fleming-Viot resampling at single sites; (2) Cannings resampling at single sites; (3) Cannings resampling and uniform reshuffling in balls around single sites (corresponding to consecutive hierarchical levels). The reshuffling adds to the migration, resulting in an \emph{effective migration} that is the superposition of the two relocation mechanisms. 

The continuum space process is a strong Markov process arising as the limit of discrete space processes on countable subgroups of $\Omega_N^\infty$ arising as level-$L$ refinements of $\Omega_N$. We show that the continuum-space process is the solution of a \emph{well-posed martingale problem} (after the latter is properly modified to deal with the special features of the uncountable geographic space $\Omega_N^\infty$) and obtain an equivalent characterisation via stochastic equations. We analyse both the finite-time and the infinite-time (= equilibrium) properties of the solution. We construct a dual measure-valued process driven by a spatial coalescent and show that there is a dichotomy between two regimes: (I) \emph{deterministic} (= locally multi-type = globally multi-type); (II) \emph{stochastic} (= locally mono-type = globally mono-type). We characterise both regimes in terms of the migration rates and the resampling rates, thereby obtaining a full \emph{classification} in term of the degree of the effective migration (subject to certain regularity conditions on the rates). This characterisation adds a \emph{new universality class} to the population genetics and ecology literature, where the block-resampling is the game changer. It confirms that the \emph{degree} of the underlying migration is the quantity that characterises the large space-time scaling behaviour. In the stochastic regime we characterise the process as the entrance law of a solution to a stochastic equation for the densities of the state w.r.t.\ the canonical reference measure on $\Omega_N^\infty \times \K$. This requires an extended analysis of the growth of monotype clusters on $\Omega_N$, where we obtain new results that are of interest in themselves. The latter in turn leads to another \emph{new universality class} in the stochastic regime for the evolution on $\Omega_N^\infty$, arising in the limit of diverging Fleming-Viot resampling.

We further consider two simultaneous scaling limits of the discrete-space process and the continuum-space process: (i) scaling of space, time and the migration rates to get stable limiting migrations; (ii) scaling of the resampling rates and the Haar measure. We show that these give limiting processes representing \emph{universality classes}. Namely, on $\Omega_N^\infty$, in the stochastic regime we get a mixture of pure Cannings processes with block-resampling and Fleming-Viot diffusion processes, while in the deterministic regime we get L\'evy-flows. 

Working on the ultrametric hierarchical group allows us to carry out a full analysis, which on Euclidean space has \emph{not} yet been achieved and remains highly challenging.
\end{abstract}

\noindent
\hspace*{0.3cm} {\footnotesize $^{1}$ 
Department Mathematik, Universit\"at Erlangen-N\"urnberg, Cauerstrasse 11, D-91058 Erlangen, Germany\\
greven@mi.uni-erlangen.de}\\
\hspace*{0.3cm} {\footnotesize $^{2}$ 
Mathematisch Instituut, Universiteit Leiden, Einsteinweg 55, 2333CC Leiden, The Netherlands\\
denholla@math.leidenuniv.nl}



\medskip\noindent
\emph{Acknowledgments:}
AG was supported by the Deutsche Forschungsgemeinschaft (through grant DFG-GR 876/16-2 of SPP-1590), FdH was supported by the Netherlands Organisation for Scientific Research (through NWO Gravitation Grant NETWORKS-024.002.003) and by the Alexander von Humboldt Foundation.
\normalsize


\small
\tableofcontents
\normalsize


\section{Introduction and main results}
\label{s.intro}

Section \ref{ss.motivation} offers the motivation of our paper, while Section \ref{ss.background} provides the background. Section \ref{ss.hierarmod} describes the hierarchical models that are the object of our investigation, Section \ref{ss.program} sketches our program, while Section \ref{ss.maintheorems} states our main theorems. Section \ref{ss.outline} provides an outline of the remainder of the paper.

All the results presented in the paper are motivated, stated and explained in Section \ref{s.intro}. The details behind these results are explained in Sections \ref{s.mi}--\ref{s.dual} and Appendices~\ref{appA}--\ref{appB}, and are technically demanding.


\subsection{Motivation}
\label{ss.motivation}


\paragraph*{Setting and construction.}

The present paper is the first in a series of three complementary papers in which we consider spatial Cannings processes on the \emph{continuum} hierarchical group $\Omega_N^\infty$ of order $N\in\N$ (= the closure of $\Omega_N^L$ over $L \in \N_0$ i.e. the direct sum over $\{-L,\ldots,-1\} \cup \N_0$ of the cyclical group $\{0,1,\ldots,N-1\}$; see Section \ref{ss.cont}). Some of these processes arise as the continuum-space \emph{rescaling limit} of spatial Cannings processes on the \emph{discrete} hierarchical group $\Omega_N$ of order $N\in\N$ (= the direct sum over $\N_0$ of the cyclical group $\{0,1,\ldots,N-1\}$), which were introduced and analysed in \cite{GHKK14}. Indeed the study of large space-time scale properties, keyword cluster formation, of these models results also in limit objects, which are a very special subclass of the continuum space processes on $\Omega_N^\infty$ (see~\cite{EF96} for this example).

The continuum-space processes take values in the set of locally finite measures on $\Omega_N^\infty \times \K$, with $\K$ a compact Polish space. Each site of $\Omega_N^\infty$, a locally compact ultrametric space, is to be thought of as a colony of individuals, each carrying a \emph{type} from $\K$, and the state of each colony describes the type distribution in that colony in the large-colony-size limit of spatial individual-based models. Each individual is being subjected to \emph{migration} and to three mechanisms of \emph{resampling}: (1) Fleming-Viot resampling at single sites; (2) Cannings resampling at single sites; (3) Cannings resampling and uniform reshuffling in balls around single sites corresponding to consecutive hierarchical levels. The \emph{initial state} is taken to be stationary and ergodic under the group action of $\Omega_N^\infty$.

We construct the continuum-space process via approximation by the processes on the countable $\Omega_N^L$ arising by $L$-truncation and are analysed by using \emph{duality} as key tool \cite{DGP24}. In interesting parameter regimes we get the process as the unique solution to a \emph{stochastic equation} with an \emph{SPDE-term}, and show that it arises as the limit of suitable discrete-space processes representing successive refinements. We provide its characterisation as the solution to a \emph{well-posed martingale problem} (which in the \emph{continuum-space} setting is not of the standard form; see \cite{EF96,DF97a,DP98} for related examples and more recently \cite{EVY19,EKLRT24}).

The continuum-space process can be either \emph{deterministic} or \emph{stochastic}, i.e., the solution is a deterministic process, respectively, a stochastic process with changes over time in the type distribution occurring in every ball around every site. Deterministic means that in the martingale problem we have a first-order operator as the generator of the underlying motion, and the solution is mass flow driven by the transition kernel of this migration. Stochastic means that in the martingale problem we have a second-order operator as the generator of a diffusion (= `diffusion-case') and/or a non-degenerate integral operator as the generator of a jump process with countably many small jumps (= `jump-case'), leading to a process with paths of infinite variation, respectively, finite variation, or a combination of both. Which of the two scenarios occurs depends on the rates of migration and resampling. We identify \emph{necessary and sufficient conditions} in terms of the transition rates, thereby obtaining a full \emph{classification} (subject to certain regularity conditions on the rates). 


\paragraph*{Structure and key properties.}

It turns out that the block-resampling is a game changer, leading to a further subdivision of the stochastic regime. Namely, an interesting feature is that, because of the block-resampling, the process can be stochastic not only for \emph{point-recurrent} effective migration on $\Omega_N^\infty$, but also for Harris-recurrent migration on $\Omega_N^\infty$, provided the block-resampling rates are large enough. This fact is important for population genetics and ecology, since Harris-recurrent migration on $\Omega_N^\infty$ is comparable to Harris-recurrent migration on $\R^2$ in the sense of potential theory.

With the help of the duality relation established for the spatial Cannings process on $\Omega_N$ respectively $\Omega_N^L$, we derive a duality relation for the spatial Cannings process on $\Omega_N^\infty$. The dual system is a measure-valued process on $\Omega^\infty_N \times \K$ with atomic states driven by a population of migrating and block-coalescing L\'evy-processes carrying types and having weights. We determine when these L\'evy-processes are not coalescing and when they are coalescing, which in the original system corresponds to \emph{locally multi-type} states, respectively, \emph{locally mono-type} states. In the first case the state of the population is a deterministic measure for $t>0$, in the second case it is an atomic \emph{random measure} that induces a countable partition of $\Omega_N^\infty$ into \emph{random clusters} of one type. 

We identify the long-term behaviour of the spatial Cannings process on $\Omega_N^\infty$ by showing that it converges to an equilibrium measure that depends on the initial type distribution and that can be classified, namely, the equilibrium measure can be either \emph{globally mono-type} (`clustering') or \emph{globally multi-type} (`co-existence'). Interestingly, all three dichotomies for the dynamics coincide, i.e., `deterministic = locally multi-type = globally multi-type' and `stochastic = locally mono-type = globally mono-type'. 

Finally, we show the existence of a special class of \emph{stable} migrations on the group $\Omega_N^\infty$ analogous to stable motions on $\R^d$, which represent the \emph{universality classes} for the large space-time behaviour of the spatial Cannings process on $\Omega_N$ \emph{and} on $\Omega_N^\infty$. The classes arise as the \emph{space-time scaling limit} of the process on $\Omega_N$, respectively, $\Omega_N^\infty$. In particular, we find that only the Cannings block-resampling gives rise to a new universality class, while the Cannings single-site resampling does not and falls in the same universality class as the Fleming-Viot resampling.

There is yet a completely different universality class of stochastic processes on $\Omega_N^\infty$, arising in the limit as the resampling rate \emph{tends to infinity} so we have not anymore a jump-diffusion. This comes from a study of the cluster formation on $\Omega_N$ for recurrent migration and is in duality with coalescing L\'evy-processes with \emph{instantaneous} coalescence. For the stable case this was already observed in~\cite{EF96,Klenke1996}, and arises from studying the cluster growth on $\Omega_N$ for strongly recurrent migration, respectively, on $\Omega_N^\infty$ for point-recurrent migration.


\subsection{Background}
\label{ss.background}


\paragraph*{Targets.}

We are interested in populations consisting of individuals carrying both a type and a geographic location (like in $\Z^d$ or the hierarchical group $\Omega_N$), evolving randomly by \emph{migration} and by \emph{change of generation} driven by Fleming-Viot and Cannings \emph{resampling}. We view the population as a family of probability measure-valued processes, living on a common type space $\K$ and being indexed by the geographic space $\Omega_N$ that labels the colonies in which the individuals are organised (referred to as configuration space). To pass to the uncountable geographic space referred to as continuum space we have to view the population as a \emph{measure on $\Omega_N \times \K$}, leading to measure-valued processes.

We focus on the hierarchical group as the geographic space for two reasons: (I) On this space it is possible to carry out explicit calculations and to analyse the large-scale behaviour of the process, which exhibits the same features as found or conjectured on $\Z^d$, respectively, $\R^d$. (II) In population genetics and in ecology, it is often \emph{more appropriate} to model the habitat as equipped with an \emph{ultrametric} classification of neighbourhoods rather than with balls in a \emph{Euclidean metric}. 

The following four problems are of general interest:
\begin{itemize}
\item[(1)]
It is challenging to construct processes describing populations evolving in \emph{continuum} geographic spaces, because standard techniques for discrete countable geographic spaces break down, compared to the standard form of martingale problems in configuration space.
\item[(2)]
\begin{itemize}
\item[(i)] 
A good understanding of continuum-space processes is lacking when the underlying migration is a \emph{jump process} rather than a Brownian motion, respectively, a diffusion process, which gives a richer panorama of spatial patterns for the types. The reason is that a Brownian motion on an ordered space as $\R$ allows for a considerable simplification of the description of the measure-valued dual process, but reduces the complexity. 
\item[(ii)] 
We must pass to a continuum space of types instead of two types, which is a significant step when we work with stochastic equations.
\item[(iii)] 
The same is true for the more volatile changes of generation that arise under Cannings resampling instead of Fleming-Viot resampling, which lead to a jump process also in type-space (with countably many jumps per finite-time interval) instead of a diffusion.
\end{itemize}
In the hierarchical setting it is possible to elucidate (i)--(iii).
\item[(3)]
It is not clear which space-time features of Cannings processes are \emph{universal} and which are not. In particular, it is necessary to distinguish between \emph{Cannings block-resampling} and \emph{Cannings single-site resampling}. This situation can only be clarified by taking scaling limits (scaling of space, time and the migration/resampling rates, i.e., continuum-space scaling limits of discrete-space processes).
\item[(4)]
The evolution of continuum-space processes needs to be described both on the level of \emph{type frequencies} (in the present paper) and on the level of the \emph{genealogies} in the population (tackled in \cite{gh2}).
\end{itemize}
In a series of three papers we contribute towards a better understanding of these problems by analysing \emph{spatial Cannings processes on the continuum hierarchical group} $\Omega_N^\infty$. 


\paragraph*{Earlier work.}

Continuum-space models are viewed as being more realistic and are expected to be easier to handle mathematically. The idea to exhibit the \emph{universality classes} for the large space-time behaviour of systems on geographic spaces such as $\Z^d$ via \emph{continuum-space limits}, arising from space-time, volatility-rate and mass scaling, was successfully applied to so-called \emph{continuum-mass branching} models. In particular, systems of interacting Feller diffusions indexed by $\Z^d$ in the continuum limit lead to Dawson-Watanabe super-processes on $\R^d$ (see \cite{DP91,D93}). The same is true for other universality classes given by systems of \emph{catalytic, self-catalytic} or \emph{mutually catalytic} branching diffusions (on $\Z^d$ in \cite{DF97a,DF97b,DP98}, on $\Omega_N$ in \cite{DGZpr}, with coloured noise in \cite{BHO16}). The structure in the mutually catalytic case and the case of stable-migration Fisher-Wright diffusions is the closest to the structure in our model, even though we have finite or countably many types (selected from a continuum) rather than two types. For these branching processes also \emph{new concepts} around the martingale problem characterisation have been introduced, which turn out to be crucial also in the class of models we discuss here. \emph{Stable-migration Fisher-Wright diffusions} have been studied on the continuum hierarchical group $\Omega_N^\infty$ in \cite{EF96}, while cluster formation has been studied on the discrete hierarchical group $\Omega_N$ in \cite{Klenke1996}. The latter gives rise to yet another model on $\Omega_N^\infty$, as we will see later. 

\begin{itemize}
\item 
In the case of Brownian migration and finite-variance branching, the continuum-space limit gives rise to an evolution that is \emph{stochastic in low dimensions} and \emph{deterministic in high dimensions}. The cutoff is model-dependent. For the above four models the deterministic regime occurs in the following cases: never (Dawson-Watanabe process), for $d>4$ (catalytic branching), for $d >2$ (mutually catalytic branching, self-catalytic branching). The same holds for spatial Fleming-Viot processes, where the cutoff is $d>1$. Similar properties for these models occur on the continuum hierarchical group $\Omega_N^\infty$ \cite{DGZpr}, where the relevant parameter is the potential-theoretic concept of the \emph{degree} of the underlying migration, instead of the \emph{dimension} of the space. Also on $\R^d$ the results are in fact based on the \emph{degree} of the migration (and translated into the dimension of the geographic space), and hence the degree is a key unifying concept. In particular, the models on $\Z^d$ have been extended to infinite-variance branching and to migration via L\'evy jump processes, providing the full spectrum of possible degrees, but this setting has so far resisted a complete understanding. Genealogical structures were studied for the Dawson-Watanabe process in \cite{DP91} (in terms of \emph{historical processes} on $\R^d$, leading to stochastic continuum-space processes).
\item
In the case of Fleming-Viot models, the focus has been on $\R$, where the stochastic dynamics is described by SPDEs \cite{KS88,MT95}, and on incorporating \emph{selection} and \emph{mutation} in \cite{EVY19,EKLRT24}. A subclass of stable-migration Fisher-Wright models was studied on $\Omega_N^\infty$ in \cite{EF96}. Also \emph{historical} and \emph{genealogical} aspects were studied. For example, on $\R$, L\'evy-webs were considered in \cite{EMS13} and \cite{MRV18}, and \emph{entrance laws} were constructed for systems of coalescing particles that start from every space-time point. For the spatial Fleming-Viot process, genealogies have been studied on $\R$ in \cite{GSW16} (by using the description of genealogies via equivalence classes of ultrametric measure spaces, leading to stochastic \emph{continuum-space genealogy} processes).
\end{itemize}


\paragraph*{Present work.}

The use of the hierarchical group as a geographic space in \emph{population genetics} and in \emph{ecology} goes back to the 1940's, in attempts to give real populations a more adapted structure than that provided by Euclidean space. For more recent examples, we refer to \cite{FG94,EF96,GdHO22,GdHO25}. It was also used successfully in \emph{chemical physics} starting in the 1980's (compare with \cite{SF83}, \cite{ABO03}). The ultrametric distance on hierarchical groups allows us to resolve questions that seem intractable on Euclidean spaces investigated in \cite{BEV10,EVY19}.

In the present paper we focus on continuum-space processes based on \emph{jump processes} as underlying migration. The main theme is the construction of the process on the level of type distribution via a \emph{well-posed martingale problem}, on the \emph{duality} with a spatial coalescent process, and on establishing the \emph{dichotomy} of a stochastic evolution versus a deterministic evolution. In our setting, the latter coincides with the dichotomy of locally mono-type states versus locally multi-type states. This coincidence is different from what happens for branching super-processes, which are stochastic evolutions but do allow for both locally mono-type and locally multi-type states. 

In the second sequel paper \cite{gh2} we will look at the \emph{genealogical structure} of the population via a graphical construction called the \emph{ancestral web} (see \cite{GKWpra, GKWprb}) and its dual the \emph{L\'evy-web}, which captures the full ancestry of the population at a given time. For Fleming-Viot processes on $\R$ this construction was done in \cite{GSW16} for Brownian motion, based on the concept of the Brownian web. For $\Omega_N^\infty$ and \emph{jump processes} as migration we expect a \emph{richer variety} of \emph{spatial patterns} that, due to the ultrametricity of the geographic space, are still manageable. 

In the third sequel paper \cite{gh3} we will look at the \emph{hierarchical mean-field limit}, respectively, the \emph{hierarchical mean-field continuum-space limit}, where the order $N$ of the hierarchical group $\Omega_N$, respectively, $\Omega_N^\infty$ tends to infinity and the \emph{degree of the migration} tends to $0$ (as for simple random walk on $\Z^2$). 

What our series of papers shows is that in the hierarchical setting we can complete a program that in the Euclidean setting remains open. 
  

\subsection{Hierarchical geographic spaces: $\Omega_N$ and $\Omega_N^\infty$}
\label{ss.hierarmod}


\paragraph*{Discrete-space model.}

\begin{figure}[htbp]
\centering
\includegraphics[width=0.85\linewidth]{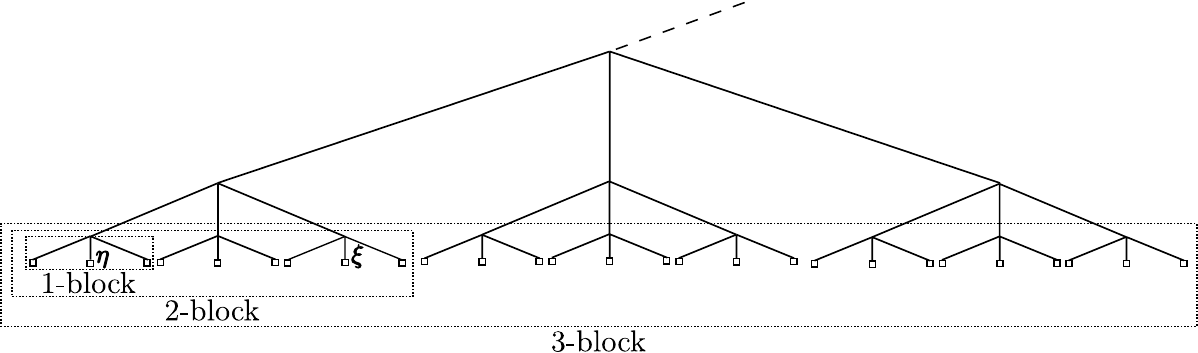}
\caption{\small Close-ups of a 1-block, a 2-block and a 3-block in the discrete hierarchical group of order $N=3$. The elements of the group are the leaves of the tree (indicted by $\Box$'s). The hierarchical distance between two elements is the graph distance to the most recent common ancestor (equal to $2$ for $\xi$ and $\eta$ in the picture) in the auxiliary infinite $N$-ary tree that is drawn on top (which represents the ultrametric distance).}
\label{fig-hierargr}
\end{figure}

In \cite{GHKK14} we introduced the following spatial multi-type population model arising as the large-colony-size limit of an \emph{individual-based} model. In the latter, individuals live in colonies, indexed by the discrete hierarchical group $\Omega_N$ of order $N \in \N$, equipped with an ultrametric that organises the colonies into balls, often referred to as blocks, at hierarchical levels $k\in\N_0$ (see Fig.~\ref{fig-hierargr}). 

Individuals carry a \emph{type}, drawn from a compact set $\K$. Individuals \emph{migrate}, \emph{resample at single sites} according to Fleming-Viot and Cannings resampling, and \emph{resample in blocks} of colonies according to Cannings resampling at all hierarchical levels $k\in\N$ simultaneously at rates that depend on $k$. In the large-colony-size limit this leads to a system of \emph{interacting measure-valued processes} indexed by $\Omega_N$, living on the state space $\CP(\K)^{\Omega_N}$, whose evolution is as follows: 
\begin{itemize}
\item
The migration is characterised by a random walk kernel on $\Omega_N \times \Omega_N$ that is built from a sequence of \emph{multi-level exchangeable migration rates} $(c_k)_{k\in\N_0}$ taking values in $(0,\infty)$. 
\item
The resampling consists of three parts:
\begin{itemize}
\item[(1)]
A \emph{volatility constant} $d \in [0,\infty)$, describing the rate of Fleming-Viot resampling in single colonies. 
\item[(2)]
A \emph{finite resampling measure} $\Lambda_0$ on $[0,1]$, describing Cannings resampling in single colonies. 
\item[(3)] 
A \emph{sequence of finite resampling measures} $(\Lambda_k)_{k \in \N}$ on $[0,1]$, describing Cannings block-resampling with block-reshuffling in $k$-balls around each colony (see Figs.~\ref{fig-canevol}--\ref{fig-reshuffle} below).
\end{itemize}
\end{itemize}
The four evolution mechanisms have different effects on the shape and the structure of the population.

The conditions that we pose on the parameters ensures that we obtain a system of measure-valued components that for all positive times are atomic measures in every component. Consequently, the dual coalescent can be started with countably many individuals at each site, which results for $t>0$ in a locally finite configuration for a fraction $1-\varepsilon$ of the dual individuals for every $\varepsilon>0$.

\begin{remark}{\bf [Cannings models]} 
{\rm 
\begin{itemize}
\item[(a)] 
The non-spatial Cannings model \cite{C74}, \cite{C75} was introduced as a \emph{neutral} population model that allows for a change of generation in which a \emph{positive fraction} of the next generation descends from a \emph{single ancestor} of the current generation, even in the rapid reproduction large-colony-size limit. This accounts for a larger volatility than in the Moran model. 
\item[(b)] 
Spatial Cannings models were introduced in population genetics to analyse the behaviour of multi-type populations under neutral but volatile evolution, with a geographic structure that allows for non-local resampling on top of migration and possibly for selection (on $\Z^d$ and $\R^d$ in \cite{BEM07}, \cite{BEV10} \cite{BEV13a}, \cite{BEV13b}, \cite{EVY19}, on $\Omega_N$ in \cite{GHKK14}, \cite{GHK18}). The new feature added in \cite{GHKK14}, respectively, \cite{BEV10} and fully executed in \cite{EVY19}, is \emph{block-resampling} (= resampling in balls), which is added to the Cannings single-site resampling and allows for additional volatility arising from \emph{local catastrophic events}.
\item[(c)] 
For Cannings models on $\Omega_N$ with appropriate choices of $(c_k)_{k\in\N_0}$, the migration on $\Omega_N$ can be made into a caricature of migration on $\Z^2$ (in the sense of potential theory), which is the most relevant geographic space for population genetics. This will be explained in more detail in \cite{gh3}, the key notion being the \emph{degree} of the migration \cite{DGW04}, \cite{DGW05}.
\end{itemize}
} \hfill$\spadesuit$
\end{remark}


\paragraph*{Continuum-space model.} 

We consider three aspects: structure, evolution, properties. 

\medskip\noindent
{\bf (I) Structure:}
Below we introduce the analogue of the spatial Canning process on the continuum hierarchical group $\Omega^\infty_N$, obtained by extending the discrete hierarchical group $\Omega_N$ downwards to all hierarchical levels $k \in \Z$, i.e., by successively splitting each site into $N$ further sites, etc.\ (see Fig.~\ref{fig-hierargr}) and choosing the topology according to an ultrametric in such a way that the new sites get closer together as the splitting level $k$ decreases to $-\infty$. 

\begin{remark}{\bf [Construction of the continuum hierarchical group]}
{\rm The continuum hierarchical group $\Omega_N^\infty$ can be viewed as a \emph{countably infinite iterative refinement} $\Omega_N^L$ of the discrete hierarchical group $\Omega_N$ in Fig.\ \ref{fig-hierargr}, in which at every level $k \in \Z_- = -\N_0$, every site is split into $N$ sites on level $k-1$ such that the distance between two sites at level $k-1$ originating from a site at level $k$ tends to zero as $k\to\infty$ (see Section~\ref{ss.cont} for further details). Think of $k\to\infty$ as the `upstairs macro-world' and $k \to -\infty$ as the `downstairs micro-world'.} \hfill$\spadesuit$
\end{remark}

\begin{figure}[htbp]
\centering
\includegraphics[width=0.85\linewidth]{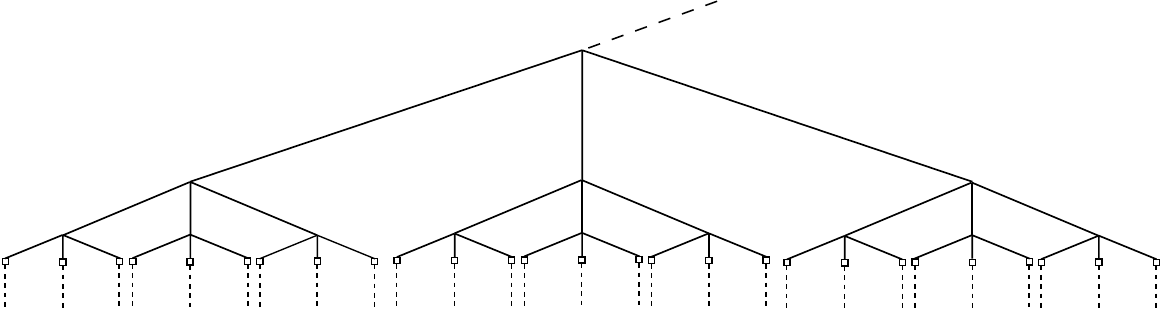}
\caption{\small Visualisation of $\Omega_N^\infty$: Each element of $\Omega_N$ in Fig.\ \ref{fig-hierargr} is refined iteratively, and is replaced by an infinite $N$-ary tree running downwards, whose vertices represent the successive refinements $\Omega^L_N$, $L \in \N$, and whose edges represent the ultrametric distance. These trees are indicated by the vertical dotted line segments attached to the leaves $\Box$ at level $0$. The hierarchical distance is the same above and below level $L=0$, representing $\Omega_N$, while the hierarchical group distance at level $-L$ shrinks to zero as $L\to\infty$ (see \eqref{hierdist}--\eqref{ag33} below). Every infinite path downwards represents a Cauchy sequence in $\Omega_N^\infty$.}
\label{fig-conthierargr}
\end{figure}

\medskip\noindent
{\bf (II) Evolution:}
Refining the \emph{evolution mechanisms} accordingly on the refined hierarchical group looks straightforward for migration and for block-resampling due to ultrametricity (even though taking an \emph{infinite} refinement needs some care). However, the combination with single-site Fleming and Cannings resampling is more subtle. We therefore construct the process as a \emph{limit} of processes running on countable successive refinements of $\Omega_N$, and characterise it as the solution of a \emph{well-posed martingale problem}. A key tool is the \emph{duality} to single-site coalescing and block-coalescing L\'evy-processes on $\Omega_N^\infty$. If we stop the refinement at a finite level $-L$, then we obtain a discrete group $\Omega_N^L$ for which we can define our process when we have single-site Cannings resampling at this level and block-resampling at levels $k \in \Z$ with $k >-L$ by well-posed classical martingale problems. Afterwards we must take the limit $L \to \infty$.

\begin{remark}{\bf [Evolution parameters]}
\label{r.417}
{\rm The continuum process has three parameters: 
\begin{itemize}
\item[$\bullet$] 
$(c_k)_{k\in\Z} =$ (muli-level) migration rates.
\item[$\bullet$] 
$(\Lambda_k)_{k\in\Z} =$ (multi-level) Cannings block-resampling measures (with reshuffling).
\item[$\bullet$] 
$d^\ast = d + |\Lambda_{-\infty}| =$ effective Fleming-Viot single-site resampling rate, arising from a Fleming-Viot single-site resampling rate $d$ and a Cannings single-site resampling measure $\Lambda_{-\infty}$.
\end{itemize}
The conditions we impose on these parameters are needed for the construction of the continuum-space process and for the identification of the universality classes capturing its space-time behaviour (see \eqref{cond1}--\eqref{regcond} for details).}\hfill$\spadesuit$  
\end{remark}  

\medskip\noindent
{\bf (III) Properties:}  
We are interested in identifying universal patterns in both the \emph{small-scale} behaviour (= `downstairs micro-world') and the \emph{large-scale} behaviour (= `upstairs macro-world'). The large space-time properties of the discrete-space model are represented by specific types of dynamics in the continuum-space model, namely, the ones corresponding to migration being a stable process in $\Omega_N^\infty$, which play the role of labels for the \emph{universality classes}. We are particularly interested in identifying whether the continuum-space process is a \emph{stochastic} evolution and not just a \emph{deterministic} mass flow coming from migration (without random changes in types). This turns out to be equivalent to the question whether the states are \emph{locally mono-type} or \emph{locally multi-type} as follows. 
\begin{itemize}
\item[(1)]
Without block-resampling, we expect the continuum-space measure-valued process, starting from an initial state that is stationary and ergodic under the group action, to be stochastic when the symmetrised migration is point-recurrent and otherwise to be a \emph{deterministic mass flow} driven by the transition kernel of the L\'evy-process on $\Omega_N^\infty$ describing the \emph{effective migration} built from migration and reshuffling. For transient symmetrised migration we expect the process to converge for time $t>0$ to the \emph{deterministic limit} $\mathfrak{h} \otimes \theta$, with $\mathfrak{h}$ the Haar measure on the topological group $\Omega_N^\infty$ (see~\eqref{2.44} for a detailed description) and $\theta$ the initial type distribution on $\K$. 
\item[(2)]
With block-resampling, we expect the following: 
\begin{itemize}
\item[(i)] 
In the regime of \emph{point-recurrent} effective migration the states projected on $\K$ are atomic, the type frequencies have a \emph{density} with respect to $\mathfrak{h}$ on $\Omega_N^\infty$ that solves a stochastic equation with an \emph{SPDE-part} driven by \emph{Brownian motions} and a \emph{jump-part} driven by block-resampling and reshuffling \emph{stochastic equation}. 
\item[(ii)] 
In the regime of \emph{Harris-recurrent} effective migration but \emph{not} point-recurrent migration with sufficiently strong block-resampling, the type frequencies exhibit $\K$-projected atomic and $\Omega_N^\infty$-projected regular states, but densities evolve \emph{only} via jump processes of the density states. 
\item[(iii)]
In the regime of \emph{transient} effective migration, we get a \emph{deterministic} limiting evolution in the Meyer-Zheng topology, since the large-scale block-resampling is instantaneously counteracted by the small-scale block-resampling, moving the system back into the state of a dynamical system, as in the case of no block-resampling.
\end{itemize}
\end{itemize}
The continuum spatial Fleming-Viot process (i.e., the model without block-resampling) is locally mono-type if and only if the migration is \emph{point-recurrent}. A detailed analysis of the size of the mono-type clusters for discrete lattice models on $\Omega_N$ has been carried out in \cite{FG94,FG96} and for the continuum hierarchical group in \cite{EF96}. Also models on $\R^d$ have been studied extensively as mentioned above. In the present paper we will show that the spatial Cannings process on $\Omega_N^\infty$ can be \emph{locally mono-type} and \emph{stochastic} also for \emph{Harris-recurrent migration}, provided the block-resampling is strong enough. This is the ultrametric analogue of stochastic evolutions on $\R^2$.


\subsection{Program}
\label{ss.program}

Before stating our main theorems, we describe the tasks of the present paper and the two sequel papers \cite{gh2}, \cite{gh3}. 


\paragraph*{$\blacktriangleright$ Construction and characterisation, measure states and genealogies.}

Our first task is to control the dynamics of the population on the continuum hierarchical group $\Omega^\infty_N$ as measure-valued process (present paper) and as genealogy-valued process in a parallel paper:
\begin{itemize}
\item 
Answer the question whether large classes of continuum-space processes exist, both for the \emph{type frequencies} and for the \emph{genealogies}, and whether the \emph{continuum-space scaling limit} of the process on $\Omega_N$ exists. In particular, we are interested in finding \emph{stochastic} continuum-space dynamics, especially in the case of \emph{critically-recurrent} migration, which is analogous to Brownian motion on $\R^2$.
\item 
Build a \emph{duality theory} linking the continuum-space process and the continuum-space limit process to measure-valued dual processes driven by \emph{"site"} or \emph{block-coalescing L\'evy-processes} on $\Omega_N^\infty$.
\item 
Characterise the regimes where the type-frequencies in the continuum limit are \emph{stochastic versus deterministic}, and show that the corresponding limits have states that are \emph{locally mono-type} (and decompose the space into countably many type-clusters), respectively are \emph{locally multi-type}. In the former case, show that the $\K$-projected states are atomic and investigate for which parameters there is a \emph{density} for the type frequencies that can be characterised via a \emph{stochastic equation} with an \emph{SPDE} that is driven by Brownian motions for the case of a \emph{diffusion} and by Poisson point processes in space-time for the case of a \emph{jump process}, or a combination of both.
\item 
Give a representation of the genealogies in terms of a \emph{graphical construction} based on the \emph{L\'evy-web} analogue of the Brownian web on $\R$ that was used to construct the state of the genealogy of the population under the Fleming-Viot process on $\R^1$ in \cite{GSW16}. This is doable on ultrametric geographic spaces. Show that the dichotomy of single-ancestor versus multi-ancestor coincides with the dichotomy of stochastic versus deterministic.
\item
Show that the genealogies are locally \emph{ultrametric-measure-space-valued} when the population is locally mono-ancestor, and are \emph{dust} when the population is locally multi-ancestor. For the former case, study the evolution of the genealogy as a Markov process~\cite{gh3}.
\end{itemize}


\paragraph*{$\blacktriangleright$ Dichotomies: microscopic versus macroscopic.}

Once the above steps have been completed, our second task is to identify the \emph{small-scale} and \emph{large-scale} properties of the limiting process in space and time, and their relation to the genealogical structure. There are two dichotomies:
\begin{itemize}
\item
On the type-frequency level, distinguish between \emph{coexistence} (= globally multi-type configurations) versus \emph{clustering} (= globally mono-type configurations) for large times, in particular, for equilibria. For the latter case, describe the growth of the mono-type clusters.
\item 
On the genealogical level, for all positive times, distinguish between \emph{locally multi-ancestor} in a sufficiently small neighbourhood of a site (= locally multi-type configurations) versus \emph{locally mono-ancestor} in a small neighbourhood of a site (= locally mono-type configurations). Identify the equilibrium for the process of marked genealogies, and show convergence to this equilibrium.
\end{itemize}


\paragraph*{$\blacktriangleright$ Hierarchical mean-field limit.}

A more detailed analysis can be achieved in the \emph{hierarchical mean-field limit} $N \to \infty$, for which the space-time scales at successive hierarchical levels separate. Our third task is: 
\begin{itemize}
\item 
Study the formation of mono-type clusters in \emph{space} and in \emph{time} as $t \to \infty$ and $N \to \infty$.
\item 
Study the \emph{hot spots} near a volatility-biased sampled site. 
\end{itemize}


\subsection{Main theorems}
\label{ss.maintheorems}

In this section we list our main theorems, seven in total, referring the reader to Sections \ref{sec:trlim}--\ref{s.results} for further details in the formulations and the introduction of duality. We organise these theorems in Sections \ref{ss.1}--\ref{ss.4} according to construction, properties for finite time, properties for infinite time, and stable migration and universality classes. Proofs will be given in Section \ref{s.dual} based on the results stated and derived in Sections \ref{sec:trlim}--\ref{s.results}. Throughout the sequel we assume that the \emph{basic conditions} on the migration and resampling rates stated in \eqref{cond1}--\eqref{cond3} are in force.

As \emph{initial condition} we assume that
\begin{equation}
\label{ag97*}
\wh X_0^\infty \in \CM_{\rm loc}(\Omega_N^\infty \times \K), 
\end{equation}
with $\CM_{\rm loc}$ the set of locally finite measures, is such that its law is \emph{stationary and ergodic} under the group action
with
\be{ag97}
\wh X^\infty_0 (\,\cdot \times \K) = \mathfrak{h}(\cdot), \qquad  
\E\big[\wh X^\infty_0 (\Omega^\infty_N \times \cdot\,)\big] = \theta(\cdot).
\ee
Here, $\mathfrak{h} \in \CM_{\rm loc}(\Omega_N^\infty)$ is the \emph{Haar measure} on $\Omega_N^\infty$, normed such that the $1$-ball has mass $1$ (see \eqref{2.44}), $\E$ stands for expectation, and $\theta\in\CP(\K)$ is the mean initial type distribution on $\K$. The latter is to be viewed as a free parameter, and we 
\begin{equation}
\label{e579}
\text{exclude the trivial case where $\theta = \delta_a$ for some $a \in \K$,}
\end{equation}
for which initially there is no randomness and the dynamics is trapped. Throughout the sequel we write
\begin{equation}
\wt\CM(\Omega_N^\infty \times \K) 
\end{equation} 
to indicate the Polish space given by the subset of $\CM_{\rm loc}(\Omega_N^\infty \times \K)$ that $\wh X^\infty_0$ lies in through assumptions \eqref{ag97*}--\eqref{ag97}.  


\subsubsection{Construction}
\label{ss.1}

Our first theorem establishes by construction existence of a continuum-space process $\wh X^\infty$ on $\Omega_N^\infty$, as limit of the corresponding (by truncation) processes $\wh X^L$ on countable subgroups $\Omega_N^L$, \emph{$L$-refinements} of the hierarchical groups $\Omega_N$. The processes $\wh X^L$ were described at the beginning of Section~\ref{ss.hierarmod} via the evolution rules. Later we show that there are two regimes, \emph{stochastic} versus \emph{deterministic} (see Definition~\ref{Ddetstoch} and~\eqref{ag6} for identification of the regime). The process $\wh X^\infty$ is \emph{characterized} as the solution of a \emph{well-posed martingale problem}.
 
\beT{T.CONV}{\bf [Continuum-space Cannings process: existence and well-posedness]}
\begin{itemize}
\item[(a)] 
There exists a process $\wh X^\infty = (\wh X_t^\infty)_{t \geq 0}$ on $\wt\CM(\Omega_N^\infty \times \K)$ such that the family of processes $(\wh X^L_t)_{t \geq 0}$, $L\in\N$, living on $\mathcal{P}(\K)^{\Omega_N^L}$ and embedded in $\wt\CM(\Omega_N^\infty \times \K)$, with $\Omega_N^L$ the $L$-refinement of $\Omega_N$, satisfies
\be{ag98}
\CL\big[(\wh X^L_t)_{t \geq 0}\big] \mathop{\Longrightarrow}\limits_{L \to \infty} \CL\big[(\wh X^\infty_t)_{t \geq 0}\big],
\ee
where the convergence is in the sense of finite-dimensional distributions. The convergence also holds in law on path space in the Skorohod topology for the stochastic regime and in the Meyer-Zheng topology for the deterministic regime.
\item[(b)] 
$\wh X^\infty$ is called the $((c_k)_{k\in\Z},d^\ast,(\Lambda_k)_{k\in\Z})$-Cannings process with block-resampling, has \cadlag\, paths, is strong Markov and has the Feller property.
\item[(c)]
For every $t>0$,
\be{ag97alt}
\wh X^\infty_t (\,\cdot \times \K) = \mathfrak{h}(\cdot) \,\, a.s., \qquad  
\E\big[\wh X^\infty_t (\Omega^\infty_N \times \cdot\,)\big] = \theta(\cdot),
\ee
and the following dichotomy holds:
\begin{equation}
\label{e527}
\begin{aligned}
&\text{either } \wh X_t^\infty (B \times \cdot\,) \in \CM_{\rm fin, atom}(\K) \, \text{for any ball } B \subseteq \Omega_N^\infty,\\
&\text{or } \wh X^\infty_t (B \times \cdot\,) = \theta(\cdot)\,\, a.s.\ \, \text{for any ball } B \subseteq \Omega_N^\infty,
\end{aligned}
\end{equation}
where $\CM_{\rm fin, atom}$ is the set of finite atomic measures. In the former case, the \emph{random} set $A_t$ of types with a positive measure somewhere in $\Omega_N^\infty$ and visible in $\wh X_t^\infty$ at time $t>0$ is non-increasing in $t$ and is random, and $A = \cup_{t>0} A_t$ is a countable random subset of $\K$.
\item[(d)] 
The process $\wh X^\infty$ is a semi-martingale measure (see \cite[Chapter 7]{D93}) and is the unique solution of a well-posed martingale problem for the sum of a mass-flow operator, a jump-process operator and a diffusion operator, with the latter two of a non-standard form. (This is detailed in Section~\ref{sec:wpmp}, in particular, in Definition~\ref{def.martiro}.)
\end{itemize}
\end{theorem}

\noindent
The process exhibits a dichotomy of being either \emph{stochastic} or \emph{deterministic} (see Theorems~\ref{T.det} and~\ref{T.STOCHDET} below). Further properties of $\wh X^\infty$ are given in Theorem~\ref{T.char} below (and depend on the particular parameter regime of interest). 

We refer to Section \ref{ss.cont} for the definition of $\Omega_N^L$ and the dynamics on $\Omega_N^L$ or $\Omega_N^\infty$. We refer to Section \ref{sss.dupro} for the new concepts appearing in Theorem \ref{T.CONV} and for an explanation why the Meyer-Zheng topology is needed (the reason being that $\wh X^\infty$ instantaneously \emph{collapses onto a subspace} of the state space, namely, the atomic measures), and to Section \ref{ss.nsmp} for an explanation why the martingale problem is non-standard.  

\begin{remark}{\bf [Explanation of the outcome]}
\label{r.548}
{\rm In part (a) of Theorem \ref{T.CONV}, the distinction between the two path topologies is important and is related to a \emph{collapse property}: the collapse to a subset of the state space where the process can diffuse via collision local times and can jump at unbounded, but effectively finite, jump rates. In part (b), the properties are standard. What happens in part (c) will become more clear via our analysis of the \emph{process} of the \emph{genealogies} in \cite{gh2}, from which we conclude, in particular, that in the stochastic regime the entrance law can be viewed as a Markov process on a Polish space, by strengthening the weak topology to the \emph{atomic weak topology} and by decomposing the population at a given time into at most countably many \emph{mono-ancestor subfamilies}, representable in geographic space by different centers and a weight. What happens in part (d) is analogous to what happens with the diffusion operator in the case of \emph{catalytic} or \emph{mutually catalytic} branching processes \cite{DF97a,DP98}, where the diffusion is also of non-standard, namely, driven by \emph{collision local times}.
}\hfill$\spadesuit$
\end{remark}

\begin{remark}{\bf [Duality]}
\label{r.duality}
{\rm A key property of the constructed Markov process $\wh X^\infty$ is the existence for every time horizon $t >0$ of a dual \emph{backward process} $Y^\infty$ and a dual \emph{measure-valued} process $\wh Y^{\ast,\infty}$ driven by a \emph{finite} or even an \emph{infinite particle system} in form of a marked partition valued spatial coalescent, namely, block-coalescence in $k$-balls $(k \in \Z \setminus \{0\})$ and single-site \emph{delayed coalescing} L\'evy-processes $Y^\infty$ on $\Omega_N^\infty$ (which also is the case for the genealogy process). This is detailed in Proposition~\ref{P.DUAL} and Proposition~\ref{P.1430} and the surrounding sections. The dual allows us to trace the \emph{genealogies of samples} of individuals taken at time $t$ all the way back to their origin at time $0$ (see~\cite{gh2}). This specific structure is behind many of the statements in the present paper. There is even a \emph{strong duality}, allowing us to represent the full measure state of $\wh X^\infty$ as an \emph{entrance law} of a \emph{measure-valued dual state} $\wh Y^{\ast,\infty}$ corresponding to an \emph{entrance law} of the spatial coalescent $Y^\infty$ driving the measure-valued dual process and starting from a locally countably infinite state. This entrance law of the dual backward process starts with \emph{countably many} dual individuals locally at each site of a \emph{countable dense subset} of $\Omega_N^\infty$. (See Section \ref{ss.dualdisc} for $\Omega_N$ and Section \ref{sss.dupro} for $\Omega_N^\infty$.)} \hfill$\spadesuit$
\end{remark}

\begin{remark}{\bf [A related model on $\Omega_N^\infty$ with infinite volatility]}
\label{r.634}
{\rm The \emph{instantaneously} coalescing L\'evy-processes $\wt Y^\infty$ arise as the dual associated with a \emph{different} stochastic process, namely, the \emph{$d^\ast \to \infty$ limit of the process $\wh X^\infty$ studied in the present paper}. The finite volatility process is of more interest due to its SPDE-representation (see Remark~\ref{rem.1968}). The infinite volatility process appears in the context of the analysis of the cluster formation for the model on $\Omega_N$ (see~\cite[Theorem 2, Section 5.1]{Klenke1996}). In the present paper we denote the infinite volatility forward process by $\wh X^{\infty,\sim}$. This will play some role in the proofs, and is of interest because it appears when we study monotype-cluster formation in models on $\Omega_N$, or for $\wh X^\infty$ on $\Omega_N^\infty$, and similarly on $\Z$, respectively, $\R$.}\hfill$\spadesuit$
\end{remark}


In special cases (see~\cite{EF96,Klenke1996}) it is known that we can let $d \to \infty$ and obtain a stochastic Markov limit process on $\wh X^{\infty,\sim}$ describing cluster formation on $\Omega_N$ in spatial Fisher-Wright models.

\begin{corollary}{\bf [Infinite volatility limit: new class of processes]}
\label{cor.779}
The processes $\wh X^\infty$ in the limit as $b \to \infty$ converge to a limiting process, denoted by $\wh X^{\infty,\sim}$, that is a stochastic strong Markov process with c\`adl\`ag paths. For this process the remaining theorems in this section hold as well, except the item in Theorem~\ref{T.char} below~\eqref{e710}.
\end{corollary} 

\noindent
We will see later in Section~\ref{sec:dual} that this amounts to passing in the dual process from \emph{delayed} coalescence to \emph{instantaneous} coalescing L\'evy-processes. This process arises in the study of the growth of monotype regions as $t \to \infty$ in the model on $\Omega_N$ (see~\cite{Klenke1996} for the stable case).

\begin{remark}{\bf [Transient migration]}
\label{r.773}
{\rm At first glance, the above observation looks paradoxical for the transient case, because due to the block-coalescence the L\'evy-processes seem not to be independent. However, the genealogies of two \emph{typical} individuals of initially countably \emph{infinite} populations in the forward and the backward picture do have paths of descend with this property. In \cite{gh2} we will see that the genealogies of the populations in the limit as $L \to \infty$ form dust, i.e., at time $t>0$ all subfamilies descending from one ancestor have \emph{frequency zero} in every ball, and hence two \emph{typical} paths of the population are independent. This leads to a \emph{law of large numbers} for the state of the process $\wh X^\infty$ as $L \to \infty$ and hence for $t>0$ to \emph{deterministic} limit states.}\hfill$\spadesuit$
\end{remark}

\subsubsection{Properties for finite time}
\label{ss.2}

The following two dichotomies characterise the properties of the orbit of the process $\wh X^\infty$ for finite time.

\begin{definition}{\bf [Deterministic versus stochastic]}
\label{Ddetstoch}
{\rm The process $\wh X^\infty$ is said to be \emph{deterministic} when the following holds: For all $t>0$,
\be{e547}
\CL[\wh X^\infty_t] = \delta_{\mu_t} \text{ for some } \mu_t \in \wt\CM(\Omega_N^\infty \times \K).
\ee
Otherwise $\wh X^\infty$ is said to be \emph{stochastic}.} \hfill$\spadesuit$  
\end{definition}

\noindent
Note that \emph{deterministic} still allows for a \emph{random} change from time $t=0$ to time $t >0$ (which is what actually happens).

\begin{definition}{\bf [Locally mono-type versus locally multi-type]}
\label{Dmonomulti}
{\rm The process $\wh X^\infty$ is said to be \emph{locally mono-type} when the following holds: For all $t>0$ and $\mathfrak{h}$-a.e.\ $x \in \Omega_N^\infty$, there exists a $k=k(t,x) \in \Z$ such that
\be{e558}
\wh X_t|_{B_k(x) \times \K} = \mathfrak{h}|_{B_k(x)} \otimes \delta_a \text{ for some } a = a(t,x,k) \in \K,
\ee
where $B_k(x)$ is the ball of radius $k$ around $x$. Otherwise $\wh X^\infty$ is said to be \emph{locally multi-type}.} 
\hfill$\spadesuit$ 
\end{definition}

Our second theorem shows that the two dichotomies coincide.
 
\begin{theorem}{\bf [Coincidence of dichotomies]}
\label{T.det}
$\wh X^\infty$ is deterministic if and only if it is locally multi-type. In particular, $\wh X^\infty$ stochastic means that it is locally mono-type.
\end{theorem}

\noindent
A similar dichotomy of stochastic versus deterministic appears in critical branching~\cite{D93}, in catalytic and mutually catalytic branching systems~\cite{DF97a,DP98}, with many unresolved issues remaining.

Our third theorem classifies the dichotomy in terms of the migration and resampling rates, and comes directly from the duality. Abbreviate 
\be{lambdakdefs}
\lambda_k = \int_{(0,1]} \Lambda_k(\dd r), \qquad \lambda^*_k = \int_{(0,1]} r^{-1} \Lambda_k(\dd r), \qquad k \in \Z.
\ee 
and assume that, on top of the \emph{basic conditions} in \eqref{cond1}--\eqref{cond3}, also the \emph{non-triviality condition} in \eqref{condextra} and the \emph{extra regularity conditions} in \eqref{regcond} are met. We say that \emph{resampling does not dominate migration} when
\begin{equation}
\label{criterion}
\sup_{k\in\Z} \frac{\lambda_k}{c_k} < \infty.
\end{equation}

\begin{theorem}{\bf [Criterion for dichotomy]}
\label{T.STOCHDET}
Assume \eqref{condextra}--\eqref{regcond}.
\begin{itemize}
\item[(a)] 
Suppose that $\lambda_k \not\equiv 0$ for some $k\in\Z$. If \eqref{criterion} holds, then
\be{ag6}
\sum_{k\in\N_0} \frac{1}{c_k} \sum_{\ell=0}^k \lambda_\ell = \infty
\quad \Longleftrightarrow \quad \wh X^\infty \mbox{ is stochastic}. 
\ee
If \eqref{criterion} fails, then $\wh X^\infty$ is stochastic. 
\item[(b)] 
Suppose that $\lambda_k \equiv 0$ for all $k\in\Z$. If $d^\ast > 0$, then
\be{ag6alt}
\sum_{k \in \N_0} \frac{1}{c_k} = \infty,\quad \sum_{k \in -\N} \frac{1}{c_k} < \infty
\Longleftrightarrow \quad \wh X^\infty \mbox{ is stochastic}.
\ee
If $d^\ast=0$, then $\wh X^\infty$ is deterministic.
\end{itemize}
\end{theorem}

\noindent
The criterion in \eqref{ag6} only looks `upstairs', while the criterion in \eqref{ag6alt} looks `upstairs and downstairs'. The condition $\lambda_k \not\equiv 0$ means that there is some $k \in \Z$ for which $\lambda_k>0$. 

The above classification is related to the behaviour of the difference process of two independent copies of the migration process. We will see in Appendix \ref{appA} that:
\begin{itemize}
\item
$\sum_{k \in \N_0} \frac{1}{c_k} = \infty$, $\sum_{k \in -\N} \frac{1}{c_k} < \infty$: the migration is point-recurrent\\ 
(and $\wh X^\infty$ is stochastic when $\lambda_k \equiv 0$ for all $k\in\Z$).
\item 
$\sum_{k \in \N_0} \frac{1}{c_k} = \infty$, $\sum_{k \in -\N} \frac{1}{c_k} = \infty$: the migration is Harris-recurrent but not point-recurrent \\ (and $\wh X^\infty$ is deterministic when $\lambda_k \equiv 0$ for all $k\in\Z$). 
\end{itemize}
The criterion in \eqref{ag6} includes right-transience (i.e., $\sum_{k\in\N_0} (1/c_k) < \infty$), can be expressed in terms of the \emph{degree} of the underlying motion (see \cite{DGW05}), is of a general nature and therefore also applies to the pair $\Z^d,\R^d$. Thus, point-recurrent migration always leads to a locally mono-type dynamics, both in the presence and in the absence of block-resampling, while Harris-recurrent migration and transient migration only do when the block-resampling is \emph{strong enough}. 

Our fourth theorem shows that key properties of $\wh X^\infty$ can be characterised for each $t>0$ in terms of a backward $(\Omega_N^\infty)^2$-process running from $t$ to $0$, $((\zeta_1^\infty (t-s),\zeta_2^\infty (t-s)))_{s \in [0,t]}$ obtained by sampling \emph{two typical individuals}, denoted by $(\zeta_1^\infty,\zeta_2^\infty)$, from the dual population of the \emph{entrance law} of the dual process starting with locally countably many dual individuals to be introduced formally later in Section~\ref{sec:dual} and denoted $\wh Y^\infty$ which is automatically (by being typical) starting from two different sites. The proof will show that the characterisation of the properties of $\wh X^\infty$ in terms of two L\'evy-processes is of a subtle nature and is not just the duality relation applied to two dual individuals. In fact, a full understanding requires an understanding of the genealogies of the population exhibited in~\cite{gh2}.

\begin{remark}{\bf [Distinction]}
{\rm Care is needed to distinguish between the dual process with two initial particles starting at two locations, i.e., a pair following the evolution of the finite spatial dual process, and two sampled particles from a dual population starting in a locally countably infinite configuration. This distinction is clarified in Lemma~\ref{l.1562}.}\hfill$\spadesuit$
\end{remark}

A \emph{deeper analysis} shows that we can classify the behaviour of $(\zeta_1^\infty,\zeta_2^\infty)$ in terms of the effective L\'evy-process. Depending on the recurrence-transience properties of the \emph{effective L\'evy-process $\zeta^\infty$}, the process is stochastic or deterministic. In the stochastic regime we can characterise $\wh X^\infty$ as the solution to a \emph{stochastic equation}. For the latter we need to keep in mind a crucial point, namely, that we have $\mfh$-populations rather than collections of populations in sites. It is only the type-structure that is discrete for $t>0$. 

\begin{theorem}{\bf [Properties and stochastic equation of continuum-space Cannings process]}
\label{T.char} 
Let $\zeta^\infty$ be the effective L\'evy-process. Then the process $\wh X^\infty$ constructed via \eqref{ag98} has the following properties in the three regimes listed below, which correspond to $\zeta^\infty$ being transient, $\zeta^\infty$ being Harris-recurrent but not point-recurrent, respectively, $\zeta^\infty$ being point-recurrent.
$\mbox{}$
\begin{itemize}
\item[(a)] 
If $(\zeta_1^\infty,\zeta_2^\infty)$ coalesce with probability $0$, then $\wh X^\infty$ is deterministic and constant equal to $\mathfrak{h} \otimes \theta$ for $t >0$.
\item[(b1)]
If $(\zeta_1^\infty,\zeta_2^\infty)$ do not hit each other but coalesce with probability $>0$, then $\wh X^\infty$ is stochastic. The coalescence can only occur via block-coalescence, $\wh X^\infty$ is a jump process (= finite variation), has countably many jumps, off the jumps has continuous paths with finite variation, and solves a stochastic equation based on a collection of Poisson point processes on $(0,\infty) \times \Omega_N^k$, $k \in \Z$. 
\item[(b2)] 
If $(\zeta_1^\infty,\zeta_2^\infty)$ hit each other with probability $1$, then $d^\ast>0$ and they coalesce eventually as $t\to \infty$, the process $\wh X^\infty$ is stochastic, and for $t >0$ has states such that:
\begin{itemize} 
\item
$\wh X_t (B \times \cdot\,)$ is atomic, with $B \subseteq \Omega_N^\infty$ a ball and a random set $A$ of atoms in $\K$.
\item
$\wh X_t^\infty (\cdot \times \{a\})$ is absolutely continuous with respect to $\mathfrak{h}$, with $a \in A \subseteq \K$ an atom and $A$ a countable random set in $\K$ that is a function of the path $\wh X^\infty$ and consists of all types occurring in $\wh X^\infty$ somewhere in space for some strictly positive time (i.e., $\wh X^\infty_t$($B \times \{a\})>0$ for some $t>0$ and some $a \in A$ and ball $B$). 
\item 
There exists for every $a \in A$ a version of the density of $\wh X_t^\infty(\cdot,a)$ the family, called 
$$
(\wh X^{\infty,\dagger}_t(\eta,\{a\})_{a \in A},
$$  
which are \emph{jointly continuous} in $(t,\eta) \in (0,\infty) \times \Omega_N^\infty$, except in block-resampling jump times in $t$, where it is c\`adl\`ag and specifies the weights of the atoms when observing the population in every $\mfh$-finite subset $B$ in $\CB(\Omega_N^\infty)$. For fixed $t$, as function of $\eta$ the values of $\wh X_t^{\infty,\dagger} (\eta,\{a\})$ are $\mfh$-a.e.\ $0$ or $1$, defining a random decomposition of $\Omega_N^\infty$ into $\mfh$-monotype balls. 
\end{itemize}
Given the above random set $A \subseteq \K$, consider the corresponding countable collection of densities of $\wh X^\infty$ w.r.t.\ $\mfh$ on $\Omega_N^\infty$:
\begin{equation}
\label{e710}
\text{$(\wh X_t^{\cdot,\ddag}(\eta,a))_{a \in A}$.}
\end{equation}
The law of this object is the entrance law of a solution of a stochastic equation (written out in \eqref{a.1917} for a countable set of components), and the equation has the following components. 
\begin{itemize}
\item 
For the continuous infinite-variation part of the path, the stochastic equation contains an \emph{SPDE-part} based on an $A$-indexed collection $W=\{W^a, a \in A\}$ of white noise measures on $\Omega_N^\infty$.
\item 
There is c\`adl\`ag jump-part of the path based on $\Pi$, a $\Z$-indexed collection of Poisson point processes on $[0,\infty) \times \Omega_N^k, k \in \Z$ (with $\Omega_N^k$ viewed as subset of $\Omega_N^\infty$). The latter generate the events of block-resampling with reshuffling in $k$-balls.
\item 
There is a continuous finite-variation part given by a deterministic measure flow via the transition semigroup of the process $\zeta^\infty$. 
\end{itemize}
Here, the noises $W,\Pi$ are independent for given $A$ and of the initial state $\wh X_0^\infty$. The conditional (on $A$) weak solution is unique for a prescribed $\wh X_0^\infty \in \wt\CM(\Omega_N^\infty \times \K)$, and hence is strongly unique.
\end{itemize}
\end{theorem}

\noindent
See Section \ref{334} for details on the form of the stochastic equation (in particular, \eqref{a.1917} with white noise) and on the Poisson point processes describing the densities associated with the types arising from the entrance law.


\subsubsection{Properties for infinite time}
\label{ss.3}

Our fifth theorem establishes convergence to equilibrium.  

\begin{theorem}{\bf [Convergence to equilibrium]}
\label{T.conveq}
\begin{equation}
\lim_{t\to\infty} \CL[\wh X_t^\infty] = \nu_\theta \text{ exists, for some equilibrium } 
\nu_\theta \in \CP\big(\wt\CM(\Omega_N^\infty \times \K)\big),
\end{equation}
and satisfies 
\be{ag97altalt}
\nu_\theta(\pi_{\Omega_N^\infty} \wh X^\infty \in \cdot\,) = \delta_\mathfrak{h}(\cdot), \qquad  
\E_{\nu_\theta} \left[\pi_\K X \in \cdot\, \right] = \theta(\cdot),
\ee
where $\pi_{\Omega_N^\infty}$ denotes projection on $\Omega_N^\infty$ and $\pi_{\K}$ denotes projection on $\K$.
\end{theorem}

The following dichotomy characterises the equilibrium behaviour of $\wh X^\infty$.

\begin{definition}{\bf [Globally mono-type versus globally multi-type]}
\label{D.glob}
{\rm The equilibrium measure $\nu_\theta$ is said to be \emph{globally mono-type} when
\begin{equation}
\nu_\theta = \int_\K \theta(\dd a)\,\delta_{\mathfrak{h} \otimes \delta_a}.
\end{equation} 
Otherwise $\nu_\theta$ is said to be \emph{globally multi-type}.}\hfill$\spadesuit$ 
\end{definition}  

\noindent
We will see that if $\nu_\theta$ is globally multi-type, then it takes the form $\nu_\theta = \delta_{\mathfrak{h} \otimes \theta}$. In the literature the dichotomy in Definition \ref{D.glob} is also referred to as \emph{clustering} versus \emph{co-existence}. 

Our sixth theorem shows that all three dichotomies coincide.

\begin{theorem}{\bf [Coincidence of all dichotomies]}
\label{T.coindichs}
The equilibrium $\nu_\theta$ is globally multi-type if and only if it is locally multi-type. 
\end{theorem}

Theorem~\ref{T.conveq} raises the question, in the mono-type regime, how fast the mono-type clusters grow in time. This will be answered in detail in Sections~\ref{ss.clgr1}--\ref{ss.clgr2}, and exhibits multiple growth schemes. 


\subsubsection{Stable migration and universality classes}
\label{ss.4}

\begin{example}{\bf [Example: stable continuum-space Cannings process]}
\label{ex.1035}
{\rm The choice
\be{ag40b}
d, \Lambda_{-\infty}, \qquad c_k = c^k, \,\lambda_k = \lambda^k, \,\Lambda_k = \lambda^k \Lambda, \quad k \in \Z,
\ee
with $d \in [0,\infty)$, $\Lambda_{-\infty} \in \CM((0,1])$, $c \in (0,N)$, $\lambda \in [0,N)$, and $\Lambda \in \CM((0,1])$ satisfying \eqref{ag7d}, represents the class of \emph{stable continuum-space Cannings processes} with \emph{block-resampling} (where the word stable refers to the migration being a \emph{stable process} on the topological group $\Omega_N^\infty$). This class is invariant under a rescaling operation (see Section \ref{ss.resc}; finer conditions are listed below \eqref{cond1}).} \hfill$\spadesuit$
\end{example}

\noindent
For stable continuum-space Cannings processes, the sum condition in \eqref{ag6} holds when either $\lambda \in (1,N)$, $c \in (0,\lambda]$ or $\lambda \in (0,1]$, $c \in (0,1]$ (corresponding to the stochastic regime, irrespective of \eqref{criterion}), and fails when either $\lambda \in (1,N]$, $c \in (\lambda,N)$ or $\lambda \in (0,1]$, $c \in (1,N)$ (corresponding to the deterministic regime, with \eqref{criterion} in force). The sum condition in \eqref{ag6alt} holds when $\lambda=0$, $c \in (0,1)$ and fails when $\lambda=0$, $c \in [1,N)$. For $c=1$ the migration is Harris-recurrent but not point-recurrent, while the sum condition in \eqref{ag6} is met for any $\lambda \in (0,N)$. Thus, for \emph{critically-recurrent migration}, turning on the block-resampling causes a switch from \emph{locally multi-type} (= deterministic) to locally mono-type (= stochastic). See Appendix \ref{appA} for further details.

Our seventh and last theorem shows that stable $\wh X^\infty$ arise as space-time scaling limits of Cannings processes with block-resampling on $\Omega_N$, and as such represent the \emph{universality classes} of the processes on $\Omega_N$ for the large space-time scaling properties. 
  
\begin{theorem}{\bf [Convergence in the stable domain of attraction and universality classes]}
\label{T.CONVRES}
Suppose that the initial state of the measure-valued process $\wh X = (\wh X_t)_{t \geq 0}$ on $\Omega_N$ is stationary and ergodic under translations, with single-site marginal type distribution $\theta \in \CP(\K)$. Consider the measure-valued process $\wh X^L = (\wh X^L_t)_{t \geq 0}$ on $\Omega^L_N$, the $L$-refinement of $\Omega_N$, obtained after an appropriate $L$-dependent scaling of space, time, migration rates and resampling rates. Then, subject to certain further regularity conditions,
\be{ag39}
\lim_{L\to\infty} \CL\big [(\wh X_t^L)_{t \geq 0}\big] = \CL\big[(\wh X_t^{\infty,\rm sta})_{t \geq 0}\big], 
\ee
where the limit is the stable version (i.e., the parameters satisfy~\eqref{ag40b}) of the process $\wh X^\infty$ in Theorem \ref{T.CONV} with scaled rates $c,\Lambda,\lambda=|\Lambda|$ and with initial state $\mathfrak{h} \otimes \theta$. Namely, the limit process has stable effective migration-reshuffling motion on the topological group $\Omega_N^\infty$. The same assertion holds for the scaling limit of Cannings processes on $\Omega_N^\infty$ constructed in Theorem \ref{T.CONV}.
\end{theorem}

\noindent
For the scaling of space, time, resampling rates and migration rates, see \eqref{ag2a}--\eqref{ag2c} in Section \ref{ss.resc}. For the further regularity conditions, see \eqref{ag7a}--\eqref{ag7c} in Section \ref{ss.resc}. For the parameters of the limiting process in terms of the original parameters, see \eqref{ag7c} and \eqref{ag13} in Section \ref{ss.resc}.

Thus, the stable continuum-space Cannings processes represent the \emph{universality classes} for the large space-time scaling behaviour of both discrete-space \emph{and} continuum-space Cannings processes, and the latter are \emph{fixed points} under the \emph{rescaling transformation} meaning the system is reproduced under application of this rescaling limit. We can identify the deterministic regime of \emph{L\'evy-mass flow} with the parameter range $1 < c \vee \lambda < N$, which represents one universality class. We see that for the case of \emph{stochastic} dynamics, which corresponds to the parameter range $c \vee \lambda \leq 1$, there are \emph{three} universality classes, namely, $c>\lambda$, $c < \lambda$, $c=\lambda$, i.e., the stable pure \emph{Fleming-Viot resampling} universality class with continuous infinite-variation paths, the pure \emph{Cannings block-resampling} universality class with c\`adl\`ag paths of finite variation, and the universality class arising from combinations of the two. There is \emph{no} single-site Cannings universality class.

\paragraph*{Application to cluster formation on $\Omega_N$} 
The proof of the above theorems leads to a generalisation and extension of the results obtained in the literature, which can be focused in Lemma~\ref{cor.2748}. It requires several pages to formulate, so that we put this in the proof section.


\subsection{Outline of the remainder of the paper}
\label{ss.outline}

The remainder of this paper builds up the tools that are needed to prove the theorems in Section \ref{ss.maintheorems} and provides a deeper understanding of the phenomena that transpire from these theorems. Section \ref{s.mi} reviews some basic facts about the discrete and the continuum hierarchical group, about random walks and L\'evy-processes running on them, and about the spatial Cannings process on the discrete hierarchical group that was introduced and analysed in \cite{GHKK14}, including duality with a spatial block-coalescent. Sections \ref{sec:trlim}--\ref{sec:wpmp} focus on the \emph{continuum-space dynamics}, Section \ref{s.results} on the \emph{universality} via a \emph{continuum-space scaling limit}, i.e., the scaling of the discrete-space dynamics towards the continuum-space dynamics when space, time and the migration and resampling rates are scaled properly. Section \ref{s.dual} provides the proofs of the main theorems based on the results in Sections \ref{sec:trlim}--\ref{s.results}. Appendix \ref{appA} collects a few basic facts about L\'evy-processes that are needed along the way. Appendix \ref{appB} mentions alternative routes to the construction of the continuum-space process. 

Sections \ref{sec:dual}--\ref{sec:wpmp} are long and demanding because they deal with new challenges surrounding the dual process and the martingale problem.


\section{Model ingredients: discrete space versus continuum space}
\label{s.mi}

In Section \ref{ss.disc} we recall the definition and the construction of the discrete hierarchical group $\Omega_N$ and the spatial Cannings process on $\Omega_N$ introduced in \cite{GHKK14}. In Section \ref{ss.dualdisc} we recall the duality of this process with a system of block-coalescing random walks, which was central in the analysis carried out in \cite{GHKK14} and is also here a key tool. In Section \ref{ss.cont} we define the continuum hierarchical group $\Omega^\infty_N$ and indicate what functions of the model parameters are needed to prepare for the construction of the spatial Cannings process on $\Omega^\infty_N$ in Sections \ref{sec:trlim}--\ref{sec:wpmp}. In Section~\ref{ss.cont} we introduce L\'evy-processes on $\Omega^\infty_N$ in a form that is needed to prepare for the duality with a system of block-coalescing L\'evy-processes. In addition, we give an outline of what is done in Sections \ref{sec:trlim}--\ref{sec:wpmp}.


\subsection{Spatial Cannings process on $\Omega_N$}
\label{ss.disc}


\subsubsection{Discrete hierarchical group $\Omega_N$}

The discrete hierarchical group is the set given by
\be{ag27}
\Omega_N = \left\{\xi = (\xi^i)_{i \in \N_0}\colon\, \xi^i \in \{0,1,\ldots,N-1\}\,\,\forall\,i\in\N_0,\, 
\sum\nolimits_{i\in\N_0} \xi^i <\infty\right\},
\ee
with component-wise addition $\mathrm{mod}(N)$ as algebraic structure and with the topology induced by the hierarchical distance
\be{ag28}
k(\xi, \eta) = \inf \{k \in \N_0 \colon\, \xi^i = \eta^i\,\, \forall\,i \geq k\}, \qquad \xi,\eta \in \Omega_N,
\ee
which is an \emph{ultrametric} and turns $\Omega_N$ into a Polish space (see Fig.~\ref{fig-hierargr}). For $k \in \N_0$ and $\eta\in\Omega_N$, write 
\be{}
B_k(\eta) = \{\xi\in\Omega_N\colon\,k(\xi,\eta) \leq k\},
\ee 
to denote the $k$-ball around $\eta$, i.e., the macro-colony consisting of all the $N^k$ colonies that are within distance $k$ of colony $\eta$.

An alternative way to define $\Omega_N$ is as
\be{}
\Omega_N = \oplus_{\N_0} Z_N,
\ee  
with $Z_N$ the cyclic group of $N$ elements. 


\subsubsection{Migration: random walks}

Let $(c_k)_{k \in \N_0} \in (0,\infty)^{\N_0}$ be a sequence of positive numbers satisfying
\be{ag35}
\sum_{k\in\N_0} \frac{c_k}{N^k} < \infty.
\ee
Define a \emph{migration rate kernel} on $\Omega_N \times \Omega_N$ by setting
\be{fdhkernel}
a^{(N)}(\xi,\eta) = \sum_{k \geq k(\xi,\eta)} \frac{c_{k-1}}{N^{k-1}}\,\frac{1}{N^k}, \quad \xi \neq \eta,
\qquad a^{(N)}(\eta,\eta) = 0.
\ee
Individuals migrate according to continuous-time random walks with kernel $a^{(N)}(\cdot,\cdot)$, i.e., for every $k\in\N$ each individual jumps at rate $c_{k-1}/N^{k-1}$ to a position that is uniformly distributed in the $k$-ball around its current position. Writing the migration rate kernel in the form \eqref{fdhkernel} is important because the sequence $(c_k)_{k \in \N_0}$ characterises the potential theory associated with the random walk. Note that the random walk is \emph{hierarchically exchangeable}.


\subsubsection{Resampling in single colonies}

We first recall the definition of the Cannings model, which describes the evolution of the types of finitely many individuals living in a single colony. Let $M \in \N$ be the number of individuals, and let $\K$ be a compact Polish space encoding the types. The evolution of the population, whose state space is $\K^M$, is as follows. We fix a volatility constant
\be{ag36a*}
d \in [0,\infty),
\ee
which gives us the rate at which \emph{Fleming-Viot resampling} occurs, i.e., each pair of individuals is chosen independently at rate $d$ and is replaced by two new individuals that are descendants of one of the parent individuals chosen at random, each inheriting the type of that parent. 

We fix a \emph{resampling measure} $\Lambda \in \CM_{\text{fin}}([0,1])$, the set of finite measures on $([0,1],\CB([0,1])$, satisfying
\be{ag36a}
\Lambda(\{0\}) = 0, \qquad \int_{(0,1]} r^{-1}\Lambda(\dd r) = \infty.
\ee
We define \emph{Cannings resampling} as follows (see Fig.~\ref{fig-canevol}):

\begin{figure}[htbp]
\vspace{0.4cm}
\centering
\includegraphics[width=0.8\linewidth]{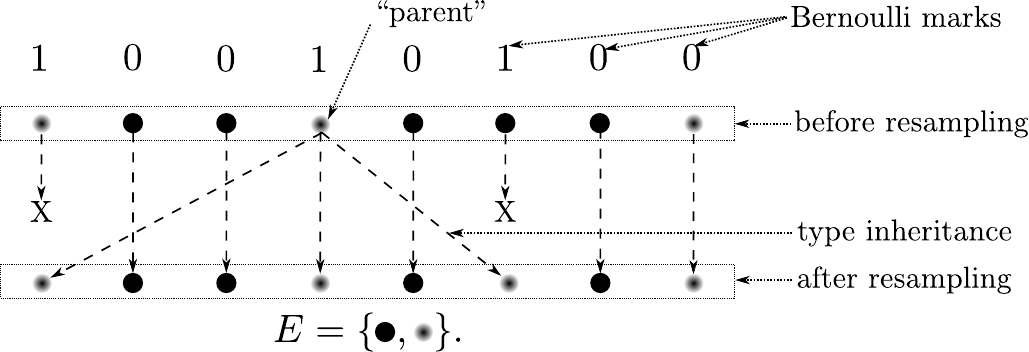}
\caption{\small Cannings resampling event in a colony of $M=8$ individuals of two types. Arrows indicate type inheritance, X indicates death.}
\label{fig-canevol}
\end{figure}

\begin{itemize}
\item
The number of individuals stays fixed at $M$ during the evolution.
\item
Initially, i.i.d.\ types are assigned to the individuals according to a given type distribution
\be{AG5}
\theta \in \CP(\K).
\ee
\item
To represent the Cannings transition as a Poisson point processes, we need to change parameters. Let $\Lambda^*$ be the $\sigma$-finite measure on $([0,1], \CB([0,1])$ defined as
\be{AG3}
\Lambda^*(\{0\}) = 0, \qquad \Lambda^*(\dd r) = r^{-2}\Lambda(\dd r), \quad r \in (0,1].
\ee
Consider an inhomogeneous Poisson point process on $[0,\infty) \times [0,1]$ with intensity measure
\be{AG4}
\dd t \otimes \Lambda^*(\dd r).
\ee
For each point $(t,r)$ in this process, we carry out the following transition at time $t$. \emph{Mark} each of the $M$ individuals independently with a $1$ or a $0$, with probability $r$, respectively, $1-r$. \emph{All} individuals marked by a $1$ are killed and are replaced by copies of a \emph{single} individual (= the parent) that is drawn randomly from all the individuals marked by a $1$. 
\end{itemize}
In this way we obtain a \emph{Markov pure-jump} process, which we call the \emph{$\Lambda$-Cannings process} with resampling measure $\Lambda$ and population size $M$ and initial state $\theta^{\otimes M}$.

\begin{remark}{\bf [Conditions on the resampling measure]}
{\rm The first condition in \eqref{ag36a} excludes Fleming-Viot resampling, which is already encoded in the parameter $d$. The second condition in \eqref{ag36a}, which is called the \emph{`no-dust condition'}, is needed to have well-defined type frequencies for the descendants of a single ancestor at positive times in the population in the large-colony-size limit $M\to\infty$ (because possibly $d=0$), and to be able to define a genealogical tree for the infinite population. Indeed, to avoid `dust formation', we need a countable number of jumps in the limit as $M \to \infty$ in every colony, so that for any $\varepsilon,\delta>0$ a fraction $\geq 1-\varepsilon$ of the individuals descend from finitely many ancestors living a lapse $\delta>0$ back in time (carrying finitely many types). We return to this issue below \eqref{ak:multi-level-generator}, in the discussion of the dual. Thus, we need either $d>0$ or $|\Lambda| >0$, and the no-dust condition for $\Lambda$ when $d=0$ to have well-defined genealogies both in the process and in its dual entrance law starting from countably many dual individuals.} \hfill$\spadesuit$
\end{remark}

Taking the type frequencies in each colony as the state and letting $M \to \infty$, we get a \emph{continuum-mass jump process} (with countably many small jumps in any positive length interval), which we call the \emph{$\Lambda$-Cannings process} with initial state $\theta$. Under \eqref{ag36a}, this process has countably many jumps in every time interval of positive length. The Fleming-Viot resampling generates, in the $M \to \infty$ limit, a diffusion process on $\CP(\K)$, the Fleming-Viot diffusion.

The $(d,\Lambda)$-resampling just described will be used in every \emph{single colony}. We next add resampling in \emph{blocks of colonies}.


\subsubsection{Resampling in blocks of colonies}

\begin{figure}[htbp]
\vspace{0.4cm}
\label{fig:reshuffling} 
\centering
\includegraphics[width=0.8\linewidth]{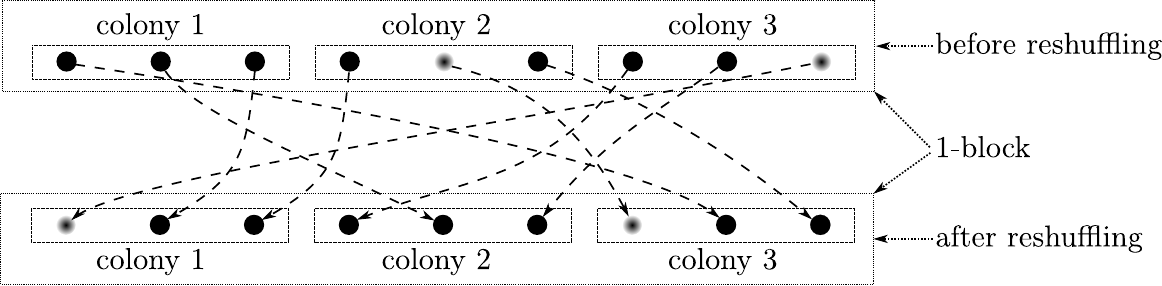} 
\caption{\small
Random reshuffling in a 1-block on the hierarchical lattice of order $N=3$ with $M=3$ individuals of two types per colony.}
\label{fig-reshuffle}
\end{figure}

Let $(\Lambda_k)_{k \in \N} \in (\CM_{\text{fin}}([0,1]))^\N$ be a sequence of finite measures on $([0,1], \CB([0,1])$. Let individuals \emph{reshuffle-resample} by carrying out the following two steps at once: 
\begin{itemize}
\item
For every $\eta\in\Omega_N$ and $k\in\N$, choose $B_k(\eta)$ at rate $1/N^{2k}$. Each individual in $B_k(\eta)$ is first moved to a uniformly chosen random location, i.e., a \emph{reshuffling} takes place inside $B_k(\eta)$ (see Fig.~\ref{fig-reshuffle}). 
\item
After the reshuffling, \emph{resample} the population in $B_k(\eta)$ with resampling measure $\Lambda_k$. 
\end{itemize}
Note that the reshuffling-resampling affects all the individuals in all the macro-colonies simultaneously and in the same manner. Because the reshuffling is done first, the resampling always acts on a uniformly distributed state (`panmictic resampling'). In what follows we will talk about `block-resampling' rather than `block-reshuffling-resampling'. Later we must remember that the reshuffling adds to the migration and that the resampling has two forms, the \emph{Cannings single-site resampling} and the \emph{Cannings block-resampling}, in particular, $\Lambda_0$ and $(\Lambda_k)_{k \in \N}$ need to be distinguished.  

\begin{remark}{\bf [Reshuffle only $1$-marked individuals]}
{\rm It will be expedient to use a \emph{modified version} of the resampling-reshuffling mechanism in this paper and would have been conceptually more clear in~\cite{GHKK14} as well. Namely, only $1$-marked individuals are reshuffled and are replaced by the type of a randomly chosen $1$-marked individual. The difference is that the \emph{change of location only applies to the $1$-marked individuals in the block.}}\hfill$\spadesuit$
\end{remark}

Note that $\Lambda_k^\ast$, the intensity of resampling and reshuffling in $k$-blocks, is not necessarily finite. Because $\Lambda_k$ is finite, the rate of change in the type remains bounded as soon as $\int_{(0,1]} r^{-1}\Lambda_k(\dd r)<\infty$, and so the latter is the only restriction needed in the block-resampling. However, for the reshuffling this condition is not enough \emph{without change of the mechanism to only reshuffle $1$-marked individuals}. We therefore assume here in this paper that we are given a \emph{single-site resampling measure $\Lambda_0$} and a sequence of \emph{$k$-block-resampling measures} $(\Lambda_k)_{k \in \N}$ satisfying
\be{e654}
\Lambda_k \in \CM_{\text{fin}}([0,1]), \quad \Lambda_k (\{0\}) = 0, \qquad k \in \N_0.
\ee
and impose the conditions 
\be{AG4bb}
\int_{(0,1]} r^{-1} \Lambda_0 (\dd r) = \infty, \qquad
\int_{(0,1]} r^{-1} \Lambda_k (\dd r) < \infty, \quad k \in \N.
\ee

\begin{remark}{\bf [Stronger condition]}
{\rm 
In~\cite{GHKK14,GHK18}, instead of the second part of \eqref{AG4bb} the \emph{stronger} condition 
\begin{equation}
\label{e574}
\int_{(0,1]} r^{-2} \Lambda_k(\dd r) < \infty, \qquad k \in \N,
\end{equation}
was imposed because \emph{all} the individuals were reshuffled even those not part of the resampling event. This introduces, compared to our version here, an additional migration, called the \emph{reshuffling}, which is of a different form because all the jump times of all the individuals happen at one time, namely, the resampling jump time. The latter condition implies that for each $k\in\N$ the number of jumps due to $k$-block-coalescence in a time interval of positive length is finite, resulting in a \emph{Markov pure-jump process}. Since in the present paper we only reshuffle individuals marked by a $1$, we can do with the \emph{weaker} condition in the second half of \eqref{AG4bb}, since this ensures that the jump rate goes to zero as $r \to 0$. The latter allows for \emph{countably many jumps} due to $k$-block-coalescence, in particular, it allows for the stable models we discussed earlier, the stable L\'evy-processes, which play an important role because they represent universality classes.}\hfill$\spadesuit$
\end{remark}
 
To get a finite rate for block-resampling, we assume that
\be{AG4c3}
\sum_{k \in\N} \frac{\lambda_k}{N^k} < \infty 
\ee
with
\be{AG4c}
\lambda_k = \Lambda_k((0,1]) = |\Lambda_k| \in [0,\infty).
\ee


\subsubsection{Process as solution to a martingale problem} 

We equip the state space
\be{e890}
E=\CP(\K)^{\Omega_N}
\ee
with the product topology to make it Polish. Let $\mathcal{F} \subset C_\mathrm{b}\big(E,\R\big)$ be the algebra of \emph{polynomials} generated by \emph{monomials}, i.e., functions of the form
\be{ak:multi-level-test-functions}
\begin{aligned}
&F(x) = \int_{\K^n} \left(\bigotimes_{m=1}^n x_{\eta_m}\big(\dd u^m\big)\,
\psi^{(n)}\big(u^1,\ldots,u^n\big)\right),
\quad x = (x_{\eta})_{\eta\in\Omega_N} \in E,\\
&n \in \N, \quad \psi^{(n)} \in C_{\mathrm{b}}(\K^n,\R), 
\quad \eta_1,\ldots,\eta_m \in \Omega_N.
\end{aligned}
\ee
The linear operator for the martingale problem
\be{ak:multi-level-generator}
L^{(\Omega_N)}\colon\,\mathcal{F} \to C_\mathrm{b}\big(E,\R\big)
\ee
has two parts:
\be{ak:multi-level-generator-decomposition}
L^{(\Omega_N)} = L^{(\Omega_N)}_{\mathrm{mig}} + L^{(\Omega_N)}_{\mathrm{res}}.
\ee
The \emph{migration operator} is given by
\be{ak:multi-level-migration}
(L^{(\Omega_N)}_{\mathrm{mig}} F)(x)
= \sum_{\eta,\zeta \in \Omega_N} a^{(N)}(\eta,\zeta)
\int_\K (x_{\zeta} - x_{\eta})(\dd a)\,\frac{\partial F(x)}{\partial x_\eta}[\delta_a]
\ee
and the \emph{reproduction} operator is given by
\bea{ak:multi-level-resampling-global}
(L^{(\Omega_N)}_{\mathrm{res}} F)(x)
&=& \sum_{k\in\N_0} N^{-2k} \int_{[0,1]} \Lambda^*_k (\dd r) \int_\K y_{\eta,k}(\dd a) \sum_{\eta \in \Omega_N} 
\left[F\left(\Phi_{r,a,B_k(\eta)}(x)\right)-F(x)\right]\\
&&+\, d \int_\K \int_\K Q_{x,\eta} (\dd u,\dd v)\, \frac{\partial^2 F}{\partial x^2_\eta} [u,v], \nonumber
\eea
where 
\be{typedist}
y_{\eta,k} = \frac{1}{N^k} \sum_{\zeta \in B_k(\eta)} x_\zeta
\ee
is the type distribution on $B_k(\eta)$, $\Phi_{r,a,B_k(\eta)}\colon\,E \to E$ is the \emph{reshuffling-resampling map} acting as
\be{ak:resampling-mapping-global}
\begin{aligned}
\Big[\big(\Phi_{r,a,B_k(\eta)}\big)(x)\Big]_\zeta 
&=
\begin{cases}
(1-r) x_\zeta + r \delta_a, &\zeta \in B_k(\eta),\\
x_{\zeta}, &\zeta \notin B_k(\eta),
\end{cases}\\
&\qquad \qquad r \in [0,1],\, a \in \K,\, k\in\N_0,\,\eta\in\Omega_N,
\end{aligned}
\ee
and
\begin{equation}
\label{e602}
Q_{x,\eta} (\dd u,\dd v) = \big[x_\eta \otimes \delta_u - x_\eta \otimes x_\eta\big](\dd u,\dd v),
\qquad x \in E,\,\eta \in \Omega_N,\,u,v \in \K.
\end{equation}
The first term in the right-hand side of \eqref{ak:multi-level-resampling-global} is the Cannings-reshuffling resampling operator, the second term is the operator of the Fleming-Viot resampling at rate $d$. 

\begin{remark}{\bf [Taylor expansion]}
{\rm The first term in the right-hand side of \eqref{ak:multi-level-resampling-global} is well-defined because of \eqref{AG4bb}--\eqref{e574}. Indeed, by Taylor-expanding the inner integral in powers of $r$ and recalling \eqref{ak:multi-level-test-functions}, we get
\be{intTay}
\int_\K y_{\eta,k}(\dd a) \left[ F(\Phi_{r,a,B_k(\eta)}(x)) - F(x) \right] = O(r^2), 
\qquad r \downarrow 0,
\ee
where we note that the term of order $r$ is zero because of \eqref{typedist}--\eqref{ak:resampling-mapping-global}. To get a well-defined reshuffling operator, the expression in \eqref{intTay} must be integrable with respect to the block-resampling measure $\Lambda^*_k$, which is guaranteed by \eqref{e574}, and summable over $k$, which is guaranteed by \eqref{AG4bb}.} \hfill $\spadesuit$
\end{remark}

\begin{remark}{\bf [Effective migration]}
\label{rem:eff}
{\rm Note that the block-coalescence with \emph{reshuffling} induces an extra migration \emph{coupled} to block-resampling. We refer to the sum of the migration and the reshuffling as the \emph{effective migration}.} \hfill $\spadesuit$
\end{remark}

The following is taken from \cite[Proposition 1.4]{GHKK14} (which treats the original model, but the proof can be adapted to the modified model). 

\begin{proposition}{\bf [Existence and uniqueness of the spatial Cannings process]}
\label{D.mvC}
$\mbox{}$\\
The $(\CF, L^{(\Omega_N)},\delta_x)$-martingale problem with $x \in E$ has a unique solution 
\be{mp1}
(X_t)_{t \geq 0}
\ee
and is called the $((c_k)_{k \in \N_0}, d,\Lambda_0,(\Lambda_k)_{k \in \N})$-Cannings process with block-resampling. This process has the property that for positive times the component processes $\{(x_\xi (t))_{t \geq 0}\colon\, \xi \in \Omega_N\}$ take values in $\CM_{\rm atomic}(\K)$ if $d$, $\Lambda_0$, $(\Lambda_k)_{k \in \N}$ are not all zero.
\end{proposition}

\begin{remark}{\bf [Representation by stochastic equations]}\label{r.repres}
{\rm It is known that under the no-dust condition the process has an at most countable \emph{random} set, say $A$, of types for the orbit for all $t>0$. This allows us to write down, based on Ito's-formula, an $(\Omega_N \times A)$-system of \emph{stochastic equations} representing the process $X$ as the unique solution to this system.}\hfill$\spadesuit$
\end{remark}


\subsection{Dual process: coalescent on $\Omega_N$}
\label{ss.dualdisc}

A key role in the analysis of the spatial Cannings process is played by its \emph{moment dual} driven by a spatial coalescent, which was introduced and analysed in \cite{GHKK14}. This coalescent takes values in the partitions of $[n]=\{1,\ldots,n\}$ (or $\N$) with partition elements carrying a location mark and has $\Lambda_0$-\emph{single-site-coalescence} and $(\Lambda_k)_{k\in\N}$-\emph{block-coalescence} with \emph{reshuffling} and \emph{migration} of partition elements. We denote the partition by $\pi=(\pi(1),\ldots,\pi(n))$, where $\pi$ is a map from $[n]$ to $\{1,\ldots,|\pi|\}$, and denote the partition elements by $\pi^{-1}(i)$, $i=1,\ldots,|\pi|$. The locations are given via a map $\xi$ from partition elements to $\Omega_N$ such that $\xi_{\pi(1)}, \ldots, \xi_{\pi(n)}$ are the respective locations of the individuals $1,\ldots,n$, where all individuals of a partition element have the same position. We also have a measure-valued dual as entrance law.


\paragraph*{Evolution of the partition valued dual process.} 

The migration of the coalescent is given by the random walk with kernel $a^{(N)}(\cdot,\cdot)$ on $\Omega_N \times \Omega_N$, and \emph{reshuffling} is induced by the block-coalescence. Coalescence in a $k$-block works as follows:
\be{e805}
\mbox{At intensity } \dd t \otimes \Lambda_k \mbox{ choose points in } [0,\infty) \times [0,1], 
\mbox{ denoted by } (t,r).
\ee
At time $t$, \emph{mark} every partition element $1$ with probability $r$ and $0$ with probability $1-r$. All partition elements marked $1$ are merged at time $t$ into a single new partition element. After that all partition elements marked $1$ are relocated uniformly.

The dual process on $\Omega_N$ is denoted by
\be{ag37a}
(Y_t)_{t \geq 0}
\ee
and is a \emph{spatial $\Lambda$-coalescent with block-resampling}, taking values in the \emph{$\Omega_N$-marked partitions} (finite or countable) and evolving as follows:
\begin{itemize}
\item 
Partition elements move independently on $\Omega_N$ according to $a^{(N)}(\cdot,\cdot)$.
\item
Pairs of partition elements on $\Omega_N$ at the same site coalesces independently at rate $d$.
\item 
Partition elements on $\Omega_N$ at the same site $\Lambda_0$-coalesce.
\item
For each $k \in \N$, partition elements located in the same $k$-block undergo $\Lambda_k$-coalescence, and the \emph{new} partition element is instantaneously relocated in the $k$-block uniformly.
\end{itemize}
The dual process $Y$ can be enriched to $Y^\ast$ by giving a partition element at time $s$ a type from $\K$ for all $s \in [0,t]$, by "sampling" this type from $X_0$ at the location at backward time $t$ and hence getting $(\Omega_N \times \K)$-marked partition elements. To define rigorously what sampling means here with a locally finite measure needs a bit of effort, we explain for the continuum space case later.


\paragraph*{Entrance law of the measure-valued dual process.}

We can also define the \emph{entrance law} of $Y$, starting at time 0 from the state with \emph{countably} many individuals at each site, as the limit of locally finitely many individuals with per site diverging occupation numbers (see~\cite{GLW05}) provided that $d>0$ or $|\Lambda_0|>0$. We can give partition elements a leading element. 
In the dual process, the no-dust condition for the coalescent allows us to define the spatial coalescents, starting in countably many partition elements and get in positive time partitions of the populations and furthermore define for partition elements at every location in $\Omega_N$ \emph{positive frequencies} summing up to one at each site. We can then assign the partition element a location type and location by sampling at the position of its leading element a type from $\K$ according to the marginal at the location. This construction allows us to define for time horizon $t$ at time $s \leftarrow [0,t]$ a measure $\wh Y^\ast$ on $\Omega_N \times \K$, with weight $1$ on sets $\{\eta\} \times \K$ for $\eta \in \Omega_N^\infty$. We call this measure $\wh Y_t^\ast$ at $s=t$. 

Note, to get this object a point of concern is, that the process \emph{does not come down from infinity} due to the migration of dual elements out of a site, leading always to countably many partition elements hence this condition is too strong. However, due to the \emph{condition} in the first part of \eqref{AG4bb}, the frequencies of the partition elements that do not belong to a finite set of partition elements can at every site in $\Omega_N$ be made arbitrarily small, so that the number of partition elements per site is \emph{finite} for positive time, \emph{except for partitions with frequency $\leq \varepsilon$} for every $\varepsilon >0$. Without this condition, `dust' would form with positive frequency as the number of dual individuals tends to infinity.


\paragraph*{Duality relation.}

We recall the \emph{duality relation} for the process on $\Omega_N$. To do so, we need some ingredients to construct the duality function $H(\cdot,\cdot)$. For $n \in \N$, $L^{(n)}$ a function on the partition elements of a partition of $\{1,\ldots,n\}$ having values in $\Omega_N^\infty$, and $\Psi_\K^{(n)} \in C_b(\K^n, \R)$, consider the function
\be{ag50}
H^{(n)}_{L^{(n)},\Psi^{(n)}} (x, \pi)
= \int_{(\Omega_N \times \K)^n} x^{\otimes n} (\dd \xi,\dd u)\, 
1_{\xi_1}(\eta_1) \times \cdots \times 1_{\xi_n}(\eta_n)
\, \Psi^{(n)} (u_{\pi(1)}, \ldots, u_{\pi(n)}),
\ee
where $x \in E = \mathcal{M}_{\mathrm{loc}}(\Omega_N \times \K)$, $\pi \colon\, [n] \to [n]$ is an $\Omega_N$-marked partition with marks $\eta_i = L^{(n)}(\pi^{-1} (i))$ for the partition element containing $i$. For $E^\prime$ the set of functions $\Psi \in \cup_{n \in \N}\, C_b(\K^n,\R)$ paired with $\Omega_N$-marked partitions $\pi \in \cup_{n \in \N} \Pi_n$ (with $\Pi_n$ the set of partitions of $[n]$), define
\be{ag51}
H (x,(\Psi, \pi)) = H^{(n)}_{L^{(n)},\Psi^{(n)}}(x,\pi), \qquad x \in E, \, (\Psi,\pi) \in E'.
\ee
The following is taken from \cite[Proposition 2.10]{GHKK14} for the original model, but can be adapted to the present setting. 
($\PP$ and $\E$ stand for probability and expectation.)

\begin{proposition}{\bf [Duality relation for the Cannings process]}
\bea{ag52}
\E\big[H (X_t,(\Psi,Y_0))\big] = \E\big[H(X_0, (\Psi, Y_t))\big], \qquad (X_0, Y_0) \in E \times E^\prime.
\eea
\end{proposition}
It is possible to strengthen this to a \emph{strong duality, a measure-valued duality} of the processes $X$ and $\wh Y^\ast$ and represent the law of $X_t$ by the law of $\wh Y_t^\ast,$
\begin{equation}\label{e1197}
 \CL [X_t] = \CL [\wh Y_t^\ast].
\end{equation}


\subsection{Spatial Cannings process on $\Omega^\infty_N$}
\label{ss.cont}

We want to approximate $\Omega_N^\infty$ through \emph{successive refinements} of $\Omega_N$ (recall Figs.~\ref{fig-hierargr}--\ref{fig-conthierargr}), namely, by splitting each site into $N$ lower-level sites, splitting each of these sites into $N$ lower-level sites, etc., and taking the closure to obtain $\Omega_N^\infty$ as a limiting set. Alongside we also want to refine the evolution mechanism based on our sequences of parameter indexed with $\Z$, namely, the migration and the resampling, with the aim to get a proper limiting processes on $\Omega_N^\infty$ from well-known processes on the hierarchical group.


\subsubsection{Continuum hierarchical group $\Omega_N^\infty$}
\label{sss.ONLdef}

For $L \in\N$, define the \emph{$L$-refined} version of $\Omega_N$ where each site in $\Omega_N$ splits into $N^L$ new sites each carrying a hierarchical structure with $L$ levels, i.e.,
\be{ag29}
\Omega^L_N = \left\{\xi=(\xi^i)_{i \geq -L}\colon\,
\xi^i \in \{0,1,\ldots,N-1\}\,\,\forall\,i \geq -L,\,\,\sum\nolimits_{i \geq -L} \xi^i < \infty\right\},
\ee
with component-wise addition $\mathrm{mod}(N)$. Define on $\Omega_N^L$ the \emph{hierarchical distance} 
\be{ag31}
k^L(\xi,\eta) = \min \{k \geq -L\colon\, \xi^i = \eta^i\,\, \forall\, i \geq k\}
\ee
and the \emph{hierarchical group distance}
\be{ag30}
d^L(\xi,\eta) = \left\{
\begin{array}{ll}
k^L(\xi,\eta), &k^L(\xi,\eta)>0,\\
\eee^{k^L(\xi,\eta)}, &-L< k^L(\xi,\eta) \leq 0,\\
\eee^{-L}, &k^L(\xi,\eta) = -L \text{ when } \xi \neq \eta, \text{ and } 0 \text{ otherwise}.
\end{array}
\right.
\ee
Define a $k$-balls on $\Omega_N^L$ as
\be{ag32}
B^L_k(\xi) = \{\eta \in\Omega^L_N \colon k^L(\xi,\eta) \leq k\}.
\ee
In this setup the new space geographic $\Omega_N^L$ can be thought of as a \emph{refinement} of degree $L$ of the old geographic space $\Omega_N$ with $N^L$-times as many single sites (recall Figs.~\ref{fig-hierargr}--\ref{fig-conthierargr}). 

Define the hierarchical distance on $\Omega_N^\infty$
\be{hierdist}
k^\infty(\xi,\eta) = \inf\big\{k \in \Z\colon\,\xi^i = \eta^i\,\,\forall\, i \geq k\big\}
\ee
and the hierarchical group distance $\Omega_N^\infty$
\be{ag30inf}
d^\infty(\xi,\eta) = \left\{
\begin{array}{ll}
k^\infty(\xi,\eta), &k^\infty(\xi,\eta)>0,\\
\eee^{k^\infty(\xi,\eta)}, &k^\infty(\xi,\eta) \leq 0,
\end{array}
\right.
\ee
obtained by extension based on $\Omega_N^L$. We can visualise $\Omega_N^\infty$ by imagining that the leaves of the tree representing $\Omega_N$ are refined to infinite $N$-ary trees running downwards (recall Fig.~\ref{fig-hierargr}). Note that sites $(\xi_i)_{i \in \Z}$ agreeing for all $i \geq -L$ have distances converging to $0$ as $ L \to \infty$. This allows us to define
\be{ag33}
\Omega^\infty_N = \left\{(\xi^i)_{i \in \Z} \colon \xi^i \in \{0,1,\ldots,N-1\}\,\,\forall\,i \in \Z,\,
\sum\nolimits_{i\in\N_0} \xi^i < \infty\right\}
\ee
equipped with the metric above to get a \emph{complete separable metric space} that is \emph{locally compact}. In fact, $\Omega_N^\infty$ is an \emph{uncountable} Abelian \emph{topological group} that forms an \emph{ultrametric} Polish space. Note that $\Omega_N^\infty$ is \emph{very different} from the Euclidean space: (i) it is \emph{totally disconnected}; (ii) it can be \emph{decomposed} into disjoint open \emph{and} closed balls of arbitrarily small positive radius; (iii) it can be seen as the \emph{set of leaves} of the completion of an infinite $N$-ary tree. This means, for example, that processes with independent increments are always jump processes: there is no Brownian motion on $\Omega_N^\infty$. There is a \emph{L\'evy-Khintchine representation} for such jump processes. Moreover, on $\Omega_N^\infty$ a \emph{Fourier analysis} is available that has tractable features due to the \emph{ultrametric} structure \cite{BCR84}. 

As a topological group, $\Omega_N^\infty$ has a Haar measure, unique up to a constant factor and we set 
\be{ag33b}
\mathfrak{h} = \mbox{\emph{normed Haar measure} on } \Omega^\infty_N,
\ee
where normed means that $\int_{\Omega_N^\infty} \mathfrak{h}(\dd\eta)\,1_{\{k^\infty(\xi,\eta) \leq k\}} = N^k$, $k\in\Z$, for all $\xi \in \Omega_N^\infty$. See for a description~\eqref{e1183} and~\eqref{2.44}

The specific maps and homomorphisms given by
\be{ag34}
\begin{array}{lll}
&(\xi^i)_{i \in \N_0} \mapsto (0, \ldots, 0, \xi^0, \xi^1, \ldots), 
&\mbox{ with } L \mbox{ zeros before } \xi^0,\\[0.2cm]
&(\xi^i)_{i \geq -L} \mapsto (\ldots, 0, \xi^{-L}, \xi^{-L+1}, \ldots),
&\mbox{ with infinitely many zeros before } \xi^{-L},
\end{array}
\ee
allow us to \emph{embed} $\Omega_N$ into $\Omega^L_N$ and $\Omega^L_N$ into $\Omega^\infty_N$ \emph{as groups}. Write $0=0^\Z$ to denote the \emph{origin} in $\Omega_N^\infty$. Note that, via this natural embedding,
\be{limHaus}
\Omega^\infty_N = \lim_{L\to\infty} \Omega_N^L,
\ee
where the convergence is in the \emph{Hausdorff metric for sets}. 

An alternative way to define $\Omega_N^\infty$ is as
\be{e1183}
\Omega_N^\infty = \big(\oplus_{\N_0} Z_N\big) \times (Z_N)^{-\N},
\ee  
with $Z_N$ the \emph{cyclic group} of $N$ elements, corresponding to the `upstairs macro-world' and the `downstairs micro-world'. Accordingly, we can write
\be{2.44}
\mathfrak{h} = \mathfrak{h}_{\N_0} \otimes \mathfrak{h}_{-\N},
\ee
where $\mathfrak{h}_{\N_0}$ is the counting measure on $\oplus_{\N_0} Z_N$ and $\mathfrak{h}_{-\N}$ is the product over $-\N$ of the uniform distribution on $Z_N$.  


\subsubsection{Migration: L\'evy-processes}
\label{Levydef}

The \emph{migration} of our population model on $\Omega_N^\infty$ will be given by a \emph{flow of mass} induced by the semigroup associated with the stochastic process $\zeta^\infty$ on $\Omega_N^\infty$. The process $\zeta^\infty$ has independent increments arising as a sum over $k\in\Z$ of jumps in blocks of size $k$, drawn according to the normed Haar measure in that block, at a rate that depends on $k$. We write this rate in the form $c_{k-1}/N^{k-1}$, $k \in \Z$, since this choice of parametrisation makes the potential-theoretic properties of the migration expressible in terms of the sequence $(c_k)_{k \in \N_0}$. To define the migration on $\Omega_N^\infty$, we need that the two-sided sequence
\begin{equation}
(c_k)_{k \in \Z} \in (0,\infty)^\Z \text{ still satisfies the one-sided condition in \eqref{ag35}.}
\label{e775}
\end{equation}
With the two-sided sequence we associate a \emph{L\'evy-measure} $m$ on $\Omega^\infty_N$ given by
\be{ag36}
\dd m = g\,\dd\mathfrak{h}, \qquad g(\xi) = \sum_{k \geq k^\infty(\xi,0)} 
\frac{c_{k-1}}{N^{k-1}}\,\frac{1}{N^k}, \qquad \xi \in \Omega^\infty_N,
\ee
which specifies the \emph{intensity of the jumps}. The intensity of the jumps in a ball are multiplied by the inverse of the volume of the ball. In this form, the recurrence and transience properties of the L\'evy-process can be expressed in terms of $(c_k)_{k \in \Z}$ simultaneously for all $N$. On that scale the L\'evy-process can explore small neighbourhoods of a site for suitable $(c_k)_{k \in \Z}$ and for all choices of $N$.

According to \cite{E89}, $m$ uniquely defines a \emph{L\'evy-process} in $\Omega^\infty_N$, which we denote by
\be{Lpr}
\zeta^\infty = (\zeta^\infty_t)_{t \geq 0}.
\ee 
This process has \emph{c\`adl\`ag} paths, \emph{stationary independent increments} and \emph{infinitely divisible marginals} that can be represented as a sum of points in a Poisson point process with intensity measure $m$. 

The generator $\mathfrak{A}$ of $\zeta^\infty$ has a \emph{domain} that is dense in the set of continuous functions and is given by
\be{e801}
\begin{aligned}
\CD(\mathfrak{A}) &= \Big\{f\in\wt\CM(\Omega_N^\infty)\colon\,\text{ there exist } B>1 \text{ and } C<\infty \text{ such that }\\
&\qquad |f(\eta+\xi)-f(\eta)| \leq CB^{k^\ast(\xi)} \text{ for all } \eta,\xi \in \Omega_N^\infty \text{ with } k^\ast(\xi) \leq 0\Big\},
\end{aligned}
\ee
where 
\be{}
k^\ast(\xi) = \min\{i\in\Z\colon\,\xi^j = 0\,\,\forall\,j>i\}.
\ee
The condition in \eqref{e801} is to be read as requiring a modulus of continuity that is uniform in space.

On $\CD(\mathfrak{A})$ the \emph{generator} $\mathfrak{A}$ of $\zeta^\infty$ is defined as (compare with~\cite[Section 7.1]{DGZpr})
\begin{equation}
\label{e857}
(\mathfrak{A}f)(\eta) = \int_{\Omega_N^\infty} [f(\eta+\xi)-f(\eta)]\,m(\dd \xi), \qquad \eta \in \Omega_N^\infty. 
\end{equation}
For a general theory about L\'evy-processes on $\R^d$ we refer the reader to \cite{PS71a}, \cite{PS71b}.


\subsubsection{Effective migration and coalescing L\'evy-processes}

In our dual model there is block-coalescence with \emph{reshuffling}, which induces extra migration \emph{coupled} to block-resampling (recall Remark~\ref{rem:eff}). Nevertheless, it turns out that consequences of this extra migration are captured by the path of a single dual individual or pairs of individuals, and therefore the \emph{effective migration} induced by the migration and reshuffling turns out through a deeper analysis to be completely controlled by the sequence of parameters 
\begin{equation}
\label{effmig}
\frac{c_{k-1}}{N^{k-1}} + \frac{\lambda^*_k}{N^k}, \qquad k \in \Z,
\end{equation}
with $\lambda^*_k = \int_{(0,1]} r^{-1} \Lambda_k(\dd r)$ (see \eqref{cond1} below), describing the motion of a \emph{single tagged} individual over time. Note that migration and reshuffling \emph{differ} when we consider $\ell$ tagged individuals with $\ell > 2$, since the latter make their jumps simultaneously.

Return to the L\'evy-process $\zeta^\infty$, which is the path in geographic space an individual follows due to the $(c_k)_{k \in \Z}$ migration process. In Appendix~\ref{appA:oneLP} in Proposition~\ref{prop.one-four-regimes}, we give criteria when this process is \emph{transient, point-recurrent}, or \emph{Harris-recurrent} but \emph{not} point-recurrent. In the point-recurrent case we identify the local time process of $\zeta$ (see Appendix~\ref{appA:pot}).

Consider the \emph{distance-from-$0$ process} $(k^\infty(\zeta^\infty_t,0))_{t \geq 0}$, which is a Markov process on $\Z$.
Define
\be{}
\sigma_0^+ = \sum_{k\in\N_0} \frac{1}{c_k}, \qquad
\sigma_0^- = \sum_{k \in -\N} \frac{1}{c_k}.
\ee
We will see in Appendix~\ref{appA:oneLP} that, under the condition in \eqref{criterion}, the distance-from-$0$ process is
\be{criteria}
\begin{array}{ll}
&\bullet \text{ recurrent when } \sigma_0^+ = \infty,\,\sigma_0^-=\infty,\\
&\bullet \text{ transient to the left when } \sigma_0^+ = \infty,\,\sigma_0^-<\infty,\\
&\bullet \text{ transient to the right when } \sigma_0^+ < \infty,\,\sigma_0^-=\infty,\\
&\bullet \text{ transient to the right and to the left when } \sigma_0^+ < \infty,\,\sigma_0^-<\infty.
\end{array}
\ee

The same considerations apply to the additional migration that we get for an individual from the reshuffling in the Cannings block-resampling, which we called the effective migration. However, the additional replacements of the marked individuals occur simultaneously and at the same time as the resampling, whereas at rate $c_k$ individuals migrate independently, and independently of their resampling events, which leads to a mass flow for our process $\wh X^\infty$ via the $(c_k)_{k \in \Z}$-transition kernel. Hence, on top of mass flow according to the continuous-time migration kernel, we also have \emph{jumps at random times}, where all individuals in a $k$-ball reshuffle. This causes a mass redistribution that takes place deterministically in blocks (= balls) by reshuffling, accompanied by a \emph{random} change of the type frequencies. This in turn creates \emph{jumps}, and so only we get c\`adl\`ag paths by applying a deterministic map at random times, maps for which uniform measures in space on $k$-blocks are preserved for the level-$k$ reshuffling part of the map. Note that no change occurs by the reshuffling when the types are uniformly distributed over $\ell$-balls with $\ell \leq k$.

With the above \emph{effective migration process} on the topological group $\Omega_N^\infty$, we have a system $Y^\infty$ of \emph{coalescing L\'evy-processes} on $\Omega_N^\infty$, with single-site coalescence of pairs at rate $d^\ast$ based on local times, independently of the migration, and block-coalescence with rate $\Lambda_k$ (coupled to a block-reshuffling). On $\Omega_N$ we also had Canning single-site resampling with rate $\Lambda_0$, requiring the additional condition in the first part of \eqref{AG4bb}, which is now replaced by the same condition for $\Lambda_{-\infty}$ for finite $L$ approximations. 


\subsubsection{Stable migration}

A special class of spatial Canning processes, which are the prototype of the possible universality classes, are the ones with \emph{stable migration} for a tagged individual.

\begin{remark}{\bf [Example: stable migration]}
\label{R.levy}
{\rm The choice 
\be{ag36b}
c_k = c^k, \quad k \in \Z, \qquad \text{ for some $c \in (0,N)$,}
\ee 
corresponds to a \emph{stable L\'evy-process} on $\Omega_N^\infty$. This process has the scaling property
\be{ag555}
\left(\pi_L\big(\zeta^\infty_{(N/c)^Lt}\big)\right)_{t \geq 0} {\stackrel{\wedge}{=}}\, 
\big(\zeta^\infty_t)_{t \geq 0},
\ee
where $\stackrel{\wedge}{=}$ denotes equality in distribution and $\pi_L$ is the shrinking operator defined by 
\be{shr}
\pi_L((i_\ell)_{\ell \in \Z}) = (i_{\ell-L})_{\ell \in \Z}.
\ee 
The latter can be written as
\be{ag36c}
\pi_L(\zeta^\infty_{(N^L)^\alpha t}) \,{\stackrel{\wedge}{=}}\, (\zeta^\infty_t)_{t \geq 0}, 
\qquad \alpha=1-\frac{\log c}{\log N} \in (0,\infty).
\ee
Because of the form of the L\'evy-measure in \eqref{ag36}, the \emph{transition kernel density} exists for $t>0$ and scales like $t^{-\frac{1}{\alpha}+o(1)}$, $t \downarrow 0$, as shown in Appendix \ref{app:stable}. This makes it analogous to an \emph{$\alpha$-stable processes} on $\R$ (see also \cite{E89,EF96,FG96}).} \hfill$\spadesuit$
\end{remark}


\subsubsection{Conditions on migration and resampling rates}

The continuum model requires the following \emph{basic conditions}:
\begin{eqnarray}
\label{cond1}
&&c_k \in (0,\infty), \quad \Lambda_k(\{0\}) = 0, \quad \lambda_k,\,\lambda^*_k \in [0,\infty), \qquad k \in \Z,\\[0.2cm]
\label{cond2}
&&\sum_{k\in\N_0} N^{-k} c_k < \infty, \qquad \sum_{k \in\N_0} N^{-k} \lambda_k < \infty,\\
\label{cond3}
&&d \in [0,\infty), \quad \Lambda_{-\infty}(\{0\}) = 0, \quad \int_{(0,1]} r^{-1}\Lambda_{-\infty}(\dd r) = \infty, \\
\label{cond4}
&& \int_{(0,1]} r^{-1} \quad \Lambda_k (\dd r) < \infty \text{ for $k \geq L$ for some $L \in \Z$.}
\end{eqnarray}
We do not need to impose a condition on the limit points of $\Lambda_k$ for $k \to -\infty$ other than on $\lambda_k$ for $k \to -\infty$.

The second condition is a requirement only for the \emph{upper tail} of the migration and block-resampling rates, and is needed to make sure that the \emph{total rate} of migration and block-resampling is finite for all transitions acting on levels $\geq -L$ for some $L \in \N$. The transitions on infinitely many lower level needs no such requirements, even though the accumulated effects on these locations still remain, and for changes of measures this has to be shown to make sense. We discuss this in the sections that introduce the operators of the martingale problems. The third condition ensures that the single-site resampling is well-defined and satisfies the \emph{no-dust condition}.

Recall at this point that $d$ is the Fleming-Viot resampling rate and that $\Lambda_{-\infty}$ is an \emph{auxiliary resampling measure} that replaces $\Lambda_{-L}$ for the approximating process on $\Omega_N^L$ \emph{before} the limit $L\to\infty$ is taken. The second half of \eqref{cond3} is effectively a no-dust condition for \emph{all} the approximation processes on $\Omega_N^L$ because of the choice of resampling mechanism at level $-L$ made in \eqref{e988} below.

We further require the following \emph{non-triviality condition}: 
\begin{equation}
\label{condextra}
\sum_{k\in\Z} N^{-k} c_k = \infty, \qquad \sum_{k \in \Z} N^{-k} \lambda_k^\ast \in \{0,\infty\}.
\end{equation}
This guarantees that the migration and the reshuffling of a typical individual accumulates countably many small jumps in every positive time interval for balls in space, and therefore differs from a pure Markov jump process (which we consider to be trivial). Furthermore, also for \emph{block} resampling (if present) the resampling is either a diffusion or a process with countably many small jumps in finite time for balls in space, or both, but not just a Markov jump process with a deterministic drift term.  

In addition, for some of the computations we require \emph{extra regularity conditions}:
\begin{equation}
\label{regcond}
\sup_{k\in\Z} \frac{c_{k+1}}{c_k} < N, \qquad 
\sup_{k\in\Z} \frac{\lambda_{k+1}}{\lambda_k} < N, \qquad
\sup_{k\in\Z} \frac{\lambda^*_k}{\lambda_k} < \infty.
\end{equation} 
The first two conditions are stronger than conditions \eqref{cond1}--\eqref{cond2}, while the third condition means that $\Lambda_k$ does not accumulate mass near zero ($\delta_0$ is no limit point) as $k\to\pm\infty$. 


\subsubsection{Hurdles along the way} 

It turns out that the generator can only be written down on a \emph{subspace} of $\wt\CM(\Omega_N^\infty \times \K)$. We have to show that we can obtain \emph{entrance laws} when starting the process in $\wt\CM(\Omega_N^\infty) \times \K$ such that in positive time it moves into this subspace. Recall the well-known fact that the classical measure-valued Fleming-Viot process for $t >0$ always has finitely many types, yet can start in general measures on the type space \cite{D93}. We have two operators to deal with, for the jumps and for the diffusion. 
\begin{itemize}
\item
When the sum over $k \in \Z_-$ of the rates is \emph{finite}, we are close to a Markov pure-jump process having finitely many jumps in finite time. However, when the sum is \emph{infinite} we get a path with countably many jumps at countably many places in a finite set, but in this case the generator contains positive and negative terms, and the sum is well defined only if the configuration around every site is eventually \emph{$\mfh$-a.e. mono-type} as $k \to -\infty$, so that only finitely many block-resampling events occur that yield a change, and the sum \emph{effectively} exists in every site of $\Omega_N^\infty$. Thus, to work out the martingale problem we have to address the issue whether the process enters specific subsets of configurations of the state space instantaneously and stays there forever, which would allow us to get existence of \emph{entrance laws} solving the martingale problem.
\item
Also the sum of the Fleming-Viot generators over the $L$-level sites decomposed into smaller and smaller sub-blocks has to be formed, and a limit object has to be formulated. The martingale problem is again given by an entrance law, but if the process starts in general elements of $\wt\CM(\Omega_N^\infty \times \K)$, then it is not necessarily atomic when projected on $\K$, or does not have densities in the space variable, or is not locally mono-type. This requires a \emph{new formulation} based on \emph{collision local times} in the continuum-space setting, because the volatility sits in \emph{rare space-time points} and only arises from migration of different types in a cluster of a given type. These issues will be addressed in Section~\ref{sec:trlim}, before we formulate the results that are needed to prove our main theorems. 
\end{itemize}


\subsubsection{Program towards the construction and characterisation of the continuum-space dynamics}

In Sections \ref{sec:trlim}--\ref{sec:wpmp} we introduce the \emph{continuum-space process} of type distribution as a \emph{measure-valued process} taking values in $E = \wt\CM(\Omega^\infty_N \times \K)$, and its \emph{dual} process on $\Omega^\infty_N$ as a \emph{measure-valued process} driven by a \emph{spatial coalescent} taking values in the \emph{$\Omega_N^\infty$-marked partition} of $[n] = \{1,\ldots, n\}$, $n \in \N$. We first pass to partitions of $\N$ via \emph{entrance laws} of the coalescent process on $\Omega_N^\infty$ starting with an individual \emph{in every site} in $\Omega_N^\infty$. For that purpose we construct a measure-valued dual obtained from the classical coalescent with finitely many individuals by starting countably many individuals in every site and taking a limit by passing to frequencies of partitions in balls as description, which specifies an element of $\wt\CM (\Omega_N^\infty)$. This is needed to be able to get a limit in a suitable topology. We introduce \emph{types from $\K$} in the state of the entrance law to get the \emph{dual process} at backward time $t$ based on the initial state $\wh X_0^\infty$ of the forward process $\wh X^\infty$, to be constructed with the purpose to end up with a dual state in $\wt\CM(\Omega_N^\infty \times \K)$. 

The duality not only helps us to get \emph{existence} and \emph{uniqueness} of the continuum-space process, it also allow us to decide when the continuum-space process is \emph{stochastic} versus \emph{deterministic} (i.e., a dynamical system), and to show that this dichotomy is equivalent to the dichotomy of \emph{locally mono-type states} versus \emph{locally multi-type states}.

We provide a construction and characterisation of the continuum-space stochastic process on $\Omega_N^\infty$ in three steps: 
\begin{itemize}
\item 
Section~\ref{sec:trlim}: We introduce a stochastic process $\wh X^\infty$ arising as the \emph{limit} of processes $\wh X^L$ on the countable geographic space $\Omega_N^L$ as $L\to\infty$. For the migration we work with \emph{continuum-space mass flow} induced by the transition kernel of the \emph{L\'evy-process $\zeta^\infty$ on $\Omega_N^\infty$}, which has been constructed in \cite{E89} and which arises as the projective limit of the laws of the L\'evy-processes $\zeta^L$ on $\Omega_N^L$. We add order-$L$ \emph{coarse-grained} block-resampling with reshuffling on the continuum-space $\Omega_N^\infty$, and also add single-site resampling in the sites of $\Omega_N^L$, both for Fleming-Viot and Cannings, leading in the limit $L\to\infty$ to a \emph{diffusion} term on $\Omega_N^\infty$. 
\item
Section~\ref{sec:dual}--\ref{s.third}: We establish existence of a limiting stochastic process $\wh X^\infty$ by characterising the limit process via the $\Omega_N^L$-duals and their limit as $L \to \infty$. This will allow us to show that the process $\wh X^\infty$ is deterministic, or is stochastic with for $t >0$ an atomic projection on $\K$ and the measures $x(\cdot \times \{a\})$, with $a$ an atom, having densities w.r.t.\ $\mathfrak{h}$. These properties define a subset of the state space containing the orbit of $\wh X^\infty$. However, eventually we want more, namely, a \emph{Markov} process and a \emph{Feller} process, which we address next.
\item
Section~\ref{sec:wpmp}: We characterise the process $\wh X^\infty$ as the solution to a \emph{well-posed martingale problem}. We use a strong duality to obtain \emph{uniqueness}. This duality allows us to also obtain the \emph{Markov property} by using results in \cite{DGP24}. The martingale problem has a different form than for discrete spaces, since the dynamics instantaneously moves into a subspace of the state space (as specified above) that is \emph{not closed} in the topology on that state space. This forces us to work with \emph{entrance laws}. Only on that subspace can we write down a \emph{classical operator acting on functions}, at least for the jump part on the subspace of measures from $\wt\CM(\Omega_N^\infty \times \K)$ viewed as state space. For the diffusion part the latter also plays a role and for the entrance law we also have to make use of a more subtle function of the measure state to obtain the \emph{volatility coefficient} in a given state in the subspace containing the orbit, which is being extracted from observations in a space-time window, resulting in a different form of the operator. This description of the volatility involves a measure on the \emph{space-time footprint} of the process (like for the local time of a classical diffusion process), and will be referred to as the \emph{different-type collision local time}. Otherwise we still have a \emph{second order operator}.
\end{itemize}

In \cite{gh2} we will provide a graphical construction via the so-called \emph{L\'evy-web}, which allows us to construct also the \emph{genealogies} of the population.


\section{Continuum-space dynamics: I. Construction of approximations}
\label{sec:trlim}

In this section we introduce the approximations $\wh X^L$, $L \in \N$, of $\wh X^\infty$, that are the backbone for the construction of the continuum-space process $\wh X^\infty$. The actual construction of $\wh X^\infty$ is carried out in Sections \ref{sec:dual}--\ref{sec:wpmp}. Section \ref{(0)} opens with some general observations. Section \ref{(i)} argues that the proper description is via \emph{locally finite Borel measures} on $\Omega^\infty_N \times \K$. Section \ref{(ii)} sets up a \emph{sequence of approximations} $X^L$ via refinement and embedding. Section \ref{(iii)} explains what hurdles towards the construction arise along the way. Section \ref{(iv)} outlines the program for the construction. Section \ref{(v)} lists the key challenges that come up in Sections \ref{sec:dual}--\ref{sec:wpmp}.


\subsection{Observations}
\label{(0)}

A continuum-space dynamics can often be justified as a limit after space-time scaling of a discrete-space dynamics, motivated by a desire to identify \emph{universality classes} (for example, Dawson-Watanabe superprocesses on general continuum geographic spaces arise from Galton-Watson branching random walks). The disadvantage of this approach to the construction of processes is that we only get special dynamics, while general dynamics are not captured. Indeed, the fixed points under the scaling represent the universality classes of the dynamics, not \emph{all} the possible dynamics. However, often it is possible to define the continuum-space dynamics analytically (for example, for superprocesses via the log-Laplace equation for the log-Laplace functional, which can be interpreted as a disguised duality relation with a deterministic evolution). In Section~\ref{s.results} we follow the former approach. In the present section we follow the latter approach, by choosing a \emph{natural approximation} with dynamics on $\Omega_N^L$, $L \in \N$, embedded in $\Omega_N^\infty$. 

To prove the approximation, we exploit \emph{duality} of the approximating processes and \emph{approximate the dual process on $\Omega_N^\infty$} by the dual process on $\Omega_N^L$ as $L \to \infty$. This duality will be our main tool to obtain a stochastic process on $\Omega_N^\infty$ matching the parameters in Remark~\ref{r.417}.

Both the block-resampling and the migration can be refined as $L$ increases, simply by adding smaller blocks, respectively, smaller jumps in these blocks, and taking the \emph{projective limit}, using ultrametricity. In some parameter regimes, we must worry about potentially diverging block-resampling rates as $L\to\infty$. Moreover, there is a more serious hurdle: both the Cannings single-site resampling and the Fleming-Viot single-site resampling are refined by \emph{shifting} the mechanism from the sites of $\Omega_N$ to the level $-L$ single sites, which are sites in $\Omega_N^L$. This is different from adding a refined evolution mechanism. If we add a mechanism one level lower and replace the mechanism at level $L$ by the lower one, then this changes the \emph{full} system and a proper \emph{scaling of rates} is needed, as explained in the next paragraph. 


\subsection{Description as measures}
\label{(i)}

For the process on the countable set $\Omega_N$, we described the type distributions in the different colonies by an element of $\CP(\K)^{\Omega_N}$. However, for the process on the uncountable set $\Omega_N^\infty$ the proper description is via a \emph{measure} on $\Omega_N^\infty \times \K$, namely, an element of $\wt\CM(\Omega_N^\infty \times \K)$, the set of locally finite Borel measures on $\Omega^\infty_N \times \K$ (note that $\Omega_N^\infty$ is locally compact and recall that $\K$ is compact). Observe we can reformulate the discrete description \emph{equivalently} like this as well. Our goal is to obtain, via approximation by processes on $\Omega_N^L$, a limiting process
\be{ag90}
\wh X^\infty = (\wh X^\infty_t)_{t \geq 0} \in D\big([0,\infty), \wt\CM(\Omega^\infty_N \times \K)\big).
\ee
Throughout the sequel we use the symbol $\,\wh{}\,\,$ to indicate that the process is \emph{measure}-valued.


\subsection{Sequence of approximations}
\label{(ii)}
 
We first define an approximating dynamics on the countable group $\Omega_N^L$, $L\in\N$. Let $X^L=(X^L_t)_{t \geq 0}$ evolve as a Cannings process on $\Omega_N^L$, which we embed in $\Omega_N^\infty$ (see Section~\ref{ss.disc}), with the following scaled parameters, evolution mechanisms and structural properties.


\paragraph*{$\blacktriangleright$ Migration and block-resampling.} 

For migration and block-resampling, use rates
\be{ag93b}
(c_k)_{k \geq -L}, \quad \quad (\Lambda_k)_{k > -L}, 
\ee
with the index set given by the unique disjoint decomposition of $\Omega_N^\infty$ into $(-L)$-balls. This set is isomorphic to $\Omega^L_N$ (recall \eqref{ag29}), and from \cite{GHKK14} we obtain the corresponding spatial Cannings process by reading the single components as densities w.r.t.\ $\mathfrak{h}^L$, the normed Haar measure on $\Omega^L_N$. 

We will prove that the $L$-truncated migration converges as $L\to \infty$ to a L\'evy-process and that the empirical distributions converge to the corresponding \emph{L\'evy mass flow} as deterministic component of the evolution. To do so, we must show that the $L$-truncated block-resampling converges to the one given by the infinite sequence $(\Lambda_k)_{k \in \Z}$ of block-resampling rates, including its contribution to the migration via the reshuffling. We must also show convergence when the migration and the block-resampling are combined, which is immediate via duality.


\paragraph*{$\blacktriangleright$ Single-site resampling.}

To include \emph{Fleming-Viot} resampling at single sites in $\Omega_N^L$, we observe that the number of terms in the Fleming-Viot generator after $L$-refinement steps increases as $N^L$ and the mass is rescaled by $N^{-L}$, which requires a scaling of the rate by $N^{L}$ to lead to something non-degenerate. More precisely, in the approximation we resample in such a way that if all the $L$-block masses are \emph{aggregated} into a single site, then the resulting scaled resampling leads to a rate-$d$ Fleming-Viot process. 

We replace the \emph{Cannings} resampling at the lowest level, i.e., $k = -L$, by a multiple of $\Lambda_{-\infty}$, where $\Lambda_{-\infty}$ is a \emph{reference resampling measure} that is needed in the construction. In summary, to compensate for the decrease of the volatility of the aggregated masses of small balls located in large balls compared to the volatility in small balls, we must rescale as we did for migration. Thus, we must choose 
\begin{equation}
\label{e988}
\begin{aligned}
&\text{single-site Fleming-Viot resampling at rate $d_L = N^L d$},\\
&\text{single-site Cannings resampling with resampling measure $\Lambda_L = N^L \Lambda_{-\infty}$.}
\end{aligned} 
\end{equation}
Note that this choice implies that the resampling aggregated in $(-L)$-balls is compatible with the migration in such balls, which allows for enough volatility to ensure that averages over the $(-L)$-balls in a $1$-ball evolve at a rate of order $1$ as $L \to \infty$. 

Note that if we were also to let $d \to \infty$ as $L \to \infty$, then we would get a \emph{degenerate process} in the point-recurrent regime, even though this process is a Markov Feller process and has a nice dual process with \emph{instantaneous coalescence} upon hitting. We do not pursue this observation further here, since we are interested in other universality classes. The instantaneous model would lead to \emph{historical} and \emph{genealogical structures}, where genealogical distances are twice the hitting times of the two corresponding backward L\'evy-processes and hence the historical process contains all genealogical information (see \cite{EF96} for a discussion of this setting). Letting $d \to \infty$ as $L \to \infty$ yields a truly stochastic process, but one that admits \emph{no} stochastic dynamics describable by SPDEs.


\paragraph*{$\blacktriangleright$ Level-$L$ parameters.} 

The \emph{level-$L$-dynamics} on the discrete space $\Omega_N^L$, arising from the migration and resampling according to the two previous items, is characterised by the parameters
\begin{equation}
\label{e1061}
(c_k)_{k \geq -L}, \quad (\Lambda_k)_{k > -L}, \quad d_L, \Lambda_L.
\end{equation}


\paragraph*{$\blacktriangleright$ Level-$L$ approximation.}

For $L\in\N$, let $\wh X^L = (\wh X^L_t)_{t \geq 0}$ be the $((c_k)_{k \geq -L},\,(\Lambda_k)_{k > -L},\,d_L,\,\Lambda_L)$-Cannings process on $\wt\CM(\Omega_N^L \times \K)$, i.e., the process corresponding to the modified mechanism on level-$L$ given by the parameters in \eqref{e1061}. Clearly, the measure $\wh X^L$ can be viewed as a measure on $\Omega_N^\infty \times \K$ via embedding in $\Omega_N^\infty$ equipped with the Haar measure $\mfh$, i.e., each site in $\Omega_N^L$ is identified with an $(-L)$-ball in $\Omega_N^\infty$ that gets the total mass of that site by taking the associated value in $\Omega_N^L$ distributed according to the normed Haar measure $\mathfrak{h}^L$ on that ball. We set
\be{fdh97}
\wh X_0^L = \wh X^\infty_0 |_{\CB(\Omega_N^L \times \K)}.
\ee
The goal is to show that there exists a process $\wh X^\infty$ such that
\be{e1768}
\CL\big[(\wh X^L_t)_{t \geq 0}\big] \mathop{\Longrightarrow}\limits_{L \to \infty} \CL\big[(\wh X^\infty_t)_{t \geq 0}\big],
\ee
where the convergence is in the sense of finite-dimensional distributions. Later we discuss \emph{convergence} in \emph{path space}, which is more subtle.


\paragraph*{Dual of approximation.} 

As dual of the process on $\Omega_N^L$, we consider a process with time horizon $t$, called $(Y_s^L)_{s \in [0,t]}$, of \emph{partitions marked with locations} in $\Omega_N^L$, embedded in $\Omega_N^\infty$ and equipped with \emph{weights}, in which the partition elements migrate and coalesce, i.e., merge and combine their weights. The state of this dual process is an atomic measure on $\Omega_N^\infty$ starting with \emph{finitely many dual individuals} at each site, which leads to a pure Markov jump process always well-defined. We need to prove convergence as $L \to \infty$ to a process on $\Omega_N^\infty$, with finitely many individuals. This is however not sufficient for our purposes. Indeed, afterwards we pass to \emph{countably} many individuals at each site in $\Omega_N^L$, for which we have to construct an entrance law. We do this by passing to initial per-site frequencies of partition elements at time $s \in (0,t]$, with $t$ the time horizon of the duality. After the process with index $L$ runs for a \emph{positive} time, it becomes a locally finite measure on $\Omega_N^L$, and by the embedding also on $\Omega_N^\infty$. In the limit as $L\to\infty$, we have to distinguish between the case where we have effectively coalescence (which turns out to lead to the stochastic regime) and the case where we do not (which turns out to lead to the deterministic regime). In the latter case we effectively only have migration and ``block-resampling with instantaneous mixing". 

In the stochastic regime, we consider a \emph{measure-valued} process $\wh Y^\infty$ on $\Omega_N^\infty$, corresponding to the initial state $\mathfrak{h}$, the normed Haar measure on $\Omega_N^\infty$, by passing to an \emph{entrance law} corresponding to starting with \emph{countably} many dual individuals in each site in $\Omega_N^L$ and letting $L \to \infty$. We have to distinguish between the two stochastic scenarios (b1) and (b2) of Theorem~\ref{T.char}. We want $Y^L = (Y_t^L)_{t \geq 0}$ to converge to a process $Y^\infty$ of which we can specify the evolution rules. To make sense of starting with countably many dual individuals at the sites of $\Omega_N^L$, we work with an \emph{entrance law} and pass to the measure-valued process $\wh Y^\infty$. This means that we pass to \emph{frequencies} of individuals in the partition elements, to get the \emph{entrance law} starting from countable many individuals per site. Namely, from $Y^L$ we pass to $\wh Y^\infty = (\wh Y_t^\infty)_{t \geq 0}$ as $L\to\infty$, with $\wh Y^\infty$ acting as the \emph{dual process} of $\wh X^\infty$ for each of the scenarios (b1) and (b2). Subsequently, we introduce \emph{types} of individuals and partition elements, with types from $\K$ at time $t$. The assignment is based on the forward process $\wh X_0$. This way we come to a measure on $\Omega_N^\infty \times \K$ as dual state of a process $Y^{\ast,\infty}$, getting a standard duality relation between $\wh X^\infty$ and $\wh Y^{\ast,\infty}$ as well. 

In the deterministic regime, where no coalescence occurs, we have populations of migrating L\'evy-processes in the entrance laws described by a measure-valued dual process, which we denote by $\wh \zeta^\infty$. This object describes populations of individuals at time $t$ by specifying the frequency of individuals arising from a certain initial value of the location. Due to the rapid local mixing in that regime, this forces the infinite populations into the state $\mfh \otimes \theta$.


\subsection{Hurdles towards the construction}
\label{(iii)}

To appreciate the hurdles that arise when trying to carry out the above program, consider the following simplification. Suppose that $d=0$, $\Lambda_{-\infty} \equiv 0$ and $\Lambda_k \equiv 0$ for $k < -L^\prime$ for some $L'<L$. Then
\be{ag94}
\CL [(X^L_t|_{\CB^{L^\prime}})_{t \geq 0}] = \CL[(X^{L^\prime}_t)_{t \geq 0}],
\ee
where $X^L_t|_{\CB^{L^\prime}}$ is the projection of $X^L_t$ onto $\CB^{L^\prime}$, the $\sigma$-algebra on $\Omega^\infty_N \times \K$ generated by the $(-L^\prime)$-balls
\be{ag92}
\big\{B_{-L^\prime }(\xi) \times K\colon\, \xi \in \Omega^\infty_N,\,K \in \CB(\K)\big\}.
\ee
The identity in \eqref{ag94} simply says that what happens downstairs does not depend on what is present upstairs. If we decrease $L^\prime$, then we simply add additional resampling in smaller balls, at rates that do not depend on what happens upstairs. Consequently, we have convergence as $L \to \infty$ to a limiting process arising as a \emph{projective limit}, provided the rate of change does \emph{not} diverge, so that the effective change remains bounded. The latter is an issue, since in the interesting parameter regimes leading to stochastic evolutions we do have diverging rates, as shown in the following. How do we get convergence of the measure-valued state on $\Omega_N^\infty \times \K$ in the general setting? If we want to describe the limiting process via an operator of a martingale problem, then we have to show that the total \emph{effective} rate remains finite as $L\to\infty$. Furthermore, dealing with measure-valued processes on $\Omega_N^\infty \times \K$ \emph{changes the shape of the duality}, as we will see below.

The lowest-level Fleming-Viot and Cannings single-site resampling at $k=-L^\prime$ \emph{perturbs} the structure in \eqref{ag94}, and therefore the existence of the limiting process is even more non-trivial than for the above simplification with block-resampling only. By compactness, there are limit points for convergence of the \emph{finite-dimensional distributions}. To obtain a limit we must show \emph{convergence} of the latter, for which \emph{duality} comes to our rescue. 


\subsection{Program for the construction}
\label{(iv)}

We know that for every $L$ the approximating process from \eqref{ag93b} and \eqref{e1061} has a finite dual process that is a block-coalescing random walk on $\Omega_N^L$, and that there is a duality function based on polynomials as described in \eqref{ag52}. We will analyse what happens with these dual processes as $L \to \infty$. After that we reformulate them as measure-valued processes. Furthermore, we pass to entrance laws of dual processes starting with countably many individuals, reformulating them as \emph{measure-valued} dual processes $\wh Y^{L,\ast}$ on $\Omega_N^L \times \K$, and show that we have limit processes $\wh Y^{\infty,\ast}$ on $\Omega_N^\infty \times \K$. This analysis will allow us to build a duality between $\wh X^\infty$ and $\wh Y^{\infty,\ast}$, and the structure of the dual will allow us to overcome the above hurdles when taking the limit $L \to \infty$. 

This approach works by embedding of $\Omega_N^L$ in $\Omega_N^\infty$ and by considering the states as densities with respect to the Haar measure in the decomposition of $\Omega_N^\infty$ into disjoint balls of hierarchical distance $\eee^{-L}$ (recall Section~\ref{ss.cont}). We normalise the mass of the measure-valued version by $N^{-L}$ to account for the increasing number of sites as $L$ increases.


\paragraph*{$\blacktriangleright$ Single-site resampling.}

We have to show that Cannings single-site resampling becomes Fleming-Viot single-site resampling observed at level $-L$ as $L \to \infty$, and that the dual of single-site $d$-Fleming-Viot resampling and $\Lambda_{-\infty}$-Cannings resampling properly scaled turn, in the limit as $L \to \infty$, into a system of coalescing L\'evy-processes such that coalescence occurs at \emph{rate $d^\ast$} (with $d^\ast = d + |\Lambda_{-\infty}|$) according to the \emph{collision local time} when two partition elements \emph{meet} at a site. Our strategy is to show that the entire dual process (with all its evolution mechanisms combined) converges in law as $L \to \infty$. 


\paragraph*{$\blacktriangleright$ Block-resampling.}

In the \emph{dual} system the block-resampling and the reshuffling have a simple structure, which is mirrored by the fact that two partition elements with paths that come within hierarchical distance $L^\prime$ of each other may or may not coalesce, with rates that depend on $L^\prime$. Consider $L$ and $L^\prime$ with $L^\prime >L$, and let $L^\prime \to \infty$ followed by $L\to\infty$. The block-resampling \emph{without} single-site Fleming-Viot resampling and Cannings resampling leads to a \emph{monotone} structure. Namely, as $L^\prime$ increases, \emph{additional} block-coalescence is activated for paths that come close to each other and therefore coalescence times get smaller. This monotonicity guarantees the existence of an $L \to \infty$ limit for the case $d=0$ and $\Lambda_{-\infty} \equiv 0$, i.e., $d^\ast=0$. 


\paragraph*{$\blacktriangleright$ Combination of single-site resampling and block-resampling.}

We need to describe the limit of the single-site Fleming-Viot coalescence and the single-site Cannings coalescence after it is \emph{combined} with migration and block-resampling with reshuffling. This task becomes simpler in the dual perspective, namely, altogether on $\Omega_N^L$ in the limit as $L \to \infty$ coalescence of paths can occur either by collision followed by coalescence at rate $d$ (with a deterministic limit evolution occurring both for $d=0$ and for $d \in (0,\infty)$, the latter with delayed coalescence, and with a stochastic limit evolution occurring for $d=\infty$), or by paths coming close and coalescing by block-coalescence. To achieve the latter, we have to show that the single-site Cannings coalescence at level $-\infty$ is either $0$, a finite positive number or $\infty$ (i.e., either no coalescence, or rate-$d^\ast$ delayed coalescence or \emph{instantaneous} coalescence), similarly as for single-site Fleming-Viot resampling. The limit for the dual can be delayed coalescence upon hitting the same site, or coalescence by coming close on infinitely many levels with an infinite accumulated rate, leading to coalescence via block-coalescence, or separation of paths leading to pure migration (and in $\wh X^\infty$ to mass flow). This happens in a way that depends on the recurrence versus transience property of the paths.


\paragraph*{$\blacktriangleright$ Completion.}

We finish the construction of the continuum-space process by showing that the $L$-\emph{truncated} block-coalescing random walk converges to the \emph{block-coalescing L\'evy-process with delayed coalescence upon meeting in a site}, and block-coalescence in $k$-balls of volume $N^k$ under the Haar measure with rate $\lambda_k$, prescribed for the corresponding forward process, for $k \in \Z$. 


\subsection{Upcoming challenges}
\label{(v)}

The arguments pointed out above proceed in two steps: 
\begin{itemize}
\item
Section~\ref{sec:dual}--\ref{s.third}: We formulate the \emph{duality} for finite $L$ as duality of \emph{two measure-valued} processes on $\Omega_N^L$, respectively, $\Omega_N^L \times \K$ (the latter gives the dual state with type marks). We show \emph{convergence} of expectations under the dual to a limit expression in terms of the resulting dual limit process $\wh Y^\infty$ on $\Omega_N^\infty$. This will give us an embedding of the $L$-processes in $\Omega_N^\infty$, for every time $t>0$, and convergence of the $L$-approximation $\wh X^L$ to a limit process $\wh X^\infty$. An important ingredient is the space-time duality, which allows us to get the finite-dimensional distributions of the process $\wh X^\infty$ in terms of the dual process $\wh Y^\infty$.
\item
Section~\ref{sec:wpmp}: The same approach as in \cite{DGP24} can be applied to define the process $\wh X^\infty$ via dual expectations, through a \emph{Markov} semigroup on $\CM (\Omega_N^\infty \times \K)$ associated with a Markov process that is in space-time duality. The measure-valued dual $\wh Y^\infty$ must therefore be equal to $\wh X^\infty$. In particular, $\wh X^\infty$ must be \emph{Markov} and, in a suitable sense, satisfy the generator criterion for duality, giving us also \emph{uniqueness} of the solution to the martingale problem. To achieve the latter we have to overcome the problem that we do \emph{not} have the standard form of the martingale problem. 
\end{itemize}
The main technical challenge is to get \emph{convergence} of the dual process $Y^L$ and its measure-valued version $\wh Y^L$ to continuum-space processes $Y^\infty$ and $\wh Y^\infty$ as $L\to\infty$, and to decide whether $\wh Y^\infty$ is the measure-valued process associated with a true coalescence process $Y^\infty$, or leads to a process $\wh \zeta^\infty$ representing independent L\'evy-processes. This in turn allows us to identify the relevant operator (more precisely, its adapted substitute in the limiting martingale problem), to show that with that operator $\wh X^\infty$ is indeed the solution to a \emph{well-posed martingale problem} and is \emph{in duality} with $\wh Y^\infty$ extended to $\wh Y^{\ast,\infty}$ by introducing types, and to identify the parameters, when this object is stochastic and when it is deterministic.


\section{Continuum-space dynamics: II. Convergence via duality}
\label{sec:dual}

In this section we construct the continuum-space process $\wh X^\infty$ as truncation limit of the discrete-space process $\wh X^L$ with the help of \emph{duality}. Section \ref{sss.pr} outlines the program for this construction. Section \ref{sss.tkp} states \emph{two key propositions} and a corollary that constitute the main route to the construction, and lists some relevant properties of $\wh X^\infty$ as a stochastic process. Section \ref{sss.dupro} introduces the measure-valued \emph{dual process} $\wh Y^{*,\infty}$ and shows how it is constructed. Section \ref{sss.dualfu} describes the \emph{duality function} in the duality relation. Section \ref{sss.sptcoal} describes the space-time coalescent needed for pathwise statements via the space-time duality~(\cite{FG94}). After some preparation on the dual process, Section \ref{sss.proofs} provides the proof of the key propositions, which uses three constituent lemmas proved in Sections \ref{sss.first}--\ref{sss.second} and Section \ref{s.third}.


\subsection{Program}
\label{sss.pr}

We reviewed the properties of the coalescent $Y$ on $\Omega_N$ in Section~\ref{ss.dualdisc}. The same properties apply to $Y^L$ on $\Omega^L_N$, which is controlled by the parameters of the continuum-space process truncated at level $-L$. By applying \eqref{ag52} to the process with these parameters, we get that $X^L$ is \emph{dual} to $Y^L$. We need a dual process with types in $\K$. To get this $Y^{\ast,L}$, we \emph{assign types} to partition elements at a backward time $t$ based on $\wh X_0^L$, respectively, $\wh X^\infty_0$ and this type is carried by all individuals in the partition elements for all times $s \in [0,t]$. The type is picked at random according to $\wh X_0$, given the location $\eta$, where the conditioning needs to be discussed (see below \eqref{ag37b}). 

We want a measure-valued dual process. We can view $Y^L$ as an atomic measure-valued process $\wh Y$ on $\Omega_N^\infty$ by assigning every dual individual weight $1$; this yields an atomic measure on $\Omega^L_N$, which we embed in $\Omega_N^\infty$. With the type assignment we have elements that induce, for each partition element, an element of $\wt\CM^{\rm atomic} (\Omega_N^\infty \times \K)$ (which is just a projection of states). This process of \emph{partitions, types} and \emph{locations} we can see as a measure on $\CI \times (\Omega_N^L \times \K)$, with $\CI$ a countable index set labelling the dual individuals. This measure we denote by
\begin{equation}
\label{e1642}
\text{$\mfY^{\ast,L},\mfY^{\ast,\infty}$, and projecting the latter on locations and types we get $\wh Y^{\ast,L},\wh Y^{\ast,\infty}$. }
\end{equation}
We can rewrite
\begin{equation}
\label{e1598}
\text{the duality between $Y^{\ast,L}$ and $X^L$ in terms of a duality between $\wh X^L$ and $\wh Y^{\ast,L}$,}
\end{equation}
which are two measure-valued processes on $\Omega_N^L \times \K$. 

To obtain a \emph{strong dual representation}, i.e., representing the state by the dual one of $\wh X^\infty$, we need to introduce 
\begin{equation}
\label{e1604}
\text{\emph{entrance} laws of $\mfY^{\ast,L}$ and projection $\wh Y^{\ast,L}$,}
\end{equation}
both starting from \emph{countably many individuals} at each site of $\Omega_N^L$. We state and establish the existence of these objects in Proposition~\ref{P.1430}, after explaining in detail all ingredients needed. We can do this by passing to the \emph{frequencies} of the individuals being in the partition elements at backward time $t$, arising from the individuals at a given site at backward time $0$. As $L\to\infty$, this gives a measure on $\Omega_N^\infty$, with the types in $\wt\CM(\Omega_N^L \times \K)$, $L \in \N \cup \{\infty\}$, both in $\Omega_N^L$ and in $\Omega_N^\infty$.

How can we make this approach rigorous? Using entrance laws, we can start with countably many individuals --- labeled by copies of $\N$ --- at every site of $\Omega_N^L$, by passing to frequencies at sites, as follows. Start with $n$ individuals at every site and take the limit $n \to \infty$, exploiting the consistency of the coalescent. When taking this limit, we must distinguish between the stochastic regime and the deterministic regime. This is where the \emph{no-dust condition} comes in, which is needed in the stochastic regime to get \emph{well-defined frequencies}, both for finite $L$ on $\Omega_N^L$ and in the limit as $L\to\infty$ on $\Omega_N^\infty$. In the deterministic regime we do have \emph{dust} formation as $L \to \infty$, and we need to proceed differently to define $\wh Y^\infty$. By adding all frequencies of partition elements with a given location and type, we get the frequencies that we need for the construction.  

Note at this point that the partition elements represent the \emph{descendants} of an \emph{ancestor} at time $0$ surviving at time $t$. In~\cite{gh2} we will endow the set of partition elements with an ultrametric structure to represent a \emph{genealogy} for the forward population. This will be carried out via a genealogy introduced in the \emph{dual} entrance law. In the stochastic regime, the processes evolve according to the entrance law in the same way as $Y^L$, except that the weights are now the \emph{frequencies} of the partitions elements, i.e., the weights of the individuals are scaled by their total number per site, which is asymptotically of order $n$. After that we can \emph{rewrite} the duality relation for $X^L$ and $Y^L$ in terms of the \emph{measure-valued} version $\wh Y^{\ast,L},\wh X^L$ of the \emph{dual} and the original process. This way we can write the \emph{strong} duality as in \eqref{e1411} below and, from the convergence of $\wh Y^{\ast,L}$ as $L \to \infty$, obtain convergence of $\wh Y^{\ast,L}$ to a process called $\wh Y^{\ast,\infty}$ and of $\wh X^L$ to a process called $\wh X^{\infty}$, in the sense of finite-dimensional distributions.

In the deterministic regime, the coalescent with its partition structure and its frequencies collapses to complete \emph{dust} as $L \to \infty$. To still get a limit measure in this situation, we consider the clusters of individuals traveling in time $t$ to specific locations in $\Omega_N^\infty$ of a specific type in $\K$ and put as weights their frequency. We have to worry about convergence, at least in the case of block-resampling, which is not visible in the limit $L\to\infty$ in the transient regime, and therefore we have to use the \emph{Meyer-Zheng topology} on the path space of $\wh X^\infty$. 


\subsection{Two key propositions and a corollary}
\label{sss.tkp}

The main results in this section are Propositions \ref{prop.1010}--\ref{P.DUAL}, together providing a \emph{unique stochastic process $\wh X^\infty$ as a limit of $\wh X^L$} as $L \to \infty$. With the help of Proposition~\ref{P.1430}, we get for $\wh X^\infty$ a measure-valued \emph{dual process} $\wh Y^\infty$ with values in $E^\prime \subseteq \wt\CM(\Omega_N^\infty \times \K)$, a process formulated and explained in detail in Section~\ref{sss.dupro}. Namely, we introduce an enriched (indicated by $\ast$) and \emph{measure-valued} (indicated by $\,\wh{}\,\,$) \emph{dual process} $\wh Y^{\ast,\infty}$ (a measure on $\Omega_N^\infty \times \K$) and write down a \emph{function-valued duality} (as we have for finite $L$) and based on a \emph{finite} partition-valued \emph{spatial Lambda-coalescent} process $Y^\infty$ on $\Omega_N^\infty$, a Markov pure-jump process, as well as a version $Y_t^{\ast,\infty}$ with \emph{types} associated to partition elements and their constituent individuals for states at time $t >0$. This is done starting with \eqref{Yhatdef} below. The process $Y^\infty$ describes, via its locations of partition elements in $\Omega_N^\infty$, the origin of the population at time $t$ in its spatial distribution, while the process $Y^{\ast,\infty}$ also describes the type of these subpopulations. Later this process will be extended to spatial coalescent with types, see Section~\ref{sss.dupro} for details.

\begin{proposition}{\bf [The process $\wh X^\infty$ on $\Omega_N^\infty$]}
\label{prop.1010}
$\mbox{}$
\begin{itemize}
\item[(a)]
Fix $L' \in \N$ and an admissible initial state of the dual process $Y^{L^\prime}$ on $\Omega_N^{L^\prime}$. Let $Y^{\ast,L}_{[0,t]} = (Y_s^{\ast,L})_{s \in [0,t]}$ be the dual process on $\Omega_N^L$ with refinement $L \geq L^\prime$ and time horizon $t$, starting from the embedding of the initial state on $\Omega_N^{L'}$ into $\Omega_N^L$. Then, for every $L^\prime$, $Y^{\ast,L}_{[0,t]}$ converges in law as $L \to \infty$ to a limit $Y_{[0,t]}^{\ast,\infty}$. 
\item[(b)] 
The weak limit of $\CL[(\wh X_t^L)_{t \geq0}]$ as $L \to \infty$ exists, in the sense of finite-dimensional distributions.
\item[(c)]
The limit is the law of a measure-valued process $\wh X^\infty$ on $\Omega_N^\infty \times \K$ that is in functional duality with $\wh Y^{\ast,\infty}$ (in the sense that \eqref{ag99} below holds).
\item[(d)] 
The process $\wh X^\infty$ has states $\wh X_t^\infty$, where
\begin{itemize}
\item
either $\wh X_t^\infty (\cdot \times \K)$ is the Haar measure on $\Omega_N^\infty$ and, for $t>0$, $\wh X^\infty (\Omega_N^\infty \times \cdot)$ is atomic and has a random density of the measure $\wh X^\infty_t (\cdot,\{a\})$ w.r.t.\ the Haar measure $\mfh$ on $\Omega_N^\infty$, which is either $0$ or $1$ $\mfh$-a.e.\ with $a \in \K$ an atom (i.e., is stochastic), 
\item
or $\wh X^\infty_t$ is identically equal to $\mfh \otimes \theta$ for strictly positive times (i.e., is deterministic).
\end{itemize}
\end{itemize}
\end{proposition}

\noindent
In Proposition~\ref{P.1430} we will see that the convergence can be extended from the enriched and recoded version of $Y^L$ to convergence on path space, for each time horizon, namely, to a measure-valued process $\wh Y^{\ast,L}_{[0,t]} = (\wh Y_s^{\ast,L})_{s \in [0,t]}$ with limit $\wh Y_{[0,t]}^{\ast,\infty}$.
 
\begin{proposition}{\bf [Duality]}
\label{P.DUAL}
For $\wh X^\infty$ with state space $E$ there exists a dual process $\wh Y^{\ast,\infty}$ with values in $E^\prime$, and a continuous duality function $\wh H$ on $E \times E^\prime$, such that $\{\wh H(\cdot,e^\prime),e^\prime \in E^\prime\}$ is measure- and convergence-determining and the following duality relation holds: 
\be{ag99}
\E_{\wh X^\infty_0}\big[\wh H(\wh X^\infty_t, \wh Y^{\ast,\infty}_0)\big] 
= \E_{\wh Y^{\ast,\infty}_0}\big[\wh H(\wh X^\infty_0, \wh Y^{\ast,\infty}_t)\big],
\qquad (\wh X^\infty_0, \wh Y^{\ast,\infty}_0) \in E \times E^\prime, \qquad t \geq 0.
\ee
The choice for $\wh H$ and $E^\prime$ is given in Section~\ref{sss.dualfu} below with $E=E^\prime=\wh E=\wt\CM (\Omega_N^\infty \times \K)$. Furthermore, the dual process is given by $(\wh Y^{\ast,\infty}_{[0,t]}(t))_{t>0}$, i.e., driven by the limit of $Y_t^{\ast,L}$ as $L \to \infty$. In particular,~\eqref{ag99} is the limiting duality relation as $L \to \infty$.
\end{proposition}

\noindent
The duality relation in \eqref{ag99} allows us to determine the \emph{finite-dimensional distributions} of a unique limiting evolution $\wh X^\infty = (\wh X_t^\infty)_{t \geq 0}$, rather than merely extract limit points through compactness of $\CL[\wh X^L]$. 

The duality relation is a classical duality between Markov processes, once we have shown that $\wh X^\infty$ is itself a Markov process. For this we check a number of conditions that are formulated in \cite{DGP24} (this approach is needed to obtain the \emph{Markov property} for $\wh X^\infty$, which we have constructed in Proposition~\ref{prop.1010} through its finite-dimensional distributions). From the latter reference we obtain the duality relation in \eqref{ag99} as a classical duality relation between two Markov processes. In Section~\ref{sec:wpmp} we strengthen Propositions \ref{prop.1010}--\ref{P.DUAL} by formulating a \emph{martingale problem} for $\wh X^\infty$ and establishing indeed the \emph{Markov property} for the limiting evolution. 

With the properties of $\wh H(\cdot,\cdot)$ and the structure of the dual process of the $\Omega_N^L$-approximations and their duality function in hand, we immediately obtain the following strengthening of the duality formulated in Corollary~\ref{C.sdual} below (see \cite{GSW16} for the case of the spatial Fleming-Viot process on $\R$). Namely, consider the \emph{particular version}, see the paragraph around~\eqref{e2009} of the \emph{entrance law} for the dual process $Y^\infty$ that we specify in Section~\ref{sss.dupro} and which amounts to starting the dual process with \emph{one} dual individual in each site of a \emph{countable skeleton and subgroup} of $\Omega_N^\infty$ here $\Omega_N^{\infty, \rm disc}=$ union over all embedded $\Omega_N^L$. This results in a strong duality for $\wh X^\infty$.

\begin{corollary}{\bf [Strong duality]}
\label{C.sdual}
For every time horizon $t$ the strong duality relation holds:
\be{e1411}
\CL \big[\wh X_t^\infty\big] = \CL \big[\wh Y_t^{\ast,\infty}\big], \qquad t > 0.
\ee
\end{corollary}

\noindent
This fact is important to establish properties of the states of $\wh X_t^\infty$, since we may just prove these properties for $\wh Y_t^{\infty,\ast}$ instead. It also implies that the dual process started in $\mathfrak{h}$ for positive time stays in a subspace of $\wt\CM (\Omega_N^\infty \times \K)$, namely, the atomic measures with atomic $\K$-marginals and locally mono-type states, so that the martingale problem yields \emph{uniqueness} of the original process, restricted to that subspace.

The strong duality can even be strengthened further as follows. 

\begin{remark}{\bf [Pathwise strong duality]}\label{r.1711}
{\rm In \cite{gh2} we will show that in the stochastic regime there is a \emph{pathwise} strong duality between the dual process driven by the $\K$-marked partitions of $\Omega_N^\infty$, i.e., $\wh Y^{\ast,\infty}$ and the measure-valued process $\wh X^\infty$ that specifies the state by the \emph{regions of a single particular type} (see \cite{GSW16} for the Fleming-Viot version on $\R$).} \hfill$\spadesuit$
\end{remark} 

\begin{corollary}{\bf [Dichotomy: dual entrance law and L\'evy-process characterisation agree]}
\label{cor.para}
Cases (a), (b1) and (b2) of Theorem \ref{T.char} on the properties of $\wh X^\infty$ correspond, respectively, to transient, Harris-recurrent, and point-recurrent behavior of the difference process of the pairs $(\zeta^\infty_1,\zeta^\infty_2)$ of two \emph{typical paths} of the \emph{entrance law} of the dual started with \emph{locally countably} many dual individuals, for example, one at each site in $\cup_{L \in \N} \wt \Omega_N^L$ with $\wt \Omega_N^\infty$ the version of $\Omega_N^L$ embedded in $\Omega_N^\infty$.
\end{corollary}

The proofs of Propositions \ref{prop.1010}--\ref{P.DUAL} are organised into three steps: in Section~\ref{sss.dupro} by defining the candidate for the \emph{dual process} on $\Omega_N^\infty$ (which leads to a further key proposition: Proposition~\ref{P.1430}), in Section~\ref{sss.dualfu} by defining the \emph{duality function} $H$, and (to obtain information about the finite-dimensional distributions) in Section~\ref{sss.sptcoal} by defining the \emph{space-time coalescent} with the space-time duality, respectively. Next, in Section~\ref{ss.prepp} we prepare for the proofs involving the more subtle object $\wh Y^{\ast,\infty}$ introduced in Section~\ref{sss.dupro}. In Section~\ref{sss.proofs} we carry out these proofs, based on three constituent lemmas that we prove in Sections~\ref{sss.first}--\ref{sss.second} and Section~\ref{s.third}, respectively. 


\subsection{Dual process: coalescent on $\Omega_N^\infty$}
\label{sss.dupro}

This section is based as starting point on the duality of the process on $\Omega_N$, reviewed in detail in Section~\ref{ss.dualdisc}. This we have to adapt to the geographic space $\Omega_n^\infty$. In this section, via Proposition~\ref{P.1430}, we define a \emph{Markov pure-jump process $Y^\infty$} called \emph{spatial coalescent} on $\Omega_N^\infty$ and introduce a version $Y^\ast$ with \emph{types} being the analogue on $\Omega_N^\infty$ of the dual process for the process on $\Omega_N^L$. We have to lift this process to a \emph{measure-valued} one to obtain a genuine dual process on $\Omega_N^\infty$. We obtain a \emph{measure-valued} version $\wh Y^\ast$ as an \emph{entrance law} from locally countably many dual individuals. Concretely, we place one individual at each site of $\Omega_N^{\infty, \rm disc}$, the dense countable skeleton of $\Omega_N^\infty$ defined as the union of all embedded copies of $\Omega_N^L$, $L \in \N$. This yields, as in Corollary \ref{C.sdual}, also a \emph{strong dual} for the process $\wh X^\infty$ that we obtained above as limit $L \to \infty$ of the truncated processes $\wh X^L, L \in \N$.


\paragraph*{Markov pure-jump process: finite spatial coalescents.}

We begin by considering the Markov pure-jump process, which is a \emph{finite} spatial coalescent on $\Omega_N^\infty$ denoted by
\be{Yhatdef}
Y^\infty = (Y_t^\infty)_{t \geq 0}.
\ee 
One mechanism is \emph{delayed} or \emph{instantaneous coalescence} of a pair of partition elements proportional to the \emph{collision local time}, respectively, the infinite volatility limit $(d \to \infty)$ when both \emph{hit} the same site. The other mechanism is \emph{effective} migration of tuples of partition elements. The latter process uses tuples of L\'evy-processes 
\be{}
\bar\zeta^L = (\bar\zeta^L_t)_{t \geq 0} \quad \text{on} \quad \Omega_N^L,
\qquad \qquad 
\bar\zeta^\infty = (\bar\zeta^\infty_t)_{t \geq 0} \quad \text{on} \quad \Omega_N^\infty,
\ee 
with migration rates given by \eqref{effmig}. As a partition-valued process with partition elements marked by locations in $\Omega_N^L$ or $\Omega_N^\infty$, these tuples are well defined for all times when they start from \emph{finitely many} dual individuals (each with their own partition element) at every site. 

More formally, consider the coalescent with finitely many initial partition elements, say $m$ with locations $\underline{\bar\zeta}^L = (\bar\zeta_1^L,\ldots,\bar\zeta_m^L)$. The partition elements perform independent \emph{$(c_k)_{k \in \Z}$-L\'evy-processes} on $\Omega^\infty_N$ that \emph{coalesce upon collision} at a rate proportional to the \emph{collision local time} of the pairs of paths of $\underline{\bar\zeta}^\infty$. In Corollary~\ref{cor.1939}, for point-recurrent effective migration, we will establish that this coalescence occurs at rate
\be{dinf}
d^\infty = \left\{\begin{array}{ll}
0, &d^\ast = 0,\\
d^\ast, & d^\ast \in (0,\infty),\\
\infty, & \text{if $d_L \to \infty$ as $L \to \infty$},
\end{array}
\right. 
\ee
where either $d^\ast = d+|\Lambda_{-\infty}|$ or $d^\ast = 0$ (the last line corresponds to \emph{instantaneous coalescence}). In addition, \emph{block-coalescence} and \emph{block-reshuffling} occur according to $\Lambda_k$ when a pair is within hierarchical distance $k$ (recall \eqref{hierdist}). Thus, we have a well-defined finite coalescent $Y^\infty$ on $\Omega_N^\infty$, as long as $d^\ast$ and $|\Lambda_k|$ are not all zero (otherwise we have independent L\'evy-processes).

\begin{remark}{\bf [Classification of coalescence]}
\label{rem.threecases}
{\rm It turns out that, under our basic conditions \eqref{cond1}--\eqref{cond3} on the migration rates and the resampling rates, we encounter the following cases: 
\begin{itemize}
\item[(a)] 
\emph{No coalescence}: This occurs when the L\'evy-processes of $Y^\infty$ on $\Omega_N^\infty$ are \emph{independent}.
\item[(b1)] \emph{Delayed/Instantaneous coalescence:} This occurs when the L\'evy-processes are \emph{point-recurrent} and the collision rate $d$ due to Fleming-Viot coalescence and/or rate $|\Lambda_{-\infty}|$ due to Cannings coalescence is positive (\emph{delayed}) or $\infty$ (\emph{instantaneous}) (for single-site coalescence it will follow from \eqref{e988} that \eqref{dinf} is the limiting coalescence rate). 
\item[(b2)] 
\emph{Block-coalescence}: This occurs when the L\'evy-processes can merge via block-coalescence rather than via delayed or instantaneous coalescence, and move also by reshuffling.
\end{itemize}
Which of these cases occurs depends on the migration rates and the resampling rates in our model. In some instances we in fact have only convergence in the \emph{Meyer-Zheng path topology}, rather than in the \emph{Skorohod path topology}, because the process may make} block-resampling jumps at random times. \hfill$\spadesuit$
\end{remark}


\paragraph*{Enrichment with types: from $Y$ to $Y^\ast$.}

To get a duality relation with $\wh X^\infty$, we can \emph{enrich} the dual $Y^\infty$ to a \emph{$\K$-marked} infinite particle system by giving, for given time horizon $t>0$, each partition element a \emph{type} from $\K$ (also attached to each constituent individual), drawn once for all times $s$ in $[0,t]$, with the state of $\wh X_0^\infty$ as the tuple of locations (see below) of the coalescent run for time $t$. This gives a \emph{path} with values in $\Omega_N^\infty \times \K$:
\be{ag37b}
Y^{\ast,\infty} = (Y^{\ast,\infty}_s)_{s \in [0,t]}.
\ee
The dual individuals carry the type from $\K$ at time $t$ of their evolution ($t$ being the time horizon). The type is sampled from $\wh X^\infty_0$ conditioned on the location at time $t$. Here, problems arise because the Haar measure is infinite and the space is continuous (conditioning $\wh X_0^\infty$ on the locations is more subtle). To resolve the former, we have to work with regular conditional distributions w.r.t.\ the Haar measure on $\Omega_N^\infty$ and ensure that the approximation arguments go through. For the latter, the conditioning is well-defined on $\Omega_N^L$ but not on $\Omega_N^\infty$: this is well-defined only if we restrict the measures to subgroups (to sites, $n$-balls) in $\Omega_N^{L,1}, \Omega_N^{\infty,n}$ with finite Haar measure. Therefore on $\Omega_N^\infty$ we consider restrictions to that set and take the limit $n \to -\infty$ of the
\begin{equation}
\label{e1697}
\text{ normalized } \wh X_0 \mid_{\CB_n(0)} \text{ conditioned on $\eta$.}
\end{equation} 
Here we must use the $L$-approximation. We first state the final result we obtain.


\paragraph*{Lifting to measure-valued processes: main results on entrance law on $\Omega_N^\infty$.}

The above finite system $Y^\infty$ can be extended to a countable system of the type of a so-called \emph{infinite particle system}, both on $\Omega_N^L$ and $\Omega_N^\infty$, starting from a \emph{locally finite} configuration in the \emph{Ligett-Spitzer space} (see~\cite{GLW05} for this concept in a dual of a spatial population genetics model) associated with the effective migration/reshuffling system, which generates a Markov process with values in the Liggett-Spitzer space.

After the preparation on $\wh Y^{\ast,L}$ and $\wh Y^{\ast,\infty}$ realising \emph{both} processes with locations in $\Omega_N^\infty$ as explained above, we can extend the convergence statement of part (a) of Proposition \ref{prop.1010} to obtain the \emph{entrance laws} from countably infinite initial states on $\Omega_N^\infty$ for $Y^{\ast,\infty}$.

\begin{proposition}{\bf [Entrance laws of the coalescent from the Haar measure]}
Let $T$ be the time horizon for the duality relation.
\label{P.1430}
\begin{itemize}
\item[(a)] Consider the dual process $\wh Y^{\ast,L}$ of $\wh X^L$. The dual process on $\Omega_N^\infty$ is given as follows: 
\begin{equation}
\label{e1357}
\wh Y_T^{\ast,\infty} = \lim_{L \to \infty} \wh Y^{\ast,L}_T, \qquad T \geq 0, \qquad 
\text{ in law on $\wt\CM(\Omega_N^\infty \times \K)$}.
\end{equation}
The left-hand side is identified as the entrance law of the process $Y^\infty$ starting with one individual at each site of $\Omega_N^{\infty, \rm disc}$ (see the last four paragraphs in Section~\ref{sss.dualfu}); it is characterised in (b) and (c) below.
\item[(b)] 
In the stochastic regime (recurrent $\zeta$), the process $(\wh Y_s^\infty)_{s \in [0,T]}$ is the unique process with the property that every \emph{finite sampled subcoalescent} of $\wh Y^\infty$ evolves according to the rates and the transitions of $Y^\infty$ on $\Omega_N^\infty$ as defined in \eqref{Yhatdef}--\eqref{dinf}. As specified in \eqref{Yhatdef} and below, every partition element is assigned the type at backward time $T$ according to $\wh X_0^\infty$.
\item[(c)] 
In the transient regime (transient $\zeta$), the finite systems of \emph{typical} $m$-tuples under the entrance law $\wh Y^\infty$ converge to $\underline{\zeta}^\infty$ (the vector of effective L\'evy-processes of order $m$). Recall the discussion below \eqref{ag37b} and note that types are assigned in the same way as in (b).
\end{itemize}
\end{proposition}

\noindent
Thus, for positive times we have a measure-valued dual process with states that are \emph{atomic measures} (which are \emph{locally finite}) and with weights that arise from frequencies of $\Omega_N^\infty$-marked partition elements starting in $\Omega_N^{\infty,\rm disc}$, which we can label by $\N$ and a certain specified location at time $0$ (in backward time). We can now define our candidate for a dual process of $\wh X^\infty$ with state at time $t$ (= time-horizon) as in \eqref{e1431}, denoted by $\wh Y_t^{\infty,\ast} \in \wt\CM (\Omega_N^\infty \times \K)$, first for a \emph{finite} number of dual individuals and afterwards for entrance laws, as shown in Sections~\ref{sss.dualfu}--\ref{sss.sptcoal}. 


\subsection{Duality function}
\label{sss.dualfu}

Consider a coalescent with a \emph{finite} number of individuals, say $n \in \N$. Introduce \emph{partitions} $\pi$ of $\{1,\ldots,n\}$, where we view $\pi$ as a map from $[n] = \{1,\ldots,n\}$ into itself such that the \emph{partition elements} consist of all $k$ that are mapped by $\pi$ into the same element. For a function $f$ of $n$ variables we can consider the function $f(u_{\pi(1)},\ldots,u_{\pi(n)})$, which we associate with the pair $(f,\pi)$. As above, we assign a \emph{location} $\xi$ in $\Omega_N^\infty$ and a \emph{type} $u$ in $\K$ to each \emph{partition element}, the former giving a $n$-tuple $\underline{\xi}$ of positions in $\Omega_N^\infty$, the latter giving a tuple of types assigned based on the initial state $\wh X_0$ of the \emph{forward} evolution. This gives us a \emph{function-valued} dual process driven by a \emph{spatial coalescent}, which is what we used on $\Omega_N^L$ together with a duality function $H$ (see \eqref{ag50}--\eqref{ag51} for $H$ and \eqref{ag52} for duality on $\Omega_N$, properly adapted to $\Omega_N^L$). 

To treat $\Omega_N^\infty$, we pass from the \emph{function-valued} dual $Y^\infty$ for $\wh X^\infty$ started in finite populations to the \emph{measure-valued} dual $\wh Y^{\ast,\infty}$ given as \emph{entrance law} from a locally countably infinite state by using a duality function $\wh H$. Returning to the \emph{function-valued duality} for processes on $\Omega_N^L$, which is based on $H$, we are now viewing the state of $X^L$ as a measure $\wh X^L$ on $\Omega_N^L$, getting a duality w.r.t.\ a function of $\wh H$ of $\wh X^L$ with a measure-valued entrance law of $Y^L$, called $\wh Y^L$, and measures indicated by the superscript $\,\,\wh{}\,\,$. This approach puts the duality in the framework of \emph{measure-valued} processes. More importantly, it allows us later to rigorously formulate the idea of starting the dual entrance law for $Y^L$ with \emph{countably} many dual individuals \emph{at every site} in $\Omega_N^L$, converging to a coalescent starting in \emph{every site of a skeleton in $\Omega_N^\infty$}, namely, in the union $\Omega_N^{\infty, \rm disc}$ of all $\Omega_N^L$, after embedding in $\Omega_N^\infty$ and after passing to frequencies corresponding to \emph{starting the dual process} $\wh Y^\infty$ in the \emph{Haar measure} $\mathfrak{h}$. 

We transfer the function-valued set up of the duality on $\Omega_N^L$ w.r.t.\ $H$ to $\Omega_N^\infty$, inserting an integration of $H(x,\cdot)$ w.r.t.\ to the Haar measure of the parameters in the function-valued dual, thereby obtaining from $H$ a new function $\wh H$ for the lifted duality of $\wh X^\infty$ to the dual $\wh Y^{\ast,\infty}$ w.r.t.\ $\wh H$. In formulas this reads as follows.

Let 
\begin{equation}
\label{e970}
\wh E = \wt\CM(\Omega^\infty_N \times \K), \quad E^\prime = \bigcup_{n \in \N} \left(\Omega_N^\infty \right)^n \times E^\prime_n \times S_n,
\end{equation}
with
\begin{equation}
\label{a1453a}
\begin{aligned}
E^\prime_n &= \text{ algebra of functions $\Psi^{(n)}$ on $(\Omega_N^\infty \times \K)^n$ generated by products }  
\Psi^{(n)}_{\Omega_N^\infty}\,\Psi^{(n)}_\K \\ 
&\qquad \text{with } \Psi^{(n)}_{\Omega_N^\infty} \in C_c ((\Omega_N^\infty)^n,\R), \Psi^{(n)}_\K \in C_b(\K^n,\R), S_n = \text{ set of partitions of } \{1,\ldots,n\},
\end{aligned}
\end{equation}
where the functions $\Psi^{(n)}_{\Omega_N^\infty}$ are in the \emph{domain} of the generator $\mathfrak{A}$ of the L\'evy-process (recall \eqref{Lpr}--\eqref{e801}). Thus, the state space of the function-valued dual $E^\prime$ is a triple of a \emph{location vector}, a \emph{partition} and a \emph{test function}.

The \emph{duality function} for the measure-valued $\wh X^L$ and $\wh Y^{\ast,L}$ on the new state space of measures is as follows. First recalling~\eqref{ag50}-~\eqref{ag52} $H\colon\,\wh E \times E^\prime \to \R$ and define
\be{e1461}
H(x,(\Psi,(\pi, \underline{\xi}))) = H_\Psi(x,(\pi,\underline{\xi})), 
\qquad x \in \wh E,\,(\Psi,\pi) \in E^\prime_n,\, n\in\N,
\ee
with
\be{e1014}
\begin{aligned}
H_{\Psi^{(n)}_{\Omega_\infty^{(n)}},\Psi^{(n)}_\K}(x,(\pi,\underline{\xi})) = \int_{(\Omega^\infty_N \times \K)^n} x^{\otimes n} (\dd \underline{\eta},\dd \underline{u})\,
\Psi^{(n)}_{\underline{\xi},\pi}(\underline{\eta},\underline{u})
\end{aligned}
\ee
and
\be{e1474}
\Psi^{(n)}_{\underline{\xi},\pi}(\underline{\eta},\underline{u}) = 1_{\xi_1}(\eta_1) \times \cdots \times 1_{\xi_n}(\eta_n)\, \Psi_{\Omega_N^\infty}(\eta_{\pi (1)},\ldots,\eta_{\pi(n)}) 
\Psi_\K(u_{\pi(1)},\ldots,u_{\pi(n)}),
\ee
where $\xi_i = L(\pi^{-1}(i))$ with $L^{(n)} = (L(1),\ldots,L(n)) \in (\Omega_N^\infty)^n$ denoting the locations of the respective partition elements. 

The duality functions $H(\cdot,e^\prime)$ evaluate \emph{finite samples} from $\wh X^\infty$ based on and parametrised by the dual state $e^\prime$, which in the duality formula is given via $Y_t^{\ast,\infty}$. Observe that $\{H(\cdot,e^\prime),e^\prime \in E^\prime \}$ is law-determining on $\wh E$. We must reinterpret the right-hand side of \eqref{e1014} in terms of a \emph{measure-valued} dual $\wh Y^{\ast,\infty}$ \emph{starting in} $\mfh$, meaning that when projected onto $\Omega_N^\infty$ it starts in the \emph{Haar measure} $\mathfrak{h}$, and at time $t$ carries a type that is assigned exactly as for $\wh Y^{\ast,\infty}_t$, as explained in Proposition~\ref{P.1430} (recall \eqref{e1357}). But what does starting in $\mathfrak{h}$ mean? What does ``reinterpreting" mean?

We have to focus on spatial averages of the components and view these averages as measures acting on the geographic space. We consider sums over the spatial components and normalise with the volume by using $\mfh$. This is effectively an integration w.r.t.\ $\mfh^L$ of the dual initial locations. In fact, we can think of $\Psi_{\Omega_N^\infty}$ as the indicator of a ball. To these objects we can apply the $H$-duality for $\wh X^L$ and $Y^L$, which then becomes a duality between $\wh X^L$ and $\wh Y^L$, based on \emph{finite coalescents}. 

To continue, recall that for the process on $\Omega_N^L$ there exists an \emph{entrance law} of the coalescent, where initially we have \emph{countably} many individuals at each site, which we can realise embedded in $\Omega_N^\infty$ by taking $\Omega_N^{\infty, \rm disc}$ with \emph{one particle per site} and taking all particles in an $L$-ball as the population in a site in $\Omega_N^L$. We saw earlier that this leads to a process that, for positive times, has locally \emph{countably} many partition elements, and we can record at each site the \emph{frequency} of a particular partition element in the state and obtain a \emph{locally finite measure} on $\CI \times \Omega_N^L$, where $\CI$ labels the partition elements and the projection on $\Omega_N^\infty$ is $\mfh^L$. 

We want to pass to the $L \to \infty$ limit. We know that $\mfh^L \Rightarrow \mfh$ as $L \to \infty$, with $\mfh$ the normalised Haar measure on $\Omega_N^\infty$. For this $L$-approximation, we can establish in Lemma \ref{l.1562} the existence of an $L\to \infty$ limit, which we call $\wh Y^{\ast, \infty}$ and which is the \emph{dual} entrance law starting in $\mfh$. The \emph{duality function} $\wh H$ is induced by the right-hand side of \eqref{e1014}. We have a pairing of $\wh Y^{\ast, \infty}$ and $\wh X^{\ast,\infty}$ and a \emph{duality function $\wh H$} on $\wh E \times \wh E$ derived from a function-valued duality function $H$ on $\wh E \times E^\prime$. Namely, $\wh H(e_1,\cdot)$ is an empirical average of $H(e_1,\cdot)$ in the location component in the $L \to \infty$ limit. Indeed, the $\wh H$-\emph{duality between $\wh X^\infty$ and $\wh Y^{\ast,\infty}$}, with $e_1,e_2 \in \wt\CM(\Omega_N^\infty \times \K)$ and $e^\prime \in E^\prime$, given by
\begin{equation}
\label{e1775}
\wh H (e_1,e_2) = \int_{E^\prime} H(e_1,e^\prime) e^{\otimes \mid \xi \mid}_2 
\big(\dd\underline{\xi},\K^{\mid\underline{\xi}\mid}\big)
\end{equation}
and, using~\eqref{e1461}, we can write this as
\begin{equation}
\label{e1775a}
\int_{(\Omega_N^\infty \times K)^n} e_1^{\otimes n} \left(\dd\underline{\eta},\dd\underline{u}\right) 
\int_{(\Omega_N^\infty)^n} \mfh^{\otimes n} \left(\dd\underline{\xi}\right) \Psi^n_{(\underline{\xi},\pi)} 
(\underline{\eta},\underline{u}).
\end{equation}


\subsection{Space-time coalescent} 
\label{sss.sptcoal}

Next, we express the \emph{finite-dimensional distributions} of the process $X^L$ in terms of the \emph{space-time dual process} $\wh Y^{\ast,\infty}_{[0,t]}$ with the help of the \emph{time-space duality} for $X^L$, using the Markov property that we have for each $L \in \N$ (see \cite{FG96} for a description and use of this technique). To that end, we \emph{enrich} the spatial type-marked coalescent by giving the partition elements (the basic individuals) the label \emph{frozen} or \emph{active}, together with a \emph{deterministic time label} indicating at what backward times the frozen elements become active. Thus, we now label initial dual individuals by their position and their \emph{activation time}: $(\xi_i,t_i)$, $i=1,\ldots,m$ and $t_m=t >t_{m-1} > \ldots > t_1>0$. 

We obtain the \emph{space-time coalescent} from the spatial coalescent $(Y_s)_{s \in [0,t]}$, with $t$ the horizon of the duality relation, by freezing the evolution of the individuals up to the time given by the first time-coordinate of their initial position, but otherwise following the rules of the spatial coalescent with block-resampling. This allows us to define a \emph{space-time dual} $\wh Y^\infty_{[0,t]}$ for the time interval $[0,t]$, where dual individuals start the evolution backwards in time, from time $t=t_m$, are activated at times $t_m > t_{m-1}>t_{m-2}> \ldots >t_1>0$, and are finally all evaluated at time $0$ by assigning the type. The dual individuals are referred to as frozen until times $t_m,\ldots,t_1$ viewed backwards. As before, the types are introduced based on the state $\wh X_0^\infty$, which gives $Y^{\ast,\infty}_{[0,t]}$.

With the help of the space-time coalescent, we can determine the \emph{finite-dimensional distributions} of $X^L$ at times $t_1,\ldots,t_m$. To do this for $\wh X^\infty$, we must use the measure-valued formulation and $\wh Y^\infty$. Take a measure-valued process on $\Omega_N^\infty \times \K$. For our process we consider \emph{polynomials} evaluated at a sequence of time points $0 < t_1 < \ldots < t_m=t$, such that the factors are monomials of the form that we use in the duality, i.e., we form a product of functions of the form of the duality function $H(\cdot,e^\prime_i)$ for some $e^\prime_i \in E^\prime$ at times $0$ with duality horizons $t_i$, $i=1,\ldots,m$, namely,
\be{ag99b}
F(X) = \prodl^m_{i=1} \langle \wh X^{\otimes n_i}_{t_i}, f_i \rangle,\quad 
n = \suml^m_{i=1} n_i, \quad n_i, m \in \N,
\ee
where $f_i \in \CB((\Omega^\infty_N \times \K)^{n_i}, \R)$ is given by
\be{e1518}
f_i(\underline{\eta}_i,\underline{u}_i) = \int_{(\Omega_N^\infty)^{n_i}} 
\mathfrak{h}^{\otimes n_i}(\dd \underline{\xi}_i) \wt f_{\underline{\xi}_i}(\underline{\eta}_i,\underline{u}_i),
\qquad \underline{\eta}_i \in (\Omega_N^\infty)^{n_i},\,\underline{u}_i \in \K^{n_i},
\ee
with
\begin{align}
\label{e1013}
\wt f_{\underline{\xi}_i}(\underline{\eta}_i,\underline{u}_i) = \left[\prod^{n_i}_{k=1} 1_{\xi_k}(\eta_k) \right] 
\psi_{\Omega_N^\infty}^{(n)}(\eta_1,\ldots,\eta_{n_i}) \; \psi_\K^{(n)}(u_1,\ldots,u_{n_i}), \\
\psi_{\Omega_N^\infty}^{(n)}(\eta_1,\ldots,\eta_{n_i}) \in C_c([\Omega_N^\infty]^{n_i},\R), \psi_\K^{(n)} \in C(\K^{n_i},\R). \nonumber
\end{align}
These \emph{monomials} generate an \emph{algebra of polynomials}, which is measure determining and convergence determining for $\wt\CM(\Omega_N^\infty \times \K)$ and the expectations, the \emph{moment measures} are law determining. Therefore a space-time duality based on these functions is very useful.

Take the space-time coalescent starting with individuals in the space-time points
\be{e1530}
\underbrace{(t_1,\eta_1^1), \ldots, (t_1,\eta_{n_1}^1)},\,
\ldots, \underbrace{(t_m,\eta_1^m), \ldots, (t_m,\eta_{n_m}^m)},
\ee
with $\underline{\eta}_i = (\eta_1^i,\ldots,\eta_{n_i}^i)$, $i=1,\ldots,m$. The initial value $\wh X_0^\infty$ determines uniquely the law of the process $\wh X^\infty_{[0,t]}$, by determining the joint marginal distributions at times $0<t_1< \ldots < t_m$ of the \emph{measure-valued state.}

Following \cite{GSW16}, from the Markov property of both $\wh X^L$, $L \in \N$, and $\wh Y^L$, $L \in \N_0 \cup \{\infty\}$, we have, as a consequence of the duality relation~\eqref{ag99}, the \emph{space-time} duality relation~\eqref{e1547} for $\wh X^\infty$.
\begin{corollary}{\bf [space-time duality for $\wh x^\infty$]}
\label{cor.1815}
\begin{equation}
\label{e1547}
\E_{\wh X_0^\infty}\big[F(\wh X^\infty_{[0,t]})\big] 
= \E_{\wh Y_0^\infty}\big[F(\wh Y^{\ast,\infty}_{[0,t]}(t))\big], \qquad t \geq 0,
\end{equation}
where $\wh Y^\infty_{[0,t]}$ is the space-time dual (just defined) on the time interval $[0,t]$, with $t$ the time horizon. 
\end{corollary}

Note that the endpoints of $\wh X^\infty_{[0,t]}, \wh Y^\infty_{[0,t]}$ are $\wh X^\infty_t$, respectively, $\wh Y^\infty_t$, with the latter being the endpoint $\wh Y_{[0,t)}^\infty(t)$ of $\wh Y^\infty_{[0,t)}$ at backward time $t$. This allows us to study the law of the process $\wh X^\infty$ via $\wh Y^{\ast,\infty}$.


\subsection{Preparations for the proof of two propositions}
\label{ss.prepp}

Before starting the proofs of Propositions~\ref{prop.1010} and~\ref{P.DUAL}, we have to construct formally and understand better the dual object $\wh Y^{\ast,\infty}$. We first present the general idea and afterwards proceed in two steps.


\paragraph*{Limiting entrance law measure-valued dual $\wh Y^{\ast,\infty}$ as truncation limit: general ideas.} 

For the strong duality and for the measure-valued process we need to pass to an \emph{entrance law} of $Y^{\ast,L}$ or $Y^{\ast,\infty}$, starting with \emph{countably} many dual individuals at some sites on $\Omega_N^L$, respectively on $\Omega_N^\infty$, where we need one individual at all sites in $\Omega_N^{\infty,\rm disc}$. We need to let $L \to \infty$ to obtain $Y^{\ast,\infty}$ from $Y_N^{\ast,L}$ for bounded configurations, and to extend this limit to the \emph{entrance laws} $\wh Y^{\ast,L}$ to obtain $\wh Y^{\ast,\infty}$, by starting from locally \emph{countably infinitely many} individuals with one individual at each site of $\Omega_N^{\infty,\rm disc}$. We need to construct these objects for $Y^L$ first. 

We define the entrance law of $Y^L$ for $L \in \N$ with countable occupancies per site by restricting both the occupation numbers (to at most $n$ per site) and the set of occupied sites (to a finite set $V_n \subseteq \Omega^L_N$ with $V_n \uparrow \Omega_N^L$), and considering the resulting \emph{locally finite} coalescent. For fixed $n$ the process is well defined as a coalescent with finitely many occupied sites. Letting $n \to \infty$, we get a Markov process on $\Omega^L$. This means that we have a well-defined infinite particle system on $\Omega_N^L$. Namely, we move to \emph{countable} local occupation per site. For that we consider first the $n$-truncated occupation per site and then let the truncation level $n \to \infty$, to get, via the projectivity of the mechanism, a limit process $Y^L$ as \emph{entrance law} from a locally countable occupation in the \emph{Liggett-Spitzer} space, as we argue in the paragraph below.

In a second step, to obtain the corresponding objects on $\Omega_N^\infty$, we use the \emph{$L$-truncation approximation} to get $Y^\infty$. Therefore we have to match the description of $Y^\infty$ as a Markov jump process and entrance law (as indicated below \eqref{Yhatdef}) with the limit process of $Y^L$. The way to achieve this is by returning to an \emph{approximation} with the embedded $Y^L$ on $\Omega_N^L$ for $Y^\infty$ on $\Omega_N^\infty$. This uses the measure-valued formulation of the entrance laws of $Y^{\ast,L}$ or $Y^{\ast,\infty}$, as formulated in Proposition \ref{P.1430}. 

However, we need not only $Y^{\ast,\infty}$ but also $\wh Y^{\ast,\infty}$ (the measure-valued \emph{entrance law} starting from \emph{countably} many dual individuals at each of the \emph{finitely} many initial sites). We get such a process $(Y_t^{\infty,\ast})_{t \geq 0}$ viewed in backward time as the entrance law from time $0$, and consisting of locally finitely many partition elements for positive backward time, if we discard a small frequency of dual individuals initially and further backward. The processes $Y^{\ast,\infty}$ and $\wh Y^{\ast,\infty}$ are more subtle to define \emph{directly} as stochastic processes on $\Omega_N^\infty$, but we will see that they give rise to a partition of $(\Omega_N^\infty,\mfh)$ as a measure space, and we work with the \emph{$L$-truncation limit} for the rigorous construction. 

The situation is as follows.
\begin{itemize}
\item
For finitely many individuals, if the block-resampling acts for $k$-values that are negative enough, then no tuples of multiple sites are in a single $k$-ball (because they are separate except at a hitting time), where on $\Omega_N^\infty$ coalescence may occur or infinitely many block-resampling events may occur. Therefore convergence as $L \to \infty$ occurs as long as the configuration has no dual partition elements in a single ball smaller than some $-k$, $k \in \N$, for times $s \in [0,t]$. We have to analyse the situation more in detail at the hitting times to get $Y^L \to Y^\infty$ as $L \to \infty$ in the \emph{Skorohod path topology}. 
\item
To obtain the \emph{strong duality} on $\Omega_N^\infty$, we need entrance laws with initially \emph{countably} many individuals locally by putting one individual at \emph{each site} in the skeleton $\Omega_N^{\infty, \rm disc}$ of $\Omega_N^\infty$, the union of all $\Omega_N^L$ (see \eqref{e2009}) 
\item
Without block-resampling and in the \emph{transient case}, we introduced measures $Y_t^\ast$ by considering the weight of paths that reach a certain location, and we assigned a specific \emph{type} from $\K$ at time $t$. With block-resampling we can do the same, but now the convergence only occurs in the \emph{Meyer-Zheng path topology}. The limit is again $\wh \zeta$, which has to be used in the duality relation since, by the dust formation, the \emph{law of large numbers} is again effectively giving deterministic states.
\end{itemize}

In summary, by using for each $L \in \N$ the embedding of $\Omega_N^L$ into $\Omega_N^\infty$ for positive times, we get in the \emph{recurrent case}
\begin{equation}
\label{e1540}
\begin{aligned}
&\text{an $\K$-atomic measure-valued process $(\wh Y_t^{L,\ast})_{t \geq 0}$ on $\CI \times \Omega_N^\infty \times \K$,}\\&\text{respectively, } (\wh Y_t^{L,\ast})_{t \geq 0} \text{ on } \Omega_N^\infty \times \K,
\end{aligned}
\end{equation}
while in the \emph{transient case} we have as dual measure a population of \emph{backward L\'evy-processes}, and in the forward language as \emph{fixed state} $\mfh \otimes \theta$ as soon as $t >0$, the \emph{mean measure} of the state at $t=0$. This setting will be developed in the proof of Proposition \ref{P.1430} in Section \ref{sss.proofs}. Thus, we need to add $\ast$ and $\wedge$ to $\Omega_N^\infty$. We now turn to the details.


\paragraph*{Existence of the entrance laws from countably infinite states on countable spaces.}
Take first the case of no block-resampling. We easily get the necessary objects on $\Omega_N$ and $\Omega_N^L$ by arguing as follows. For the Moran model the same situation was considered in the \emph{non-spatial} setting by Kingman. The \emph{entrance law} $\wh Y^{\ast}$ of the dual process for $\Lambda$-Cannings resampling is well understood in a \emph{non}-spatial Cannings process under the no-dust condition on $\Lambda$, which holds for both the Kingman-coalescent and the $\Lambda$-coalescent under our conditions on $\Lambda$. This can be verified by following the arguments introduced for \emph{non}-spatial settings, extended also to $\Omega_N^\infty$ for $\Lambda=\Lambda_{-\infty}$ or on $\Omega_N$ for $\Lambda=\Lambda_0$ (see \cite{P99,S00}). 

Moreover, the argument extends to the \emph{spatial} Cannings process on $\Omega_N^L$ under the \emph{no-dust condition} via the techniques developed in \cite{GLW05}, where this is carried out for the Moran model and it is shown that the coalescent instantly jumps to a \emph{locally finite} configuration of partition elements, and afterwards evolves in the set of atomic measures on $\Omega_N^L \times \K$, in the appropriate \emph{Liggett-Spitzer space} $\Omega_N^\infty$ (generalising this concept from $\Omega_N^L$), or any countable Abelian group of countably many locations in $\Omega_N^\infty$. The latter easily extends to the spatial case with single-site Cannings resampling under the \emph{no-dust condition} on $\Lambda$, by considering a truncated configuration on $\Omega_N^L$ that is coming down from infinity.

In the sequel we have to extend the above properties to the setting with \emph{block-resampling}, both for finite $L$, and later for $L \to \infty$ to cover $\Omega_N^\infty$. The block-resampling corresponds in the dual to a block-coalescence, which leads to a pure jump process that on $\Omega_N^L$ is easily constructed as an infinite particle system in a standard way (see~\cite{GLW05}). 


\paragraph*{Entrance law for the dual starting from the Haar measure.} 

We want next to formally define a candidate for the dual process of $\wh X^\infty$, called $\wh Y$ as function of a process $\mfY$ derived from evolving partitions. This is parallel to what we did for $\wh X^L$ by introducing an \emph{entrance law} from locally \emph{countable} states. Since $\Omega_N^\infty$ is \emph{uncountable}, this takes some care and therefore we want to start the coalescent on $\Omega_N^\infty$ with one individual at each site of some countable skeleton of $\Omega_N^\infty$ or, more conveniently, a particular and canonical \emph{skeleton $\Omega_N^{\infty, \rm disc}$} arising from the union of all $\Omega_N^L$ over $L \in \N$, each embedded in $\Omega_N^\infty$ (defined more formally in \eqref{e2009}) with one individual in \emph{every site}. 

Start the dual in \emph{finitely many balls}, but with \emph{countably many individuals} each, one in each site of $\Omega_N^{\infty,\rm disc}$ in the ball. Consider the frequencies of the partition element among the population at these sites at backward time $t$. The partition elements are located at specific sites at backward time $t$, where a type from $\K$ is sampled according to $\wh X_{0,\eta}^\infty(\cdot)$ from $\K$, with $\eta$ the \emph{vector of locations} and $\wh X_{0,\eta}$ the \emph{conditional} measure $\wh X_0^\infty$ given that the location of the partition element is $\eta$. We have finitely or countably many partition elements left at time $t$, which have a location at backward time $s \in [0,t]$ and a type at time $t$ and by construction at all times $s$ in $[0,t]$. A partition element at site $\eta$ gets as weight the \emph{frequency} of the partition element in the original population at backward time $0$ at sites $\eta^\prime$ in all of the space $\Omega_N^\infty$ and of a specific type. For each partition element this induces a measure in $\CM_{\rm fin}(\Omega_N^\infty \times \K)$ for time $s \in [0,t]$. This measure we use to define a measure on $\CI \times (\Omega_N^\infty \times \K)$, with $\CI$ the initial countable dual population.

To obtain a strong duality on $\Omega_N^\infty$, we need to start from \emph{countably many} individuals from \emph{every} bounded set in $\Omega_N^\infty$, for example, by placing i.i.d.\ sites according to the normalised Haar measure $\mfh$. Or even better, we use the $L$-approximation, which is more suitable to overcome the technical problems by starting with countably many individuals on a \emph{dense} countable subset of $\Omega_N^\infty$. A particularly good choice is $\Omega_N^{\infty, \rm disc}$, the union over $L$ of all in $\Omega_N^\infty$ embedded $\Omega_N^L$ (see \eqref{e2009}--\eqref{eq:2009}). 

On $\Omega_N^\infty$, for each partition element at backward time $t$ we consider the \emph{frequency of individuals} located initially (in backward time) at a \emph{particular site in $\Omega_N^\infty$}, as well as the type of this \emph{partition element}, to get, for every partition element and every $t>0$, an element of $\wt\CM(\Omega_N^\infty \times \K)$, as described above. We have countably many partition elements, which we can label with an index set $\CI$ (countable), labeling all initial individuals at time $0$. The joint object can be viewed as measure on $\CI \times (\Omega_N^\infty \times \K)$.

In the latter case we have a countable population, say $\CI$, marked with locations and types from $\K$, getting, for all times $t \geq 0$, $(\Omega_N^\infty \times \K)$-marked partitions of $\CI$. The partition elements we label by individuals of $\CI$, for example, the leading individual of the partition element, based on a labeling of the individuals (respectively, partition elements) with $\N$ in some per site way. The states of the resulting stochastic process give rise to a process of frequencies of a partition element in $Y^\infty$ in a bounded subset of $\Omega_N^\infty$ with types in a certain subset of $\K$, defining a measure 
\begin{equation}
\label{e1959}
\text{$\wh \mfY_t^\infty$ on $\CI \times \{\Omega_N^\infty \times \K \}$, }
\end{equation}
via frequencies of initial individuals in the partition element labeled $\iota \in \CI$, with location in $A \subseteq \Omega_N^\infty$ and type in $B \subseteq \K$ for $A,B$ measurable subsets of $\Omega_N^\infty$, respectively, $\K$. We can think of $\wh \mfY_t^\infty$ as specifying the frequency of \emph{descendants} at time $t$ of a specific \emph{ancestor} at time $0$ in the forward time-$t$ population. This description of $Y^\infty_t$ via $\wh Y_t^\infty$ leads not only to a dual process for $\wh X^\infty$, but to a dual process for the \emph{genealogy process} of our model, to be studied in~\cite{gh2}. 

Looking at types and locations of a \emph{fixed} partition element only, we obtain a measure on $\Omega_N^\infty \times \K$ (rather than on $\CI \times (\Omega_N^\infty \times \K)$). We can do this, for example, by taking a finite ball in $\Omega_N^\infty$ centered at $0$ and taking the largest one, the second largest one, etc., randomizing the ranks in ties, resulting in a collection of measures. This measure is also associated with a \emph{number of partition elements} after we have truncated the configuration, deleting an `$\varepsilon$-small' number of initial individuals elements at each site, getting what we call 
\begin{equation}
\label{e1745}
\text{a collection of locally finite \emph{occupation measures} $\wh \mfN^\varepsilon_k, k \in \N$, associated with $\wh \mfY^\infty$.}
\end{equation}
For example, we can choose to remove individuals jumping by effective migration before a suitable time $\delta(\varepsilon)$. This is the choice we work with in the present paper.

In the frequency picture we start in the Haar measure $\mfh$ on $\Omega_N^\infty$, as explained above, and in the \emph{stochastic} regime ($\zeta^\infty$ recurrent) we obtain $\wh \mfY^\infty$ with the characteristic property that a \emph{finite subcoalescent sampled} from the \emph{measure} $\wh \mfY^{\infty,\ast}_t$ is $Y^{\infty, \ast}$, as defined in \eqref{Yhatdef}--\eqref{dinf}. This \emph{sample} is defined by restricting to initial states started in countably many individuals only in some ball (or, more generally, bounded subsets of $\Omega_N^{\infty,\rm disc} \times \K$), leading for positive backward time to a finite measure that is \emph{normalised} as a \emph{probability} measure. From the state restricted to that set \emph{draw an $m$-sample}, which we follow in backward times up to the time-horizon $t$ and which evolves as the well defined Markov process $Y^\infty$ starting with these $m$ individuals. In the deterministic regime (i.e., \emph{transient} effective migration), $\wh \mfY^{\ast,\infty}$ has to be replaced by $\wh \zeta$ and $Y^\infty$ by $\underline{\zeta}^\infty$. In the stochastic regime (i.e., \emph{recurrent} effective migration), on $\Omega_N^\infty$ we also need that, at positive times and after deleting a frequency-$\varepsilon$ per site population, we get \emph{locally finite configurations} for the \emph{occupation measure} $\mfN^\varepsilon$ in the \emph{Liggett-Spitzer space}, as for the $\Omega_N^L$-version, respectively, the $\Omega_N^\infty$-version of the migration, by verifying uniformity in $L$. 

How can we get rid of the labels $\CI$ in the description? We want to identify $\CI$ with locations in the limit $L \to \infty$, and for that purpose we use $\Omega_N^{\infty, \rm disc}$ as support for the countably infinite starting entrance law already for \emph{finite $L$}. We choose $\CI$ for $L \in \N$ and $\Omega_N^{\infty, \rm disc}$ for $L=\infty$ as follows. For $\Omega_N$ we take as labels the countably many individuals in a site $\eta$ as the intersection of the $1$-ball in $\Omega_N^{\infty, \rm disc}$ around $\eta$, the $(-L+1)$ ball in $\Omega_N^L$, and $\CI=\Omega_N^{\infty, \rm disc}$ for $L=\infty$, and each site carries $1$ individual.

To treat the latter situation rigorously, we formally proceed for the case $L=\infty$ as follows. Think of starting the process in the Haar measure $\mathfrak{h}$, taking the $L \to \infty$ limit of the corresponding entrance laws on $\Omega_N^L$ with countably many dual individuals per site of $\Omega_N^L$. For that embed $\Omega_N^L$ into $\Omega_N^{\infty, \rm disc} \subseteq \Omega_N^\infty$, by putting one individual at each site of $\Omega_N^{\infty,\rm disc}$, and taking as equivalence classes all sites in an $L$ ball. Then, passing for finite $L$ to the current \emph{frequencies} of the partition elements at time $t$ in the population at time $0$ per $L$-ball, we get a \emph{measure}-valued process as entrance law for each decomposition of $\Omega_N^\infty$ in $L$-balls. The partition frequencies properly scaled give us, in the limit $L \to \infty$ of these decompositions, a random \emph{measure} on $\Omega_N^\infty$, with a.s.\ space marginal $\mathfrak{h}$. Adding types as prescribed above and projecting the measure on $\CI \times (\Omega_N^\infty \times \K)$ onto $\Omega_N^\infty \times \K$, we get from $\wh \mfY^\infty$ the \emph{$\wt\CM(\Omega_N^\infty \times \K)$-valued} dual process
\begin{equation}
\label{e1431}
\text{$\wh Y^{\infty,\ast}_T \in \wt\CM(\Omega_N^\infty \times \K)$, and associated \emph{occupation measures} $\wh \mfN^\varepsilon_k (T)$, $k\in \N$, for $\varepsilon > 0$,}
\end{equation}
where $\,\,\wh{}\,\,$ indicates measure. 


\subsection{Proof of the key propositions}
\label{sss.proofs}

\begin{proof}{[Proof of Propositions \ref{prop.1010}--\ref{P.DUAL}]}
Our task is to use \emph{tightness} of the stochastic processes $\wh X^L$ and afterwards show convergence of the dual process $\wh Y^{\ast,L}_{[0,t]}$ to a limiting Markov process $\wh Y^{\ast,\infty}_{[0,t]}$. This implies, in particular, \emph{convergence} of $\wh X^L$ to a process $\wh X^\infty$ \emph{in finite-dimensional distributions}:
\begin{equation}
\label{e1794}
\lim_{L\to\infty} \E\big[\wh H(\wh X_0^L,\wh Y_t^{\ast,L})\big] 
= \E\big[\wh H(\wh X_0^\infty,\wh Y_{[0,t]}^{\ast,\infty}(t))\big], \qquad t \geq 0.
\end{equation}
Furthermore, in the stochastic regime, for $t>0$ we need to prove that the states of $\wt Y^{\ast,\infty}$ are \emph{atomic} in the type-variable, and that the measure $\wh X_t(\cdot \times \{a\})$, with an atom $a \in \K$, has a \emph{density} taking values $0$ or $1$, while in the deterministic regime we need to prove that the state for positive time is $\mfh \otimes \theta$. 

We first show convergence of the finite dual coalescence processes, and of the expressions in the left-hand side of the duality relation, to the corresponding objects for the claimed limit of $\wh X^L$ and $\wh Y^L$ as $L \to \infty$. Afterwards we prove the other properties, by verifying them for the space-time measure-valued dual $\wh Y^{\ast,L}_{[0,t]} (t),L \in \N \cup \{\infty\}$. 

The proof comes in 4 Steps.


\subsubsection{Step $1$: Convergence via three constituent lemmas} 

We first prove Proposition \ref{prop.1010}(a). This is done by reducing the claim to three lemmas that will be proved only later, in Sections~\ref{sss.first}--\ref{sss.second} and Section~\ref{s.third}.

We must establish convergence of the \emph{finite} block-coalescing random walks $Y^L$ to the corresponding block-coalescing L\'evy-processes $Y^\infty$. To that end we consider \emph{finite coalescents} with $m$ basic dual individuals in positions such that $Y_0^L$ has a position in $\Omega_N^L$, which is embedded in $\Omega_N^\infty$ and which converges to $Y_0^\infty$ with positions in $\Omega_N^\infty$ as $L \to \infty$. This gives us a sequence of finite coalescents $(Y_t^{L,m})_{t \geq 0}$, and we want to obtain convergence as $L \to \infty$ for every $m \in \N$ to a Markov coalescent process $(Y_t^{\infty,m})_{t \geq 0}$ that has the Feller property. Note that for every $L \in \N$ the Feller property holds for every $m \in \N$.

We proceed in three layers. We begin by formulating the claim for $Y^\infty$ without coalescence, i.e., for an $m$-tuple of moving individuals $\underline{\zeta}^L$, like L\'evy-processes following the mechanism of independent \emph{migration} respectively, \emph{reshuffling} at the prescribed rates, including a single random walk. Afterwards we include the mechanism of \emph{coalescence}. Finally, we proceed with $\wh Y^\infty$, the entrance law from countably infinite initial states and the corresponding measure-valued version of the spatial coalescent.

\begin{lemma}{\bf [Convergence of the $m$-tuple of motions under migration and reshuffling to the effective motion]}
\label{RWL}
Given $\underline{\zeta}^L$ with $\underline{\zeta}^L(0) \to \underline{\zeta}^\infty(0)=(\zeta^\infty_1,\ldots,\zeta^\infty_m) \in (\Omega_N^\infty)^m$ for $m \in \N$:
\be{ag1685}
\lim_{L\to\infty} \CL[\underline{\zeta}^L] = \CL[\underline{\zeta}^\infty] \text{ in } 
D([0,\infty),(\Omega_N^\infty)^m),
\ee
where the one-dimensional marginal ($m=1$) is the \emph{effective L\'evy-process}.
\end{lemma}

\noindent
The proof of Proposition~\ref{prop.1010}(a) is based on the path convergence of the $m$-tuple "random walk". This convergence is enough when the migration is \emph{transient} and the L\'evy-processes have local time $0$ and no ``single-site" coalescence, so that we get systems of "random walks", i.e., $(Y_t^\infty)_{t \geq 0}$ is $(\zeta^\infty_t)_{t \geq 0}$, on \emph{which only the block-resampling acts}. In the transient regime we clearly have convergence at each level $k \in \Z$ of the block-resampling, for $k \geq -\ell$ and for every $\ell \in \N$, ignoring for the moment the lower block-coalescence. For $\ell \to \infty$ there is a limit path for all times, where the path $\zeta^\infty$ has hierarchical distance larger than $-\ell$ for $\ell \in \N_0$ in finite time intervals provided we start without multiple occupation of a site. In the case of multiple occupation of locations we have to use that, by our assumption on the rates for small jumps (see~\eqref{condextra}), this configuration is instantaneously dissolved. Hence, in the transience regime we are done with the convergence of $Y^L$ to $Y^\infty$ by Lemma~\ref{RWL}: $\underline{\zeta}^\infty$ started in different sites does \emph{not} have collisions in finite time, because the difference process is not point-recurrent (see~\ref{appA:oneLP}). In particular, we have the Feller property.

If the migration is \emph{Harris-recurrent} but \emph{not} point-recurrent, then the convergence follows from the path convergence of the random walk, since coalescence in the limit is then possible only via block-coalescence, and the previous argument still works because the paths starting from \emph{distinct} sites can get arbitrarily close only as $t \to \infty$. If the migration is \emph{point-recurrent}, then we still have to prove convergence to the \emph{delayed} coalescence based on the rate $d^\ast$ times the local time, to get convergence of $Y^L$ to $Y^\infty$. This goes beyond what we know from Lemma~\ref{RWL}. 
 
We condense the above discussion into the following claim on the convergence of $Y^L$ to $Y^\infty$, to be proved later.

\begin{lemma}{\bf [Convergence of limit $Y^L$ to $Y^\infty$]}
\label{RWLprec}
For parameters that satisfy the basic assumptions the following holds for finite coalescents:\\
(a) If the migration is point-recurrent, then the \emph{coalescent} $Y^L$ converges in path space as $L \to \infty$ to the true spatial coalescent $Y^\infty$ on $\Omega_N^\infty$, as defined in \eqref{Yhatdef} and the sequel. \\
(b) If the migration is Harris-recurrent but not point-recurrent, then the claim from (a) holds for $Y^L,Y^\infty$.\\
(c) If the migration is transient and there is no block-resampling, then $Y^L$ converges in path space to the $(\Omega_N^\infty)^m$-L\'evy-process $\zeta^\infty$ defined in \eqref{Yhatdef}. In the case where $(\lambda_k)_{k \in \Z}$ is non-zero the limit process has an additional block-coalescence with reshuffling that is a random functional of a version of the effective $\underline{\zeta}^\infty$. 
\end{lemma}

Next, we turn to the statements concerning the \emph{entrance} law of $Y^\infty$ and $\wh Y^\infty$, i.e., the limit $L \to \infty$ of the entrance laws of $Y^L$ starting with \emph{countably many} dual individuals at every site in $\Omega_N^{\infty,\rm disc}$ (see \eqref{e2009}), and the convergence of the \emph{measure-valued dual}, i.e., $\wh Y^{\ast,L}$ and its limit $\wh Y^{\ast,\infty}$ (recall \eqref{e1431}), which exhibits \emph{new} effects. 

\begin{definition}{\bf [Dichotomy for $\wh Y^\infty$: dust free versus complete dust]}
\label{def.dustfree}
{\rm We say that the limiting process $\wh Y^\infty$ satisfies the \emph{no-dust condition} for the corresponding entrance law when $\wh \mfY^\infty$ decomposes for $t>0$ into partition elements with a \emph{well-defined positive frequency} at sites, which sum up to $1$ at each site at each time. Otherwise the no-dust condition does \emph{not} hold. In case all frequencies are zero, we speak of \emph{complete dust formation}.}\hfill$\spadesuit$ 
\end{definition}

\noindent
The above terminology arises from the fact that the genealogies of the population at time $t$ can be described by specifying subfamilies of individuals having a common ancestor at time $s$ back, for $s \in [0,t)$, such that the frequencies of the subfamily in the population located in some bounded set in $\Omega_N^\infty$ are well defined. This object gives rise to a measure on a space $U \times (\Omega_N^\infty \times \K)$, which will be the dual state in the \emph{genealogy} process for our model, we will study in \cite{gh2}, and exists only in the case of recurrent $\zeta^\infty$. In our setup this means that a typical \emph{$m$-tuple of paths} of $\wh \mfY^{\ast,\infty}_t$ may be given by $Y^{\ast,\infty}$ with coalescence or may be given by a tuple $\zeta^\infty_t$, depending on whether we are in the recurrent or in the transient case. We have to show the \emph{alternative}: the no-dust condition holds or the complete dust formation holds for the spatial coalescent, depending on which of the two cases for the L\'evy-process $\zeta^\infty$ occurs. 

For point-recurrent migration, we can generalise what is known for the stable Fleming-Viot process to the general Cannings and Fleming-Viot model. For Harris-recurrent migration, in addition we have to show that the limit $Y^\infty$ still allows for an entrance law \emph{$\wh Y^\infty$ starting locally with countably many individuals with one particle everywhere in $\Omega_N^{\infty,\rm disc}$, satisfying as $L \to \infty$ for positive time the no-dust condition}. For transient migration without the block-resampling, we immediately have \emph{dust formation}, but we need to show that the block-resampling cannot prevent complete dust formation. The situation is therefore quite different in the \emph{recurrent} case and in the \emph{transient case}. 

More precisely, taking the measure-valued process and letting $L \to \infty$, we find the following features. In the \emph{recurrent case}, in the limit as $L \to \infty$ we have partition elements with a well-defined frequency in both for the coalescent starting with one individual in every site of $\Omega_N^{\infty,\rm disc}$, giving rise to $\wh Y^\infty$ started in $\mfh$ as the rigorous formulation of the process $Y^\infty$ starting with one individual in every site of $\Omega_N^{\infty,\rm disc}$ and hence \emph{countably} many individuals locally. In the \emph{transient case}, in the limit as $L \to \infty$ we get \emph{complete} \emph{dust}, hence we have no decomposition into partition elements with positive frequency per site, and the dual population looks like the one of $\wh \zeta^\infty$, i.e., \emph{without coalescence} but \emph{with effective migration}. In that case, $\wh Y^\infty$ is defined to represent the components of individuals starting in specific measurable subsets of $\Omega_N^\infty$. 

Finally, we need the property that the limiting process $Y^\infty$, respectively, the process under the entrance law with types $Y^{\infty,\ast}$, started with countably many individuals at each site in $\Omega_N^{\infty,\rm disc}$ (i.e., $\wh Y^{\infty,\ast}$), converges as a measure-valued process with types in $\K$ (i.e., as element of $\wt\CM(\Omega_N^\infty \times \K)$) to $\wh Y^{\infty,\ast}$, which will eventually follow from Lemma \ref{RWLprec} and Lemma \ref{l.1562}. However, according to the regime of \emph{dust} formation, respectively, of the \emph{no-dust property}, we get two different situations: 
\begin{itemize}
\item
In the \emph{recurrent regime} we have no-dust, the state of $\wh X_t^\infty$ is represented by $\wh Y_t^{\infty,\ast}$ for $t>0$, and the latter is a \emph{true} coalescent process. It is useful to note that the partition elements correspond to a decomposition of $\wh X_t^\infty$ representing the descendants at time $t$ of the different founding fathers that start to evolve at time $0$ from the initial state $\wh X_0^{\ast,\infty}$. Each subfamily corresponds to different partition elements at time $t>0$ (this property will be worked out in \cite{gh2}). It turns out that this structure of the genealogy induces a decomposition of the geographic space $\Omega_N^\infty$ into disjoint components. If $\wh X_0$ projected onto $\K$ is a \emph{diffusive} measure, then each partition element has a specific type (they are only all different if $\theta$ is a diffusive measure) and a specific frequency in bounded spatial sets. In the recurrent regime, the components in the decomposition give us the regions of $\Omega_N^\infty$ where we have a specific type from $\K$ at all the sites, a fact that we make rigorous below. 
\item
In the \emph{transient regime}, the set of partition elements forms \emph{complete dust} and so each of them has \emph{no} strictly positive weight in the case of an entrance law with initial states with space-marginal $\mfh$. Therefore we introduce the measure $\wh Y^{\ast, \infty}$, which has at backward time $s \in [0,t]$ weights of types and the geographic origin at backward time $t$ (forward time $0$) of the populations of independently evolving L\'evy-processes $\zeta^{\ast,\infty}$ with types given via $\wh X_0^{\ast,\infty}$ and with positions given via the Haar measure $\mfh$. This results in 
\begin{equation}
\label{e1980}
\text{measures $\wh \zeta^{\ast,\infty}$ on $\{\Omega_N^\infty \times \K\}$,}
\end{equation}
replacing the dual object $\wh Y^{\ast,\infty}$ \emph{with} coalescence that we had in the \emph{recurrent regime}. The measure corresponds to the frequency of individuals in the dual population in backward time $0$ in sites $\eta \in \Omega_N^\infty$ that end up at backward time $s \in [0,t]$ in some random site $\eta^\prime \in \Omega_N^\infty$. This measure at time $t$ describes the descendants of the founding fathers sharing a particular location in some subset of $\Omega_N^\infty$ at time $0$. Hence the \emph{types} are locally dispersed and are given by the measure $\theta$ independently of the location.
\end{itemize} 
Formally, the description above reads as follows, with two different definitions of the measure-valued dual $\wh Y^{\ast,\infty}$ corresponding to the two cases of recurrent or transient $\zeta^\infty$, namely, as a \emph{true coalescent}, respectively, a \emph{population of effective L\'evy-processes}, which are independent if and only if there is no block-resampling.

\begin{lemma}{\bf [Convergence of the measure-valued process $\wh Y^L$ and dichotomy: decomposition of $\Omega_N^\infty$ into mono-type regions versus local law of large numbers]}
\label{l.1562}
For all $t>0$ the following holds:
\begin{itemize}
\item[(a)] 
The $L$-approximations $Y^{\ast,L}$ of the finite dual process on $\Omega_N^L$ converge as $L \to \infty$ to the true finite continuum-space coalescent $Y^{\ast,\infty}$ or to pure random walks $\zeta^{\ast,\infty}$. The same is true for the measure-valued version, namely, for the entrance law $\wh Y^{\ast,L},\wh Y^{\ast,\infty}$ from locally countably infinite initial states. After embedding, in both cases this amounts to starting one particle at each site in $\Omega_N^{\infty,\rm disc}$, as specified in the embedding construction above~\eqref{e1431}. 

The limiting process $\wh \mfZ^{\ast,\infty}$ is either the continuum-space true coalescent $\wh Y^{\ast,\infty}$ or corresponds to independent migrating L\'evy-processes $\wh \zeta^{\ast,\infty}$. In formulas:
\begin{equation}
\label{e1585}
\CL\big[\wh Y_t^{\ast,L}\big] \Lto \CL\big[\wh Y_t^{\ast,\infty}\big] \qquad \text{or} \qquad  
\CL\big[\wh Y_t^{\ast,L}\big] \Lto \CL\big[\wh \zeta_t^{\ast,\infty}\big].
\end{equation}
The first alternative occurs for $t>0$ when the effective migration path $\zeta^\infty$ is recurrent, the second when $\zeta^\infty$ is transient. 

In particular, in all cases covered, and for $\wh H$ as in \eqref{e1775} and for all $t>0$:
\begin{equation}
\label{e1671}
\E \left[\wh H\big(\wh X^\infty_0, \wh Y^{\ast,L}_t\big)\right] \Ltoo \E \left[\wh H\big(\wh X^\infty_0, \wh \mfZ^{\ast,\infty}_t\big)\right].
\end{equation}
\item[(b)] 
There are two regimes of behaviour:
\begin{itemize}
\item[$\bullet$]
Suppose that collision occurs with probability $1$ and coalescence occurs in $Y^\infty$ at rate $d^\ast$ times the collision local time at a site in $\Omega_N^\infty$, or the hitting of arbitrary small balls by pairs occurs with probability $1$ and block-coalescence occurs with positive rate (recall the assumption in~\eqref{condextra}). 

Then for every $\eta \in \Omega_N^\infty$ there exists a $k=k(t,\eta) \in \Z$ such that $\wh Y^{\ast,\infty}_t$ restricted to $B_k(\eta)$ is equal to $\mathfrak{h} \otimes \delta_a$ for some $a \in A$, with $A$ a random countable subset of $\K$, which is a function of the path in time $[0,\delta)$, where every $\delta >0$ can be chosen. In particular, the state of $\wh Y^{\infty,t}(\cdot,\{a\})$ has a density w.r.t.\ $\mfh$.
\item[$\bullet$]
Suppose that no single-site or block-coalescence occurs. Then $Y^{\ast,\infty}$ is a collection of independent L\'evy-processes carrying a type and $\wh Y_t^{\ast,\infty}=\mathfrak{h} \otimes \theta$ for $t>0$. If block-coalescence does occur (recall~\eqref{condextra}), then the latter still holds in the \emph{Meyer-Zheng topology} on path space. 
\end{itemize}
\end{itemize}
\end{lemma}

Our goal in Step 1 was to prove Proposition \ref{prop.1010}(a). This is implied by Lemmas \ref{RWL}, \ref{RWLprec} and \ref{l.1562}, which are proved in Sections~\ref{sss.first}--\ref{sss.second} and Section~\ref{s.third}. In the remaining Steps 2--4 we prove Proposition \ref{prop.1010}(b)--(d) and Proposition \ref{P.DUAL}. Proposition \ref{P.1430} is implied by Lemma \ref{l.1562}. 


\subsubsection{Step $2$: Finite-dimensional distributions}
 
By the duality relation for finite $L$ and the convergence in Proposition~\ref{P.1430}, we have convergence of the expectations of the duality function for every $t>0$. Since the polynomials from the set $\{H(\cdot,Y_0^\infty)\colon\,Y_0^\infty \text{ finite coalescent}\}$ are law- and convergence-determining, we obtain convergence of the $t$-marginal distributions from Proposition \ref{prop.1010}(a) and Proposition \ref{P.1430}. By the space-time duality from \eqref{e1547}, we can extend this conclusion from the marginal to the \emph{finite-dimensional distributions} of the path of $\wh X^\infty$, which proves Proposition \ref{prop.1010}(b).


\subsubsection{Step $3$: Duality}

Proposition \ref{prop.1010}(c) and Proposition \ref{P.DUAL} are an immediate consequence of the previous result.


\subsubsection{Step $4$: Deterministic versus stochastic}

Proposition \ref{prop.1010}(d) is a consequence of Corollary \ref{C.sdual} and Lemma \ref{l.1562}(b), which establishes this property for the dual process. In particular, we have proved Proposition~\ref{prop.1010} once we prove the three lemmas in Sections~\ref{sss.first}--\ref{sss.second} and Section~\ref{s.third}.
\end{proof}


\subsection{Proof of the first constituent lemma}
\label{sss.first}

Lemma \ref{RWL} follows from a short argument of stochastic analysis, since its basis is the behaviour of finitely many spatial L\'evy-processes.

\begin{proof}{[Proof of Lemma \ref{RWL}]}\mbox{}\\
For the case where the change of location comes from the $(c_k)_{k \in \Z}$-L\'evy-process only, the claim follows from the construction as a \emph{projective limit}: the accumulated jumps added as $L \to \infty$ tend to $0$, and the $L^\prime$-truncated position of $\zeta$ remains constant for $L >L^\prime$ (i.e., above $-L^\prime$). We therefore focus on the motion coming from \emph{reshuffling}. The law of a \emph{single} observed walk is of the same type as above, and hence converges. The new aspect is that the migration steps of a tuple of, say, $m$ walks occur all at the \emph{same time} and not at independent waiting times. However, we can use that these steps are independently uniformly distributed in $-L$ balls, which have a projective structure, and converge uniformly in space to $0$ as $L \to \infty$ and $m$ is fixed. 

Namely, with the reshuffling, the $m$-tuple is a random walk in $(\Omega_N^\infty)^m$ such that the jump times of \emph{all} the marked walks in a $k$-ball chosen for resampling are \emph{coupled}, in the sense that all jumps to new locations occur at the \emph{same time} when the jump is taken as a reshuffling step in the $k$-ball. The jumps for the different components, however, are still independently taken according to the uniform distribution in the $k$-ball, for all $k$. Recall that for the \emph{migration} step this is different, because each component is chosen with equal probability, the jump rate is given by $m(c_k/N^{k-1})$ for size-$k$ jumps, and the jumps are drawn according to the uniform distribution. However, both have a projective structure and the jump sizes for $k \geq L$ tend to zero as $L \to \infty$. The position projected onto $\Omega_N^k$ only changes for jumps of size $k^\prime \geq k$. Therefore convergence of the $L$-approximation holds in both cases and we get convergence of $m$-tuples to the \emph{effective migration} for $m$-tuples.
\end{proof}


\subsection{Proof of the second constituent lemma}
\label{sss.second}
 
Lemma \ref{RWLprec} follows from a standard argument of stochastic analysis, since its basis is the behaviour of spatial coalescents.

\begin{proof}{[Proof of Lemma \ref{RWLprec}]}\mbox{}\\
Since the previous lemma gave us the convergence of the processes corresponding to pure migration/reshuffling, it remains to prove  
for the process on $\Omega_N^L$ that the vector of time points for the successive \emph{coalescence events} $T_i^L$, with $L \in \N$, $1 \leq i \leq k$, $1 \leq k \leq n$, and with $n$ the initial number of individuals up to the time horizon $t$, converge to the coalescence times of the process, which we defined on $\Omega_N^\infty$ and denoted by $T_i^\infty$, $i=1,\ldots,k$. 

For point-recurrent migration, the latter objects are defined in terms of the \emph{hitting times} $H_i^L, i \in \N$, respectively, $H_i^\infty$, $i \in I$, with $I$ a countable set with accumulation points of two or more dual elements together with the joint occupation times accumulated up to times $t$. Since in the limit the hitting points are immediately left we have now the \emph{collision local times} associated with the tuple during which a coalescence event may occur. We may construct this object for a given realization of the walks with exponential clocks for the coalescence events. In case of delayed coalescence without block-resampling this is the same as marking the time of the actual coalescence. For block-resampling and for Harris-recurrent or transient migration the picture is more subtle.

We go through these three cases for each of the claims in part (a) and (b) of the lemma.

\medskip\noindent
\textbf{Part (a).} 
Using the results derived in Appendix~\ref{appA} and work done in the literature, we have the following result for the coalescent process $Y^\infty$ and for the underlying migration path $\zeta$. 

\begin{corollary}{\bf [Local times]}
Consider the case of point-recurrent effective migration. Then the following hold.
\label{cor.1939}
\begin{itemize}
\item[(a)]
For the effective migration path $\zeta = (\zeta_s)_{s \geq 0}$ starting from $0$ there exist jointly continuous \emph{local times}
\begin{equation}
\label{e1926}
\{(L^0(s,\xi))_{s\geq 0}\}_{\xi \in \Omega_N^\infty}
\end{equation}
such that the accumulated local times up to time $t$ of $\zeta$ are given by the $\wt\CM([0,\infty) \times \Omega_N^\infty)$-valued \emph{local time process} $(L^0_t)_{t \geq 0}$ characterised by the property
\begin{equation}
\label{e1930}
L^0_t(\varphi) = \int^t_0 \dd s \int_{\Omega_N^\infty} \mfh(\dd\xi)\, L^0(s,\xi)\, \varphi(s,\xi), \quad t \geq 0,
\end{equation}
where $\varphi$ is any bounded measurable function on $[0,t] \times \Omega_N^\infty$. In particular, $L^0(s,\xi)$ is the density of $L^0_t$ w.r.t.\ the product of the Lebesgue measure on $[0,t]$ and the Haar measure on $\Omega_N^\infty$. Alternatively, $L^0_t$ is the space-time occupation measure up to time $t$ of $\zeta$ starting from $0$, written
\begin{equation}
\label{e2004}
L^0_t(\varphi) = \int_0^t \dd s \int_{\Omega_N^\infty} \mfh(\dd\xi)\,1_{\{\zeta_s=\xi\}}\,\varphi(s,\xi), \quad t \geq 0.
\end{equation}
\item[(b)]
The local time $L^0(\cdot,\cdot)$ arises as the limiting \emph{space-time average} of the occupation time of $\zeta$ in balls,
\begin{equation}
\label{e1936}
N^L \int^t_0 \dd s \int_{\Omega_N^\infty} \,\mfh(d \eta)\, 1_{\{\zeta_s = \eta \in \CB_{-L}(\xi)\}}
{\mathop{\longrightarrow}_{L \to \infty}} L^0(t,\xi) \text{ a.s.\ with $\zeta_0 = \xi$}.
\end{equation}
As a consequence, the rate of coalescence of pairs is $d^\ast L(t,0)$, associated with the difference process started in $\zeta_0^1-\zeta_0^2$, and the first hitting time of $\xi$ when $\zeta$ starts from $\eta$ is given by 
\begin{equation}
\label{e1940}
\inf \{t>0\colon\, L^\eta(t,\xi)>0\}.
\end{equation}
\end{itemize}
\end{corollary}

\begin{proof}{[Proof of Corollary~\ref{cor.1939}]}
\label{pr.cor1939}
The main tool to establish these properties is the \emph{local-time} theory of L\'evy-processes on \emph{ultrametric topological groups} worked out in Appendix~\ref{appA:IPLT}. In \cite{EF96}, the local time at a site in $\Omega_N^\infty$ is worked out for one specific case of our model class, namely, the case of \emph{stable migration}, i.e., $c_k=c^k$ for some $c \in (0,1)$, and \emph{no} block-resampling. The arguments in that paper are based on \emph{explicit formulas} for the \emph{transition kernel} and the \emph{characteristic function} in terms of $c$. 

In Appendix~\ref{appA:pot}, we look at the evolution in a new way, allowing us to obtain a general formula for the \emph{transition kernel density} of the migration and for the \emph{characteristic function}, valid for general parameters $(c_k)_{k \in \Z}$ and $(\lambda_k)_{k \in \Z}$. With these formulas we can still use the approach in \cite{EF96} to get, in Lemma~\ref{lemmA4}, the \emph{existence} of the local times of the underlying \emph{effective} L\'evy-process $\zeta$ on $\Omega_N^\infty$, as scaling limit of the joint pairwise occupation times of the $L$-approximations of $\zeta$. This approach works because we can replace the \emph{explicit} formulas in terms of (c) in \cite{EF96} by new formulas, with the same four key properties for the resulting expressions, as outlined in the proof of Lemma~\ref{lemmA4}.
\end{proof}

To complete the proof of part (a), we distinguish between two cases, (i) and (ii), dealing with the setting without, respectively, with block-resampling.

\medskip\noindent
\textbf{(i)}
\emph{Case: without block-resampling}. Note that \eqref{e1936} implies that the coalescence rate per collision time of the limiting pair of coalescing individuals is finite and equal to $d^\ast$, due to the choice of the resampling rates and the built-in scaling. This allows us to identify the new migration mechanism of $Y^\infty$ in the absence of block-resampling as the \emph{spatial $(c_k)_{k \in \Z}$-L\'evy-process} with \emph{delayed coalescence} at rate $d^\ast$ times the collision local time of a pair of partition elements. We can combine this structure with the space-time duality of $\wh X^\infty$ and $\wh Y^\infty$ developed earlier, to obtain as $L \to \infty$ convergence of $\wh X^L$ to $\wh X^\infty$, as explained above.

Namely, note that, due to the form of the rates $c_k/N^k$, $\lambda_k/N^k$, $k \in \Z$, and the factors $dN^L$ and $|\Lambda_{-\infty}| N^L$ (recall \eqref{e988}), we indeed obtain a \emph{limiting migration/reshuffling structure with local times}, and with coalescence rate $d^\ast$ per local time for a pair of two dual elements, since we get local times for the limiting \emph{difference L\'evy-processes} in $0$, or any other site in $\Omega_N^\infty$, due to the scaling of the parameters with the factor $N^k$. More precisely, we have that the rates $(c_k)_{k \geq -L}$, $(\lambda_k)_{k >-L}$, give us the existence of a \emph{collision local time} of pairs of effective migration (see Corollary~\ref{cor.1939}), and $d_L$, $\Lambda_{-L}$ actually give us a contribution to the total coalescence rate for pairs in macroscopic space-time regions, in the limit as $L \to \infty$, in the form of an \emph{occupation measure} of the \emph{local time process} in space-time sets. In other words, we have a \emph{delayed} coalescent \emph{based on the local time} of collision of pairs of paths and a finite rate $d^\ast$, in the limit as $L \to \infty$, \emph{provided} we can establish that the $\Lambda$-coalescence in the limit gives \emph{only coalescence for pairs}. The latter means that the \emph{higher} joint collision local time of $m$-tuples for $m>2$ as $L \to \infty$ are $o(1)$ and zero in the limit, and hence equals the sum of the local times of all possible pairs plus a correction term $O(N^{-L})$.

Given our scaling, we must show that the local times on $\Omega_N^L$ are \emph{in mean} strictly increasing over positive time intervals uniformly in $L$, and that the local times on $\Omega_N^\infty$ are strictly increasing over finite time intervals as well. This is a consequence of the moment equations.

To establish convergence of the \emph{first hitting times} for point-recurrent migration, we note that this holds because we have point-recurrence in the limit $\zeta^\infty$, and therefore the path convergence implies that convergence of the hitting times to the one of the limit. By~\eqref{e1936}, the accumulated occupation times converges, after the $N^L$ scaling, to the \emph{collision local time} for a given pair, and the difference random walk has a strictly positive local time at site $0$ after \emph{first} hitting that site, even for the $L \to \infty$ limit $\zeta^\infty$. We therefore have the proof of part (a) for the case where $(\Lambda_k)_{k \in \Z}$ are all $0$.

\medskip\noindent
\textbf{(ii)} \emph{Case: with block-resampling}. 
We need to show that \emph{adding the block-resampling} for point-recurrent migration preserves the convergence property $Y^L \to Y^\infty$ as $L \to \infty$ that we have shown in the model with only single-site resampling. The process of pure block-resampling with reshuffling is well-defined for \emph{finite} coalescents for every $L \in \N$ and for $L=\infty$. However, we must combine both the mechanisms of resampling for these finite coalescents. 

To that end we use the \emph{Trotter formula} that intertwines two groups of mechanisms, to get the Markov process corresponding to the sum of all operators for the dual process. The first group contains the combined other mechanisms, the second group contains the various block-coalescences. Note that the process of finite spatial $Y^\infty$ is a Feller process on subsets of a Polish space that are closed under the dynamics. The subsets are defined by $|Y^\infty| \leq m$, $m \in \N$, and are Polish spaces. Thus, we see that we have defined a Markov process $Y^\infty$ with c\`adl\`ag paths on the \emph{finite configurations}, with right-continuous paths (in backward time) \emph{on the space of finite coalescents}, which is the \emph{limit of the truncated processes $Y^L$}.

If we were to choose the model where the rate $d=d_L$ diverges with $L$, then we would get $d^\ast = \infty$, corresponding to \emph{instantaneous} coalescence. The first case above with $\Lambda_L=d \in (0,\infty)$ is the more interesting one, as we argue next.  

\begin{remark}{\bf[Cluster process as infinite volatility limit]}
\label{rem.1968}
{\rm If we let $d^\ast \to \infty$ as $L \to \infty$, then we get an \emph{instantaneously} coalescing dual process $\wt Y$, and an entrance law $\wh Y^{\infty,\sim}$. The forward process $\wh X^{\infty,\sim}$ is constructed similarly as $\wh X^\infty$. The latter process is degenerate, in the sense that we do \emph{not} obtain an \emph{SPDE representation} as we did for $d^\ast<\infty$, with a \emph{second-order diffusion operator} for the new process, whose \emph{quadratic variation diverges}. It is, nevertheless, a valid strong Markov and Feller random process. It does \emph{also} occur naturally as a \emph{continuum space rescaling limit} $\wh X^{\infty, \rm sta, \sim}$ and therefore does represent a \emph{new} universality class for our models on countable geographic spaces which are jump-diffusions.}\hfill$\spadesuit$ 
\end{remark}

\medskip\noindent
\textbf{Part (b).} First look at $Y^L,Y^\infty$. For Harris-recurrent but not point-recurrent migration, we have to prove that the coalescence times converge to the ones we have in the claimed limit system, and that there is no limiting single-site coalescence but \emph{only block-resampling} with $(\Lambda_k)_{k \in \Z}$. The former follows from \ref{appA:oneLP}, implying that in the absence of point-recurrence we get convergence of the local times to zero for the difference walk as $L \to \infty$. 

Next look at the block-resampling. We use the obvious \emph{coupling} for our models with different $L$. The coalescence times arising in the $L$-approximation of the coalescent will change as $L \to \infty$, since the added lower-level block-coalescing events lead to a coalescence \emph{before} upper-level block-coalescing events can act. This decreases the coalescence times successively. A strict decrease can only occur if partition elements are closer than $\eee^{-L}$ and a $k$-block-coalescence for some $k< -L$ is activated. This means that we need to consider the time spent by pairs of the tuple in the respective \emph{$k$-sausages, weighted by the coalescence rate}.

The block-coalescent events introduce changes of order $N^{-\mid k \mid}$, but at a rate that increases with a factor $N^{|k|}$. However, this factor is modified by $\bar \lambda_k$ and spreads with the mass flow of the effective migration, which is Harris-recurrent. The monotone structure leads to convergence of $Y^L$ to the claimed $Y^\infty$, since the paths of the tuple have a minimal distance if we start in different locations over a finite time horizon. However, if we start with multiplicity, then we immediately have for $T>0$ different locations by our assumptions.

\medskip\noindent
\textbf{Part (c).} In the case where $\zeta^\infty$ is transient, we have $d^\ast=0$ in the limit, there is no ``single-site coalescence", and the claim follows from Lemma \ref{RWL}, since we have a system of effective random walks. For transient migration the random walks with block-resampling have a minimal distance over finite time since the difference walk does not hit $0$, and hence the block-resampling rate has a \emph{finite} total block-resampling rate and convergence is immediate.
\end{proof}


\section{Continuum-space dynamics: III. Cluster formation, proof of the third constituent lemma}
\label{s.third}

Lemma \ref{l.1562} requires a long argument, since serious issues concerning the \emph{entrance} laws of initially \emph{locally infinite} initial states of coalescents from that into locally finite states come into play.

\begin{proof}{[Proof of Lemma \ref{l.1562}]}
The proof splits into two parts for each claim (a) respectively (b) of the Lemma, and are called: Part (a) and Part (b). 

Both are long because the details depend on the \emph{regime of effective migration under consideration} and the process is an \emph{entrance law} from locally \emph{infinite} states. In each of the three migration regimes
\begin{itemize}
\item[(i)] point-recurrent,
\item[(ii)] Harris-recurrent but not point-recurrent,
\item[(iii)] transient,
\end{itemize}
a completely different strategy is required. This means that the two parts each split further into \emph{three} cases. Part (b) is the hardest, since we must enter into considerable detail on cluster formation.

The proof relies on the results derived in Appendix \ref{appA} on L\'evy-processes on $\Omega_N^\infty$, and on results about \emph{cluster formation} for system on $\Omega_N^L$ for the case of effective recurrent migration we will obtain.  


\subsection{Part (a)}

The first claim in \eqref{e1585} follows from Lemmas \ref{RWL} and \ref{RWLprec}. The second claim in \eqref{e1585} we get from tightness and convergence of the finite coalescents $Y^L$ in path space to $Y^\infty$ as defined below~\ref{Yhatdef}. 

To obtain convergence as $L \to \infty$ of the \emph{measure-valued entrance laws} of $\wh Y^L$ to a measure-valued entrance law $\wh Y^\infty$ --- which arises as the entrance law for $Y^\infty$ starting from countably many dual individuals in $\Omega_N^{\infty,\rm disc}$, the dense countable subset of starting sites --- we need to show three things:
\begin{itemize}
\item[(I)]  
$\wh Y^{\ast,\infty}$ exists as defined in and above~\eqref{e1431}.
\item[(II)] 
$\wh Y^L$ converges as $L \to \infty$ to the unique claimed limit. 
\item[(III)] 
The necessary \emph{tightness properties} for the frequencies of partition elements in the various locations hold in the process $\wh Y^L$ as $L \to \infty$ (i.e., the \emph{no-dust condition} of Definition~\ref{def.dustfree} holds, so that in $\wh Y^{\ast,\infty}$ \emph{these frequencies are well-defined}). 
\end{itemize}
The term tightness is used above since this property is, first of all, tightness in $(\ell_1)^{\Omega_N^L}$, but also is the tightness condition in $\U_{\rm loc}$, the space of \emph{genealogies} in the model. This issue is further studied in~\cite{gh2}.


\subsubsection{(I)--(II) Existence and convergence}
\label{sss.ect}

Assume for the moment \emph{tightness} (which is the hard part and which will be proved below). We can think of a probability measure $\nu$ on a Polish space $E$ as the limit of the empirical measures of i.i.d.\ samples drawn from $\nu$ indexed by $\N$. We can do this first for finite measures and then for locally finite measures. Here, we consider locally finite random measures, using that the law is \emph{uniquely} determined by the \emph{moment measures}. It suffices to consider the entrance laws $\wh Y^L$ of the level-$L$ process from countably many individuals at each site of $\Omega_N^L$, to restrict the observation to some chosen ball, and to obtain a finite measure which, after normalisation, becomes a \emph{probability} measure. 

To show convergence in law of $\wh Y^L$ on $\wt\CM(\Omega_N^\infty \times \K)$, we need to know that the collection of all the $m$-coalescents sampled independently evolve according to the rules used to define $Y^L$, where we suppress further parameters from the notation for brevity. To sample the $m$-tuples, we consider dual individuals starting in a specific ball of $\Omega_N^\infty$ restricted to $\Omega_N^L$. Then, for every $m \in \N$, the corresponding \emph{empirical measure} of $m$-tuples (according to the normed Haar measure on $\Omega_N^L$ restricted to some ball) converges in law as $L \to \infty$ to the \emph{empirical measure of the $m$-samples of the coalescent} $Y^\infty$. The latter empirical measure is $\wh Y^\infty$ by the definition of $\wh Y^\infty$ and the fact that the polynomials are measure- and convergence-determining, and so the \emph{moment measures} converge to the moment measures of the limit. 

Indeed, the first observation is that, in case of tightness of partitions as introduced above, $\wh Y^\infty$ has the following \emph{finite subcoalescent property}: 
\begin{align}
&\text{subcoalescents sampled according to the Haar measure in a ball in $\Omega_N^\infty$ are given}\\
&\text{by the process $Y^\infty$ of \eqref{Yhatdef} started in that sampled \emph{initial} configuration.}
\nonumber
\end{align}
Here, $Y^\infty$ is the coalescent that we described at the beginning of Section \ref{sss.dupro}, specifying the jump rates for the transitions of the state of this Markov pure-jump process. 

The above claim on the $L\to \infty$ convergence $\wh Y^L$ to $\wh Y^\infty$ follows by using Lemma~\ref{RWLprec}, giving the convergence of every $m$-tuple as $L\to \infty$ to the limiting tuple, using a Skorohod embedding, and combining the Feller property of $Y^\infty$, starting from a fixed number $m$ of individuals as initial state, with the fact that the embedded $\mfh^L$ are converging weakly to $\mfh$ on $\Omega_N^\infty$. The Feller property implies that when the initial sies converge to a limit site in $(\Omega_N^\infty)^m$, then the paths of the processes converge in law as well.


\subsubsection{(III) Tightness respectively dust}

It remains to show the above assumed \emph{tightness property}, reading it as a \emph{no-dust condition} respectively show \emph{complete dust} if the no dust condition does not hold. Therefore we have to distinguish between the \emph{three migration regimes} (i)-(iii) for the effective L\'evy-process $\zeta$.


\paragraph{(i) Point-recurrent regime}

$\mbox{}$\\[0.2cm]
The proof comes in 7 Steps: Steps 1--2 are preparatory, Steps 3--7 execute the argument.

\medskip\noindent
{\bf 1.}
The \emph{non}-spatial $\Lambda$-coalescent may come down from infinity. This holds for migration on a finite set, but \emph{not} for the \emph{spatial} $\Lambda$-coalescent when the space is \emph{infinite}. The reason is that individuals may move out of a finite set before coalescence and not return in a given time interval with positive probability, which produces \emph{countably} many partition elements after positive time for the entrance law. However, contrary to the coming down from infinity property, the \emph{no-dust condition} on a countable geographic space is persistent under the addition of migration. This shows that some care is needed.

\medskip\noindent
{\bf 2.}
Define the countably infinite (even locally) dense \emph{subgroup} $\Omega_N^{\infty,\rm disc}$ (a special \emph{skeleton}, which is a \emph{subgroup}) of $\Omega_N^\infty$ given by
\begin{equation}
\label{e2009}
\Omega_N^{\infty,\rm disc} = \bigoplus_\Z Z_N.
\end{equation}
Recall that the entrance laws for $Y^L$ starting from \emph{countably many} individuals at \emph{each site} can be embedded in $\Omega_N^{\infty, \rm disc}$ starting from \emph{one} individual at each site of $\Omega_N^{\infty,\rm disc}$ (see above~\eqref{e1431}). To get the tightness property, we need Corollary~\ref{cor.1939}. This implies, in particular, that the coalescent $Y^\infty$, respectively, the $\K$-marked version $Y^{\ast, \infty}$, started in countably \emph{infinite} many individuals, \emph{one} at each site of $\Omega_N^\infty$ has the following properties.


\paragraph*{Properties of $\wh Y^\infty, \wh Y^{\ast,\infty}$ entrance laws.} 

The measure-valued dual $\wh Y^{\ast,\infty}$ in the model with $d+|\Lambda_0| >0$ (and hence $d^\ast >0$), with $Y^\infty$ started at a time horizon $T>0$ and observed at backward time $t \in (0,T]$, satisfies:
\begin{equation}
\begin{aligned}
&- \text{$\wh Y^\infty_t$ is started in $\mfh$ is \emph{locally finite} on $\Omega_N^{\infty}$}.\\
&- \text{$\wh Y_t^{\ast,\infty}(B,\cdot)$ with $B$ a bounded measurable subset of $\Omega_N^\infty$ is an \emph{atomic} measure}.\\
&- \text{$\wh Y_t^{\ast, \infty}(\cdot,B)$, $B \in \CB(\K)$, has an \emph{occupation measure} $\mfN^\varepsilon$ on $\Omega_N^\infty$ (recall~\eqref{e1745})}\\ 
&\quad\, \text{that is \emph{locally finite} in $\Omega_N^\infty$}.\\ 
&- \text{For $\wh \mfY_t^\infty$ there are countably many locations of partition elements in $\Omega_N^\infty$ 
and a fraction}\\
&\quad\, \text{$1-\varepsilon$ of the weight in a ball whose $\mfh$-weight is finite sits in \emph{finitely many sites}},\\
&\quad\,\, \text{i.e., the $\mfN^\varepsilon$ are locally finite}.\\ 
&- \text{The \emph{full frequency} $\mfh$ is carried by all these countably many partition}\\
&\quad\,\, \text{elements of $\wh \mfY_t^{\ast,\infty}$ in any bounded measurable subset of $\Omega_N^\infty$ (= \emph{no dust}).} 
\end{aligned}
\label{eq:2009}
\end{equation} 
The five properties in \eqref{eq:2009} are settled next in Steps 3--7.

\medskip\noindent
{\bf 3.}
On $\Omega_N^L$, the above properties follow from the \emph{no-dust property} of the \emph{on-site resampling}, which means that $1-\delta$ of the mass comes down from infinity for every $\delta > 0$ in positive time. To see why, we assume that $\Lambda_{-\infty}$ \emph{comes down from infinity}, and consider the individuals not moving before time $\delta=\delta(\varepsilon)$ chosen such that they are at least of frequency $1-\varepsilon$. They come down from infinity via coalescence at a site by time $t>0$. We let $\varepsilon \downarrow 0$, and with it $\delta(\varepsilon) \downarrow 0$, to argue with the decreasing mass of migrants that what we do not cover eventually has $0$-frequency. Hence, we have a \emph{$\delta$-truncated locally finite occupation measure $\mfN^\delta$} for the $\delta$-truncated coalescent at time $t$. In fact, we consider the dual population with their location at time $s \in [0,t)$ as a configuration of weighted locations in the \emph{Liggett-Spitzer space} (for this concept in the context of coalescents, see \cite{GLW05}). This property generalises when we assume only the \emph{no-dust condition} of $\Lambda_{-\infty}$. In that case we immediately have the no-dust property of the partition elements of the dual $L$-process (and hence the no-dust property of the local type frequencies contained in the $L$-forward process). 

\medskip\noindent
{\bf 4.}
To get the claim for $\Omega_N^\infty$, we have to transfer the argument from $\Omega_N^L$ to $\Omega_N^\infty$. First note that after the $L \to \infty$ limit is taken the $\Lambda_{-\infty}$-Cannings single site resampling turns into Moran resampling at rate $|\Lambda_{-\infty}|$). The technical point for the transfer is that now the effective migration path $\zeta^\infty$ leaves its starting site \emph{instantaneously}, but with high probability only moves into a small neighbourhood in the limiting process $Y^\infty$, and so does the joint return of two paths to the diagonal in other sites. As a consequence, for $L=\infty$, instead of all the migration steps considered in Step 3 we have to consider those steps taken by dual individuals leading \emph{out} of the chosen \emph{small} ball (a $k$-ball), but allow for steps \emph{within} the small ball. Hence the pair $\varepsilon,\delta (\varepsilon)$ from the argument in Step 3 in $\Omega_N^L$ now become $k,\delta(k)$ with $k \in \Z$, and we have to pass to frequencies of the partition elements in this ball and to the local times at all sites down from the ball size $k$, to obtain the effective coalescence time for the particular local dual subpopulation. Using this approach, we argue that the process comes down from infinity in small balls. 

Why does this line of reasoning work? The reason is, first, that on $\Omega_N^\infty$ the collision local time of pairs exists, and second, that a pair of paths on $\Omega_N^L$ coalesces at a strictly positive constant rate $d^\ast$ times the \emph{collision time} of a pair of partition elements. Hence, the rate for the first coalescence to appear increases by $d^\ast\,(\#^L_n)^2$, with $\#^L_n$ the number of pairs in a site, with $n$ dual individuals at level $-L$ in all the selected sites when \emph{no} migration event occurs. 

Turn to $\Omega_N^\infty$. If we consider a single ball in $\Omega_N^{\infty,\rm disc}$ and put on $n$ of the sites a dual individual, then we have $d^\ast$ times the collision \emph{local} time as the coalescence rate. By the same argument as before, we get that the total configuration comes down from infinity when we start with one dual individual at each site of $\Omega_N^\infty$ intersected with the ball. This follows by considering the \emph{local time functionals} of all the pairs of paths from the distinguished site, which in law is the same as the local time of $\zeta^\infty$ in $0$, started in the small ball with two particles. This observation extends to finitely many small balls. We have to argue now that this amounts effectively to a collapse of spatially distributed \emph{locally infinite} initial conditions to finite ones. Thus, we consider all sites in a ball, say a $k$-ball restricted to $\Omega_N^{\infty,\rm disc}$ with countably many initial individuals. 

\medskip\noindent
{\bf 5.}
We pass to a \emph{measure-valued entrance law} on the state space as the building block for our argument. Here, note that as $L \to \infty$ the number of possible paths \emph{increases}, but this is compensated because the time scale is \emph{speeded up by $N^L$} at $(-L)$-level. The fact that pairs are now chosen from a \emph{ball} at level $k \in \Z$ instead of a single site requires a correction factor, since the distance between sites results in a loss of coalescence rate up to the hitting time, which depends on the \emph{locations of the two individuals in the ball.} However, this spatial effect is already taken care of by using the \emph{measure-valued formulation} of the martingale problem, where states are \emph{measures} representing \emph{averaged components}. For the martingale problem of an increasing process averaged over the number of pairs from a ball, we need to show that this leads to coalescence rates of the form ``a \emph{positive} constant times $d^\ast$'', which is determined by the local time of $\zeta$ at $0$ as a function of the distance from $0$ in the average w.r.t.\ $\mfh$, since $d^\ast$ is claimed to be the factor that we get for the density of this quantity w.r.t.\ $\mfh$ in the $L \to \infty$ limit.

The approximation of the entrance law from the countable initial state, as limit of level-$L$ approximations of the \emph{initial} state, is obtained by restricting $\Omega_N^{\infty, \rm disc}$ to $\Omega_N^L$. We know that for finite $L$ the entrance law from countably many individuals per site of $\Omega_N^L$ results in a positive time $t$ in locally finite states. By using for $|\Lambda_{-\infty}| \neq 0$ (due to the no-dust condition for $\Lambda_{-\infty}$) the nature of the migration, it also results in states with a frequency $\geq 1-\varepsilon$ of \emph{finitely many partition elements per bounded set}. Note we can realise the latter by starting with one individual at each site of $\Omega_N^{\infty,\rm disc}$ and using the embedding of $\Omega_N^L$ into $\Omega_N^\infty$. This means that if we can prove for the $L$-approximations that we have a measure $\mfN^\varepsilon$, with the coming down from infinity property \emph{uniformly in $L$}, then we are guaranteed a \emph{locally finite occupation measure $\mfN^\varepsilon$} associated with the state for the entrance law of $Y^\infty$ starting everywhere in $\Omega_N^{\infty,\rm disc}$ with one particle (i.e., \emph{locally countably infinite} in $\Omega_N^\infty$), or at least such a property up to a frequency $\varepsilon$ for every $\varepsilon>0$. This condition holds because, as $L \to -\infty$, all rates are scaled up by $N^{|L|}$, compensating for the splitting of old sites into $N^{|L|}$ new ones. We need that $d^\ast$ is still the correct factor, because the average of components is governed by \emph{averages of local times for pairs}, \emph{averaged} over a sequence of \emph{shrinking balls}, which give the claimed local rate $d^\ast$ local time at zero. This is where the continuity of the local time at $0$ starting in $\eta \in \Omega_n^\infty$ comes into play. We now give the details.

\medskip\noindent
{\bf 6.}
We proceed by considering the dual dynamics on $\Omega_N^\infty$ (which we know exists as $L \to \infty$ for \emph{finite} initial states), and recalling that initial states are countably infinite only in sites from $\Omega_N^L$ and become locally countable in case $L=\infty$. We saw earlier that these entrance laws exist. We need to show \emph{tightness} of the local frequencies of partition elements in $L$ for these $L$-approximations of the initial state. This means that at each step we add $N-1$ further sites to each single site, where we have already countably many dual individuals, so that we increase the number of partition elements by ``increasing" the original number of basic individuals by $N^L$. But, since we speed up the rates by an \emph{additional factor $N$ at each splitting}, there is compensation of these effects by increased coalescence rates. We want to show that this leads in the limit as $L \to \infty$ to locally finite states in the sense of a finite $\mfN^\varepsilon$ , with \emph{well-defined frequencies for partition elements} in the initial population for each $t>0$. This can also be viewed as a tightness property in the weak-atomic topology as $L \to \infty$.

The strategy to prove the latter is to use the \emph{Liggett-Spitzer norm} for each finite $L$, based on the norm that we have in a single component, which is given by the number of partition elements, and afterwards let $L \to \infty$, applied to all partition elements that are making up most of the mass. Indeed, in the case where we only have the no-dust property and no coming down from infinity, we can apply this approach only to most of the frequency (per site). For the first case, we take out those dual individuals trying to jump before some time $\delta=\delta (\varepsilon)>0$ with jump size $\geq -\ell$ for some $\ell \in \N$, and treat them separately. The rest behaves for time $\delta$ like a non-spatial coalescent, during which time it may come down from infinity (or at least to a countable number of partition elements in the second case), as we saw below \eqref{eq:2009}. We argue this for the case where the process comes down from infinity, and afterwards argue that in the limit as $L\to\infty$ we are in fact in that specific case. 

The coming down from infinity property of $\Lambda_{-\infty}$ allows us to obtain that $1-\varepsilon$ of the frequency in a ball comes down from infinity, for every $\varepsilon >0$. This means that there are finitely many ancestors of families among individuals starting in a ball and surviving up to time $t$ somewhere as leading element.

In the dual backward model the instantaneously emigrating dual individuals from the time horizon $t$ back are moving still in the small neighbourhood, where they coalesce also with other migrants from that ball according to delayed coalescence at rate $d^\ast$ times \emph{their} collision local time. We need to guarantee that this holds in all balls. Here, we use that the total rates within a $k$-ball are infinite compared to the rates to move out (at least in the parameter regime that we assume here). This provides a mixing of partition elements in the $k$-ball or $k^\prime$-balls with $k^\prime<k$. The resampling arises at rate $d^\ast$ during the local time of two paths in the same location before the observation time $t$. Hence, as $k \to \infty$, we have delayed coalescence and the whole population is \emph{coming down from infinity}, provided we can \emph{stochastically bound} the local time of paths (which all start in a small ball) \emph{from below} by a factor times the local time starting from the same site and run for a strictly positive time. As a consequence, the dual population comes down from infinity locally and does \emph{not} form dust.

\medskip\noindent
{\bf 7.}
To get a locally strictly positive coalescence rate equal to $d^\ast$ times the local time at $0$ for $\zeta^\infty$ starting at $0$ and run for a positive time, we need to consider the first hitting time $H$ of a pair of effective L\'evy-processes starting from two sites in a ball. We claim that, because the effective L\'evy-process $\zeta^\infty$ is point-recurrent,
\begin{equation}
\label{e2104}
\lim_{k \to -\infty}\, \inf_{ {\eta,\eta^\prime \in \Omega_N^\infty} \atop {d^\infty(\eta,\eta^\prime) \leq k} } 
\PP_{\eta,\eta^\prime}\big(H< t\big) = 1
\quad \forall\,\,t>0
\end{equation}
(recall \eqref{ag30inf}). Indeed, let $T_\ell$, $\ell \in \Z$, denote the time for the pair of effective L\'evy-processes to reach distance $\ell-1$ when they start from distance $\ell$. Then $\PP_{\eta,\eta^\prime}(H< t) = \PP(\sum_{\ell \leq k} T_\ell<t)$ when $d^\infty(\eta,\eta^\prime) = k$. But in Proposition \ref{prop.two-four-regimes-K}(a) in Appendix \ref{app:digression} we show that the distance between a pair of \emph{point-recurrent} effective L\'evy-processes is transient to the left and with probability 1 hits $-\infty$ in finite time, i.e., $\PP(\sum_{\ell \leq k} T_\ell<\infty) = 1$ for every $k \in \Z$. Hence $\lim_{k\to-\infty} \PP(\sum_{\ell \leq k} T_\ell < t) = 1$ for all $t>0$, and the claim follows. In view of \eqref{e2104}, starting from $\eta,\eta^\prime$ in a $k$-ball with $k\to-\infty$, the collision local time is stochastically equal to the local time at $0$ for $\zeta^\infty$ starting in $0$, up to an error given by min (hitting time, $t$) for all $t>0$. Because $d^\ast t $ is the quadratic variation of the evaluations of types by a moment equation, we have therefore completed the proof in the point-recurrent regime. 

\begin{remark}{\bf [A different approach to dust-free verification]}
\label{rem.2190}
{\rm Alternatively, we may show in the \emph{forward} model that we only have finitely many ancestors from the subpopulation starting at time $0$ in a ball (i.e., a subset of $\Omega_N^\infty$ with finite weight) surviving up to a positive time. If we also have migration, then we need to show that, starting from a population of finite mass, there are finitely many surviving \emph{immigrants} in the ball before time $t$ that immigrated in that ball recently and that survive until time $t$. This property is well-known for the branching model and we need to argue that it persists here, because in small times and small populations the resampling model behaves like branching. Thus, altogether, we have \emph{finitely} many ancestors of the population originating in some ball surviving in arbitrary large balls. Hence a $(1-\varepsilon)$-frequency has finitely many ancestors for all $\varepsilon >0$, while the others immigrated and can have \emph{countably} many ancestors. This is a \emph{no-dust condition} for the local population.}\hfill$\spadesuit$
\end{remark}


\paragraph{(ii) Harris-recurrent but not point-recurrent regime}

$\mbox{}$\\[0.2cm]
First we give a reformulation, and later we conclude the proof under Part (b). The reformulation comes in 3 Steps.

\medskip\noindent
{\bf 1.}
What happens with the process $\wh Y^L$ when $\zeta$ is only Harris-recurrent? We need to show \emph{tightness} of the frequencies of the partition elements in the entrance law starting from locally countably many dual individuals by putting one dual individual at every site in $\Omega_N^{\infty, \rm disc}$, this time for the case where we have, in the $L \to \infty$ limit, only \emph{block}-resampling. Furthermore, we need to show in Part (b) that the process is locally mono-type because it is locally occupied by one partition element. Recall that in the \emph{recurrent} but \emph{not} point-recurrent regime we have $d^\ast=0$ in the limiting process $Y^\infty$, and hence \emph{no} on-site coalescence. For that reason we need to use a condition guaranteeing that, asymptotically as $L \to \infty$, the no-dust condition holds in a ``single site", i.e., zooming in on $k$-balls, allowing us to get \emph{locally} non-degenerate and dust-free partition frequencies, i.e., locally finitely many partition elements for $1-\varepsilon$ of the frequency, and for Part (b) in the limit $k \to -\infty$, we need to cover small enough balls. For this we need that
\begin{align}
\label{a2231}
&\text{eventually in small balls the block-resampling guarantees that the coalescent locally}\\
&\text{has at least $(1-\varepsilon)$-frequency for \emph{finitely many partition elements} for every $\varepsilon >0$}\nonumber\\
&\text{and these partition elements cover small enough balls completely.}\nonumber 
\end{align}
This means that eventually the entrance law collapses to a single partition element in any sufficiently small ball; in other words, coalescence takes place before the ball is left.

To prove the first part of \eqref{a2231}, note that in the example of a \emph{stable} process we have $\Lambda_{-\infty}^{\rm block} = \Lambda = \lim_{k\to\infty} \Lambda_k$, which does \emph{not} satisfy the no-dust condition. We want to replace this property by the requirement that the \emph{normalised $\Lambda_k$} converges for $k \to -\infty$ to a limit, called $\Lambda_{-\infty}^{\rm block}$ (which in the limit replaces a whole sequence of mechanisms as we zoom into a site), with the goal to achieve that the coalescent at least effectively satisfies the \emph{no-dust condition}, having an infinite sequence of transitions that forces the very small balls to collapse to one partition element. This must be shown to ensure that stochasticity is generated via the infinite rate of block-resampling after passing through \emph{countably} many \emph{small jumps} right after starting (in backward time), in case the partition elements have not yet been trapped locally in one partition element. This means that it is the total mass of $\Lambda_k$, $k \in \Z$, that matters for that property. For that reason we do not need to make any extra assumption on $\Lambda_{-\infty}^{\rm block}$, as in \eqref{cond3}, but we need the second part of condition~\ref{condextra}. 

Note that many dual particles in a $k$-ball forces that many smaller balls must be occupied, leading to a growing number of transitions in the $L$-truncated system. At the same time, many partition elements are removed by the smaller ball block-resampling events collapsing small balls. We have to show that this brings us down to \emph{finitely} many partition elements in \emph{positive time} for most of the frequency. 

\medskip\noindent
{\bf 2.}
Look at the forward evolution in time. Since the single-site resampling vanishes in the limit as $L \to \infty$, we need to establish that the block-resampling on $\Omega_N^\infty$ is still well defined, since it involves countably many jumps that are locally small as long as small balls are not covered by one partition element. These jumps affect only small balls shrinking to zero, but there are countably many small balls in a finite ball, each of which has its cascade of small jumps, which are all affected by transitions as long as they are not part of one. Recall that if we are in a locally mono-type state, then the effective transition rate is locally finite and we can define the transitions taking place as in a Markov pure-jump process. Therefore we have to establish that we reach such a state \emph{immediately}. Below we show the \emph{no-dust condition} and the \emph{mono-type property in the forward process} by looking for the formal proof \emph{backward} and using the dual. 

The no-dust property of the state that is reached in strictly positive time comes from the fact that we have \emph{infinitely many jumps} below any size $k$, which is perfectly well-defined in locally mono-type states, and we need to argue that these force any state to \emph{instantaneously} enter a \emph{locally mono-type state}. Note that, according to our conditions, over finite time there is always a largest jump affecting a chosen ball. Namely, \emph{looking forward}, the jumps of the process are given by the \emph{deterministic resampling map} combined with an independent choice of type depending on the present state, which a.s.\ is the identity \emph{if and only if} the state of the $k$-ball concerned is \emph{mono-type}, and the same holds for $k^\prime <k$. Next, in the forward model we look in more detail at a sequence of random mappings generating the jumps. The jump removes a local fraction of frequencies and adds a piece of \emph{frequency of one type uniformly in the $k$-ball concerned}, and chosen at random according to the present frequency in the ball. This means that local peaks of a single type loose proportionally and are replaced by multiple types that are more uniformly distributed. The \emph{unique fixed state} of the mapping restricted to the $k$-ball is the mono-type configuration in that $k$-ball. 

Thus, the sequence of jumps acting on a $k$-ball or a $k^\prime$-ball with $k^\prime <k$ change the state via a step that concerns a countable selection of $k^\prime$ to perform a $k^\prime$-jump. Each such jump introduces a difference that in space is a fraction of $\mfh$ in that small ball. We need to show that there is always a $\wt k$ such that for all $k^\prime <\wt k$ the state is already \emph{locally mono-type}, and is so as well for all $k$-transitions for $k \leq \wt k$. For that we must show that every limit point of the state of a ball in the $L \to \infty$ limit is locally mono-type, in the sense that the normalised measure restricted to $B_k(\eta)$ for some $\eta$ converges to $\mfh \otimes \delta_a$ for some $a \in \K$. We must also use the existence and continuity of the density of the state of $\wh X_t^\infty$ w.r.t.\ $\mfh$, to obtain the claim that this type prevails in a small enough ball around $\eta$. The latter is based on the fact that \emph{countably} many small jumps occur in a finite time interval, by our assumptions on the parameters. We need to argue that therefore the countably many small jumps occurring just after the time horizon $t$ amount to a \emph{jump into that locally mono-type state instantaneously} after time $t$ via an \emph{entrance law}. For that we return to the dual process again.

\medskip\noindent
{\bf 3.}
In the backward picture for an \emph{entrance law} of \emph{locally countable initial states} this means that, before the $k$-jump occurs, the smaller $k^\prime$-jumps have brought the state locally into just one partition element, which is subsequently fixed in a $k$-ball under the $k^\prime$-coalescence events for $k^\prime \leq k$. At a fixed positive backward time $t$, we therefore seek coalescent events visible for balls of a minimal size $L \in \Z$.

We use the convergence of the truncated systems on $\Omega_N^L$ to verify the two properties of the orbit of the map formulated above for the entrance law after a jump time forward, and here we look backward and use the duality. In the Harris-recurrent regime, the first property (locally finite) and the second property (mono-type) are intertwined, which is different from the case $d^\ast > 0$. In the strongly recurrent regime, to get both properties (forward the mono-type property and backward the locally single partition element) we use a strategy that we develop later in the proof of Part (b). We complete the argument there.


\paragraph{(iii) Transient regime}

$\mbox{}$\\[0.2cm]
The proof comes in 4 Steps.

\medskip\noindent
{\bf 1.} 
For an entrance law starting from locally countably many dual individuals, \emph{tightness} of the sequence of partition elements at a certain geographic location works in the recurrent case, as we have established above, \emph{but fails in the transient case}, where we do \emph{not} have tightness, and in fact can show that there is the formation of \emph{complete dust}. Namely, starting with locally countable duals, for the \emph{population of dual paths} we just get independent effective L\'evy-processes with values in $\Omega_N^\infty$ for $m$-tuples of \emph{typical} individuals drawn from the entrance laws on $\Omega_N^L$ in the $L \to \infty$ limit, which is \emph{different} from the limit $Y^\infty$ of the finite dual processes we get from $Y^L$. This means that we have a well-defined $Y^\infty$ in all cases, but $\wh Y^\infty$ exhibits \emph{two} possible regimes: \emph{coalescing} motions versus \emph{independent} motions, the dual based on the coalescent with block-resampling versus the dual $\wh \zeta^\infty$ of independent effective L\'evy-processes. We must verify next that in the limit as $L \to \infty$ the no-dust condition does \emph{not} hold in the \emph{transient} regime and we get \emph{complete} dust formation. First we consider the case without block-resmapling, and afterwards add block-resampling.

\medskip\noindent
{\bf 2.}
We already dealt with the case \emph{without} block-resampling, since we get single-site coalescence at rate $d^\ast=0$ in the transient regime, hence \emph{independence of the dual paths} in any $m$-tuple, and hence \emph{complete dust} in the $L \to \infty$ limit. Therefore we defined $\wh Y^{\ast,\infty}$ as $\wh \zeta^\ast$ in \eqref{e1980}. For the empirical measures in the $L \to \infty$ limit we have a law of large numbers, namely, duality immediately gives that the types are in the state $\theta$, as follows from an explicit computation, because the mean measure is $\theta$ and the variance of the evolved state is zero for all times. So, we henceforth focus on the transient regime \emph{with} block-resampling. 

\medskip\noindent
{\bf 3.}
\emph{With} block-resampling the situation is more subtle. What happens in the \emph{forward} process is that the averaging of many small balls, by its construction constituting a state in a macroscopic ball, leads to a \emph{law of large numbers} due to spatial independence of the ``noise" due to the small jumps. Namely, the action of the higher-level block-resampling is instantaneously ``destroyed" by the \emph{countably many} low-level block-resampling events due to the infinite total rate valid, off the constant state, in each of the countably many balls smaller than the $k$-ball, together with the mixing through migration and/or reshuffling, which is transient. This generates a noise that acts ``instantaneously" after a higher-level transition, so that again a law of large numbers takes effect and we have only finite-dimensional distribution convergence to the state $\mfh \otimes \theta$. Another point of view is that the $\Omega_N^L$ systems converge for $t\to \infty$ to shift-ergodic equilibria. Therefore the averages over many sites converge to the mean state, which is $\theta \otimes \mfh$. 

Here are the formal details. If block-resampling is present, then we need to show in the forward process that the mono-type component created at a resampling/reshuffling jump (after taking some fraction away from every type in the $k$-ball and adding it to one of the types present in the $k$-ball, with $k$ some level) \emph{is dispersed by the infinite sequence of the smaller block-resampling events in smaller balls} contained in the $k$-ball. To make this formal, we pass to the \emph{backward} picture and the dual process. We have to show that two \emph{typical} dual individuals, drawn from a ball under the measure $\wh Y^{\ast,\infty}$, do not coalesce in finite (backward) time in the $L \to \infty$ limit, and hence in the forward picture the correlations between two disjoint balls tend to zero, so that the state in a given ball is the mean measure $\mfh \otimes \theta$ by the law of large numbers. Back to the dual, the probability of two \emph{typical} dual individuals from a $k$-ball to be at some time $t$ in the same $k^\prime$ ball with $k^\prime <k$ and coalescing tends to zero as we pass from finite coalescents to the entrance law. Therefore in the transient case we have no tightness of the partition element frequencies in spatial subsets of $Y^{\ast,\infty}$, leading to the entrance law of \emph{trivial partition frequencies} for the countable dual population $\wh Y^{\ast,\infty}$, which we have defined above as $\wh \zeta^{\ast,\infty}$ (complete dust) rather than as a spatial coalescent leading to non-degenerate partition elements.  

Observe that the jump through $\Lambda_k$-coalescence with block reshuffling is $0$ if and only if the $k$-ball concerned is in one partition element and the largest partition element has a chance equal to its size not to loose at coalescence. The idea is to establish that, as we zoom in on a site $\eta$ through $k$-balls, we see that the piece of mass of one partition element in a state that after, say a $k$-jump, is created is removed instantaneously by fragmented pieces as the time horizon $t$ changes to $>t$ (recall the \emph{coalescent with rebirth} construction). Namely, the rates of migration are such that $\zeta^\infty$ is transient, and the local time in $0$ is $0$. Therefore the neighbouring sites do \emph{not} anymore quickly get in the same partition element, and so we have an \emph{effect of separation of partition elements} locally. To get this rigorously, we have to show that the coalescent probabilities between two typical individuals from different small balls tends to zero as we shrink the ball sizes $k$, i.e., $k \to -\infty$. We have to see what happens in the limit $L \to \infty$ to the corresponding probabilities, since for fixed $L$ we have non-degenerate partitions. This means that the aggregated populations on $\Omega_N^L$ have as $L \to \infty$ vanishing coalescence probabilities for all time horizons $t>0$. The populations of $\wh Y^{\ast,L}$ as $L \to \infty$ form dust and the measure state equals $\mfh \otimes \theta$.

More precisely, consider $k,k^\prime \in \Z$ with $k > k^\prime$ and $\eta^1,\eta^2$ sites sampled from $\wh Y^{\ast,\infty}_t$ for $t>0$ in two distinct $k^\prime$-balls in the $k$-ball. Then the probability to coalesce for the two typical paths starting in $\eta^1,\eta^2$ converges to zero as $k^\prime \to -\infty$ uniformly in $k$, and the measure state of $Y_t^\ast$ follows a law of large numbers and equals its mean $\mfh \otimes \theta$.

\medskip\noindent
{\bf 4.}
We need to get \emph{convergence in path space}. Looking at the path, we see in the \emph{forward view} that the block-resampling events lead to jumps if the block is multi-type at the corresponding random jump times and a state from $\mfh \otimes \theta$ appears, which, however, does \emph{not persist} and is \emph{instantaneously} left. Note here from the duality relation at time $t$ that with probability $1$ no such jump occurs. Therefore the convergence as $L \to \infty$ cannot happen in the Skorohod topology in path space $D([0,\infty),E)$, since we have to account for exceptional (random) times and the instantaneous removal of the state therein. This is remedied by using the \emph{Meyer-Zheng topology} on the set $D([0,\infty), \wt\CM (\Omega_N^\infty \times \K))$ instead of the stronger Skorohod topology used for c\`adl\`ag paths. In the Meyer-Zheng topology a sequence of paths $x_n$ from $D([0,\infty),E)$ converges to an element $x$ in that same space when
\begin{equation}
\label{e2151}
\int^t_0 \dd s\,\,d(x_n(s),x(s)) \ntoo 0 \quad \forall\,\, t>0, \text{ with } d(\cdot,\cdot) \text{ the metric in $E$},
\end{equation}
and \emph{two paths $x$ and $x^\prime$ are equivalent} in the Meyer-Zheng topology (i.e., form one site) when the above expression for $x$ and $x^\prime$ is identically $0$. Therefore we have to verify the convergence of paths as $L \to \infty$ in the Meyer-Zheng topology.

To make this task precise we use the \emph{backward view} in the entrance law for the dual process starting from the Haar measure, and proceed as follows. Couple the different block-resamplings for the different $L$ by adding successively new transitions at levels $k \downarrow -\infty$, to construct a sequence of transitions of level $k^\prime < \ell$ backward in time after a jump of a certain level $k=\ell$, which leads in forward time at the target time to a state $h \otimes \theta$ in finite windows of observations. Note that we have countably many small jumps ($<\ell$) occurring arbitrarily close to the jump times that we consider. Since the jumps change the path in finite time slots \emph{uniformly} at most by the jump size, convergence in the Meyer-Zheng topology follows.

\medskip\noindent
This completes the proof of Part (a) in Lemma \ref{l.1562} except for the Harris-recurrent regime, where a gap still has to be closed.

 
\subsection{Part (b): (i) Point-recurrent regime}

Part (b) is concerned with the key qualitative features of the continuum system on $\Omega_N^\infty$, and involves completely different arguments in the three migration regimes. Hence our exposition again has three corresponding parts. All lemmas assume the migration regime without further notice. We need to address the same items (I)--(III) as in Part (a).


\subsubsection{(I)--(II) Existence and convergence}

By the same argument as in Section~\ref{sss.ect}, existence and convergence follows once we establish \emph{tightness}.


\subsubsection{(III) Tightness: program}

\paragraph{Preparation}
We need to show that if we start $Y^\infty$ from locally countably many dual individuals by putting \emph{one} individual at \emph{each site of $\Omega_N^{\infty, \rm disc}$}, then the coalescent entrance law in a ball $B_\ell(\wt \eta)$ intersected with $\Omega_N^{\infty,\rm disc}$ has the following property. If we zoom in on a site $\eta \in \Omega_N^{\infty,\rm disc} \cap B_k(\wt \eta)$ at time $t >0$, and restrict the dual population to individuals that \emph{stayed in} $\CB_k(\eta)$ until the positive backward time $\delta$, and then in this subpopulation observe for $k< \ell$ the number $N^\delta_{k,t}(\eta)$ of partition elements at time $t$, then this number is non-increasing in $k$, is \emph{finite} for every $k \in \Z$ a.s.\ for $d>0$ or $|\Lambda_{-\infty}|>0$, and is \emph{coming down from infinity} (otherwise, under the no-dust condition, this holds only with frequency $1-\varepsilon$ for all $\varepsilon > 0$). The reason is that, as $L \to \infty$, the limiting measure $\wh Y_t^{\ast,\infty}$ in $\Omega_N^\infty$ is locally finite due to the fact that the scaling holds uniformly. 

If we let $\delta \downarrow 0$, then we obtain the entrance law from the set of countable states that are built with the frequencies of the partition elements restricted to individuals remaining a strictly positive backward time in the ball. In this way we get the random variable $N^0_{k,t}(\eta)$, the number of partition elements in the coalescent starting originally in the ball $\CB_k(\eta)$ at backward time $t$ in the entrance law. Note that, as $t \downarrow 0$, this number diverges at every site in $\Omega_N^{\infty, \rm disc}$. 

What happens when the ball in which we observe \emph{shrinks to a site}? The claim we prove next is that
\begin{equation}
\label{e1687}
\lim_{k \to -\infty} N^0_{k,t}(\eta) = 1 \quad \forall\,t >0,
\end{equation}
so that there is a $k\in\Z$ such that $N^0_{k^\prime}(\eta)=1$ for all $k^\prime \leq k$ for $\mfh$ every $\eta$. The proof tackles \emph{two different features} we formalise next.


\paragraph{Feature 1: Sites are mono-type}

$\mbox{}$\\[0.2cm]
The reason for \eqref{e1687} is the fast time scale for the very negative levels, making the behaviour look like the longtime behaviour of systems on $\Omega_N^L$. We proceed as follows. We know from \eqref{e2104} that in the \emph{point-recurrent regime} the hitting probability by a fixed time $t$ of a pair of sites from $\CB_k(\eta)$ converges to $1$ as $k \to -\infty$. If $N^0_k(\eta) \geq 2$, then there exist two populations of different types with positive frequencies each, say $(\delta^k_1,\delta^k_2)$, with the two $\limsup$'s strictly larger than zero, denoted by $\delta_1^\infty,\delta_2^\infty$ (later we discuss the possibility that one is zero exactly at time $t$). Consider the pair of paths that have not coalesced by time $t$ and start in $\CB_k(\eta)$. View these paths backward in time. We use the existence of the local time for the difference path. Namely, being at the site $\eta$ allows us to conclude that the local time a.s.\ increases already at time $0$, and up to time $t$. Positivity of both frequencies $\delta_1^\infty,\delta_2^\infty$ therefore means that there exist two \emph{typical} individuals in the initial configuration of the process $Y^\infty$ whose paths do not meet and coalesce by time $t$ starting in arbitrary small balls around $\eta$ at time $0$. The paths are typical paths starting in $B_k(\eta)$. Therefore the two paths must converge to each other as $k \to -\infty$. In the limit $k \to -\infty$ this means that, viewed backward, they have \emph{both} started from $\eta$ without coalescing by time $t$. 

To continue we argue by contradiction. If $d^\ast>0$, then the paths hit by time $t$ (forward) must coalesce \emph{instantaneously} once the system starts from \emph{countably} many individuals at a site, and the paths must be typical. This contradicts the assumption of starting with one individual at each site in $\Omega_N^{\infty,\rm disc}$, so the case $\delta^\infty_1,\delta^\infty_2>0$ is not possible. What happens if $\delta_1\delta_2=0$ at time $0$? Then the population starting in $\eta$ forms one partition element in the limit $k \to -\infty$, with a particular type. 

How to obtain this property for a small ball around $\eta$ with the same type? Here, we use the \emph{continuity of the density} in combination with it being $0$ or $1$ in a type, as we saw before, which implies that there is a neighbourhood in which the $\K$-value at $\eta$ is also taken. Therefore, to prove \eqref{e1687}, we must show that we have a \emph{continuous} density w.r.t.\ $\mfh$. More precisely, we need to show that the measure is \emph{absolutely continuous} w.r.t.\ $\mfh$ and that the density has a \emph{continuous version}.


\paragraph{Feature 2: Existence of a continuous density}

$\mbox{}$\\[0.2cm]
We have to argue that $\wh X_t^{(\infty)}$, as limit of $\wh X_t^L$ projected onto a fixed type $a \in A$, has a density w.r.t.\ the Haar measure $\mfh$, where the $\wh X^L$ embedded in $\Omega_N^\infty$ have a density w.r.t.\ $\mfh^L$. We must argue that this density can be represented as a density function on $\Omega_N^\infty$ by identifying all the sites in an $L$-ball when we decompose $\Omega_N^\infty$ into disjoint open and closed $L$-balls. Note that, for every $L$, we have \emph{continuous functions}, since in $\Omega_N^\infty$ they are \emph{constant in the $L$-ball}. We have to argue that these densities satisfy
\begin{equation}
\label{a2338}
\begin{aligned}
&\text{$\wh X^L(\cdot,\{\mfa\})$ converges as $L\to\infty$ to a \emph{limit density} in $L_1(\mfh)$}\\
&\text{whose measure equals $\wh X^\infty(\cdot,\{\mfa\})$ on $\Omega_N^\infty \times \K$}
\end{aligned}
\end{equation}
when restricted to a bounded geographic set and a random $\mfa \in A$, with $A$ a random set depending on the path $\wh X^\infty$ in the time interval $(0,\delta)$, where $\delta$ can be any strictly positive number. 

Furthermore, we want to show convergence in law as $L \to \infty$ to a random measure having a version of a \emph{density} w.r.t.\ $\mfh$, $\mfh$-a.e., that is a 
\begin{equation}\label{a2586}
\text{\emph{$\{0,1\}$-valued} and \emph{continuous} function. }
\end{equation}
This means we are locally in points \emph{monotype in some neighborhood}. We know on $\Omega_N,\Omega_N^L$ for $t \to \infty$ we get such a behaviour. 

These properties of the process $\wh X^\infty$ suggest that they represent the cluster growth seen in countable geographic spaces in the large-time limit. The situation is more subtle, however, since the corresponding effect in $\Omega_N^L$ leads to $\wh X^{\infty,\rm sta}$, which is different from $\wh X^\infty$. The proof we construct is based on an extended and refined analysis of the cluster formation on $\Omega_N,\Omega_N^L$ looked at the right way to see $\wh X^\infty$ and relate this to $\wh X^{\infty,\rm sta}$.

The reason why the above works is that for every fixed $L$ the discrete model converges, as $t \to \infty$ in finite balls of any size, to a \emph{mono-type state}, where the latter are growing in time as $t\to \infty$ in a certain manner. This is what we need because, by construction, for the $\Omega_N^\infty$-dynamics we speed up the resampling time on the $(-L)$-level by $N^L$, which should lead to a local monotype structure. The former mentioned behaviour of the discrete model on $\Omega_N$ as $t \to \infty$ is called \emph{cluster formation}. By shifting the active level from $-L$ to $-\infty$, we turn this effect into something different: it now happens \emph{instantaneously} for positive time. For various models on countable groups as geographic space, the \emph{growth of clusters} has been studied in detail. Different forms of cluster growth have been identified depending on model parameters. To describe the cluster growth in space and time, we use the \emph{volume of balls} "covered" by one type as a function of time. This has been extensively studied on $\Z^d$ and on $\R^1$,~\cite{A79,FG94,Winter2002,D93}.

The regime of \emph{fast clustering} is what we expect in the point-recurrent case, where the volume of the ball covered by the mono-type state grows like the range of the underlying effective random walk. This has been studied on $\Omega_N$ in a special case (see \cite[Theorem 2]{Klenke1996}, which treats the \emph{stable case} only and has \emph{no} block-resampling). We need to allow for \emph{general} parameters, without the self-similarity under the truncation and the shift transformation in \cite{Klenke1996}, and we need to \emph{add block-resampling}. The cluster formation has been taken up in~\cite{GHKK14} in the context of spatial Cannings models on $\Omega_N$ \emph{with} block-resampling and our choice of parameters, but only in the $N \to \infty$ limit. Here we need the same for \emph{fixed $N$}.


\paragraph*{Outline for the proofs.}

The argument to get tightness (item (III)) for the point-recurrent regime is organised into 7 Steps: Step 1 formulates the \emph{proof strategy}, Steps 2--7 give the details, with Steps 2--5 containing the \emph{concepts} and the \emph{necessary lemmas} and Steps 6--7 the \emph{details of proof} of these lemmas.


\subsubsection{(III) Step 1.\ Explanation and strategy}

We need to show that the rate to move out of a $(-L+\ell)$-ball for $\ell\in\N$ becomes small as $L \to \infty$ (compared to the motion \emph{inside} this ball), which means that the individuals become uniformly distributed on small enough balls. The consequence is that the individuals in the $(-L+\ell)$-ball resample rather than emigrate far away, and end up in other partition elements, or they started at another site and immigrated from that site at the observation time $t$, as we zoom in on the site. We want to argue that, as we shrink the ball around $\eta$, the densities restricted to the $k$-ball around $\eta$ successively converge to a constant as $k \to - \infty$, which is $1$ for one particular $a \in A$ and $0$ otherwise. We want to do this by using duality and by showing that the population $\mfh$-a.e. in the $k$-ball is in \emph{one} partition element for $k$ below some random value $K$, for $\mfh$-a.e.\ $\eta \in \Omega_N^\infty$. In the backward picture the dual population from that ball is given by partition elements that have locations at time $0$ somewhere in $\Omega_N^\infty$, where the type is sampled for the partition element and \emph{all} its members. Hence the claim is that $B_k(\eta)$ for $t>0$ is covered by one partition element for small enough $k$. More precisely, the claim is that 
\begin{equation}
\label{e2197}
\begin{aligned}
&\text{$\mfh$-a.e.\ a site of $B_k(\eta)$ is covered by \emph{one partition element}}\\ 
&\text{when $k \leq K$ for some random $K \in \Z$ for $\mfh$-a.e. point $\eta$,}
\end{aligned}
\end{equation}
where we need to \emph{define even more precisely} what we mean by being covered. We do this in Step 2.

\begin{remark}{\bf [A feature of genealogy]}
\label{r.2355}
{\rm Note that, in the forward process in a given ball $B_k(\eta)$, the population in forward time descends at a given time $t$ from ancestors at time $0$, which have last immigrated to the ball at some time $s$ and stayed there until time $t$. The claim is that for $k \leq K$ the ball $B_k(\eta)$ is covered by one subfamily of full frequency for \emph{large $t$ on $\Omega_N^L$} and, due to the limit $L \to \infty$, also \emph{for every $t>0$} for the \emph{corresponding process on $\Omega_N^\infty$}. Namely, as $L \to \infty$, the resampling is speeded up by $N^L$, and therefore the question is whether this \emph{over compensates} the splitting into $N^L$ many active sites compared to the $L=0$ model, so that for $L=\infty$ we see the cluster formation \emph{for every $t>0$}. The fact that it compensates is standard, but we need to make sure that our conditions on the parameters and the fact we zoom in on small balls lead to a further speeding up, so that at time $t$ on $\Omega_N^\infty$ we see where macroscopic cluster formation occurs, leading to the local clusters for $L=\infty$ (as in discrete geographic spaces at late times there are large almost mono-type balls).}\hfill$\spadesuit$
\end{remark}

We need the following fact based on the \emph{$L$-approximation}, combined with what we know for the discrete geographic space $\Omega_N^L$ about the behaviour as $t \to \infty$. Namely, on $\Omega_N^L$ for balls that grow with time $t$ in a specific way, we can for a sequence of $n$-balls around a site determine exactly such a sequence of times $t_n$, with $\lim_{n \to \infty} t_n = \infty$, for which the state converges to a \emph{mono-type state} in the weak topology. For every \emph{fixed} ball, this follows from the ergodic theorem as $t \to \infty$ on $\Omega_N^L$. In particular, finite sets get covered by the mono-type clusters in the strongly recurrent case. We can therefore ask at what scale do we have to increase time to cover \emph{growing} balls, and with what probability. In many models it is known that there is a \emph{deterministic scale} as a function of $t$ at which clusters grow, and \emph{random multiples} of that scale of cluster sizes are observed on that scale. 

More precisely, in the strongly recurrent regime, we can consider a set of scales that depends on $d^\ast$ and on the effective migration coefficients, such that the \emph{critical ball grows like a random variable times a deterministic function of the volume calculated from the recurrent Green function} of the effective migration. The law of the random variable is determined via duality as a specific stochastic process on $\Omega_N^\infty$. We use duality to verify some key properties of this random variable. The critical ball is the ball around a site such that the limiting empirical measure in that ball and all smaller balls is asymptotically a $\delta$-measure on one type as we show later generalizing~\cite{EF96}. This is explained in detail in Steps 2--4 and is then \emph{modified} so that it can be applied to $\wh X^\infty$.


\subsubsection{(III) Step 2.\ Description of cluster growth: basic objects of analysis}

The results we need on cluster formation in population models \emph{refine} earlier work on this topic (see, for example, \cite{FG94,EF96,FG96,Klenke1996,Winter2002,GHKK14}). We make~\eqref{e2197} rigorous as follows. Note that the spatial size of a cluster in the geometry of $\Omega_N$ is in comparison with $n$-balls of Haar measure $N^n$. To describe the \emph{cluster formation} in the $\Omega_N^L$-approximation, for a fixed but arbitrary $L \in \N_0$, we consider the configuration of $\CM_1(\K)$-elements in $k$-balls for $k \in \Z$ around some fixed $\eta \in \Omega_N^L \subseteq \Omega_N^\infty$, and in this $k$-ball with $N^{L+k}$ sites we form
\begin{equation}
\label{e2199}
\text{the \emph{normalised empirical measures:} $E_{t,k}^L$ of $\{\wh X_t^L(\eta^\prime,\cdot),\, 
\eta^\prime \in \CB_k(\eta)\}$},\,k \in \{-L,-L+1,\ldots\}
\end{equation}
normalised by the Haar measure $\mfh$ of the ball, i.e., by $N^k$. In this way we get a \emph{probability measure on $\K$} for every $\eta \in \Omega_N^L$, i.e., an empirical transition kernel $E^L_{t,k}$ from $\Omega_N^L$ to $\K$ for every $k \in \{-L,-L+1,\ldots\}$. By the ergodic theorem for $\wh X_t^L$, this object converges as $t\to \infty$ to $\delta_\mfa$ for some random $a \in \K$, provided we observe $\wh X_t^L$ only in $k$-balls with $k \in \N$. This means we can ask: for a site $\eta$, at what time $t_k$ is the $k$-ball around it \emph{covered} by one type? We say that \emph{at times $(t_k)_{k \in \N}$ there is a cluster of size $k$ around $\eta$}. We clarify the term ``covered'' rigorously by fitting it to what we need for the dual state, following the program in~\eqref{e2197} and below.

\medskip\noindent
\textbf{A problem.} 
The situation is actually more subtle. Namely, we cannot just work with the largest ball around $\eta \in \Omega_N^\infty$ that is \emph{completely covered} by $\delta_a$ in each \emph{geographic site}, since we have immigration in a site and emigration from a site, generating small bubbles of disturbances of the simple picture. This is incorporated when we are working, firstly, with \emph{empirical measures} for the description of the state in this ball, as introduced in \eqref{e2199} and, secondly, with \emph{$\varepsilon$-almost covered} balls as follows. We have to fix an $\varepsilon >0$, ask for type-$a$ frequency $\geq 1-\varepsilon$, and call this an \emph{$\varepsilon$-cluster of $a$}, since for $\varepsilon$ small enough there is only one such $a$.


\paragraph*{Cluster size description.}
To describe the cluster \emph{size} of the growing mono-type clusters as $t\to\infty$ at a fixed site $\eta \in \Omega_N^\infty$ (w.l.o.g.\ $\eta=0$), we choose the ball size $k=k(t)$ in such a way that, in the limit as $t \to \infty$, we see the very same $\delta$-measure limit and do not see this limit for any \emph{asymptotically larger} choice of $k(t)$. The resulting cluster in space has a \emph{random} volume (under the Haar measure $\mfh$) $N^{k(t)}$, and a corresponding \emph{random factor} of the volume arises by taking $k(t) + \ell$ for some random $\ell$. We want to describe the \emph{asymptotic size} of this random variable and its distribution by writing the function $k(\cdot)$ as (recall that $k(\cdot)$ is the logarithm with basis $N$ of the volume)
\begin{equation}
\label{e2204}
k(t) = f^L(t) + Z^L
\end{equation}
for some \emph{random variable} $Z^L$ with values in $\Z$ and some \emph{deterministic} real-valued scale function $t \to f^L(t)$. In the limit as $t \to \infty$, we seek $f^L(t)$ to be asymptotically equal to $f^L + g(t)$ for some \emph{deterministic} \emph{scale function} $g\colon\,[0,\infty) \to \Z$. We need an expression for the \emph{scale $g$} at which the clusters grow, in the sense that the volumes of the corresponding clusters are equal. This leaves us with the task to identify the triple
\begin{equation}
\label{e2427}
(f^L,\,g,\,Z^L)
\end{equation}
and later get limits $(f^\infty,g,Z^\infty)$ of these objects as $L \to \infty$, where $g$ is deterministic and $f^\infty$ and $Z^\infty$ are random variables with $f^\infty= g+Z^\infty$. 

The problem explained above means that in \eqref{e2427} we obtain $(f^{L,\varepsilon},g^\varepsilon,Z^{L,\varepsilon})$ for all $\varepsilon>0$. Afterwards we let $\varepsilon \downarrow 0$ and work with this limit. We want to identify some $(f^L,g,Z^L)$ as asymptotically equivalent in the limit as $t \to \infty$ for \emph{all small enough $\varepsilon$} and hence get a cover for all $\varepsilon>0$ in this limit. In this way $\Omega_N^\infty$ is decomposed into disjoint random balls labeled by a random $a \in \K$ giving the \emph{asymptotic cluster field} of balls, which are disjoint and are of a \emph{deterministic} asymptotic \emph{order $g$} (with a correction), giving the \emph{random size} called \emph{$Z^L$}. 

The \emph{existence} of these objects we verify below together with a refinement needed to get the existence.


\subsubsection{(III) Step 3.\ Cluster growth, part 1: cluster around a fixed site}
\label{ss.clgr1}

We can rephrase the above approach by reversing the role of time and size, picking $n$, looking for the time-scale we need to work with as $n \to \infty$, and choosing the time scale $(t_n)_{n \in \N}$ of observation times by requiring that 
\begin{equation}
\label{e2212}
k(t_n) = n.
\end{equation}
We can then use the perspective in \cite[Theorem 2]{Klenke1996} on cluster formation in systems of \emph{stable} interacting \emph{Fisher-Wright diffusions} on the hierarchical group with strongly recurrent migration, and employ the results for this special model (which is also in our model class) to identify the only possible \emph{candidate} for the triple in \eqref{e2427}. What we define here as a \emph{cluster} is a more restrictive object than that in~\cite{Klenke1996}.


\begin{example}
{\rm To get the proper size of the clusters, we consider \emph{two dual individuals at maximal distance} in a $k$-ball. These individuals have to meet by time $t_n$ with a probability tending to $1$ as $n\to\infty$. For our model this time is determined by explicit calculation for the \emph{stable case} with parameters $c_k=c^k$, $k \in \N_0$, without block-resampling. Namely, if $L=0$ and $c_k=c^k$ with $c<1$ and \emph{no} block-resampling taking place, then the sequence $(t_n)_{n\in\N}$ is given by (see~\cite{Klenke1996}) 
\begin{equation}
\label{e2222}
t_n=\left( \frac{N}{c} \right)^{n +1}, \text{ with } f^0 = -\frac{1}{\log (N/c)}, \quad g(t) = \frac{\log t}{\log(N/c)}, 
\quad Z^0 = E^0_{\infty,(\cdot)},
\end{equation}
with the argument $(\cdot)$ running over $\Z$.} \hfill$\spadesuit$
\end{example}


\paragraph*{Cluster growth.}

For the \emph{cluster growth} on countable spaces $\Omega_N^L, L \in \N$, we have the following result, which is the cluster theorem in~\cite{Klenke1996,EF96} for the particular case of stable migration. This result is too weak for our purposes; the necessary strengthening is provided below. For the notion of a sequence \emph{martingale measures} $(m_k)_{k \in \Z}$ we refer to \cite[Section 7.1]{D93}, which introduces a class of ``measure-valued martingales".

\begin{lemma}{\bf [Fast clustering for $\wh X^L$ in point-recurrent regime]}
\label{l.2386}
For the $L$-truncated version of the process $\wh X^\infty$ the cluster growth is as follows. For balls of size $n \in \N$ (= volume $N^n$) there exists a sequence $(t^L_n)_{n \in \N}$, determined in Lemma~\ref{l.2420} below, such that
\begin{equation}
\label{e2218}
\CL \left[\left(^n E^L_{t^L_n, n-m}\right)_{m \in \{-L,-L+1,\ldots\}}\right] 
\nto \CL \left[\left(^\infty E^L_{\infty,m}\right)_{m \in \Z}\right],
\end{equation}
where
\begin{equation}
\label{a2392}
\begin{aligned}
&\text{$(^\infty E^L_{\infty,m})_{m \in \Z}$ is the realisation of an \emph{entrance law} from time $\infty$ to time $-\infty$}\\
&\text{\emph{of a martingale measure on} $\K$ with $|m| \to \infty$ limits $E^L_{\infty,-\infty} = \delta_\mfa$ and $E^L_{\infty, \infty} = \theta$,}\\
&\text{where $\mfa$ is a $\K$-valued random variable distributed according to $\theta$.}
\end{aligned}
\end{equation}
The right-hand side in~\eqref{e2218} does \emph{not} depend on $L$. It arises as a functional of the process $\wh X^{\infty,\rm sta,\sim}$, the latter being the limit as $d^\ast \to \infty$ of the process arising as the continuum-space limit $\wh X_t^{\infty,\rm sta}$ of the system on $\Omega_N$, respectively, $\Omega_N^\infty$ of Theorem~\ref{T.CONVRES}. (The superscript $\rm sta$ stands for stable, and marks the specific choice of the parameters in the general $\wh X^\infty$.) This functional is obtained by taking all $\mfh$-averages over balls of sizes $k \in \N$ around $0$. 

The relations in ~\eqref{e2218} and \eqref{a2392} hold for $L=\infty$ as well, i.e., for the process $\wh X^\infty$ on $\Omega_N^\infty$ constructed in Theorem~\ref{T.CONVRES} via a rescaling limit.
\end{lemma}

The statement in \eqref{a2392} does not allow us to say whether the entrance law \emph{hits} $\delta_\mfa$ or only \emph{converges}. One might think that our concept of a cluster would be verified if we would have the \emph{actual hitting} and continuity of the limit, a bit of this is correct, but in the end this is however in fact more. Note, furthermore, that the formulation of the statement is not strong enough. What we want to capture is that we get large regions of almost one type as time progresses. These regions we call clusters of a given type and we want a deterministic scale in which we observe the random size of these objects. This is not yet contained in the above statement, which talks about the process observed around a \emph{fixed} site at a \emph{fixed} time. The required strengthening of~\eqref{e2218} and~\eqref{a2392} is given in Step 4.

Note that the last statement of the Lemma~\ref{l.2386} is per se \emph{not} allowing us to consider $n,L$ tending to infinity \emph{simultaneously}. However, if we can establish enough uniformity in $L$ as $n \to \infty$, then this will give us that the local behaviour of $\wh X^\infty$ is representing the large space-time scale behaviour of the model given by the $L$-truncation. The limiting $\Z$-indexed process is \emph{not} Markov, as this property holds \emph{only} in the limit as $N \to \infty$, which is the case treated in \cite{GHKK14}, where the \emph{transition kernel can be identified} as the equilibrium of the McKean-Vlasov or Fleming-Viot process with rate-$c_m$ and immigration-emigration from a source $\theta$, with $\varrho$ the state at $m-1$ associated with $(c_m,\theta)$. This equilibrium can be given explicitly via the \emph{Poisson-Dirichlet distribution} (see \cite{D93}, where this is carried out for the Fleming-Viot model with immigration-emigration for the choice $(c_m,\theta)$ of the parameters).


\paragraph*{Cluster size and time.}

The formula in \eqref{e2222} has to be replaced by a different sequence for our more general model, via the transition rates for $\zeta^\infty$. We have to identify the maximal radius $n$ of a ball such that the \emph{probability of coalescence} for two copies of $Y^\infty$ starting anywhere in that ball tends to $1$ as $t\to\infty$. Therefore we have to calculate the probability for two copies of $\zeta^\infty$ to hit $0$ by time $t$ starting from sites at distance $n$, and derive asymptotic formulas as $n\to \infty$. In \eqref{e2222} we have an explicit formula for the \emph{stable case}. For the general case we need to use the asymptotics as $t\to\infty$ of $G_t(\cdot,\cdot)$, the time-$t$ Green function on $\Omega_N^\infty$, and from that distill the scale function $g$. The Green function depends on $(c_k,\lambda_k)_{k \in \Z}$, but we only need its \emph{asymptotic behaviour} as $k \to \infty$ and show that this controls the cluster size. In other words, we need to understand how the cluster growth can be read off from the behaviour of $c_k, \lambda_k$ as $k \to \infty$.

\begin{lemma}{\bf [Cluster size and time on $\Omega_N,\Omega_N^L,\Omega_N^\infty$]}
\label{l.2420}
Both \eqref{e2218} and \eqref{a2392} hold in the case without block-resampling for the sequence of time points $(t_n)_{n\in\N}$ with asymptotics 
\begin{equation}
\label{e2321}
t_n \sim (N-1)\,N^{n-2} \sum_{\ell=0}^{n-1} \frac{1}{c_{\ell-1}}, \qquad n\to\infty,
\end{equation}
when $(c_k)_{k\in\N_0}$ is regularly varying at infinity. In the case with block-resampling, the formula for the effective difference random walk needs to be used, which leads to
\begin{equation}
t_n \sim (N-1) N^{n-2} \suml^{n-1}_{\ell=0} \dfrac{1}{\bar c_{\ell-1,N}}, \qquad n\to\infty,
\label{e2572}
\end{equation}
with $\bar c_{\ell-1,N} = 2c_{\ell-1} + N^{-1}(\lambda^*_\ell-\lambda_\ell)$ (recall \eqref{lambdakdefs} and see \eqref{1901pre}).
\end{lemma}

\begin{proof}
\label{pr2573}
We distinguish between without and with block-resampling.

\medskip\noindent
$\bullet$ \emph{Without block-resampling:}
The proof follows \cite{Klenke1996}. The key object is the \emph{recurrent potential kernel} $A$ of the process $\zeta$, defined as
\begin{equation}
\label{e2575}
A(n) = G(0,0) - G(\eta,0) \text{ when } k(\eta,0) = n, 
\end{equation}
which equals (see \eqref{e2553L})
\begin{equation}
\label{Greenapprox}
A(n) = \int_0^\infty \dd u \sum_{\ell=0}^{n-1} \frac{1}{N^\ell} \left[\eee^{-s_\ell u} - \eee^{-s_{\ell-1}u}\right]
= \sum_{\ell = 0}^{n-1} \frac{1}{N^\ell} \left(\frac{1}{s_\ell} - \frac{1}{s_{\ell-1}}\right)
= \sum_{\ell = 0}^{n-1} \frac{1}{N^\ell} \frac{s_{\ell-1} - s_\ell}{s_{\ell-1}s_\ell},
\end{equation}
where we define (see~\eqref{e2561})
\begin{equation}
\label{e2586}
s_\ell = \sum_{k \geq \ell} \frac{c_k}{N^k}. 
\end{equation}
According to \cite{Klenke1996}, in the \emph{stable} case the crossover time equals (without block-resampling)
\begin{equation}
\label{e2590}
t_n = N^nA(n).
\end{equation}  
But this formula also holds for the general class of $(c_k)_{k \in \N}$ we consider here, by our results on the transition kernel.

We want estimates for the right-hand side of~\eqref{e2590}. A rough estimate is $s_\ell \asymp c_\ell/N^\ell$, $\ell \to \infty$ (see \eqref{sapprox}), and hence $A(n) \asymp \frac{1}{N} \sum_{\ell=0}^{n-1} (1/c_{\ell-1})$, $n\to\infty$. By recurrence, the sum diverges as $n\to\infty$. Finer estimates are possible. For $(c_k)_{k\in\N_0}$ \emph{regularly varying at infinity}, we have $s_\ell \sim \frac{N}{N-1}\,c_\ell/N^\ell$, $\ell \to \infty$, and hence 
\begin{equation}
\label{e2594}
A(n) \sim \frac{N-1}{N^2} \sum_{\ell=0}^{n-1} \frac{1}{c_{\ell-1}}, \qquad n\to\infty.
\end{equation} 
For $(c_k)_{k\in\N_0}$ \emph{stable}, i.e., $c_k = c^k$ with $c<1$, we have $s_\ell = \frac{N}{N-c}(c/N)^\ell$, and hence we get
\begin{equation}
\label{e2598}
A(n) = \left(\frac{N-c}{N}\right)^2 \sum_{\ell=0}^{n-1} \frac{1}{c^\ell} 
\sim \left(\frac{N-c}{N}\right)^2\,\frac{c}{1-c}\, \frac{1}{c^n}, \qquad n \to \infty.
\end{equation}

\medskip\noindent
$\bullet$ \emph{With block-resampling:}
To incorporate the \emph{block-resampling}, we have to replace the $(c_k)_{k \in \Z}$-L\'evy process by the $(c_k,\lambda_k)_{k \in \Z}$-\emph{effective random walk}, in particular, the effective \emph{difference} random walk relevant for the above considerations. In addition, we have to analyse the cluster formation created by the block-resampling and the effect that this has on the \emph{cluster growth}. In particular, we have to verify whether we can \emph{reduce} the analysis to properties of the \emph{difference random walk of a single pair of individuals}. We return to this point in Steps 5--6 and show that this reduction is possible. This allows us to work out the change in the formulas for the growth properties, i.e., for $(t_n)_{n \in \N}$.

The formula in~\eqref{e2590} still holds, but a change in the argument is needed. This change concerns, in particular,~\eqref{e2594}-\eqref{e2598}, which have to be written out for $(c_k,\lambda_k)_{k \in \Z}$. In the transition rate we need to replace $c_k$ in~\eqref{e2321} by $\bar c_{k,N}$, $k \in \N_0$, in~\eqref{1901pre}, and use Section~\ref{sss.block}. This gives us~\eqref{e2572}, provided the Green function of the difference random walk is the correct quantity, which we will show in Steps 5--7. 
\end{proof}


\subsubsection{(III) Step 4.\ Cluster growth, part $2$: space-time cluster fields}
\label{ss.clgr2}

As explained above, we have to strengthen Lemma~\ref{l.2386}. We proceed in three substeps, two focussing on the new concepts and one on the main result.


\paragraph{Step 4.1: Space-time description of cluster formation}

$\mbox{}$\\[0.2cm]
So far we have analysed the growth of the cluster size around a specific but arbitrary \emph{single site} and at a \emph{fixed time}. We need a geographically more \emph{global} perspective on the clusters and on their growth viewed from varying systems of $n$-balls decomposing $\Omega_N^\infty$, and we need this for various large times $t>0$. In other words, we need to consider the asymptotics in an enriched form, namely, in \emph{its space-time global spatial form}.

We decompose $\Omega_N$ into \emph{disjoint} $n$-balls, so that mutual distances between sites in different balls are always $>n$, and each ball has $N^n$ sites. This means that we consider countably many \emph{spread-out} locations $\eta \in \Omega_N$ around each of which we consider growing spatial macroscopic windows, namely $n$-balls, where we choose $n$ and time suitably. In each ball we observe the restricted measure state, so that we obtain a new countable configuration of measures on $n$-balls. For $L \in \N \cup \{\infty\}$ this leads to the following sequence of limit fields:
\begin{equation}
\label{e2445}
\text{\emph{space-time field} $\left(^n E_t^L(\eta)\right)_{\eta \in \Omega_N^\infty, t \geq 0}$, \qquad $n \in \N$}.
\end{equation}
In this way we get a \emph{spatial field} indexed by a countable \emph{$n$-coarsening} of the set $\Omega_N$ (identifying sites within hierarchical distance $n$, indicated by $\uparrow$), which we call
\begin{equation}
\label{e2449}
\Omega_N^{\uparrow n},
\end{equation}
and which, in its structure as set and an Abelian group, is a copy of $\Omega_N$. These fields are \emph{continuous in $\eta$} and \emph{c\`adl\`ag in $t$} for $t>0$, and we can carry out the scaling analysis with $t \to t_n t$ and $n\to \infty$. After introducing all the new ingredients, we state in Lemma~\ref{l.2512} the necessary extended results. 

We consider this object in the limit as $t \to \infty$ and $n \to \infty$. This description of the cluster growth contains the dependence structure between \emph{different} and \emph{separated} large clusters, when viewed spatially, and goes beyond the work in \cite{Klenke1996}, but uses information in \cite{FG96}. This leads, for every $n \in \N$, to a \emph{field of measures} on $\K$, allowing us to describe the different clusters present in different large scale locations and at different large times in the infinite geographic space for each $n,t_n$, and in this way couple the limits $n \to \infty$ and $t\to \infty$. This yields the random fields (use $\dagger$ to indicate the field)
\begin{equation}
\label{e2453}
^nE_t^L= \left(^n E_{t,n-m}^{L,\dagger}(\eta)\right)_{\eta \in \Omega_N^{\uparrow n}} \qquad m \in \{0,\ldots,n\}, \,t>0.
\end{equation}
Here, each component describes the state within the $n$-balls of a \emph{disjoint decomposition} of $\Omega_N$ into $n$-balls. This leads to objects that are \emph{$\CM_1(\K)$-valued fields indexed by $\Omega_N$}. The key idea now is to embed all these fields on countable sets into the uncountable set $\Omega_N^\infty$ in such as way as to obtain fields indexed by the \emph{continuum-cardinality} space $\Omega_N^\infty$. For that purpose we relabel all $\Omega_N^{\uparrow n}$, $n \in \N$, by embedding them in $\Omega_N^\infty$ and afterwards shift them downward by the scaling parameter $n$. In this way we get, in the limit as $n \to \infty$, an \emph{$\Omega_N^\infty$-indexed field} of measures on $\K$. This in turn leads to an extension of the large clusters around a single site explained in~\eqref{e2218} to a new object describing the configurations of \emph{$k$-tuples of \emph{$n$}-clusters} in \emph{ensembles of large windows}, both on $\Omega_N$ and on $\Omega_N^L$. This in turn gives us fields of empirical measures in collections of expanding (with $g(t_n)$) windows \emph{indexed by $\Omega_N^\infty$}, that are \emph{constant in $(-n)$-balls}. 

If, to screen for \emph{different size} large clusters, we consider the objects just constructed for $n+\ell$ with $\ell \in \{-n,-n+1,\ldots\}$, then in the limit as $n \to \infty$ (with $t=t_n$) we obtain a hierarchy of objects that we can label by $\Omega_N^\infty$ in a natural way, shifting $\Omega_N^{\uparrow n}$ by $n$ downwards. In the limit as $n\to \infty$ we get \emph{$\Omega_N^\infty$-indexed spatial field of measures on $\K$}. This in turn gives us, for every $L \in \N$, \emph{sequences} in $n \in \N_0$ of \emph{$\Z$-indexed} \emph{processes} of \emph{fields of empirical measures} \emph{on $\K$}. If we include the time structure by adding $t$, then we get \emph{$\Z$-indexed processes} of the constructed space-time fields:
\begin{equation}\label{e2628}
\left(^n E^L_{t,n-m}(\eta), t>0, \eta \in \Omega_N^\infty\right)_{m \in \Z}.
\end{equation}
We look at this object as $t \to \infty$ and $n \to \infty$, where we see the pattern in the \emph{spatial} formation of clusters when we couple $t$ and $n$ by setting $t=t_n$, with $t_n$ the order of time in which the cluster grows. 

We turn now to the macroscopic \emph{time} structure of cluster formation, i.e., considering different values for $t$ exhibiting the \emph{time aspect} of cluster growth. Clustering occurs if individuals in an $n$-ball coalesce with a positive, but smaller than $1$, probability in the $n \to \infty$ limit. By considering, for $t$ in~\eqref{e2628}, $st_n$ for the macroscopic time of a process in $s>0$, we obtain the \emph{time structure of the cluster formation}, and for each $n$ a path of a Markov process converging as $n \to \infty$ to an entrance law from time $0$. Therefore Lemma~\ref{l.2386} has to be extended to \emph{macroscopic} space-time cluster growth, which is given in Lemma~\ref{l.2512} below. 


\paragraph{Step 4.2: Geometric view: from empirical measures to cluster configurations}
\label{wat}

$\mbox{}$\\[0.2cm]
So far we have focused on \emph{empirical measures} of sites in nested balls. Next, we want to shift to a \emph{geometric} point of view by observing \emph{clumps} or \emph{clusters of monotype sites} or almost monotype sites. Therefore we turn to the description of \emph{clusters of monotype regions} in space and their \emph{growth with time}. The latter in turn lead us to \emph{local properties of the $\Omega_N^\infty$-dynamics}, by representing $\wh X^\infty$ via the field of continuous and locally constant densities w.r.t.\ $\mfh$ and exhibiting \emph{macroscopic locally mono-type} regions.

More specifically, we consider grids of sites in space that are \emph{$n$-spread-out} at the \emph{cluster-growth scale given by $g$} at time $t_n$. For every $n \in \N$, this defines subgroups of elements at distance $g(t_n)$ by decomposing $\Omega_N$ into disjoint $k$-balls with $k=g(t_n)$. These balls can be decomposed into $(k+\ell)$-balls, with $\ell \in \Z$. Each such $(k+\ell)$-ball can for some $\varepsilon >0$ be an $\varepsilon$-cluster of a certain type $\mfa \in \K$, meaning that $1-\varepsilon$ of the mass of the ball is of type $\mfa$, or there is no such type, in which case the ball is labeled by $\varnothing$. We look for the largest constant $\ell$ for which we get such a cluster $(\mfa,\mfc)$. We continue with the rest of the $k$, until everything is labeled. This gives us a \emph{field of $(\mfa,\mfc)$ indexed by} $^{n\downarrow}\Omega_N$, which in the limit $n\to \infty$ will be $\Omega_N^\infty$, and which is describing what happens in the limit as $n\to \infty$ for the process on the geographic space given by the discrete space $\Omega_N$ or $\Omega_N^L$. For $\varepsilon \downarrow 0$, if a limit field exists, then it describes the cluster formation, and we call the corresponding random variable $^n \mfC^L$ when this limit exists:
\begin{align}
\label{e2567}
&\text{$^n \mfC^L$ is a $\Z$-indexed collection of \emph{$\K \cup \{\varnothing\}$-marked decompositions} into balls of various}\\ 
&\text{\emph{random sizes $k$} with $k\in \Z$. The \emph{mark} function $\Omega_N^\infty \to \K$ is 
called $\mfA \in L_1(\Omega_N^\infty,\K,\mfh)$.}\nonumber
\end{align}
It suffices that, as $n\to \infty$, the $\varnothing$-labeled balls vanish for every $\varepsilon >0$, and indeed we can take $\varepsilon \downarrow 0$ to get a countable decomposition of the limit state into monotype balls. 

More precisely, we can at stage $n$ embed the system of averages corresponding to the rescaling in the cluster analysis as a hierarchy of averages over $k$-balls with $k=0,1,\ldots,n,n+1,\ldots,n+\ell$, each represented as a field with $^k \Omega_N$ as index set, and afterwards shift the levels down by $n$, in particular, from level $n$ to level $0$. This gives us a version of $\Omega_N^n$ and, in the limit as $n\to \infty$, via the embedding of $\Omega_N^k$ in $\Omega_N^\infty$ for every $k \geq 2$, the set $\Omega_N^\infty$ as index set for the limiting object $^\infty \mfC$. This procedure also works when we move from $\Omega_N$ to $\Omega_N^L$, and look at $(^\infty \mfC^L)_{L \in \N}$ and its behaviour as $L \to \infty$.

We will have to introduce a topology on these families of marked balls. This topology should be compatible with the description via empirical measures given above. The \emph{topology} for the space of the cluster state, induced by the process with geographic space $\Omega_N$ or $\Omega_N^L$ in the limit $n\to \infty$, is based on encoding $^\infty \mfC^L$ via \emph{atomic measures} on $\Omega_N^\infty \times \K$. The atoms are given by the \emph{location} of the cluster and its \emph{type} and the $\mfh$-\emph{weight} of the ball. We equip the state space with the \emph{atomic weak topology of measures}. The latter reflects the idea that \emph{weight, location} and \emph{type} of a cluster converge to the same quantities of the limit object. Here, the location of the ball is its center, which is the site $\eta \in \Omega_N$ in the ball with \emph{zero} in the expansion w.r.t.\ $N$ at the $k$-th position when $k$ is the exponent in the size of the ball.

Now, recalling the paragraph below~\eqref{e2427}, we also need a topology for $\varepsilon$-clusters when $n$ is \emph{finite}. Here we proceed as above, assigning an atomic measure with an atom representing at least a fraction $1-\varepsilon$ of the total local mass in the ball in one site, as above, but we keep the restricted measure on the exceptional part of the balls, so that these must converge weakly. We take the topology known as of the \emph{weak atomic topology} of measures on a Polish space, the space of all not necessarily atomic measures on $\Omega_N^\infty \times \K$. This way we obtain a \emph{Polish state space} for the random clusters $\mfC$.

The limiting object as $n\to \infty$, when it exists, is a random variable $^\infty \mfC_\varepsilon$ taking values in the \emph{$\{\K \cup \varnothing \}$-marked disjoint decomposition of $\Omega_N^\infty$} in asymptotic $\varepsilon$-mono-type clusters, i.e., balls of a certain radius $k$ and a certain location, with the site $\eta \in \Omega_N^\infty$ in that ball. This has the form of a subset of $\{(\eta,k,a) \in \Omega_N^\infty \times \Z \times \K\}$ when we have a cluster size of asymptotic order $g$. For given $t$, the constants $k \in \Z$ in the spreading parameters correspond to the \emph{macroscopic sizes} on the $g$-scale of the different \emph{clusters}, with $(\eta,a)$ corresponding to the different \emph{types} and the \emph{location} of their corresponding typed geographic \emph{clusters} on $\Omega_N$. In the limit as $\varepsilon \downarrow 0$, when it exists, we obtain the asymptotic geometric description of the cluster growth of $\wh X$ on $\Omega_N^L$ on the random variable $^\infty \mfC^L$. The latter object contains a transition kernel $\mfA$ from $(\Omega_N^\infty,\mfh)$ to $\K$, where the actual values are given by a random set $A^L$.


\paragraph{Step 4.3: Results for cluster growth}

$\mbox{}$\\[0.2cm]
For the present extended description of \emph{cluster growth} on $\Omega_N^L$, introduced as key object in Steps 4.1--4.2, we have the following result, which we will prove in Steps 6--7, but will \emph{adapt} and then use in Step 5 to get the properties of the process $\wh X^\infty$ on $\Omega_N^\infty$.

\begin{lemma}{\bf [Fast clustering: spatially global cluster growth for a space field on $\Omega_N^L$]}
\label{l.2512}
In Lemma~\ref{l.2386} replace $E^L_{t,n-m}$ by $^n E^{L,\dagger}_{t,s,n-m}$ with $t=st_n$ for $s>0$, and modify $\delta_\mfa$ by $\delta_{(\mfA_\eta)_{\eta \in \Omega_N^\infty}}$, with $\mfA=(\mfA_\eta)_{\eta \in \Omega_N^\infty}$ a $\K^{\otimes \Omega_N^\infty}$-valued field, and modify $\theta$ by $\theta^{\otimes\Omega_N^{\uparrow \infty}}$. Then
\begin{itemize}
\item The claims in \eqref{e2218} and \eqref{a2392} hold.
\item The associated cluster filed $^n \mfC^L$ converges as $n\to \infty$ to the cluster field $^\infty \mfC^L$ without the label $\varnothing$ a.s.
\item The law of the latter does not depend on $L$.
\end{itemize}
\end{lemma}

\noindent
We will strengthen this statement further in Step 5, giving more specifics on the properties of the limiting field $\mfA$ (see Lemma~\ref{l.2423}).

We still have to address the problem described below~\eqref{e2427} and show that we have really obtained the \emph{asymptotic cluster size} from a geometric point of view with our choice for $g^L$ describing the order of magnitude. The point is that the fraction of sites that do not have the highest frequency type as $t \to \infty$ tends to zero as $m \to -\infty$, and hence the fraction of these sites is $o(N^n)$, which is referred to as the regime of \emph{fast growing clusters}, which allows us to take $\varepsilon \downarrow 0$, as $n \to \infty$. (This is the second bullet in Lemma~\ref{l.2512}.) 

In the proof we have to relate the objects $^n \mfC^L, ^\infty \mfC^\infty$ with the coarser description via $^nE^{L,\infty},E^\infty$, where we work without explicitly using the \emph{atomic} measure structure and the \emph{geographic compact support structure} of the densities of the empirical measures w.r.t.\ $\mfh$ to describe the states. Here, we have to use the weak atomic topology, which is stronger than the weak topology. Hence an additional \emph{tightness argument} is needed.


\paragraph*{Clusters.}

Note that on $\Omega_N^L$ for fixed $L$ and $t \to \infty$ we obtain \emph{limiting space-time clusters}. But for finite $t$ we need not be exactly at the boundary, i.e., be mono-type (the ``pure" state arises only in the limit where the \emph{locally constant} state appears). This is leading to an ``asymptotic" decomposition of $\Omega_N^L$ into $\mfh$-equivalence classes of balls with a $\mfh$-null set of exceptional sites after we have scaled space and letting time tend to infinity, provided these \emph{macroscopic limit clusters} occur already from some random $m \in L_1 (\mfh)$ downwards all the way down to $-\infty$. This reflects the \emph{picture of growing clusters covering the space} as we evolve with $t$, with for each cluster a \emph{random size} for the given \emph{deterministic scale} denoted with $g$. Thus, we want that \emph{the limit field $E^L_{\infty,m^\prime}$ in a site $\eta$ is $\mfh$-mono-type} in an $m^\prime$-ball around $\eta$ for $\mfh$-a.e.\ $\eta \in \Omega_N^\infty (\eta)$ for $m^\prime \leq m$ and some random function $m \in L_1(\mfh)$. Of course, the levels $0,1,\ldots,k$ have this property up to zero weight, but the question is whether it holds all the way up to some $m_n$ that is asymptotic to $n+m$ as $n \to \infty$ for some $m$, or whether the volume of the cluster could be $o(N^n)$ for $n \to \infty$, or not reach some volume comparable to $N^n$. This is \emph{excluded} when we have the property called \emph{fast growing clusters}. This scenario presents some technical problems, since the field lives on a continuum space: we cannot directly speak of values at a site and must work with evaluations against test functions instead. What we must prove is the following multi-dimensional version in time on the process level of the previous result. 

\begin{lemma}{\bf [Fast clustering on $\Omega_N$ and limit properties in global space-time field]}
\label{cor.2748}
In the regime where $\wh X^\infty$ is point-recurrent, the $L$-truncations on $\Omega_N^L$ have fast growing clusters, i.e.,
\begin{equation}
\label{e2750}
\left\lbrace \left(^\infty E^L_{\infty,s,m^\prime}, (\eta)\right)_{s>0}, \, m^\prime \in \Z \right\rbrace
\end{equation}
is continuous (or c\`adl\`ag in the presence of block-resampling) and has a version of the density w.r.t.\ $\mfh$ that is continuous in $\eta$ and $\{0,1\}$-valued. This limiting object can be identified as the field of empirical measures $E$ obtained from $(\wh X_s^{\infty,\rm sta,\sim})_{s>0}$ for all $\wh X^\infty$ in the domain of attraction of $\wh X_s^{\infty,\rm sta}$ under the continuum-space scaling limit according to Theorem~\ref{T.CONVRES}, with all the details of the scalings being specified in Section~\ref{ss.resc}.
\end{lemma}

\noindent
We give the proof of this lemma in Step 5.

Thus, we need to show that the \emph{fast growing clusters} corresponds to the \emph{point-recurrent regime} on $\Omega_N^\infty$, keeping in mind the discussion of the problem below \eqref{e2204}.

\begin{lemma}{\bf [Fast growing clusters in point recurrent regime]}
\label{l.2421}
Consider $\wh X^\infty$ \emph{without} block-resampling and its $L$-truncation $\wh X^L$. The regime of fast growing clusters for each $\wh X^L$ corresponds to point-recurrent migration associated with $\wh X^\infty$.
\end{lemma}

\noindent
This property has been derived in \cite{GHKK14} in the $N \to \infty$ limit. Here, we need to use the potential theory of $\zeta^\infty$ to obtain that we have fast growing clusters \emph{exactly} in the \emph{point-recurrent} regime for $\zeta^\infty$. For the strongly recurrent migration this is known from~\cite{Klenke1996} for the case of \emph{stable} migration with parameter $c<1$. We have to deal here with the case of a strongly recurrent migration that is not necessarily stable. We give the proof in the \emph{general migration} case in Steps 6--7. For the critically recurrent regime and in the stable case with $c=1$, this is known from \cite{Klenke1996} and clusters are \emph{not} fast growing. This means that we need to show that those parts of space that are \emph{not} contained in large clusters have volume $o(N^n)$ in the point-recurrent case, via the results on the empirical measures. In the Harris-recurrent but not point-recurrent regime the clusters have volume $o(N^n)$). Thus, we have to generalise the result in \cite{Klenke1996} for the stable case to the whole Harris-recurrent regime, and show that clusters are $o(n)$ for other recurrent sequences $(c_k)_{k \in \N}$. This is addressed in the section on Harris-recurrent but not point-recurrent migration
 

\paragraph*{Block-resampling.}
 
Next, we have to allow for the presence of \emph{block-resampling}. Here, as worked out in Steps 6--7, the point is to show that the asymptotics of finite coalescents starting in $n$-spread-out sites is carrying over to the case with added block-resampling and reshuffling. The reshuffling is easily taken care of, since we are using the \emph{effective} migration. The important point is that \emph{block}-resampling also does not change the qualitative picture: in a sequence of coalescent events, the final step involving the last pair of partition elements only involves pairs with adapted rates (summing over all rates for $k$-blocks with $k \geq n$). Then, indeed, Lemma~\ref{l.2512} has the following corollary.

\begin{corollary}{\bf [Global cluster growth with block-resampling]}
\label{cor.2604}
Consider $\wh X^\infty$ and its $L$ truncation $\wh X^L$. In the case of \emph{point-recurrent} effective migration underlying $\wh X^\infty$, Lemmas~\ref{l.2512}-~\ref{l.2421} continue to hold in the presence of block-resampling.
\end{corollary}

\begin{proof}
\label{pr.2797}
We present the list of what is needed later.

First, we need to show that the behaviour in the point-recurrent regime continuous to hold with an additional block-resampling mechanism in the effective migration. Second, we have to show that block-resampling leads to \emph{fast cluster growth} for recurrent but not point-recurrent migration. The latter is treated in the section devoted to Harris recurrent but not point-recurrent
migration.

In transferring the argument to the case with block-resampling, the key is to work with the \emph{dual representation} of the state of $\wh X^L$ embedded in $\Omega_N^\infty$ via the strong duality together with the entrance law. In the dual entrance law we get an extra term inducing the creation of an additional uniformly distributed piece of weight on the ball of change to one partition element, and a reduction of all other elements by a fraction of their elements.

We first examine whether the behaviour of specific finite coalescents (those which are spread-out initially on growing boxes of radius $n$ with $n \to \infty$) carries over to the presence of block-resampling. The point is that we still want to establish the \emph{fast cluster growth}. We now have a \emph{jump} process component in the forward process, instead of a diffusion. On the other hand, the rate of the jumps taking place in a $j$-ball, does depend \emph{only on $j$} and \emph{not} on the state of the population, which is different from the effect of a diffusion. We can analyse the block-resampling first without the single-site Fleming-Viot or $\Lambda$-Cannings part, and then take on the combination. The details of the proof of Lemma~\ref{l.2386} and Lemma~\ref{l.2512} will be given in Steps 6--7.


\subsubsection{(III) Step 5.\ Cluster formation and application to local properties}

Corollary~\ref{cor.2604} claims that the $L$-truncations of $\wh X^\infty$ each have \emph{fast growing clusters} in the point-recurrent regime, and \emph{not} just in the Harris-recurrent regime, when we do \emph{not} have block-resampling. The proof of the results on the properties of the state of $\wh X^\infty$ seem to require us to answer the question whether in \eqref{a2392} the limit process $(^\infty E^L_{\infty,m})_{m \in \Z}$ (recall \eqref{e2199}) only converges a.s.\ to a state $\delta_\mfa$, or whether it \emph{hits} $\delta_\mfa$ a.s.\ for some \emph{random $m$}, which would be the \emph{asymptotic random cluster size} in the \emph{deterministic scale $g^L$}. As we pointed out, the formulation is naive and also not correct, since we condition on observing a \emph{whole field} on a continuum set at a \emph{single site}, and also too strong to get the property of a cluster scale. The true picture we need for our purposes is more subtle and developed now.


\paragraph*{Field of empirical measures.}

We have to think about the $\Omega_N^\infty$-indexed \emph{field} of empirical measures and the random level function specifying the ball of a single type as an equivalence class of a function $m(\cdot)\colon\,\Omega_N^\infty \to \Z$, such that in the ball $B_{m(\eta)}$ the type is $\mfh$-a.e.\ constant equal to some $a \in A \subseteq \K$ for $\mfh$-a.e.\ $\eta \in \Omega_N^k$ with $A$ a \emph{countable} random set. Furthermore, $m(\cdot)$ is locally $\mfh$-constant and $m(\cdot)$ is continuous $\mfh$-a.e. The latter scenario would justify the term \emph{cluster scale} for $g^L$ and \emph{random cluster size} for $m(\cdot)$. 

The scenario gives us a limiting object \emph{$m(\cdot) \in L_1(h,\K)$ for the cluster growth on $\Omega_N^L$}, which gives a \emph{limit field} $\wh X^{\infty,\sim}$ with having similar \emph{local properties} as we have claimed for our $\wh X^\infty$. This limit process arises as the rescaling limit of a process with a special choice of parameters, namely, the limit as $d^\ast \to \infty$. This allows us to get the properties that we have claimed for the process on $\Omega_N^\infty$ via the \emph{$L$-truncation approximation} of the clustering, accounting for the fact that we do no shift by $L$, but add something at level$-L$, as we explain below.

With the four lemmas from Steps 2--4 in hand, we next argue that the results on the \emph{cluster formation} for $\wh X_t^L$, $L \in \N_0$, on $\Omega_N^L$, which for $t \to \infty$ generate a limit field $\wh X^{\infty,\rm sta,\sim}$, carry over to $\wh X^\infty_t$, for which we claim the \emph{decomposition} of $\Omega_N^\infty$ into \emph{mono-type balls} for the state of $\wh X_t^\infty$ for $t>0$. To that end, we need the sharper \emph{macroscopic space-time} version of the cluster growth analysis in Lemma~\ref{l.2512}, which we gave in Lemma~\ref{cor.2748}-~\ref{l.2421} and which provides a sharper picture of the cluster growth and, indeed, in the \emph{$t\to \infty$ limit} gives the formation of real mono-type clusters in a suitable sense applicable in Lemma~\ref{l.2423} below. Below we show that this implies the local properties of $\wh X^\infty$ that we formulated in the beginning of the paper. At the same time it gives a \emph{new} result on the cluster growth on $\Omega_N$, extending and refining~\cite{Klenke1996,EF96}. 


\paragraph*{Cluster growth versus truncation limit: two different limits.}

We have to relate the results on the cluster growth limit for $\Omega_N^L$ and $t \to \infty$ to the truncation limit $L \to \infty$ at strictly positive times $t$. To that end, we first recall that, by the convergence results we have shown above, namely, $\wh X^L \to \wh X^\infty$, we get, for every $n \in \N$,
\begin{equation}
\label{e2461}
\left(\left(^n E^{L,\dagger}_{t,n-m}\right)_{m \in \{-L,-L+1,\ldots\}}\right)_{n \in \N} \text{ converges to } 
\left(\left(^n E^{\infty,\dagger}_{t,n-m}\right)_{m \in \Z}\right)_{n \in \N} \text{ as } L \to \infty,
\end{equation}
where the right-hand side is the functional of $\wh X_t^\infty$ that we get when we apply the functionals of the cluster analysis that are given via the empirical measures to $\wh X^\infty$ at time $t$. If we put $t=t_n$ and let $n \to \infty$, then we get $\wt X^{\infty,\rm sta, \sim}$.

Furthermore, the system on $\Omega_N^L$ behaves for the \emph{cluster growth} as $t \to \infty$ like a system on $\Omega_N$, where \emph{time is accelerated by $A(n) N^n$}, whereas for the truncation limit $L \to \infty$ we speed up resampling by $N^L$ only. The space for the cluster growth analysis is the \emph{space scaled down}, because cluster sizes are \emph{shrinked by $n$}, i.e., volumes in geographic space are scaled down by $N^{-n}$. The parameters $(\bar c_k,\bar \lambda_k)_{k \in \N}$ and $\bar d,\bar \Lambda$ arise as $\bar c_k=c_{k-n}, \bar \lambda_k=\lambda_{k-n}$, $\bar d=d,\bar \Lambda=\Lambda_{-\infty}$. On $\Omega_N^L$ our cluster growth analysis in Steps 2--4 with $t_n^L=t_n$ gives:
\begin{equation}
\label{e2466}
\text{$(^n E^{L,\dagger}_{t_m, n-m})_{m \in \{-n,-n+1,\ldots\}}$ converge as $n\to \infty$ to a limiting 
object $(^\infty E^{L,\dagger}_{\infty,m})_{m \in \Z}$.}
\end{equation}
The right-hand side does not depend on $L$, since the parameters agree for every index larger than $L$. 

The limit in the right-hand side arises as a functional of the $d^\ast \to \infty$ limit of the \emph{continuum space limit}, the one called $\wh X^{\infty,\rm sta}$ and treated in Theorem~\ref{T.CONVRES} on $\Omega_N$, respectively, $\Omega_N^L$, which we denoted by $\wh X^{\infty,\rm sta,\sim}$ associated with \emph{instantaneously coalescing} L\'evy-processes. We therefore see that, for the application of the \emph{cluster growth analysis} on $\Omega_N^L$ to the analysis of local properties of the $\Omega_N^\infty$-process $\wh X^\infty$ an ingredient is still missing. 

We are now dealing with the situation of an \emph{$(n,L)$-indexed array} of empirical measure processes (see~\eqref{e2628}), where $n$ and $L$ tend to $(\infty,\infty)$ in a very specific way, since we are interested in times $t N^L$ for $t>0$ and $L\to \infty$, instead of $L$ fixed. We had an observation time $t_n \to \infty$ as $n\to \infty$ in the cluster growth analysis on $\Omega_N$, where $t_n = A(n)N^n$. Furthermore, our \emph{active level} for the population is \emph{not} level \emph{$0$ shifted to $-n$}, as in the cluster formation analysis. Instead, with a level \emph{$-L \to -\infty$} added, with \emph{rates of migration} being \emph{preserved} and with the \emph{resampling rate up-scaled} at the ground level, the reference level of the \emph{observation of the state is $0$} in both cases but only after the shift in the cluster formation analysis. Since we have a measure state on an $(-L)$-ball, respectively, $\Omega_N$, it therefore is an $n$-average of active components with $n\to\infty$. We match this here by putting $n=L$.

Since the former role of $n$ is now played by $L$ and the active level is $-L$, we get that the resampling rates at the active level at $-L$ is asymptotically $d^\ast N^L$ for the Fleming-Viot resampling, rather than the speeding up of time in the clustering regime corresponding to $d^\ast=\infty$ and $\mu_{-L} N^{L}$ for the $\Lambda$-block-resampling and reshuffling and $N^L c_{-L}$ for migration. Hence we have the new \emph{cluster observation time replacing} ``$t_n$" of the \emph{cluster growth analysis} is, for the \emph{truncation limit}, \emph{replaced} by $tN^L$. Migration rates at the \emph{observation level} remain \emph{fixed}, and hence are $O(1)$ rather than approaching limiting rates tending to $0$ and compensated by the factor $A(n)$. In particular, they both converge, or in fact are equal, to the prescribed rates for $\wh X^\infty$, rather than those of the continuum space rescaling limit $\wh X$ on $\Omega_N$, denoted by $\wh X^{\ast,\rm sta,\sim}$. 

For $t>0$ this produces an entrance law that equals $\wh X^{\ast,\infty}$, as we proved for a \emph{locally clustered} state around the \emph{$0$-level at $O(1)$-times $t$}. The situation we have here is a backward view from action level-$(-\infty)$ to level-$(0)$ as the level of observation, instead of the situation in the cluster growth analysis on $\Omega_N,\Omega_N^L$ with the action level $0$ shifted to $-n$, the observation level eventually $\infty$, and shifting by $-n$ yields a \emph{different} limit object. Namely, we get $\wh X^\infty$, whereas the cluster growth analysis on $\Omega_N$ is asymptotically described by $X^{\infty,\rm sta, \sim}$, the $d^\ast \to $ limit of the former. Nevertheless, a key point is that some features like \emph{monotype clusters} are similar, and this is crucial, \emph{so that the proof can be built via cluster-growth results} already existing.


\paragraph*{Local cluster formation.}

For verification, consider the \emph{empirical measure} $E$, which again contains a \emph{macroscopic time variable $t$}, so that $(\infty,t)$ instead of $\infty$ appears, representing the fact that now $t^L_n$ is generated by a rate $\sim N^L$ in regular cases and is observed at a time points $t\,t_n^L$ in the rescaling limit, where $n=L$ and $L \to \infty$, so that the running times of the process give a path in $s$, where $tN^L+s$ is the \emph{current} time of observation. Hence $^n E_t^{L,\dagger}$ becomes $^n E^{L,\dagger}_{t,s,m}$ with $t \to t\, t_n$, $s \in \R$, and $n=L$. This is a process in the macroscopic time $t$, and for fixed $t$ in $s$. 

We know that, exactly in the regime of space-time scaling described above, we get as a limiting field states that are in law equal to $\wh X^\infty_t$ \emph{and} are locally monotype, whereas on a smaller time scale we get deterministic state and on a larger time scale we get one cluster. So, we first need to prove that exactly the proposed coupling of $(n,L)$ as $n,L \to \infty$ (obtained by putting $n=L$, with the reference level of observation fixed at $0$, and $-L$ the level of activity) gives us a \emph{limit}, namely, the \emph{state of $\wh X_t^\infty$ for $t>0$}. Recall that, for the cluster growth analysis for $t \to \infty$ with \emph{fixed} $L$, we get a \emph{different} process, namely, the continuum space limit of the system on $\Omega_N$ taken to the limit $d^\ast \to \infty$ called $\wh X^{\infty,\rm sta,\sim}$ (see~\eqref{ag39}).

The limit $E^{\infty,\dagger}_{\infty,t,m}$ has a version of a \emph{density} w.r.t.\ $\mfh$, which is $0$ or $1$ and is \emph{continuous} as a function of $(t,\eta)$ except at block-resampling jumps. Hence, for some \emph{random negative level $m$ we have}, for $\mfh$-a.e.\ all sites $\eta \in \Omega_N^\infty$, the same types based on the properties of $E_{\infty,t,\cdot}^{\infty,\dagger}$. Furthermore, we want to show that the latter scenario has a corresponding structure of limiting clusters that, for the $L \to \infty$, is a limit $^{L} \mfC^{L}$ specified by $^\infty \mfC^\infty$. 

\begin{lemma}{\bf [Local cluster formation of $\wh X^\infty$]}
\label{l.2423}
Consider $\wh X^\infty$ in the point-recurrent regime, together with the corresponding $L$-truncations on $\Omega_N^L$ all embedded in $\Omega_N^\infty$.
\begin{itemize}
\item[(a)] 
The process $E=(^{L} E^{L,\dagger}_{\infty,\cdot})$ has the property
\begin{equation}
\label{e2426}
\left(^{L} E^{L,\dagger}_{\infty,\cdot}\right) \Lto \left(^{\infty} E^{\infty,\dagger}_{\infty,\cdot}\right) \text{ in path space.}
\end{equation}
\item[(b)] 
$E^\infty=^{\infty} E^{\infty,\dagger}_{\infty,\cdot}$, when restricted to a fixed type $\mfa \in \K$, has a density w.r.t.\ $\mfh$ that for every $t>0$ and $\eta \in \Omega_N^\infty$ takes values in $\{0,1\}$ and is \emph{continuous} in these variables $(t,\eta)$, except for $t$ at the block-resampling jump times.
\item[(c)] 
The field $^{L} \mfC^L$ converges as $L \to \infty$ to $^{\infty} \mfC^\infty$. The latter equals the decomposition of the support of $\wh X^\infty$ into macroscopic $\mfh$-mono-type balls covering $\mfh$-a.e.\ sites in $\Omega_N^\infty$.
\end{itemize}
\end{lemma}

We can also say something about the \emph{global} cluster formation for $\wh X^\infty$ for $t \to \infty$ in the strongly recurrent case. In that case another process on $\Omega_N^\infty$ appears, namely, $\wh X^{\infty,\rm sta,\sim}$.

\begin{remark}{\bf [Global cluster formation for $\wh X^\infty$ via $\wh X^{\infty,\rm sta,\sim}$]}
\label{r.2876}
{\rm There is a cluster \emph{growth analysis of} $\wh X_t^\infty$ on $\Omega_N^\infty$, analogous to the one for $X$ or $X^L$ on $\Omega_N$ or $\Omega_N^L$. Namely, for $t\to \infty$, follow the same principle applied to the density field corresponding to $\wh X^\infty$ at every given time w.r.t.\ a decomposition of $\Omega_N^\infty$ into $\mfh$-monotype $k$-balls for $k \in \Z$ around $\mfh$-a.e.\ $\eta \in \Omega_N^\infty$. Then scale time and space, pass to $n$ averages and shift by $-n$ down, and take the limit $t_n \to \infty$. This limit has a law specific to the universality classes of the $\wh X^\infty$. Namely, take the \emph{continuum space rescaling limit} $\wh X^{\infty,\rm sta}$ of $\wh X^\infty$ in our class of processes on $\Omega_N^\infty$ (see Theorem~\ref{T.CONVRES}). Let $d^\ast \to \infty$ to get the functional of $\wh X^{\infty,\rm sta,\sim}$ as \emph{\emph{limiting cluster growth description}} of $\wh X^\infty$ for $t \to \infty$.}
\end{remark}

In Steps 6 and Step 7 we verify Lemma~\ref{l.2386} up to Lemma~\ref{l.2423}.


\subsubsection{(III) Step 6.\ Proof of four key lemmas}

We give the proof of Lemmas~\ref{l.2386},~\ref{l.2512},~\ref{cor.2748} and~\ref{l.2423}. We are in the \emph{point-recurrent} regime and have to prove the convergence in~\eqref{e2218} and carry out \emph{its application to $\wh X^\infty$, leading to the extension \eqref{e2426}}, and verify that the limit field in both lemmas has locally \emph{mono-type} states in a proper sense. The latter requires that we show that the state is locally mono-type (i.e., at a site) and that this type prevails in a suitable small ball. We start with the \emph{convergence} arguments for the two lemmas, and postpone the proof of the local mono-type property, which is in fact a \emph{continuity property}, to Step 7.

The convergence of the $n$-indexed \emph{sequence of empirical measures $E$} in~\eqref{e2218}, and later the $L$-indexed sequences in~\eqref{e2426}, is proved by using the \emph{method of moment measures}, namely, by showing that the \emph{space-time moment measures} of these empirical measures converge to the space-time moment measure of the claimed limit state as $L \to \infty$. For that purpose we calculate the moment measures for the respective situations, in the \emph{two groups} of lemmas.

\paragraph*{(I)}
We focus on the first group of lemmas and explain the tasks to be completed, before diving into the execution of the proofs.

\paragraph*{(i)} 
We begin with Lemma~\ref{l.2386} on the \emph{cluster formation} of $\wh X^\ast,\wh X^{L,\ast}$ on $\Omega_N$, respectively, on $\Omega_N^L$ viewed around a fixed site. The moment measures can be calculated by using the \emph{duality relation}, respectively, \emph{space-time duality} for $\wh X^L$ and $\wh X^\infty$, in terms of the \emph{initially $n$-spread-out spatial coalescent }$\wh Y^\ast$, respectively, $\wh Y^{\ast,L}$ as $n \to \infty$, to some limiting coalescent (to be specified), called $\wh Y^{\ast,\dagger}$, which is different from $Y^{\ast,\infty}$ for $\wh X^\infty$. Let $t$ be the macroscopic time horizon for the duality relation and $s$ the running backward time (from $0$ to $t$). As explained earlier, we have to show convergence of finite spatial coalescents $Y^L$ \emph{embedded in $\Omega_N^\infty$} to the claimed limit coalescent $Y^{\ast,\dagger}$. 

Note that we have to deal with a \emph{new situation}, namely \emph{initial states} that are \emph{$n$-dependent} and start in $n$-spread-out locations. A new process arises, since we observe at a \emph{late time} $t_n$ the measure state (on $\K$) in a \emph{large ball} of volume $N^n$, such that \emph{typical} $k$-tuples of dual individuals in such a ball are at distance of order $n$, i.e., they are \emph{$n$-spread-out}, and therefore, on the large time scale where they meet first, subsequent re-hitting is on smaller time scales and so ``not observable" on the large scale. First, we determine the \emph{limiting dual process} $Y^{\ast,\dagger}$ on $\Omega_N^\infty$, which --- unlike $\wh Y^{\ast,\infty}$ --- consists of \emph{instantaneously} coalescing L\'evy-processes on $\Omega_N^\infty$. Second, we show \emph{convergence} of the appropriate dual processes on $\Omega_N^L$ (embedded in $\Omega_N^\infty$) to the dual process on $\Omega_N^\infty$. 

If we focus on clusters around one site, then in the limit as $N \to \infty$ the process $Y^{\ast,\dagger}$ looks like a \emph{coalescent with emigration} on a site plus a cemetery, which also for fixed $N$ is a good caricature. We return to this later. 

\paragraph*{(ii)} 
A similar strategy applies for Lemma~\ref{l.2512}--\ref{cor.2748}, by looking at clusters spatially globally, respectively, via the space-time process, where we get \emph{instantaneously} coalescing L\'evy-processes, respectively, space-time processes.  

\paragraph*{(II)}
Next turn to Lemma~\ref{l.2423} for which there are two tasks similarly as above. Now, we have $L$-indexed sequences and $L \to \infty$, with $L$-spread-out initial states that now really are a \emph{$(-L)$-shrinking} of the distances and an $N^{-L}$-shrinking of the $\mfh$-measure, although with an $N^L$-expansion of the number of sites in a ball (by the definitions on level $-L$). The $Y^{\ast,L}$ converge to a different spatial coalescent on $\Omega_N^\infty$, namely, to $Y^{\ast,\infty}$. To identify the latter limit coalescent, we argue in the appendix that a pair of coalescing random walks on $\Omega_N^L$ embedded in $\Omega_N^\infty$ converges to a pair of \emph{delayed coalescing L\'evy-processes}, and coalescence occurs at rate equal to $d^\ast$ \emph{collision local time} (also in our more general setting). We have to prove this for $k$-tuples with $k>2$, where the coalescence appears at a rate that arises as the sum over the rates of all ${k \choose 2}$ pairs. We need that in the spread-out case all these attempts occur independently. This fact has been shown already in~\cite{Klenke1996}.

The proof proceeds as outlined above, now in \emph{two parts} for each of the \emph{two groups of lemmas} separately.

\begin{proof}{\bf [(I) Proof of Lemma~\ref{l.2386}]}
\label{pr.2803}
We treat the three cases (without block coalescence, only block coalescence and the combination) subsequently in (i)-(iii). 

\medskip\noindent
\textbf{(i)} Consider the case \emph{without} block-resampling. We want to study the ball averages in increasing balls at large times and establish their asymptotic behaviour. It arises as a rescaling limit on a copy of an $n$-ball of $\Omega_N$, where time and space are rescaled with $n$, where $n \to\infty$ and we use the duality relation to verify that the moment measures converge. Therefore we need to show that the finite coalescent with $n$-spread-out initial locations on $\Omega_N$ in some $n$-ball, where pairs may meet in an $n$-ball before leaving it and are observed at the \emph{large time $t_n$} as $n \to \infty$, converge to a limit process after rescaling time to a coalescent with \emph{instantaneous coalescence} upon hitting or to \emph{escape}. More precisely, we show that these converge for $n \to \infty$ to a \emph{Kingman-coalescent} at a single site with migration to a cemetery (\emph{emigration}). 

The latter amounts to showing that asymptotically no \emph{multiple} coalescence occurs and every pair acts independently. Here we can use results on coalescents from~\cite{Klenke1996} and extend these results to the more general model considered in our paper. Recall Theorem~\ref{T.CONVRES}, Remark~\ref{r.634}, Remark~\ref{rem.1968} and the explanation of the ingredients. The argument in the general case is built on a similar scaling result for spatial coalescents. 

The proof follows the steps in~\cite{Klenke1996}, replacing the \emph{explicit calculations} therein for the case $c_k=c^k$, $c < 1$, ($c$ is called $c^{-1}$ in \cite{Klenke1996}) and using our \emph{general} formula for the migration transition kernel, which leads to the formula for $t_n$ in \eqref{e2572}, which involves the recurrent potential kernel $A(\cdot,\cdot)$ of the difference walk of a pair and the asymptotic behaviour of $A(n)$ defined in \eqref{e2575} as $n \to \infty$. The details are left to the reader.

\medskip\noindent
\textbf{(ii)} We next look at the case where we have \emph{only block-resampling}, and determine that again the $(t_n)_{n \in \N}$ are given by the formulas that we have obtained by considering the \emph{effective} migration of \emph{pairs} under block-resampling and reshuffling, which we study in Appendix~\ref{appA} and which yield the dichotomy. This means that we need to take the scaling limit of the effective difference random walk and its recurrent potential $A(\cdot,\cdot)$ on the $n$-ball at a time of order $t_n$ as $n \to \infty$, and determine the $t_n$ via $A(\cdot,\cdot)$, $N^n$ and the $\mfh$ measure of the macroscopic unit ball. The effective difference random walk started at hierarchical distance $n$ and observed at time $t_n$ has to be shown to coalesce, after the scaling of time and space, to a limit \emph{instantaneous coalescence} with \emph{escape} as above. We need that $k$-tuples, in the limit as $n\to\infty$, again act as pairs \emph{independently}, since the block-resampling rate is summable. Hence, pairs among $k$-tuples become pairs at distance $n$, coalescing upon hitting or escaping on that scale. The rate of escape is given in terms of the transition kernel, in particular, via the recurrent potential kernel $A(\cdot,\cdot)$ and its asymptotics of $A(n)$ as $n \to \infty$.

At this point the \emph{coupling} of the \emph{time-points for the coalescences} in the $k$-tuple has to be handled. But this implies that here the asymptotics for $k$-tuples is the same, since the coalescence appears for a whole random sub-tuple at the same time, namely, on those elements of the tuple in that $j$-ball for which the $j$-block-resampling rings, but where many pairs at once coalesce, namely, all those formed by individuals marked $1$ in that $j$-ball. The limit coalescent is now a different object and is a death process with different coalescence rates than those that arose in (i) due to the effective independence of pairs to meet in a site for coalescence with rate $k \choose 2$ for a $k$-tuple to coalesce. Now these pairs have to meet in $j$-balls to be exposed to $j$-ball $\Lambda_j$-coalescence, whose rate is independent of the number of partitions elements. But they also have to be marked $1$, which happens independently leading to a rate proportional to $k$. Finally, we have to sum over $j$.

In the block-resampling we have, with $k_i$-partition elements in the $i$-th occupied $j$-ball, that $\ell$ partition elements merge at rate 
\begin{equation}
\sum_{j \in \N} \sum_{i \in \N} \sum_{\ell=0}^{k_i} \lambda_j \int_0^1 x^\ell (1-x)^{k_i-\ell} \Lambda_j (\dd x)
\end{equation}
into one partition element. We have to describe the limiting spatial configurations, which are richer than the single-sites we had in (i). Nevertheless, all this decoration is, in the limit $n \to \infty$, \emph{completely removed}: due to the spread-out initial conditions there are \emph{no} multiple collisions in the limit, and everything is \emph{reduced to pairs} and hence to the effective difference random walk. 

To conclude the proof, it suffices that the limiting process has only block-coalescences in \emph{pairs}, since otherwise $k$-tuples with $k \geq 3$ would have to be in the same ball after a large time and the probability for that decays more rapidly as $n \to \infty$ than for a pair. This is not compensated by the growth of the number of $k$-tuples in the $n$-box, and therefore a typical coalescence is asymptotically a pair coalescence. (This was already noted in \cite{Klenke1996,LS06} for specific models and without block-resampling, and we have it for our general model class as well.)

\begin{proof}{\bf [Proof of Lemma~\ref{l.2512}]}
\label{pr.29362512}
We need to get the multi-dimensional distribution of $n$ ball empirical measures from tuples of separated balls, used as windows of observation, and need to show that this multidimensional distribution converges to a limit. We use the corresponding \emph{mixed moment measures}, which we express via the duality relation by studying finite coalescents starting in the set of locations corresponding to the moment measures.

In doing so, we have to address the behaviour of groups of $k_1,\ldots,k_\ell$-tuples of coalescing individuals in which each tuple is at distance of order $n$, but are from different windows at distances that are multiples of order $n$, so that they are in $n$-balls, but these balls are at distances $a_1 + n,\ldots,a_\ell + n$ with $a_i \geq 1$, so that they are well separated. The analysis is parallel to that above, but now the pairs are only exchangeable within an $n$-ball, even though the hitting times are of the same order. Therefore we have pairs within distance $n$ and pairs at distance $a_i+n$. No new ideas are required to handle this case.
\end{proof}

The convergence as $n \to \infty$ of empirical measures has to be strengthened to convergence of clusters $^n \mfC^L, ^\infty \mfC^L$, $L \in \Z$, in a stronger topology than the weak topology on $\wt\CM (\Omega_N^\infty \times \K)$, namely, the weak-atomic topology. This gap will be addressed in Step 7.

\begin{proof}{\bf [Proof of Lemma~\ref{cor.2748}]}
\label{pr.2923}
The main points here are the convergence and the existence, and the properties of the density for the limit process $\wh Y^{\ast,\dagger}$ and its space-time version. The first point is built on the \emph{space-time function-valued duality} and runs parallel to the above replacing spatial coalescents with the \emph{space-time coalescent} we introduced. The second point uses the \emph{strong dual representation} by the entrance law of the spatial coalescent from countably infinite states, but now with entrance configurations that are countable not only at time horizon $t$, but also for the individuals that become unfrozen at earlier times $0 < s_0 < s_1 < \ldots < s_\ell < t$. The task is to lift the proofs of Lemmas~\ref{l.2386}--\ref{l.2512} to the \emph{space-time} context by adding the multiple time points $s_0,s_1,\ldots,s_\ell$. The point here is that we need to work with times $s_0,s_1,\ldots,s_\ell$ that are of the order of $t_n$, but differ by factors less than $1$ from $t_n$, i.e., that we observe at times $s_0 = b_0 t_n,s_1 = b_1 t_n,\ldots,s_\ell = b_\ell t_n$ with $0 < b_0 < b_1 < \ldots < b_\ell < 1$. 

The analysis of the space-time duals for $n \to \infty$ can be carried out with the same tools as in the previous lemma. Hence we get convergence to the claimed limit object. Again, these objects do not depend on $L$, since the limit objects only depend on the tails, namely, $(c_k,\lambda_k,)_{k \geq K}$ for any $K \in \N$.

We also have to prove the existence and identify the properties of the \emph{density} in the limiting measure state $\wh Y^{\ast,\dagger}$ on $\Omega_N^\infty \times \K$ w.r.t.\ $\mfh$. Also this gap will be addressed in Step 7.
\end{proof}

\begin{proof}{\bf [(II) Proof of Lemma~\ref{l.2423}]}
\label{pr.2812}
We focus on part (a) of the lemma: Part (b) and (c) will be postponed again to Step 7. Again, we treat separately the cases without or with block-resampling and the combination of both mechanisms. Note that the underlying L\'evy-process of the limit is the $(c_k,\lambda_k)_{k \in \Z}$-L\'evy process $\zeta^\infty$. We did already show that the underlying effective migration of a $k$-tuple on $\Omega_N^L$ embedded in $\Omega_N^\infty$ converges to the $k$-tuple $\underline{\zeta}^\infty$ of such effective L\'evy-processes $\underline{\zeta}^\infty$, and hence next we have to add the resampling operators as well. 

\medskip\noindent
\textbf{(i)} Consider the case \emph{without} block-resampling. We have to show that the finite coalescent on $\Omega_N^L$, started in $k$-spread-out $(k \geq -L)$ initial states and observed at time $t_L=tN^L$, converges as $L \to \infty$ to the state of a \emph{spatial coalescent on $\Omega_N^\infty$}, based on a migration mechanism given by a system of delayed coalescing L\'evy-processes on $\Omega_N^\infty$. The delayed coalescing is of the form $d^\ast$ times the aggregated collision local time of pairs of partition elements. The convergence statement for the migrating population remains true after adding the \emph{delayed coalescence} to the migrating population, since we showed that the scaled joint occupation times of the L\'evy-processes on $\Omega_N$ or $\Omega_N^L$ converge to the collision local times of the pairs of the $\underline{\zeta}^\infty$. Therefore the dual spatial coalescents converge to the claimed limits.

\medskip\noindent
\textbf{(ii)}
The case of block-resampling only again uses that we have the same as in the treatment of pairs in Appendix~\ref{appA}, and we have to address how to get that no $k$-tuples meet in the $n \to \infty$ limit to coalesce for $k >2$. However, we showed already for the previous lemma.

\medskip\noindent
\textbf{(iii)} Again the combination of mechanisms in the $k$-tuple has to be addressed, where the distribution of the pairs is just the replacement of the effective walk by the difference walk, as above. We can use the \emph{Trotter formula} for the sum of two operators of a martingale problem, corresponding to the spatial Cannings process operators and as second summand the block-resampling with reshuffling operators. The combination of the mechanisms works due to the fact that the $\mfh$-locally monotype balls property is preserved in the jump case or at least is restored in positive time, under \emph{both} mechanisms alone as resampling mechanism. Therefore we can represent the semigroup of the process as a sum over convolution products acting over small time intervals, acting alternatively with either one, and the claim follows.
\end{proof}


\subsubsection{(III) Step 7.\ The limit process has continuous locally $\{0,1\}$-valued densities in space}

This step consists of 3 substeps.


\paragraph{Step 7.1: Closing a gap in the argument}

$\mbox{}$\\[0.2cm]
We close the gap that was left in the argument in Step 6. We already have the limit fields in the cluster growth analysis on $\Omega_N$ or $\Omega_N^L$, namely, $\wh X^{\infty,\sim}$, and for the limit state of the $L$-approximations as $L \to \infty$, the process $\wh X^\infty$. For \emph{both} we still have to show the \emph{local monotype property}, meaning that $\wh X^\infty$ restricted to any occurring type $\mfa,a$ has a density w.r.t.\ to $\mfh$ with a version that is \emph{$\{0,1\}$-valued} and \emph{continuous} in space, c\`adl\`ag in time, such that there is continuity in $t>0$ \emph{without} block-resampling at that time and point in space. We view the objects obtained in part (a) and (c) of Lemma~\ref{l.2423}, respectively, Remark~\ref{r.2876}, as the limit as $L \to \infty$, and the statements about cluster formation in Lemma~\ref{cor.2748} as $n \to \infty$. We distinguish between the case without and with block-resampling.

\medskip\noindent
\textbf{(i) Cluster formation.} 
For the analysis of the \emph{cluster formation} of $\wh X$ on $\Omega_N$, respectively, $\wh X^L$ on $\Omega_N^L$ we have to argue similarly, but with $t \to \infty$, to get the result on the limiting object arising there in the scaled space embedded in a finite ball in $\Omega_N^\infty$, namely, to be covered by a collection of $\mfh$-(single type) balls. The idea is that \emph{finite} sets in the geographic space $\Omega_N^L$ have a monotype $t \to \infty$ limit, and that this is the case for sets that \emph{increase} with time. We can identify how fast the sets must grow from small monotype sets within the scaled limiting ball. In the limit field on $\Omega_N^\infty$, we get a disjoint decomposition into balls carrying $\mfh$-a.e.\ a specific type. 

We start with a weaker property first, where we use the \emph{strong duality} applied at time horizon $t_n$, with initial configurations with countably many initial dual particles at each site of $\Omega_N$ (or $\Omega_N^L$) and embedded in $\Omega_N^\infty$. We observe the process by looking at the subcoalescent of a ball of size $n$ as a typical macroscopic unit ball. The observed initial configuration in that ball grows like $N^n$, and we want to run the dynamics for time $t_n$ and show that the dual population is, in the limit as $n \to \infty$, \emph{locally finite} in $\Omega_N^\infty$, for every backward time $s>0$ from the time horizon $t$. This is not straightforward due to the migration, even when the resampling is coming down from infinity, since there are possibly countably many small partition elements that emigrated, but their \emph{frequency} is small, and tends to zero in the limit as $n \to \infty$ so that the limit $\wh X^{\infty,\sim}$ has a continuous $\{0,1\}$-valued densities and is locally monotype. 

\medskip\noindent
\textbf{(ii) $L$-approximation.} 
 For the purpose of $L$-approximation, we use the strong dual representation of $\wh X^L,\wh X^\infty$ by $\wh Y^{L,\ast}, \wh Y^{\infty,\ast}$ and consider the dual for finite $L$, started with countably many dual individuals in each site of $\Omega_N^L$. We use versions of $\Omega_N^L$ embedded in $\Omega_N^\infty$, with \emph{one} particle in every site of the \emph{skeleton} and \emph{subgroup} $\Omega_N^{\infty,\rm disc}$ of $\Omega_N^\infty$. We show that, for our $L$-approximation limit, in the chosen $k(t)$-ball around $\eta$ there is coalescence to one partition element in the proper sense, which means that we have to show for every $a \in A$ that the level-sets of $1$ for the density $\wh Y^{\infty,\dagger} (\cdot,a)$ contain $m$-balls $\mfh$-a.e.\ for some $m \in \Z$, where the value is $1$ $\mfh$-a.e.\ in this ball and the union over these balls has full $\mfh$-measure when we take the union over $a \in A$. To use the results on cluster formation, we have to verify the claim for its limit object on $\Omega_N^\infty$. 


\paragraph{Step 7.2: Intermezzo on infinite particle systems}

$\mbox{}$\\[0.2cm]
We view the spatial coalescent as an \emph{infinite particle system} on $\Omega_N$ and similarly on $\Omega_N^L$ and even on $\Omega_N^\infty$. It is important to embed this population in $\Omega_N^\infty$ to get a system on this geographic space and get the limit object on $\Omega_N^\infty$. The goal is to show that with frequency $1-\varepsilon$ for every $\varepsilon >0$ we have \emph{configurations} in a \emph{Liggett-Spitzer space}. The latter quantifies and generalizes the requirement "finite" configurations on infinite geographic space and on $\Omega_N^\infty$ of \emph{locally finite} configurations. In the case of the Fleming-Viot resampling the process on \emph{finite} geographic spaces comes down from infinity. On \emph{infinite} spaces this holds only for an $(1-\varepsilon)\mfh$ part of the population because of migration. However, we can approximate the population up to $\varepsilon \mfh$ by elements in a \emph{Liggett-Spitzer space}, which is a Banach space of configurations in $\N^{\Omega_N}$ (compare~\cite{GLW05}). This can be generalised to Cannings resampling with $\Lambda$-measures satisfying the \emph{no dust} condition.

We want to show that the property persists $\mfh$-a.e.\ in the $n\to \infty$ limit $\wt\CM(\Omega_N^\infty)$ on the geographic space $\Omega_N^\infty$. Once we have this, we have to argue that $\mfh$-a.e.\ we are locally finite on a ball around a site up to $\mfh$-mass at least $(1-\varepsilon)\mfh$ for \emph{all} $\varepsilon>0$. In other words, we try to work with a \emph{subpopulation} of the spatial coalescent process on $\Omega_N^L$ driving the dual process, respectively, in the limit the $\CM(\Omega_N^\infty)$-valued dual coalescent process, loosing a fraction $\varepsilon \mfh$, but such that this is a \emph{locally finite} infinite particle system, which evolves for every positive backward time $s$ like a spatial \emph{instantaneous coalescent} for the analysis of cluster growth, respectively, is \emph{delayed coalescing} for the analysis of $L$-truncation limits. Furthermore, this extends to $n=\infty$, respectively, $L=\infty$ and so even the $\varepsilon$-truncation can be dropped in the point-recurrent case.


\paragraph*{Strong duality for $L$-truncation.}

For our proofs we use a combination of the \emph{strong dual representation} combined with a particular construction of the dual entrance laws $\wh Y^{\ast,L}$ (or $\wh Y^{\ast,n,\sim}$) in $\Omega_N^\infty$ for $n,L \in \N \cup \{\infty\}$ starting with one particle in every point of $\Omega_N^{\infty,\rm disc}$, via an \emph{embedding} in $\Omega_N^\infty$. Furthermore, studying the $L$-truncation we use a \emph{coupling} of terms of the finite-$L$ case and the infinite-$L$ case. For this purpose we can represent the $\underline{\zeta}^L,\underline{\zeta}^\infty$ in a \emph{coupled way}, working just with the $L$-truncation approximation property of $\zeta^\infty$, where we suppress jumps in "size" $<-L$, so that all these processes $\wh Y^{\ast,L}$ arise by suppressing jumps of size smaller $-L$ are constructed on \emph{one probability space} and on \emph{one geographic space} $\Omega_N^\infty$. 

To see the $\{0,1\}$-valued and continuity property for $L=\infty$ of the density $\wh Y^{\infty,\dagger}(\cdot,\{a\})$, which is equal to $\wh X^{\infty,\dagger}(\cdot,\{a\})$, we have to turn from the empirical measures $E$ to the cluster description of states given by $\mfC$. We need that $^{\infty} \mfC^{\infty}$ provides a decomposition into \emph{monotype balls}, and similarly for $n=\infty$ and $\wh Y^{\infty,\dagger,\sim}$. For both the cluster growth analysis on $\Omega_N$ (or $\Omega_N^L$) or the $L$-truncation, the \emph{strong dual representation} is important. We explain this both for $\wh X$ with $X^{\infty,\sim}$, respectively, $\wh X^\infty$ and the entrance laws $Y^\sim$ and $Y^{\infty,\sim}$. The following property is needed in both cases for the cluster growth analysis and the $L$-truncation.

Recall the embedding of the countable initial state on $\Omega_N^L$ into $\Omega_N^\infty$. For $E^L$, $L \in \N$, we have a representation of the density in $\eta$ by taking the limit $m \to -\infty$ of the \emph{$m$-ball empirical measures} around a site $\eta \in \Omega_N^\infty$, which yields a random $\bar \mfa \in \CM_1( \K)$ such that $\bar \mfa$ is the limit at $\eta \in \Omega_N^\infty$, and we already showed that $\bar \mfa=\delta_\mfa$ for $L=\infty$. Namely, on the countable dense set $\Omega_N^{\infty,\rm disc}$ by the strong dual representation we have that there exists a function $\mfa\colon\, \Omega_N^\infty \to\K$ ($\dagger$ indicates the density of a measure)  
\begin{align}
\label{e2869}
&\wh X_t^{\infty,\dagger} (\eta,a) = \wh Y_t^{\infty,\dagger}(\eta,a), \qquad a \in A_t, \mfh\text{-a.e. }\eta,\\ 
&\text{with } \wh Y_t^{L,\dagger}(\eta,a) = 1_{\{\mfa (\eta)=a\}}, \qquad (\eta,a) \in \Omega_N^{\infty,\rm disc}  
\times \K, \, L \in \N \cup \{\infty\}, \nonumber
\end{align}  
where $\mfa(\eta)$ is the type of the dual path of $Y^L$ starting in $\eta$ for the time horizon $t$, and $Y^L$ starts with one individual in every site of $\Omega_N^{\infty,\rm disc}$. 

Similarly this representation holds for studying the cluster formation in $\Omega_N,\Omega_N^L$ for $\wh X^\sim,\wh X^{L,\sim}$ and $\wh Y^\sim,\wh Y^{L,\sim}$.


\paragraph*{How about strong duality for $L=\infty$?}

Suppose that we know that the density is \emph{continuous} also for $L=\infty$. Then, from the approximation property as $L \to \infty$ and by the \emph{continuity} of the density for $a \in A$ and $\eta \in \Omega_N^\infty$, we have the property 
\begin{equation}
\label{e2887}
\wh Y^{\infty,\dagger} (\eta,a) = \liml_{m \to -\infty} \wh N^{-m} Y^{\infty,\ast} 
\left(B_{-m} (\eta) \times \{a\}\right) \quad \mfh\text{-a.e.} \quad \text{a.s.}
\end{equation}
when the dual starts in $\mfh$. Hence, $\mfh$-a.e.,
\begin{equation}
\label{e2866}
\wh X_t^{\infty,\dagger} (\eta,a) = 1_{\{\mfa (\eta)=a\}}, \quad (\eta,a) \in \Omega_N^\infty \times \K.
\end{equation}  
We already showed that the density has the value $0$ or $1$ on $\Omega_N^\infty$. The level sets of all $\mfa \in A$ give the decomposition of $\Omega_N^\infty$ into disjoint sets. The two representations are equal in law. Hence, to conclude the argument, for $\wh Y^{n,\dagger,\sim} (\cdot,a)$, respectively, $\wh Y_t^{L,\dagger}(\cdot,a)$ on $\Omega_N^\infty$, we need to prove the \emph{continuity} for $n=\infty$, respectively, $L = \infty$, which is the same as the decomposition property in ``balls" of a single type.

This requires us, for the dual entrance law $\wh Y^\infty$ or $\wh Y^{\infty,\sim}$ starting with one dual individual in the sites of $\Omega_N^{\infty,\rm dis}$ restricted to an $m$-ball containing $\eta$, to obtain the probability that this $m$-ball collapses into \emph{one} single partition element up to a \emph{frequency}-zero subset by time $t>0$ is $0$ for $m$ below a threshold. We can cover a given finite ball $B \in \CB (\Omega_N^\infty)$ with an union of finitely many balls at time $t$ with $1-\varepsilon$ of its frequency having collapsed by time $t$ for \emph{every $\varepsilon>0$}. Once we have that, we find a sequence of balls whose union has full $\mfh$-measure on the collapse property. This can be phrased by saying that the necessary continuity is produced $\mfh$-a.e.\ instantaneously in backward time $s>0$ in $^t \wh Y_s^\infty$ for the time horizon $t$. Consequently, for $\wh X^\infty$ or $\wh X^{\infty,\sim}$ we get \emph{locally monotype} for $t>0$ by the strong duality.

Viewed forward the above means that in a genealogical view the descendants of single founders at a certain time back ending in some ball contributing $\mfh$-measure $(1-\varepsilon)\mfh$ at least, have immigrated first at times $\delta_0,\delta_1,\ldots$ back, evolve over time $t-T_{\delta_k}$, and occupy at time $t$ an $m_k$-ball $B$ around $\eta$, for $k \in \Z$ small enough, $\eta \in B$ (see~\cite{gh2}). 
\qed


\paragraph{Step 7.3: Continuation of the proof for cluster formation}
\label{c73}

$\mbox{}$\\[0.2cm]
Next we return to the cluster formation on $\Omega_N$ and try to verify the property~\eqref{e2869}--\eqref{e2866} above to complete the argument carried out in Step 6 to then modify and use it for the $L$-truncation. We have to control the number of partition elements in space, as a function of $n$ and the number of points at a site, say $\ell$, to force the number of individuals in one partition element to be large enough. Recall that we have $N^n$ sites in a macroscopic ``unit"-ball, each possibly carrying partition elements that must coalesce. Recall that we deal with an entrance law (from a locally countably infinite population) in which diverging sets of individuals may behave differently but have frequency $0$ and are therefore \emph{not} observable. 

We have to show, with $1-\varepsilon$ of the local frequency at positive backward time, we only have \emph{finitely many} partition elements, how does this number changes with $n$ and what can we conclude for $n=\infty$, i.e., the \emph{limiting state on $\Omega_N^\infty$}. This was resolved in~\cite{GLW05} for the Moran model here called the Fleming-Viot model. We mark the individuals moving before time $\delta$. The remaining population has frequency $<\varepsilon$ when we choose $\delta=\delta(\varepsilon)$ properly. This can be done for all the dynamics considered here. By construction, the population not moving is \emph{coming down from infinity} if $d^\ast>0$ or $|\Lambda_0|>0$. Therefore we have \emph{finitely many} partition elements in finite sets. The population moving before time $\delta$ is coalescing at the new location and altogether can be controlled as a function of $t$, so that the approximation property for $\varepsilon \downarrow 0$ holds \emph{uniformly in $t \geq \delta$}. Consequently, the spatial coalescents are \emph{locally finite infinite particle systems} on $\Omega_N$. In fact, the paths are for positive backward time in a \emph{Liggett-Spitzer} space. If we only have the no-dust property for $\Lambda_0$, then we have to make a further $\varepsilon^\prime$-approximation, taking finitely many per site with frequency $1-\varepsilon^\prime$ locally.

\medskip\noindent
\textbf{(I)} \emph{The $\varepsilon$-locally finite property for $n \to \infty$ limit.}
We want to conclude that, in the limit as $n \to \infty$, we have as entrance law, starting from countably many individuals at each site, and turns into a \emph{locally finite coalescent} $Y_s^{\infty,\sim}$ for $s>0$, $s \in t$, with $t>0$ the time horizon for the backward process. This means that we have to show that this property is preserved in the limit $n=\infty$. The number of sites in $\Omega_N$ embedded in $\Omega_N^\infty$ grows like $N^n$ in the scaled object, but the accumulated (in time) rate of resampling in that volume, which is active up to time $t_n$, increases by $A(n)N^n$, resulting in tightness of the local occupation numbers in $n$.

To verify this property at $n=\infty$, we approximate $\Omega_N^{\infty,\rm disc}$ (restricted to a finite ball) by finite sets $\bigcup_{\ell=0}^m \Omega_N^\ell$, and then let $m \to \infty$. Since in the limit we have \emph{instantaneous} coalescence upon hitting and the difference process of a pair is point-recurrent, it turns out that the entrance law must come down from infinity locally in positive time. For that we note that around every site we have for every ball a replacement of position by the uniform distribution at infinite aggregated rate from $k$ to $-\infty$, for $k \in \Z$, by the assumption of point-recurrence. This means that dual individuals keep coalescing and coming down down from infinity locally as follows for $\mfh$-almost all $\eta \in \Omega_N^\infty$. 

Consider a ball of size $k$ and all decompositions into $k-1,k-2,\ldots$ balls. Along every such string of refinement, the rate of coalescence arises from paths hitting and hence is based on the number of pairs within a $k^\prime <k$ ball and the migration rate per site aggregated along the string is \emph{infinite}. If there would be countably many partition elements locally, then we could form a sequence of pairs converging to a single site, but these pairs could not have survived for $\mfh$-almost all $\eta \in \Omega_N^\infty$ for positive time as we see immediately by contradiction.

\medskip\noindent
\textbf{(II)} \emph{Monotype decomposition property for the $n \to \infty$ limit.} We have to \emph{strengthen} property (I) above to get geometrically a cluster decomposition, requiring that clusters containing \emph{full balls} exist for the $n=\infty$ limit model and their union has full $\mfh$-measure. For \emph{finite $n$} due to migration we do not have this property, since always we can escape from a ball and hence do not collapse in a ball in finite time. But we do for the \emph{(infinite) refinement} of a finite ball in $\Omega_N^\infty$, via the internal structure in connection with \emph{infinite} rates, for \emph{small transitions}. 


\paragraph*{Background 1: \emph{Genealogy}.}

We need to conclude that the local finiteness of the configuration of the partition elements can be strengthened to the decomposition property. We consider the genealogical structure of the dual population and its relation to the structure of the corresponding locations. The dual population carries a \emph{genealogical structure} $\CU_t$ at backward time $t$ for time horizon $t$ in the form of a $(\Omega_N^\infty \times \K)$-marked ultrametric measure space $\CU_t = (U_t,r_t,\kappa_t,\mu_t)$, with $U_t$ the \emph{individuals}, $r_t$ the \emph{genealogical distance}, $\mu_t$ the \emph{sampling measure} and $\kappa_t$ the \emph{mark function} assigning type and locations to the individuals, i.e., $\kappa_t\colon\,U_t \to \Omega_N^\infty \times \K$. This will be studied in great detail in~\cite{gh2}. The idea to continue is to use the ``genealogical structure" of the geographic space $\Omega_N^\infty$, which is also an \emph{ultrametric measure space} $(\Omega_N^\infty, d^\infty,\mfh)$ (see Section~\ref{sss.ONLdef}). 

For the population of \emph{backward} independent L\'evy-processes starting in a site in $\Omega_N^{\infty,\rm disc}$, we can define for every pair of distinct sites the first backward hitting time $H(\eta,\eta^\prime)$. Turning on the coalescence we prune this population of paths to get the ones that are \emph{$\Omega_N^\infty$-marked geodesics} connecting the starting point with the end point, by requiring that they are first in one partition element after $H(\eta,\eta^\prime)$. This defines an ultrametric measure space by taking its completion $\CU_t$, in which $r_t(\eta,\eta^\prime) \geq 2 H(\eta,\eta^\prime)$. 

The strategy to relate the two ultrametric structures on $\Omega_N^\infty$ and on $U_t,t \geq 0$ is based on the following observation. We can identify the initial individuals $U_0$ in the dual entrance law by the sites in $\Omega_N^{\infty,\rm disc}$ and their position at backward time $s=0$ with that value. In this way we obtain a population of backward paths with time index $s$ running up to time $t$, the \emph{time horizon} of the duality. This gives us the \emph{marked geodesics} of the \emph{genealogy} of the population at time $t$. This population can be decomposed into the $2t$-balls corresponding to the different founding fathers of the population. These components always carry the same type. This means that for given $\eta$ we look for $\eta^\prime$ where $d^\infty(\eta,\eta^\prime)$ is such that the $d^\infty(\eta,\eta^\prime)$-ball is $\mfh$-a.e.\ contained in the $2r_t(\eta,\eta^\prime)$-ball.


\paragraph*{Background 2: \emph{Genealogical and geographic neighbourhood}.}

Thus, we need to find a \emph{genealogical neighbourhood} that contains a \emph{geographic neighbourhood}, so that a cloud of close, i.e., at distance $r_t(\cdot,\cdot) < 2t$, relatives covers a ball of size $k^\ast \in \Z$. The same holds also for genealogical balls $r_s$ for $0<s<t$. This is a property of the population of $\Omega_N^\infty$-marked geodesics. This property holds because under the entrance law large proportions of the coalescent \emph{collapse instantaneously} and the elements in the partition element do not have enough time to dislocate far and hence almost all cover a full geographic ball around a site at time $t$. In other words, \emph{positive fractions} of the population adding up to $1$ coalesce in \emph{shorter} time scales than migration of a certain \emph{size} takes place, due to the form of the rates proportional to the number of pairs and the \emph{rates of coalescence} acting with \emph{delay} during \emph{collision local times} after hitting first, in the $L$-truncation limit, respectively, \emph{instantaneously} in the application to the cluster growth analysis. For fixed $n$ we \emph{always have $\varepsilon$-escape} of the population, but as $n \to \infty$ (similarly for $L \to \infty$) this becomes \emph{negligible} below a level $m \in \Z$ under the \emph{assumption} of \emph{point-recurrence}.

The important structure of the dual at $Y_t^{\infty,\ast}$ for time horizon $t$ is described by the behaviour of the probability to hit $\eta$ from $\eta^\prime$, with $\eta,\eta^\prime$ in a ball of size $k$, before leaving the $k$-ball, denoted by $q_k(\eta,\eta^\prime)$, namely, its convergence to $1$ as $k \to -\infty$, such that even the ball intersected with $\Omega_N^{\infty,\rm disc}$ collapses with \emph{frequency $1$} for some \emph{random} $k^\ast \in \Z$ sufficiently small. This is derived from a function $m_t(\cdot)$ on $\Omega_N^\infty$ ($t$ the time horizon) before a substantial number of partition elements leave this ball. It means that the frequency of the elements of $\Omega_N^{\infty,\rm disc}$, jumping outside a ball \emph{before} coalescing with another partition element starting in this ball, is $O(N^{-L})$, so that instantaneously pairs coalesce at positive rates, and there is a first (in backward time) jump of a partition element of $\eta^\prime$ out of a ball around $\eta$ with full $\mfh$-measure of such a $\eta^\prime$. Here, we order the elements in the ball and $\Omega_N^L$ successively, adding the ones for $L=0,1,\ldots,$ etc. In the end, via our \emph{coupling} of the objects for different $L$, respectively, $n$, this means that $\mfh$-a.e.\ site $\eta$ and $\mfh$-a.e.\ site $\eta^\prime \in B_{k^\ast}(\eta)$ have the same common ancestor, whose descendants cover $B_{k^\ast}(\eta)$ at time $t$. Since, for finite $L$, we always have $\varepsilon$-escape for every $\varepsilon >0$ from balls, this structure occurs \emph{only} in the \emph{limit} $L \to \infty$, respectively, $n \to \infty$. Note that, viewed forward, this is \emph{immigration} of new subfamilies into the ball from the outside, perturbing the leading subfamily. The fact that the $\varepsilon$-escapes can collapse arises from the fact that for fixed $\varepsilon$-error the property holds from some finite $n_0(\varepsilon$) (respectively, $L_0(\varepsilon)$) onwards. The issue is to have $n_0(\varepsilon)$ (respectively, $L_0 (\varepsilon)$) bounded from above on $(0,\varepsilon^\prime)$ for every $\varepsilon^\prime >0$.

The next step is to obtain that the whole space decomposes into \emph{countably many disjoint monotype balls}. For that we have to exclude ``condensation" of arbitrarily many small balls in sets of \emph{positive} Haar measure. 


\medskip\noindent
\paragraph*{(III) Transfer from cluster growth to $L$-truncation.} 
Suppose that we have settled everything for the setup for cluster formation on countable spaces like $\Omega_N$ or $\Omega_N^L$. How do we pass from these results to the local properties of $\wh X^\infty$? We want to use the \emph{cluster growth results} for $\wh X^L$ and settle the \emph{local properties} of $\wh X^\infty$ by using its truncation approximation via $\wh X^L$ for $L \to \infty$. It suffices to prove that the limiting $\wh Y_t^{\infty,\ast}$ restricted to a particular type has \emph{$\{0,1\}$-valued} and \emph{continuous} densities w.r.t.\ $\mfh$ in \emph{space and time} on a collection of disjoint balls of \emph{full} $\mfh$-\emph{measure}.

We need to show that in $\wh X^L$ for large time, i.e., $tN^L$, we have cluster growth forming in the $L \to \infty$ limit, namely, local $\mfh$-monotype regions $\mfh$-a.e.\ in a suitable topology for \emph{all} $t>0$. In the cluster growth analysis for $X$ on $\Omega_N$ or $\wh X^L$ on $\Omega_N^L$, we observe at level $n$, up to observation time $A(n) N^n$, shift and scale in space by shifting levels by $-n$, and get the limit $\wh X^{\infty,\rm sta,\sim}$, i.e., $d^\ast = \infty$. However, here we do \emph{not} shift levels by $n$ downward, but \emph{add} $L$-levels below level $0$, and observe for large $L$ at levels of $O(1)$ corresponding to \emph{averaging} over sets of \emph{finite} $\mfh$-measure. Hence we have $A(O(1))$ instead of $A(L)$ as factor in the formula for $t_n$ the cluster observation time (recall~\eqref{e2218}, and therefore the local cluster observation time is $tN^L$ the time to observe \emph{monotype clusters} in the $L \to \infty$ limit on level $0$ replacing the $A(L)\,N^L$ for the cluster growth analysis. 

We have to argue that the above is due to increasing $L$, since we get faster rates and scaled changes at the \emph{active} level as in the cluster growth analysis which is $-L$, and are studying cluster formation at level-$0$ in sets of size of order $1$, which leads to monotype cluster formation applied to $\wh X^L$. Furthermore, the state is observed after \emph{aggregating over $L$ steps upwards} from the \emph{active level}, and we are still getting a diffusion with volatility $O(1)$. Since now the state is observed on levels $O(1)$, which leads to $\wh X_t^\infty$ for $t>0$ rather than the \emph{infinite volatility process} $\wh X^{\infty,\rm sta,\sim}$, arising from the \emph{cluster formation} on countable groups observed at large level-$n$ at late time. 

The volatility of balls of size $-L$ increases to $\infty$ when the ball does not become monotype, and aggregated over unit volume is still $O(1)$. Then eventually the values in grid points become $0$ or $1$ for all $a$-densities. For this to be true in an array of processes indexed by $L$ and $n$, with the aggregation of states by $n$ levels, requires that we take:
\begin{equation}\label{e3224}
\text{the limit as $n$ and $L$ tend to infinity such that $n=L$. }
\end{equation}
The issue is whether this scale is the one we found when studying the \emph{empirical measure} in Lemma~\ref{l.2512}, when we take into account that there is no shift by $n$ and we observe at levels $O(1)$, leading to the cluster behaviour described by $^n \mfC^L$ in that lemma. 

The idea to make this rigorous is to transfer this property on cluster formation from $\wh X^L$ for $t \to \infty$ and viewed as a measure on $\Omega_N^\infty \times \K$, via the described embedding to an \emph{instantaneous} effect for macroscopic time $t>0$ for $\wh X^\infty$, by a truncation approximation in a sharp enough \emph{strengthening} of the \emph{weak topology of measures} on $\CM(\Omega_N^\infty \times \K)$ viewed as processes in macroscopic times $t$ bounded away from $0$, which guarantees convergence to monotype clusters, i.e., the object $^n \mfC^L$ to our limiting $(n,L)$-arrays of \emph{processes} along the $n=L$ axis called $^\infty \mfC^\infty$.


\paragraph*{Two remaining tasks for the proof.}

In analogy with the cluster growth analysis above, this leaves us with two tasks for the entrance law $\wh Y^\infty$ which are intertwined:
\begin{itemize}
\item[\textbf{(A)}] The convergence property we proved for part (a) remains true in the \emph{stronger topology} applying it to $^n \mfC^L$ in the $L \to \infty$ and $n \to \infty$ limit as $^\infty \mfC^\infty$ fixed above. 
\item[\textbf{(B)}] Both the \emph{locally finite} property and \emph{local monotype-ball decomposition} property hold for $\wh Y^\infty$, the limit of the array as $(n,L) \to \infty$ with $n = L$ and for the cluster growth $\wh Y^{\infty,\ast,\sim}$.
\end{itemize}
We explain these task now in more detail.

\medskip\noindent
\underline{Part (A).} The stronger topology we want to use for the convergence is based on the atomic structure of the measures on $\K \times \Omega_N^\infty$ arising for finite $L$ from the atomic structure of the geographic space and the collapse to countably many types and for $L=\infty$ from the locally monotype property of the limit $L \to \infty$ in space and the locally finite structure of the infinite particle system, combined, allowing to define an atomic structure in space putting the mass of the ball in one point in the ball and its type.

The question is how the structure described in the intermezzo in Step 7.2 leads to monotype clusters containing full balls for the $L$-truncation limit $^\infty \mfC^\infty$ at all macroscopic times $t>0$. We proceed as in the analysis of the cluster growth on $\Omega_N$ above, and bring into play the \emph{atomic} features of the measure states. Note first that $\wh X^L_t$ has a countable random set $A_L$ of types for $t>0$, which we can describe as an \emph{atomic measure} on $\K \times \Omega_N^\infty$, with weights on $A_L \times \Omega_N^L$. We want to \emph{sharpen} the $L \to \infty$ convergence in \emph{law} based on the weak topology on $\wt\CM(\Omega_N^\infty \times \K)$ obtained earlier, by using the weak \emph{atomic} topology. We can do the same for the dual processes $\wh Y^{\ast,L}$ and obtain partition elements with positive frequencies at sites that specify atoms in $\CI \times \Omega_N^\infty \times \K$, and need to verify the convergence in the \emph{weak atomic topology}. For that we look at the number of partition elements for a subpopulation with sufficiently large accumulated frequencies at sites. Return now to the case of finite $L$.

We know for the embedded $\wh Y^L$ that the density on $\Omega_N^\infty$ w.r.t.\ $\mfh$ for every $L \in \N$ exists by construction and is uniformly continuous in $\eta$ on balls and $t$ in intervals in $(0,\infty)$ (and only with block-resampling would be c\`adl\`ag in $t$) and clusters in finite sets and large time. We must prove that this \emph{continuity in space} is \emph{uniform} in $L \in \N$ and therefore the limit $L \to \infty$ is continuous. 

To do so, the first step is to make use of the cluster growth analysis, showing that the densities w.r.t.\ $\mfh$ \emph{and} in the limit exist, i.e., $^{\infty}E^{\infty, \dagger}_{\infty,m} (t,\eta)$ from~\eqref{e2426} is \emph{continuous} in the variable $\eta \in \Omega_N^\infty$, \emph{uniformly} on \emph{finite} subsets (under $\mfh$) of $\Omega_N^\infty$ and for macroscopic times $t \geq \varepsilon$ with $\varepsilon >0$. 

This scenario is based on the corresponding continuity of the joint \emph{collision local times} of the pairs of the limiting duals, respectively, the limiting \emph{local times} for the underlying difference L\'evy-processes $\zeta^L,\zeta^\infty$, whose \emph{semigroups} generate the deterministic \emph{mass flow part} of $\wh X^L,\wh X^\infty$. Therefore we have to control the modulus of continuity as a function of $L$, for which we need tools.


\paragraph*{Lyapunov function.}

For $s \in [0,t]$, consider the function
\begin{align}\label{e3077}
&\text{$h^t_s \colon\,\Omega_N^\infty \times \Omega_N^\infty \to [0,\infty)$ such that $h^t_s (\eta,\eta) = 0$ 
for all $\eta \in \Omega_N^\infty$}\\
&\text{and $h(\eta,\eta^\prime) \in (0,\infty)$ for all $\eta,\eta^\prime \in \Omega_N^\infty$ 
with $\eta \neq \eta^\prime$.} \nonumber
\end{align}
We next introduce a Lyapunov function that measures the collision local time of pairs of paths of a population of paths of descent at a backward time $s \in [0,t]$ accumulated up to time $s$ and aggregated over the space with $\mfh$, with $t$ the current time horizon of the population. This increases with $s \to \infty$ as $s \uparrow t$. This quantity gives the total rate of coalescences in the population after the backward time $0$. We compare this function for $L \in \N \cup \{\infty\}$, increasing $L$ step by step by one. 

This comparison has to be organised in a nested scheme of successively refined systems of balls around $\eta \in \Omega_N^\infty$, each equipped with the restricted $\mfh$-measure, with $\eta$ also governed by the $\mfh$-measure. If we increase $L$ by one, then we replace the stirring in each ball by $N$ stirrings in smaller balls, and \emph{upscale} the coalescence rates to \emph{higher} rates of coalescence. In this way, as we increase $n$ or $L$, we respectively decrease or increase the integral over a ball. We need to show that, as $L \to \infty$, there is a limit function corresponding to locally monotype structures, i.e., locally equal to $0$, because in positive time the coalescent collapses in a whole ball. The argument runs as follows.

Fix $t>0$. We claim that the Lyapunov function $h^t_s$ converges to $0$ locally in $L^1(\mfh)$ as $L \to \infty$, uniformly in $s \in [\delta,t]$ for every $\delta > 0$.

\medskip\noindent
\textit{Monotone comparison.} Recall the nested scheme described above: increasing $L$ by one replaces the stirring on each $L$-ball by $N$ stirrings on $(L-1)$-sub-balls and upscales the coalescence rates. Hence $h^t_s$ is \emph{non-increasing} in $L$ for every fixed $s \in (0,t]$ and $\mfh$-a.e.\ $(\eta,\eta^\prime)$. Since $h^t_s(\eta,\eta') \geq 0$, the monotone limit $h^t_{s,\infty}(\eta,\eta') = \lim_{L\to\infty} h^t_s(\eta,\eta')$ exists pointwise. We must show that this limit is $0$ locally in $L^1(\mfh\otimes\mfh)$. The existence of the monotone limit is established. The fact that the limit equals zero follows from the point-recurrence argument given below.

\medskip\noindent
\textit{The limit is zero.} In the point-recurrent regime the difference process of any pair of dual paths hits $0$ (i.e., the two paths occupy the same site) with probability $1$, and the accumulated local time at $0$ grows to $\infty$. At level $L$ this collision local time is exactly what $h^t_s$ measures, aggregated over pairs and integrated against $\mfh\otimes\mfh$ over a ball. Because the collision rate for a pair of paths in an $m$-ball is bounded below by $c_m^{-1}$ (see the estimates in Appendix~\ref{appA}), and $\sum_{m} c_m^{-1} N^m = \infty$ in the point-recurrent regime, the accumulated coalescence probability in any finite ball converges to $1$ as $L\to\infty$. Hence $h^t_{s,\infty} = 0$ $\mfh\otimes\mfh$-a.e. By dominated convergence (using the uniform $L$-bound on $h^t_s$ derived from the total resampling mass), the convergence holds in $L^1(\mfh\otimes\mfh)$ over every finite ball.

\medskip\noindent
\textit{Uniformity in $s \in [\delta,t]$.} The function $s \mapsto h^t_s$ is non-decreasing, and the bound on $h^t_\delta$ is uniform in $L$ by the cluster-growth estimates (Lemma~\ref{l.2386}). Therefore the $L^1$-convergence of $h^t_s \to 0$ holds uniformly in $s \in [\delta,t]$.

\medskip\noindent
\textit{Tightness in the weak atomic topology.} Fix a ball $B \in \CB(\Omega_N^\infty)$. We verify tightness of $(X^L_t)_{L \in \N}$ in the weak atomic topology on $\CM(\Omega_N^\infty \times \K)$ by checking the two standard conditions of \cite[Section~2]{EK94}:

\medskip\noindent
\textit{(i) Compact containment (uniform mass bound).}
The total atom-mass of $X^L_t$ in $B \times \K$ equals $\mfh(B)$ for every $L$ (since $X^L_t$ is a probability kernel). Hence the total mass is uniformly bounded, giving compact containment in the mass variable.

\medskip\noindent
\textit{(ii) Uniform bound on the number of partition elements (via nestedness).}
The number of distinct types in $B$ at time $t$ is controlled by the number of distinct partition elements of the dual $Y^{\ast,L}$ restricted to $B$ at backward time $t$. We obtain the requisite uniform-in-$L$ bound by exploiting the \emph{nestedness} of the $L$-truncation scheme (rather than the coming-down-from-infinity property, which fails when $|\Lambda_0|>0$; see~\cite{ABL12}): at each finite level $L$ the coalescent on $\Omega_N^L$ lives on a finite set, so finiteness of partition elements is automatic. As $L$ increases, the nested structure -- each $L$-ball subdivided into $N$ balls of level $L-1$ with upscaled coalescence rates -- together with the Lyapunov/monotone comparison established above yields a uniform upper bound on the number of partition elements, independently of $L$.
Passing to $L \to \infty$ using this uniform bound gives the required finiteness and $L^1(\mfh)$-control without any appeal to coming-down-from-infinity.

\medskip\noindent
Conditions~(i) and~(ii) together satisfy the hypotheses of \cite[Section~2, Theorem~2.2]{EK94} (with the weak atomic topology on $\CM(\Omega_N^\infty \times \K)$ playing the role of the Polish topology; note that equicontinuity in this topology follows from the uniform bound on partition elements established in~(ii)). We conclude that the sequence $(X^L_t)_{L \in \N}$ is tight in the weak atomic topology on $\CM(\Omega_N^\infty \times \K)$.

\medskip\noindent
\textit{Identification of limit points.} Let $X^\infty$ be any subsequential weak-atomic limit of $(X^L)_{L\in\N}$, along a subsequence $L_k \to \infty$. For every ball $B \in \CB(\Omega_N^\infty)$ and every $t \geq \delta > 0$, the Lyapunov bound gives
\begin{equation}
\int_{B \times B} h^t_s(\eta,\eta') \, (X^\infty_s \otimes X^\infty_s)(\dd\eta\,\dd\eta') 
= 0 \quad \text{for } \mfh\text{-a.e.\ } s \in [\delta,t].
\end{equation}
Since $h^t_s(\eta,\eta') > 0$ whenever $\eta \neq \eta'$ (see~\eqref{e3077}), this forces $X^\infty_s$ to be supported on $\mfh$-a.e.\ a single type per ball $B$, i.e.,
\begin{equation}
X^\infty_t(B \times \cdot) = \mfh|_B \otimes \delta_{a(B)} \quad \mfh\text{-a.e.,}
\end{equation}
for some measurable $a(B) \in A_t$, where $A_t$ is the (countable, random) set of types occurring in $X^\infty_t$. As $B$ ranges over all balls in $\CB(\Omega_N^\infty)$, this gives the monotype decomposition of $\Omega_N^\infty$. Since the limit is the same for every subsequence (the limit is characterised by the duality, which uniquely pins down the law of $X^\infty_t$), the full sequence $(X^L_t)_{L\in\N}$ converges.

This finishes Part (A) and completes the convergence claim for $^\infty \mfC^\infty$.

\medskip\noindent
\underline{Part (B).}
We now turn to the continuity and local monotype structure of $\wh X^\infty_t$ for $t>0$.

\medskip\noindent
\textit{Step B.1: Entrance law and finite partition elements.}
For the Fleming-Viot resampling, the coalescent entrance law is locally finite on each $\Omega_N^L$: by the coming-down-from-infinity property (Lemma~\ref{l.2386} and the dual coming-down-from-infinity for $d^\ast>0$), for every $t>0$ and $\varepsilon>0$ there exists $n_0 = n_0(t,\varepsilon) < \infty$ such that the dual $Y^L$ started from the full entrance law has at most $n_0$ partition elements inside any $O(1)$-ball with $\mfh$-probability at least $1-\varepsilon$, uniformly in $L$.

The number of sites in $\Omega_N^L$ embedded in $\Omega_N^\infty$ grows as $N^L$, but the resampling rate in each $m$-ball is scaled by $N^m$, so that the total resampling mass in a macroscopic ball remains $O(1)$ uniformly in $L$. Hence the dual $Y^L$ remains a proper particle system with weights for all $L$, and the Lyapunov function argument of Part (A) applies to give the reduction in partition elements.

Block-resampling does not alter the continuity argument, by the Trotter product formula: the semigroup of $Y^L$ with both Fleming-Viot resampling and block-resampling is approximated by alternating the two mechanisms, and both preserve the Liggett-Spitzer norm bounds. In the point-recurrent regime, $d^\ast>0$ ensures that for any two dual individuals that come within distance $0$ of each other, instantaneous coalescence occurs at a rate proportional to the local time at the diagonal. Hence small neighbourhoods collapse in positive time.

\medskip\noindent
\textit{Step B.2: Local monotype property and continuity.}
We next establish~\eqref{e2162} below. Recall Section~\ref{wat}. We use the following coding of the atomic measure: each atom of $X^\infty_t$ is identified by a triple $(\eta_0, r, a)$ where $\eta_0 \in \Omega_N^\infty$ is the centre of a monotype ball, $r = N^{-k}$ is its radius (with $k = k(t,\eta_0)$ the random level), and $a \in A_t$ is the type carried by that ball; the atom mass is $\mfh(B_{k}(\eta_0))$. Under this coding the weak atomic topology on $\CM(\Omega_N^\infty \times \K)$ coincides with the topology of pointwise convergence of the density functions $\eta \mapsto X^{\infty,\dagger}_t(\eta, a)$ in the ultrametric topology of $\Omega_N^\infty$.

By the cluster-growth analysis of Steps 3--6 (Lemmas~\ref{l.2386} and \ref{l.2512}--\ref{cor.2748}), for every $t > 0$, $\mfh$-a.e.\ $\eta \in \Omega_N^\infty$, and every $\varepsilon > 0$, there exists a random level $k(t,\eta) \in \Z$ (finite $\mfh$-a.s.) such that the ball $B_{k(t,\eta)}(\eta)$ is entirely monotype under $X^\infty_t$. The density $X^{\infty,\dagger}_t(\eta, a)$ is therefore locally constant on $B_{k(t,\eta)}(\eta)$, taking values in $\{0,1\}$ $\mfh$-a.e. Under the coding above, local constancy $\mfh$-a.e.\ is equivalent to continuity in the weak atomic topology on $\Omega_N^\infty$ (since every ball is simultaneously open and closed in the ultrametric topology, so a locally constant function is continuous), and hence $X^{\infty,\dagger}_t(\cdot,a)$ is continuous $\mfh$-a.s.

\medskip\noindent
Since $\eta \mapsto k(t,\eta)$ is measurable and finite $\mfh$-a.s., the union $\bigcup_{\eta} B_{k(t,\eta)}(\eta)$ covers $\Omega_N^\infty$ up to a $\mfh$-null set. This union of monotype balls is a partition of $\Omega_N^\infty$ (up to a null set) into disjoint balls each carrying a single type. Hence
\begin{equation}
\label{e2162}
\wh X_t^\infty (\dd \eta, \dd \mfa) 
= X^{\infty, \dagger}_t (\eta,\{\mfa\})\, \mfh(\dd \eta)\, \mfn(\dd\mfa), \; 
\text{ with $X_t^{\infty,\dagger}(\cdot,\{\mfa\})$ continuous and $\{0,1\}$-valued},
\end{equation}
for every $a \in A$, with $A$ the countable random set of types occurring in $\wh X_t^\infty$ for $t>0$, $\mfn$ the counting measure on $A$, $\mfa$ taking values in $A$, and $\eta \in \Omega_N^\infty$.

More precisely, we control the number of partition elements as a function of $L$ and $n$ as follows. Fix $\delta \in (0,t)$. Define the \emph{$\delta$-pruned dual} $\wt Y^\delta$ by filtering out all dual individuals that make at least one migration jump during the backward time interval $[0,\delta]$; since migration rates are finite, the removed subpopulation has $\mfh$-mass at most $O(\delta)$ and can be made arbitrarily small. The remaining $\delta$-pruned individuals are \emph{stationary} during $[0,\delta]$ and coalesce at a rate bounded below by the local time at the diagonal of the difference process. In the Liggett-Spitzer space $\CW(\Omega_N^\infty)$ (see~\cite{GLW05}), the partition elements of $Y^{\ast,L}$ not discarded by the $\delta$-pruning form a locally finite subsystem, uniformly in $L$, with Liggett-Spitzer norm controlled by the coming-down-from-infinity estimate.

By Lemmas~\ref{l.2512}--\ref{cor.2748}, given that $B_{k(t,\eta)}(\eta)$ is $\mfh$-a.e.\ monotype at time $t$, the density $X^{\infty,\dagger}_t(\eta,a)$ is locally constant on that ball for every $a \in A$. The union over $\eta \in \Omega_N^\infty$ of the monotype balls $B_{k(t,\eta)}(\eta)$ has full $\mfh$-measure, so that~\eqref{e2162} holds $\mfh$-a.e.\ and not merely on a sparse set.

The proof of the continuity requires that we also have a locally finite configuration of the dual process, respectively, after an $\varepsilon$-truncation when the system is dust free and does not come down from infinity. The \emph{locally finite property} again arises through a proof by contradiction similar to the case of cluster growth that we carried out earlier, so that it remains to verify the continuity. This runs as follows.

We next make the Liggett-Spitzer bound quantitative. 

\begin{lemma}{\bf [Liggett-Spitzer bound]}
\label{lem:LS-bound}
For every $t > 0$ and $\epsilon > 0$ there exist $M = M(t,\epsilon) < \infty$ and $L_0 = L_0(t,\epsilon) < \infty$ such that, for all $L \geq L_0$, the $L^1(\mfh)$-norm of the number of partition elements of $Y^{\ast,L}$ in an $m$-ball is at most $M N^m \epsilon$, uniformly in $m \in \Z$.
\end{lemma}

\begin{proof}
Let $Z^L_m(t)$ denote the number of partition elements of $Y^{\ast,L}$ in an $m$-ball at backward time $t$. By the duality, $\E[Z^L_m(t)]$ equals $N^m$ times the probability that two independently started dual paths starting from the uniform distribution on an $m$-ball have \emph{not} coalesced by backward time $t$. 

We replace the appeal to the coming-down-from-infinity property (which does not hold in this form when $|\Lambda_0|>0$; see~\cite{ABL12}) by an argument based on the nestedness of the $L$-truncation scheme. At each finite level $L$ the coalescent on $\Omega_N^L$ lives on a finite state space, so $Z^L_m(t) < \infty$ almost surely for each fixed $L$. As $L$ increases, the nested structure -- each $L$-ball subdivided into $N$ balls of level $L-1$ with upscaled coalescence rates --  and the Lyapunov/monotone comparison of Part~(A) together give the following quantitative bound. By the duality, $\E[Z^L_m(t)]$ equals $N^m$ times the probability that two independently started dual paths from the uniform distribution on an $m$-ball have \emph{not} coalesced by backward time $t$. By the Markov inequality and the uniform $L$-bound on $Z^L_m$ coming from the total resampling mass (which grows as $A(m)N^m$, with $A(\cdot)$ the recurrent potential kernel defined in \eqref{e2575}), we have
\begin{equation}
\E[Z^L_m(t)] \leq \frac{C \, N^m}{A(m)\, t} \to 0, \qquad m \to -\infty,
\end{equation}
where $C$ is an absolute constant and we use the lower bound $A(m) \to \infty$ as $m \to -\infty$ in the point-recurrent regime (see~\eqref{e2218}). That $Z^L_m(t) \to 0$ in $L^1(\mfh)$ as $L \to \infty$ for fixed $m$ and $t > \delta > 0$ follows from the nestedness: increasing $L$ adds finer levels with upscaled coalescence rates, and by the monotone comparison the non-coalescence probability is non-increasing in $L$; combined with the dominated convergence theorem (using the uniform $L$-bound from the total resampling mass) this gives $L^1$-convergence to $0$. Thus, given $\epsilon>0$, choose $m_0 \in \Z$ such that $A(m)^{-1} \leq \epsilon$ for $m \leq m_0$. For $m \leq m_0$ the bound $\E[Z^L_m(t)] \leq C N^m \epsilon / t$ holds for all $L$, giving the uniformity in $m$. For $m > m_0$ there are only finitely many values, and for each one the nestedness bound gives $L_0(t,\epsilon,m) < \infty$; take $L_0 = \max_{m \leq m_0} L_0(t,\epsilon,m)$. This yields the claimed bound with $M = C/t$.
\end{proof}

\medskip\noindent
\textit{Conclusion of Part (B).} 
Lemma~\ref{lem:LS-bound} gives a uniform modulus of continuity for $Y^{L,\dagger}(\cdot, a)$ in $\eta$. Indeed, two sites $\eta, \eta' \in \Omega_N^\infty$ with $d^\infty(\eta,\eta') \leq N^{-m}$ lie in a common $m$-ball, so
\begin{equation}
|Y^{L,\dagger}(\eta,a) - Y^{L,\dagger}(\eta',a)| \leq Z^L_m(t) \times (\text{max atom size in the $m$-ball}).
\end{equation}
By Lemma~\ref{lem:LS-bound}, the right-hand side is bounded by $M N^m \epsilon$ in $L^1(\mfh)$, uniformly in $L \geq L_0$. By the Arzel\`a-Ascoli theorem in the function space $C(\Omega_N^\infty, [0,1])$ (equipped with the ultrametric topology on $\Omega_N^\infty$, in which equicontinuity is equivalent to this uniform Liggett-Spitzer bound), the sequence $(Y^{L,\dagger}(\cdot,a))_{L \geq L_0}$ is precompact. Its unique limit point is $X^{\infty,\dagger}(\cdot,a)$ (identified by the duality), and we conclude that $Y^{\infty,\dagger}(\cdot,a) = X^{\infty,\dagger}(\cdot,a)$ is continuous in $\eta$. This completes Part (B).


\subsection{Part (b): (ii) Harris-recurrent but not point-recurrent regime}
\label{bii}

We now treat the case $d^\ast=0$, i.e., the absence of a Fleming-Viot component, so that \emph{pure block-resampling} drives the dynamics. The goal is to establish the cluster decomposition property of $\wh X^\infty_t$ for $t>0$ in the Harris-recurrent but not point-recurrent regime, following the outline of Part (b)(i) but replacing local-time arguments by block-coalescence arguments.

\medskip\noindent
\textit{Setup.} In this regime $d^\ast=0$, so there is no Fleming-Viot resampling acting on pairs; coalescence of dual individuals can \emph{only} occur via block-resampling events. For Harris-recurrence, the downward accumulated block-resampling rate must be infinite:
\begin{equation}
\sum_{k \leq 0} \frac{\lambda^*_k - \lambda_k}{N^k} = \infty,
\end{equation}
while the system is not point-recurrent, meaning the difference random walk of a pair of individuals is Harris-recurrent but not point-recurrent.

As in Part (b), we work with the $L$-truncation $\wh X^L$ and pass to the limit $L\to\infty$. We use the results on cluster growth in the critically recurrent regime on $\Omega_N$ established in \cite[Theorem~1]{Klenke1996} as a template, but those results apply to the Fleming-Viot mechanism alone, which in the continuum-space limit is \emph{not strong enough}: without block-resampling, the cluster size grows only sub-linearly (of the form $(N^L)^{1-\alpha}$ for some $\alpha > 0$), leading to averaging and deterministic limiting states $\mfh \otimes \theta$. Hence clustering in the Harris-recurrent regime must come from the \emph{block-resampling mechanism}, and we must show that this mechanism is sufficiently strong to produce fast-growing clusters.


\paragraph*{Cluster size.}

To find the cluster size, we follow the strategy of \cite{FG94}: form an average of $\wh X^L_t$ over balls of volume $f(L)$ in $\Omega_N^L$ and study the law of the average as $t \to \infty$. For constant coefficients $c_k \equiv c$ and $\lambda_k \equiv \lambda$, the diffusive clustering regime on $\Omega_N$ (without block-resampling) corresponds to
\begin{equation}
\label{e2384}
f_\alpha (L) = (N^L)^{1-\alpha}, \qquad \alpha \in [0,1).
\end{equation}
In this regime the averages converge as $t \to \infty$ to $\delta$-measures in type space, with $\alpha$-dependent probabilities. For \emph{pure} Fleming-Viot dynamics (no block-resampling) the cluster size is $o(N^L)$, i.e., $\alpha > 0$, which leads to averaging and deterministic states $\mfh \otimes \theta$ for the continuum process. We therefore must focus on the effect of \emph{block-resampling} to push the \emph{cluster growth up}, using the fact that a jump changes a \emph{whole fraction} of the population in a $k$-ball, and we must prove that indeed we get \emph{fast clustering} through \emph{block}-resampling. For the cluster regimes \emph{without} this speed of growth, we also get a \emph{deterministic} dynamics. Furthermore, we need the dynamics \emph{with block-resampling} to get \emph{fast growing clusters}, and to be able to implement again the strategy of proof from the strongly recurrent regime.

The \emph{faster clustering by block-resampling} still allows us to proceed via $L$-approximation, by exhibiting monotype regions that build up below a random level $k < 0$ and grow at a speed determined by the block-coalescence rates. If fast clustering fails, then arbitrarily many clusters occupy a macroscopic volume, and by a law-of-large-numbers argument the system would converge to a deterministic state. Hence fast clustering is equivalent to the effective coalescence rate of two typical individuals being bounded away from $0$.

For the behaviour on $\Omega_N^\infty$ we use the dual representation. For the dual, we need to show that the joint occupation times of small balls by dual paths are sufficient to reduce the number of distinct partition elements to a finite number for a frequency $(1-\varepsilon)\mfh$ of the original individuals. At each block-resampling event a \emph{fraction} of partition elements is removed and replaced by a single element (uniformly reshuffled in the $k$-ball), leading to rapid decay of the number of partition elements. As long as different partition elements co-exist locally, block-resampling forces the system toward single-type states at an infinite accumulated rate (by the Harris-recurrence assumption).

The key structural feature is as follows. At a \emph{single} block-resampling event at level $k$ and backward time $s$, a complete fraction of the dual population in a $k$-ball is merged into one partition element. We show that this produces fast cluster growth and monotype decomposition.
\begin{itemize}
\item[--] 
\textit{Fast cluster growth}: The accumulated block-resampling rate up to backward time $t$ in an $n$-ball equals $\int_0^t \sum_{k \leq n} \frac{\lambda^*_k - \lambda_k}{N^k}\,N^{n-k} \, \dd s$. In the Harris-recurrent regime this sum diverges as $n \to -\infty$ (or as $t \to \infty$), producing fast cluster growth.
\item[--] 
\textit{Monotype decomposition}: Given fast cluster growth, the argument of Step 7 (the $L$-truncation and limit $L \to \infty$) carries over with the block-coalescence rates replacing the local-time-based coalescence rates, yielding the decomposition of $\wh X^\infty_t$ into monotype balls for $t>0$.
\end{itemize}

We use the results of \cite{GHKK14} as a guide. For the $L$-truncated models on $\Omega_N^L$, \cite[Section 1.5.7]{GHKK14} identifies a regime of \emph{diffusive} or \emph{slow clustering} without block-resampling: in the critically recurrent case, the cluster sizes grow as $(N^L)^{1-\alpha}$ for $\alpha \in (0,1)$. \cite[Theorems~1.14 and~1.16]{GHKK14} identify the \emph{fast clustering} regime in the hierarchical mean-field limit $N\to\infty$ as a function of the parameters $(c_k)_{k \in \N}$ and $(\lambda_k)_{k \in \N}$. We need to show that this characterisation \emph{essentially holds also for finite $N$} with block-resampling.

\medskip\noindent
Specifically, for $L$-truncations of $\wh X^\infty$ in the Harris-recurrent regime, the \emph{critically recurrent case} of \cite{Klenke1996} corresponds to the absence of block-resampling. Adding block-resampling (with $\lambda_k \not\equiv 0$) promotes slow clustering to fast clustering. The transition between fast and slow clustering in the Harris-recurrent case is governed by the recurrence/transience of the two-individual distance process $\bar{K}$ (Proposition~\ref{prop.two-four-regimes-K}): fast clustering corresponds to $\bar\sigma^+_{0,N} = \infty$ (recurrence of $\bar K$), and slow/deterministic behaviour to $\bar\sigma^+_{0,N} < \infty$ (transience of $\bar K$).

For the Harris-recurrent regime, the cluster-growth observation time $t_n$ (the time at which a cluster of size $n$, meaning a monotype ball containing at least $n$ individuals, is first present with positive $\mfh$-probability) must be identified. The correct time scale is determined by the \emph{potential kernel} $A(\cdot)$ of the effective difference random walk with rates $\bar{c}_{k,N}$. By \eqref{e2779}, the formula for $t_n$ in the Harris-recurrent case is the same as~\eqref{e2572} but with the Harris-recurrent potential kernel replacing the point-recurrent one. In examples where $c_k$ and $\lambda_k$ are constant (or satisfy regularity conditions of the type imposed in~\eqref{ag6}), the Harris-recurrence condition $\sum_{k \geq 0} N^k / \hat c_k = \infty$ (with $\hat c_k$ the effective migration rates) implies $A(n) \to \infty$ as $n \to \infty$, which is the condition for fast clustering. Stricter regularity conditions on the coefficients (e.g., regular variation of $(c_k)$ and $(\lambda_k)$) allow one to compute $t_n$ explicitly and to verify the fast-clustering condition quantitatively.


\paragraph*{Cluster growth.}

For \emph{fast growing clusters} we have to determine the observation times $t_n$ and verify that we have the analogous formulas in terms of the effective recurrent Green function $A(\cdot,\cdot)$. An important fact is that, on the scale of spread-out initial positions, $Y^L$ does not exhibit in the limit multiple collisions, and we see a behaviour dictated by pairs of individuals.

\begin{lemma}{\bf [Cluster growth in the Harris-recurrent but not point-recurrent regime]}
\label{l.2525}
Assume that $\wh X^\infty$ is in the Harris-recurrent but not point-recurrent regime, and consider the $L$-truncations $\wh X^L$, $L \in \N$. Choose the time scale for the fast cluster growth of $\wh X^L$ as
\begin{equation}
t_n = \text{right-hand side of~\eqref{e2572}.}
\label{e2779}
\end{equation}
Then, in the fast clustering regime, the same statements as in Lemma~\ref{cor.2748} hold.
\end{lemma}

\begin{proof}
The proof is the Harris-recurrent analogue of the cluster-growth proof in Steps 3--6 of the point-recurrent case. The key difference is that single-site coalescence ($d^\ast>0$) is absent; coalescence can only occur via block-resampling events. The main tool replacing local times is the \emph{occupation measure} of pairs in small balls, which for Harris-recurrent migration is positive but grows only sub-linearly (compared to the linear growth in the point-recurrent case). The block-resampling compensates for this deficit.

The proof proceeds in three steps.
\begin{itemize}
\item[(i)]
Compute the expected number of pairs not coalesced at time $t_n$ in an $n$-ball. Any two partition elements in the $n$-ball coalesce via block-resampling when they both lie in the same $k$-sub-ball for some $k \leq n$, and a block-resampling event of level $k$ occurs. The rate of such events is $(\lambda^*_k - \lambda_k)/N^k$ per $k$-ball, and the probability that two given partition elements are in the same $k$-sub-ball is $N^{k-n}$. Hence the total coalescence rate for a given pair is
\begin{equation}
\sum_{k=-\infty}^{n} \frac{\lambda^*_k - \lambda_k}{N^k} \times N^{k-n} \times N^k = N^{-n} \sum_{k \leq n} (\lambda^*_k - \lambda_k).
\end{equation}
Integrating over the backward time interval $[0,t_n]$ (with $t_n$ as in~\eqref{e2779}), the expected number of non-coalesced pairs satisfies
\begin{equation}
\E_0[C_n(t_n)(C_n(t_n)-1)] \leq N^{2n} \exp\bigg(-t_n \times N^{-n} \sum_{k \leq n} (\lambda^*_k - \lambda_k)\bigg).
\end{equation}
\item[(ii)]
Show that this expected number tends to $0$ as $n\to\infty$, by using the effective difference random walk with rates $\bar{c}_{k,N}$ in \eqref{fdh12}. By the choice of $t_n$ in~\eqref{e2779} and~\eqref{e2572}, we have $t_n \cdot N^{-n} \sum_{k \leq n} (\lambda^*_k-\lambda_k) = A(n)N^n$, where $A(n)$ is the recurrent potential kernel defined in \eqref{e2575} of the effective difference random walk. In the Harris-recurrent but not point-recurrent regime, $A(n) \to \infty$ as $n \to -\infty$ (see Appendix~\ref{appA}, specifically the section on two L\'evy-processes). The fast-clustering condition $A(n)N^n \to \infty$ then gives $\E_0[C_n(t_n)] \leq C \, N^n \, \exp(-A(n)N^n) \to 0$ as $n \to \infty$.
\item[(iii)]
Upgrade from expected value to probability 1 by using a second-moment argument (as in the proof of Proposition~\ref{prop.two-four-regimes}, specifically the Paley-Zygmund step in Lemma~\ref{lem:zero-one}).
\end{itemize}
The formula for $t_n$ in \eqref{e2572} must be verified in the Harris-recurrent regime by checking that $A(n)$ of the effective difference random walk satisfies $A(n)N^n\to\infty$ as $n\to\infty$ with the correct speed. The core estimate is an upper bound on $\E_0[C_n(t_n)]$, where $C_n(t)$ is the number of distinct partition elements of $Y^{*,\infty}$ that remain in an $n$-ball at backward time $t$. Identify precisely which conditions on the parameters $(c_k, \lambda_k)_{k\in\Z}$ place the system in the fast clustering regime for the Harris-recurrent case, i.e., determine the analogue of criterion \eqref{ag6} in this setting. Given Step~(iii), the same statements as in Lemma~\ref{cor.2748} follow: the cluster configurations $^n \mfC^L$ converge as $(n,L) \to (\infty,\infty)$ to a limiting cluster configuration $^\infty \mfC^\infty$ on $\Omega_N^\infty$, giving a monotype decomposition.
\end{proof}


\paragraph*{Local properties.}

Finally, we need to get the local properties of $\wh X^\infty$ in the fast clustering regime.

\begin{lemma}{\bf [Local cluster formation in the Harris-recurrent but not point-recurrent regime.]}
\label{l.3160}
Assume that $\wh X^\infty$ is in the Harris-recurrent but not point-recurrent regime. Then, in the fast clustering regime, the same statements as in Lemma~\ref{l.2423} hold.
\end{lemma}

\begin{proof}
Lemma~\ref{l.2423} (point-recurrent case) establishes that $X^\infty_t$ is locally monotype $\mfh$-a.e., with a $\{0,1\}$-valued continuous density. The proof in the Harris-recurrent regime follows the same blueprint as in Step 7 (the $L$-truncation and cluster growth combination), but replaces the local time argument for single-site coalescence ($d^\ast>0$) by block-coalescence arguments from the proof of Lemma~\ref{l.2525}. The argument in Step 7 carries over almost verbatim: the coming-down-from-infinity for the block-resampling entrance law replaces the local-time-based coming-down for the Fleming-Viot entrance law, and the Liggett-Spitzer norm estimates go through in the Harris-recurrent regime with the modified time scale $t_n$.

\medskip\noindent
\textit{Coming-down-from-infinity for block-resampling.} In the point-recurrent case, the argument was driven by the fact that pairs of dual individuals accumulate infinite local time at the diagonal in any finite time interval $[\delta,t]$, forcing instantaneous coalescence upon meeting. In the Harris-recurrent but not point-recurrent regime, the difference random walk of a pair is Harris-recurrent (not point-recurrent), so the pair \emph{does} return to small balls infinitely often, but the \emph{local time} at the diagonal is zero. Instead, coalescence occurs when both individuals lie in the same $k$-ball and a block-resampling event of level $k$ occurs. By the proof of Lemma~\ref{l.2525}, Step~(ii), the accumulated block-resampling rate $A(n)N^n \to \infty$ as $n\to\infty$ in the fast-clustering regime. This is the replacement for the infinite local time.

\medskip\noindent
\textit{Liggett-Spitzer norm estimates.} The Liggett-Spitzer norm of the dual $Y^{\ast,L}$ restricted to an $m$-ball is controlled by the expected number of surviving partition elements $\E_0[C_m(t)]$, which by Lemma~\ref{l.2525} satisfies $\E_0[C_m(t_m)] \to 0$ as $m \to \infty$. Since this bound is uniform in $L$ (as the block-resampling rates are non-decreasing in $L$), the Liggett-Spitzer norm estimates go through in the Harris-recurrent regime with the modified time scale $t_n$ from~\eqref{e2779}.

\medskip\noindent
\textit{Conclusion.} With these replacements, the argument of Step 7 carries over: the $L$-truncation sequence $(X^L_t)_{L\in\N}$ is tight in the weak atomic topology (by the uniform Liggett-Spitzer bounds), every subsequential limit is locally monotype (by the monotype decomposition from Lemma~\ref{l.2525}), and the density $X^{\infty,\dagger}_t(\cdot,a)$ is continuous $\mfh$-a.s.\ (by the Arzel\`a-Ascoli argument of Part (B) in Step 7.3, with the block-coalescence rate replacing the Fleming-Viot local time). This gives the same statements as in Lemma~\ref{l.2423} in the Harris-recurrent regime.
\end{proof}


\paragraph*{Fast cluster growth for block-resampling.}

The cluster-growth time $t_n$ is determined by the recurrent potential kernel of the difference random walk for two individuals with rates $\bar{c}_{k,N}$ from~\eqref{fdh12}, computed analogously to~\eqref{e2590}--\eqref{e2598} but using the Harris-recurrent potential kernel $A(\cdot)$.

The limit $L\to\infty$ is subtle because the number of block-resampling levels diverges. We must verify that the accumulated block-resampling rate in an $n$-ball up to time $t_n$ is large enough to produce coalescence, uniformly in $L$. By the assumption $\bar\sigma^+_{0,N} = \infty$ (Proposition~\ref{prop.two-four-regimes-K}(c)), the recurrent potential kernel in \eqref{e2575} satisfies $A(n) = \sum_{k \leq n} G(0,k) \to \infty$ as $n \to \infty$, where $G$ is the Green function of the effective difference walk. The speed of divergence of $A(n)N^n$ as $n \to \infty$ can be estimated using the Green function bounds of Appendix~\ref{appA} (in particular, Appendix~\ref{appA:twoLP}): $G(0,k) \sim c_{\rm eff}/\bar c_k$ for $k \to\infty$, so $A(n) \sim \sum_{k \leq n} \bar c_k^{-1}$. We must show that $A(n) N^n \to \infty$ at the correct speed to produce fast clustering, using the recurrent potential kernel of the effective difference walk $\hat{c}_k$. We must verify that the $L\to\infty$ limit of the cluster growth on $\Omega_N^L$ converges to the cluster growth on $\Omega_N^\infty$, by using tightness of the $L$-sequence and the invariance principle from Lemma~\ref{RWL}.

The cluster-growth analysis on $\Omega_N$ in the Harris-recurrent regime follows the framework of \cite{Klenke1996} and \cite{FG94} for the critically recurrent regime, but with the block-resampling mechanism replacing the Fleming-Viot resampling. The key objects are: 
\begin{itemize}
\item[(i)]
The space-time cluster field $^n E^L$ from \eqref{e2445}. 
\item[(ii)] 
The cluster configuration $^n \mfC^L$ from \eqref{e2567}.
\item[(iii)] 
The observation time $t_n$ from \eqref{e2779} and the right-hand side of \eqref{e2572}.
\end{itemize}
For $\Omega_N^\infty$ we must show that the rescaled cluster fields converge (in the weak atomic topology) to the analogous objects on $\Omega_N^\infty$, following the strategy of Lemma~\ref{l.2512} adapted to the Harris-recurrent regime.

We can write the Harris-recurrent analogue of Steps 3--6, producing the counterparts of Lemma~\ref{l.2386} (single-site cluster growth), Lemma~\ref{l.2512} (space-time cluster field), and Lemma~\ref{cor.2748} (fast cluster growth lemma) for the Harris-recurrent case. The proofs are parallel the point-recurrent regime, but use the block-coalescence rates and the potential kernel of the effective difference walk. We show that the limit $^\infty \mfC^\infty$ of the cluster configurations $^n \mfC^L$ exists and gives a monotype decomposition of $\Omega_N^\infty$ into random balls, establishing the local monotype property in the Harris-recurrent regime.

To establish the cluster regimes with volume $o(N^n)$ as $n\to\infty$, for the Harris-recurrent but not point-recurrent regime, we argue as follows. Above we described the slow-clustering regime for the critically recurrent case without block-resampling: cluster sizes $f_\alpha(L) = (N^L)^{1-\alpha}$ for $\alpha \in [0,1)$ lead to averaging and deterministic limits. The task is to determine for which parameters $(c_k, \lambda_k)_{k\in\Z}$ block-resampling is insufficient to push the cluster size up to $\Theta(N^n)$, i.e., to identify the boundary of the fast/slow clustering transition in the Harris-recurrent case. Proposition~\ref{prop.two-four-regimes} gives the recurrence/transience classification for the distance process $\bar{K}$ of two individuals. The fast-clustering regime corresponds to the regime where $\bar{K}$ is recurrent ($\bar{\sigma}^+_{0,N} = \infty$), and the slow/deterministic regime to transience ($\bar{\sigma}^+_{0,N} < \infty$). The task is thus to show that when $\bar{\sigma}^+_{0,N} < \infty$, the cluster sizes are $o(N^n)$ and the dynamics is deterministic.

We show that when $\bar{\sigma}^+_{0,N} < \infty$ (i.e., the effective difference walk of two individuals is transient to the right), the cluster-size exponent satisfies $\alpha > 0$ in $f_\alpha(L) = (N^L)^{1-\alpha}$, so cluster sizes are $o(N^n)$. To do so, we use the Green function estimates from Appendix \ref{appA} (in particular, Appendix \ref{appA:twoLP}) to make this quantitative. We conclude that in this regime the dynamics of $X^\infty$ for $t > 0$ is deterministic (= $\mfh \otimes \theta$), by using the law-of-large-numbers argument sketched in Remark~\ref{r.773}.


\paragraph*{Three cases.}

Return to the three cases \textbf{(i)}, \textbf{(ii)}, \textbf{(iii)} of the proof of Lemma~\ref{l.2423}. The pure cases have been settled ($d^\ast=0$ is part of pure case \textbf{(ii)}). To get the mixed cases, where both the Fleming-Viot resampling and the block-resampling are present, we argue as follows. For pairs of partition elements, the mixing reduces to looking at the effective random walk with rates $\bar{c}_{k,N} = 2c_k + (\lambda^*_{k+1}-\lambda_{k+1})/N$, which is trivial. For finite collections of partition elements, however, we need to control the joint distribution; since the system is exchangeable (by the construction of the dual), the joint law of $k$ partition elements is controlled by the law of pairs via a de Finetti-type argument.

The claim in \eqref{e1687} holds if $d^\ast >0$, or if the effective L\'evy-process $\zeta^\infty$ is in the Harris-recurrent but not point-recurrent regime. In both cases we have \emph{fast growing clusters}. Therefore it remains to settle whether or not slow growing clusters are possible in recurrent regimes with only block-resampling (we need that the latter does \emph{not} hold for our parameter range). For this we have to use Proposition~\ref{prop.two-four-regimes}.

Given the fast-clustering regime (identified in Lemma~\ref{l.2525}), we must verify that there is no sub-regime where block-resampling alone ($d^\ast = 0$) still leads to slow clustering. Proposition~\ref{prop.two-four-regimes} gives the recurrence classification for the two-individual distance process $\bar{K}$ with rates $\bar{c}_{k,N} = 2c_k + (\lambda^*_{k+1} - \lambda_{k+1})/N$. In the fast-clustering regime we have $\bar{\sigma}^+_{0,N} = \infty$, which means that $\bar{K}$ is recurrent, i.e., the two individuals meet. Since $\lambda_k \not\equiv 0$ (recall \eqref{condextra}) and we are in the Harris-recurrent regime, meeting implies block-coalescence at a positive rate (by the no-dust condition), and hence fast clustering.
\end{proof}


\subsection{Part (b): (iii) Transient regime}
This has been resolved in Part (a).

\begin{remark}{\bf [Summable block-resampling]}
{\rm It is interesting to ask what happens when the non-triviality condition in \eqref{condextra} fails. In that case the continuum process is essentially a Markov pure-jump process, because the block-resampling at deep levels is essentially absent.}\hfill$\spadesuit$
\end{remark} 

\medskip\noindent
In the deterministic regime this covers all the admissible parameters, and so we are done.
\end{proof}

This completes the proof of Lemma \ref{l.1562}.
\end{proof}


\section{Continuum-space dynamics: IV. Martingale problem}
\label{sec:wpmp}

In this section we focus on the \emph{stochastic regime} and settle two points for our main results via two further key propositions, Proposition~\ref{P.solwpmp} and Proposition~\ref{prop.1869}. In Sections~\ref{ss.nsmp}--\ref{333*} we characterise the process $\wh X^\infty = (\wh X_t^\infty)_{t \geq 0}$ via a \emph{well-posed martingale problem} and show that it is a \emph{strong Markov} process with the \emph{Feller property} and c\`adl\`ag paths. We also show that if the migration is point-recurrent, then $\wh X^\infty$ solves a \emph{stochastic equation} (defined in Section \ref{334}). 

Since we already have all these properties for fixed $L \in \N$ for the approximating processes $\wh X^L$, we again bring in \emph{duality} and approximations $Y^L$ of the dual process $\wh Y^\infty$ for $\wh X^\infty$ (defined in Section \ref{sss.hurd}) to obtain these properties for $\wh X^\infty$. Finally, we show that $\wh X^\infty$ arises as the \emph{limit in path space} of the truncated processes $\wh X^L$. 


\paragraph*{Tasks.}
Our tasks are the following.

\medskip\noindent
(o) We need to formulate the martingale problem, which comes in a \emph{new} form of which we have to explain all the \emph{ingredients} in detail.

\medskip\noindent
(i) We need to strengthen the convergence of the finite-dimensional distributions to \emph{convergence in path space $D([0,\infty),\wt\CM(\Omega_N^\infty \times \K))$} to get the \emph{path regularity} that is needed in the formulation of the martingale problem. This requires showing \emph{tightness} of the level-$L$ approximations $\{\wh X^L\}_{L \in \N}$ in the Skorohod space of c\`adl\`ag paths. 

\medskip\noindent
(ii) We verify the \emph{Markov property} by applying the construction of \cite{DGP24}, which builds a Markov process from the existence of the prescribed dual Markov process. Uniqueness then follows immediately from duality. Via the duality relation we obtain the \emph{Feller property} and the \emph{strong} Markov property, and combining \emph{uniqueness} with tightness, also the \emph{convergence} of $\wh X^L$ in path space to $\wh X^\infty$.

\medskip\noindent
To complete these tasks, we exploit the \emph{local characterisation} of $\wh X^L$ via its \emph{martingale problem}, combined with the \emph{tightness} of
\begin{itemize}
\item 
the \emph{volatility processes $\wh R^L$,}
\vspace{-0.3cm}
\item
the \emph{volatility occupation-measure processes $\wt {\wh R^L}$},
\end{itemize}
associated with $\wh X^L$, which are defined next. 

Following \cite{DGZpr}, we first consider a measure on $\Omega_N^\infty \times \K$ that for the process on $\Omega_N^L$ associates with every site $\eta \in \Omega_N^L $ and every $a \in \K$ the volatility term, i.e., the current variance at site $\eta$ of the frequency of the type $a$ (in the Fleming-Viot case this equals $x_\eta (a) (1-x_\eta (a))$), which gives us a $(\Omega_N^L \times \K)$-indexed configuration. To prepare for the model on $\Omega_N^\infty$, we consider the induced \emph{volatility measure} $\wh R(x) \in \wt\CM(\Omega_N^L \times \K)$, defined by recording the volatility for every type $a$ at each site $\eta$, as well as the \emph{volatility occupation measure} up to some time $t$. As key ingredient this gives us for the process $(\wh X_t^L)_{t \geq 0}$ a \emph{triple} of \emph{stochastic processes}
\begin{equation}
\label{e1645}
\big(\wh X_t^L,\wh R^L_t, \wt {\wh R^L_t}\big)_{t \geq 0}  
\text{ with } \wt {\wh R^L_t} = \int_0^t \, \wh R^L_s \, \dd s.
\end{equation}
Based on the above, we formulate a \emph{stochastic equation} for the solution of the martingale problem.

\paragraph*{Outline of the remainder of this section.}
The remainder of this section is organised around two proofs, worked out in Section~\ref{ss.nsmp}--\ref{333*}, respectively, Section~\ref{334}, as follows. Section \ref{ss.nsmp} provides some background on \emph{non-standard martingale problems}. Section \ref{sss.keyp} states a key proposition that constitutes our main result for the well-posedness of the martingale problem. Section \ref{sss.hurd} lists the hurdles towards the description of the operators in the martingale problem. The proof uses three key lemmas, which are stated in Sections \ref{331}--\ref{332} and are proved at the end of Section \ref{333}. Section \ref{331} focusses on \emph{block-resampling}, Section \ref{332} focusses on \emph{mass flow} and \emph{single-site resampling}, while Section \ref{333} combines the different evolution mechanisms. Section \ref{333*} provides the proof of the key proposition. Section \ref{334} formulates and proves a \emph{stochastic equation} for the continuum-space process that holds when the migration is point-recurrent, which constitutes a further key proposition. 


\subsection{Non-standard martingale problems}
\label{ss.nsmp}

Properties \eqref{ag97alt}--\eqref{e527} say that the \emph{orbit of the path} in positive time resides in a \emph{small subset} of the state space, in which it subsequently evolves forever. This subset of states is \emph{not} topologically \emph{closed}, so that we need to use the concept of \emph{entrance laws} from time $0$ to describe the complete evolution. Indeed, by introducing a stronger topology on the subset we can in principle get a proper Markov process on this subset, but this process cannot start in an arbitrary initial state anymore. Furthermore, the volatility for the diffusion part is observable only in space-time \emph{windows} and \emph{not} in space-time \emph{points}, which makes the formulation of the martingale problem different from the standard form, as is explained in more detail next. 
\begin{itemize}
\item
The \emph{first reason} for the change of setting is that specifying the operator must involve an instantaneous transition into a dynamically closed subset of unions of countably many locally mono-type regions (see \eqref{e1664} for a formal definition), which turns out \emph{not} to be closed in the \emph{standard} topology on the state space $\wt\CM (\Omega^\infty_N \times \K)$. Closedness holds only in the \emph{stronger} so-called \emph{atomic} weak topology on $\wt\CM (\Omega_N^\infty \times \K)$, taken on the $\K$-marginal as the usual topology and strengthened on the $\Omega_N^\infty$-marginal in a way that is explained later (see Section \ref{333}). Here, the key word is \emph{type-cluster-decomposition} of $\Omega_N^\infty$ and a related atomic measure on $\Omega_N^\infty$. 
\item
The \emph{second reason} for the change of setting is that the volatility of the diffusion part of the path is no longer given by an operator acting on functions on $\Omega_N^\infty$, respectively, $\Omega_N^\infty \times \K$. Namely, the \emph{volatility coefficient} is specified as a function of the state by using a measure-valued \emph{space-time functional} via a \emph{bivariate measure} on $[0,\infty) \times (\Omega_N^\infty \times \K)$: the \emph{different-type collision local time} up to the current time $t$ of observation. The latter arises as a \emph{local space-time average} of Fleming-Viot operators acting on the $\K$-variables functions on $\Omega_N^\infty \times \K$ resulting in functions on $\Omega_N^\infty$ to be averaged. However, for \emph{fixed time} $t$ this expression is zero on mono-type regions and hence (in certain parameter regimes) almost everywhere on the orbit. It therefore does not catch the volatility in the system, which \emph{increases} as a \emph{local time} at \emph{random times} and \emph{random locations} in space. This is due to the fact that, in the approximation with finite grids as $L \to \infty$, volatility at sites tends to $0$ or $\infty$ as $L \to \infty$, and only as a volatility \emph{measure} in space-time has a non-trivial well-defined limit. 
\end{itemize}
The above hurdles will be explained and formalised in detail in Section \ref{332}. They are common in the theory of continuum-space evolutions, namely, for the case of Fisher-Wright diffusions, catalytic, respectively, mutually catalytic branching with \emph{two} types (see \cite{EF96,DF97a,DP98}). However, here we have a \emph{more complicated form}, because we have a richer type-structure due to the fact that $\K$ is countable or even \emph{uncountable}. The reason is (as in the two-type case) that the current state, even when locally mono-type, in \emph{space-time} windows is globally erratic enough to produce, through the \emph{combination} of migration and resampling, enough volatility in spatial balls via the \emph{collision local time of different types}, so as to change the type decomposition when averaged over space and time. Namely, due to migration mono-type balls can be invaded by other types, which start excursions in small spatial balls in the form of mono-type balls, and expand into other areas.

The fact that migration in the original individual-based model is based on jump processes rather than on diffusions results in a more complex pattern of the mono-type regions than for Brownian motions in Euclidean space (as in \cite{A79,GSW16}) or for branching models with only two types (as in \cite{EF96,DF97a,DP98}).


\subsection{Key proposition}
\label{sss.keyp}

Define the following subspace of the state space of $\wh X^\infty$
\begin{equation}
\label{e1664}
\begin{aligned}
&\wt\CM^\ast (\Omega_N^\infty \times \K)\\ 
&= \{ x \in \wt\CM (\Omega_N^\infty \times \K)\colon\,
x \text{ is \emph{locally mono-type} and has an \emph{atomic} $\K$-marginal measure}\}.
\end{aligned}
\end{equation}
Note that the elements in this set have a density w.r.t.\ $\mathfrak{h}$ in the form of a density sequence indexed by the atoms in $\K$. This will be the subspace in which the orbit lies in the recurrent regime.


\paragraph*{Martingale problem.}
 
The key result for tasks (o) and (i)--(ii) is the following proposition (part (c) of which is needed because the geographic space $\Omega_N^\infty$ in uncountable). The parts are stated in the order in which they will be proved. Several subtle ingredients have to be introduced rigorously. Afterwards we give guidance.

\begin{proposition}{\bf [$\wh X^\infty$ is the solution to a well-posed martingale problem]}
\label{P.solwpmp}
\begin{itemize}
\item[(a)] 
$\wh X^\infty$ is a strong Markov process with values in $\wt\CM(\Omega_N^\infty \times \K)$, with the Feller property, and with c\`adl\`ag paths.
\item[(b)] 
$\wh X^\infty$ solves a well-posed $(\Pi,\mfL)$-martingale problem on $D([0,\infty),\wt\CM (\Omega_N^\infty \times \K))$ for an \emph{entrance law} at time $0$ into the subset $D([0,\infty),\wt\CM^\ast (\Omega_N^\infty \times \K))$ of the state space, i.e., the $\CL(\wh X^\infty)$ satisfies for all $f \in \Pi$:
\begin{align}
\label{ag1548}
&\Big(F(\wh X^\infty_t) - F(\wh X^\infty_0) - \int_0^t \dd s\,(\mathfrak{L} F)(\wh X^\infty_s)\Big)_{t \geq 0} \text{is a martingale, for all $F \in \Pi$} \\
&\wh X_0^\infty = x_0 \in \wt\CM (\Omega_N^\infty \times \K) \text{ and } \wh X_t^\infty \in \wt\CM^\ast (\Omega_N^\infty \times \K) \text{ for } t >0. \nonumber
\end{align}
Here, $F$ is a test function in the form of a polynomial in $\Pi$ on $\wt\CM(\Omega_N^\infty \times \K)$ (see~\eqref{e1013}), and $\mathfrak{L}$ is a martingale-problem operator with a domain consisting of functions that are non-zero only on $\wt\CM^\ast(\Omega_N^\infty \times \K)$, and is of the form $\mathfrak{L}^{\rm jump}+\mathfrak{L}^{\rm diff}$, where 
\begin{itemize}
\item $\mfL^{\rm jump}$ is an integral operator for a jump process that is well defined for $F$ restricted to $\wt\CM^\ast(\Omega_N^\infty \times \K)$,
\item $\mfL^{\rm diff}$ is a (generalised) second-order operator of the (continuous, infinite variation) path part of the process with values in $\wt\CM^\ast(\Omega_N^\infty \times \K)$ for positive time.
\end{itemize}
Both operators are given in explicit form on $\wt\CM^\ast(\Omega_N^\infty \times \K)$, the latter with coefficients that are given via an occupation measure of a volatility functional of the process $\wh X^\infty$ (for the concept of ``generalized" more detail is given in Definition~\ref{def.martiro} in Section~\ref{sss.mart}). 
\item[(c)] 
\begin{itemize}
\item[(c1)] 
The process $\wh X^L$ converges to $\wh X^\infty$ in path space.
\item[(c2)] 
The \emph{volatility process} $(R^L_t)_{t \geq 0}$ and the \emph{volatility occupation process} $(\wt {\wh R^L_t})_{t \geq 0}$ are tight in $L$, and converge as $L\to\infty$ (recall~\eqref{e1645}).
\item[(c3)] 
Also the triple process $((\wh X_t^L)_{t \geq 0},(\wh R_t^L)_{t \geq 0},(\wt {\wh R^L_t})_{t \geq 0})$ is tight in path space, and converges as $L\to\infty$, with the limit triple process (= process, different-type collision local time, the latter accumulated) denoted by
\begin{equation}
\label{e1670}
\big(\wh X_t^\infty, L_{\wh X^\infty}(t), L^\ast_{\wh X^\infty} (t)\big)_{t \geq 0}.
\end{equation}
\end{itemize}
\end{itemize}
This triple is characterised by a martingale problem for an entrance law at time $0$ in (b).
\end{proposition}

Two observations concerning this key proposition are in order.


\paragraph*{New ingredients.}
Several new ingredients appear in the theorem above, which are somewhat complex and are made rigorous as we go along. The test function $F \in \Pi$ is described in \eqref{ag99b} and its sequel. The definition of $\mathfrak{L}$ turns out to be non-standard for both parts of the operator, as explained in Sections \ref{331}--\ref{333}. In short, $\mfL^{\rm jump}$ is only defined in $\wt\CM^\ast(\Omega_N^\infty \times \K)$, while $\mathfrak{L}^{\rm diff}$ needs a new concept to describe the \emph{local volatility coefficient} appearing in the formulas for $\mfL^{\rm diff}$, which is based on 
\begin{equation}
\text{$\wt {\wh R}$ = the $\Omega_N^\infty$-version of the process $\wh{\wt R}^L$ from \eqref{e1645} in the limit as $L \to \infty$.}\label{e3333}
\end{equation}
For the precise formulation of the martingale problem indicated in \eqref{ag1548}, see Definition \ref{def.martiro}, Lemma \ref{EffFV} and Corollary \ref{cor.Ito}. The triple in \eqref{e1670} is identified in \eqref{e1754} and is given explicitly in terms of the process $\wh X^\infty$ in \eqref{e1793}.


\paragraph*{Different-type collision local time.}

Our starting point is that the discrete-space process on $\Omega_N$ has the properties stated in Proposition~\ref{P.solwpmp}(a)--(b), and hence so do the processes on $\Omega_N^L$, $L \in \N$. In particular, we can write down explicit expressions for the operator $\mathfrak{L}^L$ analogous to \eqref{ag1548}, based on the formulas below \eqref{e890}. We refrain from explicitly doing so and refer to \cite{GHKK14}. 

To get the corresponding formulas for the \emph{limit $L \to \infty$}, new ideas are needed to get the properties in Proposition~\ref{P.solwpmp}. For (a) we use the construction of a Markov process, i.e., a transition kernel $(P_t)_{t \geq 0}$ on the state space, based on the \emph{duality} relation, as laid out in \cite{DGP24}. For (b) we must specify the \emph{operator} in the martingale problem defined on an appropriate class of test functions based on the limiting formulas for the processes $\wh X^L$. 

There are two \emph{new effects} in \emph{continuum-cardinality} geographic space. 
\begin{itemize}
\item[(I)] 
Because of the interaction between sampling and migration, we may have that in a ball of type $a$ a sub-ball of type $a^\prime$ develops, starting \emph{in a single site} and afterwards being boosted by resampling events. Therefore the description \emph{cannot} be expressed by a classical operator $\mathfrak{L}$, since such an operator would vanish on every locally mono-type state at any fixed time. The reason for this is that the volatility coefficient in the form as we had it for $L \in \N$ degenerates (otherwise the form of the operator would indeed be a classical second-order operator, as in the discrete case). Instead, the volatility coefficient must be based on the \emph{space-time picture} through the \emph{different-type collision local time measure} (a measure on $[0,\infty) \times (\Omega_N^\infty \times \K)$). The latter leads to a second-order operator $\mathfrak{L}^{\rm diff}$ of different form by giving the \emph{quadratic variation process} of the processes $\{(F(\wh X_t^\infty))_{t \geq 0}\}_{F \in \CF}$, with $\CF$ a measure- and convergence-determining set of test functions on the state space. 
\item[(II)]
The migration and resampling evolve in \emph{countably many} but \emph{small} jumps.
\begin{itemize}
\item
For the Fleming-Viot resampling, to capture the \emph{quadratic variation} for paths of evaluation processes of $\wh X^\infty$ we introduce as key ingredient 
\begin{equation}\label{e3355}
\text{the \emph{different-type collision local time} $L_{\wh X^{\ast,\infty}}^{\infty,\ast}(t)$}, 
\end{equation}
which is a bivariate measure on $[0,\infty) \times (\Omega_N^\infty \times \K)$ specifying the \emph{local volatility}, but now \emph{accumulated} up to time $t$, of the continuum-space process in \emph{space-time} windows, which is driven by a stochastic process $(L_{\wh X^\infty}(t))_{t \geq 0}$ with values in $\CM(\Omega_N^\infty \times \K)$, replacing the $\wh R^L$ in the limiting object (this works also for the limit of the Cannings resampling with resampling measure $\Lambda_{-\infty}$). Hence, for the Fleming-Viot resampling the form of the compensator in the martingale problem changes (and this is true also for the single-site Cannings resampling on the lowest hierarchical level, where in the limit $L \to \infty$ we obtain an evolution with an \emph{infinite-variation} continuous path, or a zero-path). 
\item
Added to these mechanisms is a \emph{finite-variation} part of the path driven by \emph{jumps}, coming from the Cannings \emph{block-resampling} given by $(\Lambda_k)_{k \in \Z}$, which is of a different nature and, depending on the parameters, either has the classical form or requires that we work with \emph{entrance} laws from time $0$, because of \emph{diverging rates} as $k \to -\infty$ on states that are \emph{not} locally mono-type. In other words, the operator is \emph{not} defined on the test functions defined on the \emph{full given Polish state space}, as in the standard martingale problem, but only on the restriction to $\wt\CM^\ast(\Omega_N^\infty \times \K)$, which is not a \emph{closed} subspace in the weak topology. To obtain a Polish space we need to strengthen the topology to the \emph{atomic} weak topology.
\end{itemize}
\end{itemize}


\subsection{Hurdles towards the formulation of the martingale problem}
\label{sss.hurd}

To formulate the martingale problem stated in Proposition \ref{P.solwpmp}, there are \emph{four} hurdles, which we deal with in Sections \ref{331}--\ref{334}. Along the way, three lemmas and a corollary are needed to state the most important technical facts.
\begin{itemize}
\item
Section \ref{331}: Formulate the \emph{$k \in -\N$} part of the \emph{block-coalescence operator} in the case where the rates \emph{sum up to infinity}. Make sure to guarantee that for all parameters we still have a meaningful expression for the operator $\mathfrak{L}^{\rm jump}$. A priori the operator in the stochastic regime is only defined in sites where the state is locally equal to a $\delta$-measure at a single type in some neighbourhood of the site, or in the deterministic regime, where it equals a constant-in-time state, namely, the measure $\delta_{\mfh \otimes\theta}$. Both form a set into which the process $\wh X^\infty$ instantaneously moves, via an entrance law. To get a well-defined operator for the further evolution, we need, for certain choices of the block-resampling parameters, an instantaneous transition of the process from a point in the state space into the set of states, which in the stochastic regime are locally mono-type and in the deterministic regime are equal to a specific state. 
\item
Section \ref{332}: Formulate the martingale problem for \emph{migration mass flow} and for \emph{single-site resampling}. For the latter there are three tasks:
\begin{itemize}
\item
Define a stochastic process $L_{\wh X^\infty}(t)(\dd\eta,\dd u)$ on $\Omega_N^\infty \times \K$, replacing the volatility in the classical operator on $\Omega_N^L$ based on the Fleming-Viot operator on $\CP_1 (\K)$. This is needed to cope with the instantaneous entry of the process into a small subset, and because the operator $\mathfrak{L}^{\rm diff}$ (more precisely, the volatility constant it contains) only has a replacement via the \emph{space-time} functional $L^\ast_{\wh X^\infty}$, the \emph{different-type collision local time}, which is a measure on $[0,\infty) \times (\Omega_N^\infty \times \K)$ that has a simple form on this subset but whose support is \emph{not} closed in the given topology. 
\item
Extend the solution to a setting where the state is a priori \emph{not atomic} in its $\K$-marginal and \emph{not absolutely continuous} in its $\Omega_N^\infty$-marginal, and therefore has to move \emph{instantaneously} into the above subset of the state space, i.e., we have to define an \emph{entrance law}. 
\item
Calculate the volatility constant $d^\ast$ in front of $L_{\wh X^\infty}(t)$.
\end{itemize}
\item 
Sections \ref{333}--\ref{333*}: Formulate rigorously the \emph{combined martingale problem} for the two resampling mechanisms, one for the continuous-path part with \emph{infinite variation} and one for the c\`adl\`ag jump-path part with \emph{finite variation}, and also for the continuous and finite-variation part of the migration. (This is relevant only for the case of stochastic limit evolutions as $L \to \infty$.)
\item 
Section \ref{334}: Formulate a \emph{stochastic equation} (SPDE) for the entrance law from time $0$ for the martingale problem solution, for the \emph{densities} of the state w.r.t.\ $\mathfrak{h}$ corresponding to the different atoms in type space $\K$, when two independent L\'evy-processes have hitting probability $1$. In that case we get not only atomic $\K$-marginals of the measure, but also absolutely continuous states with respect to the Haar measure $\mathfrak{h}$, and the state is locally mono-type. Do the same for the occupation measure on $[0,\infty) \times (\Omega_N^\infty \times \K)$. We obtain a different-type collision local time on the set of measures having these two properties, denoted by $\CM^\ast(\Omega_N^\infty \times \K)$, and a well-posed martingale problem on that set. After that we need to obtain an \emph{entrance law} by showing that other initial measure states \emph{jump} or move \emph{continuously} and \emph{instantaneously} into this subset.
\end{itemize}


\paragraph*{Outline.}

In what follows, we consider the pure-jump part (i.e., block-resampling) in Section~\ref{331}, the diffusion part in Section~\ref{332}, and the combined mechanisms in Section~\ref{333}. We prove Proposition~\ref{P.solwpmp} in Section~\ref{333*}.


\subsection{Block-resampling: the pure-jump part} 
\label{331}

We first consider the dynamics with \emph{just block-resampling}, giving us a jump process $\wh X^\infty$ with c\`adl\`ag paths, for which we formulate a martingale problem. 

We specify an algebra of \emph{test functions} $\Pi$, the \emph{polynomials}, generated by the collection of monomials
\begin{equation}
\label{e.1013}
F(X) = \int_{(\Omega^\infty_N \times \K)^n}
X^{\otimes n} (\dd \underline{\eta}, \dd \underline{u})\,\varphi^{(n)} (\underline{\eta},\underline{u}) , 
\qquad \varphi^{(n)} \in C_b \left((\Omega_N^\infty \times \K)^n,\R \right), n \in \N_0,
\end{equation}
which defines a collection of \emph{measure} and \emph{convergence-determining} functions on $\wt\CM(\Omega_N \times \K)$. It suffices to consider the algebra of functions generated by monomials, functions of the form $\varphi^{(n)}=\varphi_{(\Omega_N^\infty)^n}\, \varphi_{(\K)^n}$. 

Observe that $\wh X^\infty$ is a finite measure on finite measurable subsets of $\Omega_N^\infty$. Hence, normalising the restrictions of $x \in \CM(\Omega_N^\infty \times \K)$ to disjoint balls and using a disjoint decomposition of $\Omega_N^\infty$ obtained by choosing some radius for the balls characterised by a level $k \in \Z$, we can define the \emph{regular conditional measure} given $\eta$ as 
\begin{equation}
\label{e1800}
x_\eta(\dd u) = x(\dd u \mid \eta),
\end{equation}
which is an element of $\CP(\K)$.

Next, we come to generating paths of a \emph{jump} process via the change induced by the \emph{reshuffling} (i.e., a migration where all individuals jump simultaneously, with the jump sizes chosen independently and uniformly) and by the \emph{block-resampling}. For these transitions the compensator is \emph{expected} to be given by the operator
\begin{equation}
\label{e1236}
\left(\mathfrak{L}^{\rm Can} F \right) (x) 
= \sum_{k \in \Z} \int_{\Omega_N^\infty} \mathfrak{h}(\dd\eta) \int_{[0,1]} \Lambda_k^\ast(\dd r) \int_{\K} x_\eta(\dd a)
\left[F\left(\Phi_{a,r,\eta}(x)\right)-F(x)\right] 
\end{equation}
with (recall \eqref{ak:resampling-mapping-global})
\begin{equation}
\label{e1240}
\left[\Phi_{a,r,\eta}(x)\right]_{\zeta} = 
\left\{\begin{array}{ll}
(1-r) x_{\zeta} + r \delta_a, &\zeta \in B_k(\eta),\\
x_\zeta, &\zeta \notin B_k(\eta).
\end{array}
\right.
\end{equation}
Note that this term \emph{vanishes} for states that, for $\mathfrak{h}$-almost all $\eta$, are mono-type in $B_k(\eta)$ for some $k \in \Z$. 

However, the total block-resampling rate for $k \leq 0$ is \emph{infinite}, and so the operator is \emph{not defined} on \emph{all of $\wt\CM(\Omega_N^\infty \times \K)$}, but only a subset, namely, $\wt\CM^\ast(\Omega_N^\infty \times \K)$ from \eqref{e1664}. We need to show that the operator allows us to define the martingale problem for any \emph{entrance law} on the (admissible) initial states, and that for positive times the process moves into the subset of $\wt\CM(\Omega_N^\infty \times \K)$ consisting of \emph{locally mono-type states} with an \emph{atomic} $\K$-marginal. This is achieved through the following lemma, which will be proved in Section \ref{par.2070}.

\begin{lemma}{\bf [Entrance law]}
\label{l.1622}
\begin{itemize}
\item[(a)] 
If the block-resampling rates diverge as $k \to -\infty$ or $d^\ast >0$, then in the stochastic regime the process $(\wh X^\infty_t)_{t \geq 0}$ has, for positive times, states that are absolutely continuous with respect to $\mathfrak{h}$ and locally mono-type, and instantaneously enters such states for all initial laws. The entrance law is determined by the duality relation. 
\item[(b)] 
If $d=0$ and $\Lambda_{-\infty}=0$, then the entrance law $\wh X^\infty$ solves the $(\mathfrak{L}^{\rm Can},\Pi)$-martingale problem. In particular, $\mathfrak{L}^{\rm jump}=\mathfrak{L}^{\rm Can}$. 
\end{itemize}
\end{lemma}

\noindent
The above means that we cannot define the operator on the state space $\wt\CM(\Omega_N^\infty \times \K)$, but have an \emph{entrance law} from that state space into the subset $\CM^\ast(\Omega_N^\infty \times \K)$ where $\mathfrak{L}^{\rm jump}$ is defined. However, this space is not a closed subset, and hence is \emph{not} a Polish space, so that the restricted martingale problem is not of standard form. We have to fix this by \emph{strengthening} the weak topology on the state space and pass to a version of the \emph{atomic} weak topology to obtain a Polish state space on which we can use standard tools from stochastic analysis.

 
\subsection{Mass flow and single-site resampling: the diffusion part}
\label{332}

Next we consider the case \emph{without} block-resampling and $d^\ast >0$, i.e., the stochastic regime, but the \emph{diffusion} case. In that case $\wh X^\infty$ has continuous paths, which for $t>0$ immediately enter $\CM^\ast(\Omega_N^\infty \times \K)$. Note that the drift part has no problems and comes from the \emph{migration}, which in $\wh X^\infty$ corresponds to mass flow given by the semigroup of the $(c_k)_{k \in \Z}$-L\'evy-process. Let $\mathfrak{A}$ be the associated generator of the martingale problem of the migration, respectively, reshuffling on $\Omega_N^\infty$ defined in \eqref{e857}--\eqref{effmig}, and recall that $\mathcal{D}(\mathfrak{A})$ denotes its domain. This operator defines a unique semigroup for itself, which is a deterministic mass flow that can be added to the pure-diffusion martingale part.

Section \ref{sss.mart} describes the martingale problem, Section \ref{sss.dtclt} builds the different-type collision local times (Lemma \ref{l.1786}), Section \ref{sss.density-dtclt} identifies the density of the latter (Lemma \ref{EffFV}), while Section \ref{par.2070} provides the proof of Lemmas \ref{l.1622}, \ref{l.1786} and \ref{EffFV}.


\subsubsection{Martingale problem}
\label{sss.mart}

The martingale problem arising from the operator defined in \eqref{e857} for migration requires a \emph{new form} for the resampling part (which is the diffusion part), in particular, it requires a \emph{new condition}, namely, a specification of the \emph{increasing processes} of the moment functions $F(\wh X_t^\infty), F \in \Pi$ by a space-time functional of the path up to the present time. Furthermore we need to define this Markov process as an \emph{entrance law} at time $0$ on $\CM(\Omega_N^\infty \times \K)$ into $\CM^\ast(\Omega_N^\infty \times \K)$, on which we have a martingale problem for the process.

We can use the general alternative standard characterisation of a stochastic process via a \emph{martingale problem for processes with continuous paths}, by specifying a \emph{drift} and a \emph{quadratic variation}, where the latter does not have the standard form, but is given in the form of a \emph{collision-time functional} of different types in space-time, since it has the form of a \emph{local time}. This functional arises from the singular form that the operator $\mathfrak{L}^L$ develops in the limit as $L \to \infty$ due to the scaling of the rates for the single-site Fleming-Viot and Cannings resampling, which can be captured only in the form of a \emph{space-time} and \emph{occupation} integral. We can reformulate the martingale problem after replacing the operator $\mathfrak{L}$ in \eqref{ag1548} by the following. 

\begin{definition}{\bf [Martingale problem for the continuum-space diffusion process]}
\label{def.martiro}
{\rm We say that the law of the stochastic process $\wh X^\infty = (\wh X_t^\infty)_{t \geq 0}$ with values in $\wt\CM(\Omega_N^\infty \times \K)$ has an entrance law at time $0$ with values in $\wt\CM^\ast(\Omega_N^\infty \times \K)$ and is a solution to the $(\mathfrak{L}, \Pi,\CL[X_0^{\ast,\infty}])$ martingale problem when the following holds. 
\begin{itemize}
\item 
There is a \emph{drift operator} $\mfA$ defined on $\CD(\mfA) \subseteq C_b(\Omega_N^\infty \times \K)$; this is a first-order operator that represents the migration and generates the semigroup of the L\'evy-process. (Compare with \eqref{ak:multi-level-migration} for $\Omega_N$.)
\item 
There is an algebra of functions on $\Omega_N^\infty \times \K$, generated by functions of the form $\varphi(\eta,u) = \varphi_{\Omega_N^\infty}(\eta)\,\varphi_\K (u)$ with $\varphi_{\Omega_N^\infty}\colon\,\mathcal{D}(\mathfrak{A}) \to \R$ measurable, continuous and integrable with respect to $\mathfrak{h}$ and with $\varphi_\K\colon\,\K\to\R$ continuous. 
\item Two conditions for $\CL [\wh X^\infty] \in \CP(D([0,\infty), \wt\CM(\Omega_N^\infty \times \K)))$ are satisfied. 
\begin{itemize}
\item
The \emph{initial state} at time $0$ is the given random $\wh X_0^\infty \in \wt\CM(\Omega_N^\infty \times \K)$, whose law is translation invariant on $\Omega_N^\infty$, and 
\begin{equation}
\label{e1737}
(M^\varphi_t)_{t \geq 0} \text{ is a \emph{martingale}},
\end{equation} 
with
\begin{equation}
\label{e949}
M^\varphi_t = \langle \wh X^\infty_t,\varphi\rangle - \langle \wh X^\infty_0,\varphi\rangle - \int_0^t \dd s\, 
\big\langle \mathfrak{L}\varphi,\wh X^\infty_s\big\rangle
\end{equation}
and $\mathfrak{L}\varphi = (\mathfrak{A} \varphi_{\Omega_N^\infty}) \varphi_\K$. 
\item
For the solution $\wh X^\infty = (\wh X_t^\infty)_{t \geq 0}$
\begin{equation}
\label{ar1749}
\text{there exists an $\wt\CM^\ast(\Omega_N^\infty \times \K)$-valued stochastic process $(L_{\wh X^\infty}(t))_{t > 0}$},
\end{equation}
inducing a \emph{continuous non-decreasing} process $(L^\ast_{\wh X^\infty}(t))_{t >0}$ with values in $\CM(\Omega_N^\infty \times \K)$ given by 
\begin{equation}
\label{e1766}
L^\ast_{\wh X^\infty}(t)(\dd \eta,\dd u) = \int^t_0 \dd s\,L_{\wh X^\infty}(s)(\dd \eta,\dd u),
\end{equation}
defining a bivariate measure 
\begin{equation}\label{e2838}
\text{$L^\ast_{\wh X^\infty}$ on $[0,\infty) \times (\Omega_N^\infty \times \K)$,}
\end{equation}
called the \emph{different-type collision local time up to time $t$}, such that the \emph{quadratic variation process} of the \emph{evaluation process} $M^\varphi_t = \langle \wh X^\infty_t,\varphi\rangle$ is given by
\begin{equation}
\label{e959}
\langle\langle M_t^\varphi \rangle\rangle_{t \geq 0} 
= \left(d^\ast \big\langle L^\ast_{\wh X^\infty}(t),\varphi^2 \big\rangle\right)_{t \geq 0}
\end{equation}
with 
\be{e1750}
\big\langle L^\ast_{\wh X^\infty}(t),\varphi^2 \big\rangle 
= \int_{\Omega_N^\infty \times \K} L^\ast_{\wh X^\infty}(t) (\dd \eta,\dd u)\, \varphi(\eta,u)^2.
\ee
(Compare with \eqref{ak:multi-level-resampling-global}--\eqref{intTay} on $\Omega_N$.)
\end{itemize}
\end{itemize}
}\hfill$\spadesuit$
\end{definition}

\noindent
We will see that we can find a path $(L_{\wh X^\infty} (t))_{t >0}$, taking values in the set $\wt\CM^\ast (\Omega_N^\infty \times \K)$, such that in $D([0,\infty),\wt\CM(\Omega_N^\infty \times \K))$ we obtain a \emph{unique entrance law} $\CL [(\wh X_t^{\ast,\infty})_{t \geq 0}]$ on the subset $\wt\CM^\ast(\Omega_N^\infty \times \K)$. So there is no need to define the functional $L_{\wh X^\infty}$ on the whole state space, since the solution does not involve this functional elsewhere for times $t>0$. Note, however, that we need the path regularity of the solution at $t=0$.

 
\subsubsection{Building the different-type collision local time}
\label{sss.dtclt}

How do we construct $L_{\wh X^\infty}$? Crucial is the combination of resampling and mass flow, which keeps producing volatility locally as long as the process is not in a trap. Indeed, $L_X$ is constructed via an approximation, namely, via an \emph{averaging} produced by the migration transition kernel of the states of the solution (an idea that is based on the Trotter-formula for generators that are a sum of two terms), respectively, via a more accessible version based on the specific structure of $\Omega_N^\infty$ and the L\'evy-processes.

Note that the state $\wh X^\infty$ is transported by the flow given by the semigroup $(a_s)_{s \in [0,\infty)}$, which for $s >0$ has a density with respect to $\mathfrak{h}$, because the underlying effective migration process $\zeta^\infty$ has a \emph{transition density}. This allows us to introduce the different-type collision local time on the \emph{subset $\wt\CM^\ast(\Omega_N^\infty \times \K)$} of the state space as follows. 

Let $a^{\otimes 2}(\cdot,\cdot)$ be the \emph{transition density} on $(\Omega_N^\infty \times \Omega_N^\infty) \times (\Omega_N^\infty \times \Omega_N^\infty)$ for two independent $a$-L\'evy-processes on $\Omega_N^\infty$. We need to replace $a^{\otimes 2}(\cdot,\cdot)$ by an \emph{effective transition density}, accounting for the fact that with a certain probability the jumps of a pair of L\'evy-processes occurs at the same time due to the \emph{reshuffling}. We denote this transition density by $\wt{a}^{(2)}(\cdot,\cdot)$.

For $y \in \CP(\K)$, we set
\begin{equation}
\label{e1080}
Q_y(\dd u,\dd v) = y(\dd u) \delta_u(\dd v) - y(\dd u)y(\dd v)
\end{equation}
and look at $Q_{x_\eta}(\dd u,\dd v)$ (recall \eqref{e1800} for the definition of $x_\eta$). We use the fact that $\wh X^\infty$ restricted to balls and projected on $\K$ is atomic. Along the path the sets of atoms that are formed are \emph{decreasing} countable sets $A_t$, all contained in a \emph{countable} set $A$. 

The next lemma identifies the different-type collision local time for $t>0$ on the \emph{orbit} of the path of $\wh X^{\ast,\infty}$, and will be proved in Section \ref{par.2070}.

\begin{lemma}{\bf [Different-type collision local time]}
\label{l.1786}
In the point-recurrent regime with either $d>0$ or $|\Lambda_{-\infty}|>0$, for every $t>0$ the limit
\begin{equation}
\label{e1754}
L^\ast_{\wh X^\infty}(t) (\dd\eta,\dd u) 
= \left[\int_0^t \dd s \,\,\lim_{\delta \downarrow 0} \frac{1}{\delta} \int_0^\delta \dd s^\prime \,A_{s,s^\prime}(\eta,u)\right]
\mathfrak{h}(\dd\eta) \, \mathfrak{n}(\dd u)
\end{equation}
exists in probability, where
\begin{equation}
\label{e1754alt}
A_{s,s^\prime}(\eta,u)
= \int_{\Omega^\infty_N} \mfh(\dd\eta^\prime)\,
\wt{a}^{(2)}_{s^\prime}\big((\eta,\eta),(\eta^\prime,\eta^\prime)\big) \int_{\K} Q^u_{\wh X_{s+s^\prime,\eta^\prime}^\infty}(\dd v) 
\end{equation}
with 
\begin{equation}
\label{mfrndef}
\mathfrak{n}(\dd u)\, Q^u_{\wh X_{s+s^\prime,\eta^\prime}^\infty}(\dd v) = Q_{\wh X_{s+s^\prime,\eta^\prime}^\infty}(\dd u,\dd v),
\end{equation}
$\mathfrak{n}$ the extension to $\K$ of the counting measure on the countable set of types that occur on the path $\wh X^\infty$, and $\wh X_{s+s^\prime,\eta^\prime}^\infty \in \CP(\K)$ the density of the state in $\K$ at the time-space point $(s+s^\prime,\eta^\prime)$ (recall \eqref{e1800}), which is a measurable function of $[0,\infty) \times \Omega_N^\infty$.  
\end{lemma}

\noindent
In fact, using the monotonicity and the continuity in $t$ of the different-type collision local time, we obtain \emph{path convergence} of the expressions in \eqref{e1766} and \eqref{e1754} for $t >0$, so that we have a candidate for the compensator process for $(M_t^{\varphi^2})_{t \geq 0}$, and hence for the quadratic variation for $\langle M_t^\varphi,\varphi \rangle$. 

Note that we can use the expressions for the $L$-approximations and their increasing process also for the same quantity that we get from \eqref{e959}, which are more explicit because they give densities for the states via the definition of these objects. If we have \emph{densities}, then we can carry out the approximation limit in \eqref{e1754} and get a formula in terms of a \emph{function on $[0,\infty) \times \Omega_N^\infty$}. In the \emph{point-recurrent} regime we can therefore write down an explicit formula, as shown is in the next section. 


\subsubsection{Identifying the density of the different-type collision local time}
\label{sss.density-dtclt}
 
Set
\begin{equation}
\label{e1077}
L_x(\dd\eta,\dd u) =x(\dd\eta,\K) \int_\K Q_{x_\eta}(\dd u,\dd v) = \mathfrak{h} (\dd\eta) \int_\K Q_{x_\eta}(\dd u,\dd v).
\end{equation}
The following quantity is the candidate for the stochastic process $(L^\ast_{\wh X^\infty}(t))_{t > 0}$ in \eqref{ar1749}--\eqref{e2838} associated the solution $\wh X^\infty=(\wh X^\infty_t)_{t \geq 0}$ of the martingale problem in Definition \ref{def.martiro}:
\begin{equation}
\label{e1793}
L^\ast_{\wh X^\infty}(t) (\dd\eta,\dd u) = \int_0^t \dd s\, L_{\wh X_s^\infty}(\dd\eta,\dd u)
= \mathfrak{h}(\dd\eta) \int_0^t \dd s \int_\K Q_{\wh X_{s,\eta}^\infty}(\dd u,\dd v).
\end{equation}
An important observation is that the path in \eqref{e1793} is \emph{rough} as a function of $s$, so that the integral over $s$ is non-trivial. We claim that
\begin{equation}
\begin{aligned}
&- \text{the expression in \eqref{e1754} is given by \eqref{e1793}},\\ 
&- \text{the resampling rate $d^\ast$ is given by $d + |\Lambda_{-\infty}|$.}
\label{e1821}
\end{aligned}
\end{equation}
We will prove this claim in Section \ref{par.2070}. From the following lemma, which will also be proved in Section \ref{par.2070}, we get that $L_{\wh X^\infty}$ exists on the positive orbit of $\wh X^\infty$. The quantity $d^\ast$ has been identified as the single-site coalescence rate of the dual process $Y^{\ast,\infty}$ with entrance law $\wh Y^{\ast,\infty}$, and we bring this now in connection with the volatility rate in the forward process $\wh X^{\ast,\infty}$.

\begin{lemma}{\bf [Effective Fleming-Viot resampling rate]}
\label{EffFV}
For point-recurrent migration,
\be{e1811}
d^\ast = d + |\Lambda_{-\infty}|.
\ee
If $d^\ast >0$, then the process $(\wh X_t^\infty)_{t>0}$ has, when localised in balls in $\Omega_N^\infty$, atomic $\K$-marginals and in $\Omega_N^\infty$ locally mono-type states, such that the quadratic variation of the infinite-variation part of the path $(\langle \wh X_t^\infty,\varphi \rangle)_{t \geq 0}$ is given by the stochastic process $(L_{\wh X^\infty}(t))_{t \geq 0}$ via Definition \ref{def.martiro}. Therefore the time-inhomogeneous drift operator (see Corollary \ref{cor.Ito})
\begin{equation}
\label{e2041*}
\mathfrak{L}^{\rm diff} = \mathfrak{L}^{\rm diff}_{t,\wh X^\infty}
\end{equation} 
is well defined for polynomials, with their domain restricted to the set $\wt\CM^\ast(\Omega_N^\infty \times \K)$ containing the positive orbit of the process $\wh X^\infty$, in the sense that a formula for the derivative of the expectation of the polynomials in $t$ holds via Ito's formula.
\end{lemma} 

From the formula above we see that the operator almost has a standard form, since it acts on \emph{polynomials} like a \emph{second-order differential operator}. Only the coefficient in front of the monomial is not like the normal measure-valued Fleming-Viot operator. It is a monomial of order $2$, but has a more complicated form, involving \eqref{ar1749}--\eqref{e1750}. We have to prove that this measure is strictly positive and that $\langle \wh L_{X_t^\infty},\varphi \rangle$ is increasing over every positive interval in $t$ when $\varphi>0$, and is integrable on configurations with different types with positive space-time weight. For this we use the moment measures $E[\langle \wh X_t^\infty,\varphi \rangle^n]$, $n \in \N$. The evaluations of the moment measures satisfy the following.

\begin{corollary}{\bf [Operator $\mathfrak{L}^{\rm diff}$: evolution of the moment measure]}
\label{cor.Ito}
In the point-recurrent stochastic regime, the moment measures of $\wh X^\infty$ without migration and block-resampling satisfy
\begin{equation}
\label{e1872}
\frac{\dd}{\dd t} \E \big[\langle \wh X_t^\infty,\varphi \rangle^n\big]
= \tfrac12 d^\ast n (n-1)\, \E\Big[\big\langle L_{\wh X_t^\infty},\varphi^2\big\rangle 
\langle \wh X_t^\infty,\varphi\rangle^{n-2}\Big], \qquad n \in \N,\, t>0. 
\end{equation}
In particular (recall \eqref{e2041*}),
\begin{equation}
\label{e2946}
(\mfL^{\rm diff} F)(X) = \mathfrak{A}F + \tfrac{1}{2}\,d^\ast n(n-1) \langle L_{\wh X_t^\infty},\varphi^2 \rangle \dfrac{\dd F}{\dd X}(X),
\end{equation}
where $F(X)=\langle X,\varphi\rangle^n$.
\end{corollary}

\noindent
The moment equations for the migration are well known and for block-resampling can be written down as a Markov pure-jump process analogous to \eqref{ak:multi-level-migration} and \eqref{ak:multi-level-resampling-global}, i.e., we can write down an analogue of \eqref{e1872} that includes the effective migration and block-resampling, by adding the action of the corresponding operators on polynomials. These relations allow us to read off for the single-site resampling process its \emph{duality with a coalescent} on $\N$, via an operator identifying two variables in the test function defining the polynomial of degree $n \in \N$, respectively, a spatial coalescent with migration and block-coalescence on $\Omega_N^\infty$.  

The following lemma is an alternative to Lemma~\ref{l.1786}.

\begin{lemma}
\label{l.intermediate}
For every $t>0$,
\begin{equation}
\label{e1781}
\langle L^\ast_{\wh X^\infty_t}, \varphi^2 \rangle = \lim_{L \to \infty} \sum_{\eta \in \Omega_N^L} 
\int_0^t \dd s \int_{\K \times \K} Q_{\wh X_{s}^\infty[B_{-L}(\eta)]}(\dd u,\dd v)\,\,
N^L \int_{B_{-L}(\eta)} \mathfrak{h}(\dd\eta')\, \varphi^2(\eta',u),
\end{equation}
where $\wh X^\infty_s[B_{-L}(\eta)] \in \CP(\K)$ is given by
\begin{equation}
\label{e0000}
\wh X^\infty_s[B_{-L}(\eta)](\dd u) = \frac{\wh X^\infty_s(B_{-L}(\eta),\dd u)}{\wh X^\infty_s(B_{-L}(\eta),\K)} 
\end{equation}
and $\mfh(B_{-L}(\eta)) = N^{-L}$. 
\end{lemma}

\begin{proof}
The expression under the limit equals $\langle L^\ast_{\wh X^L_t}, \varphi^2 \rangle$. By passing to the limit $L\to\infty$ we get $\langle L^\ast_{\wh X^\infty_t}, \varphi^2 \rangle$.
\end{proof}


\subsubsection{Proof of the three constituent lemmas}
\label{par.2070}

In this section we give the proof of Lemmas \ref{l.1622}, \ref{l.1786} and \ref{EffFV}, and of the claim in \eqref{e1821}.

\begin{proof}{[Proof of Lemma \ref{l.1622}]}
The transition kernel has been determined by the dual process and the transition semi-group was constructed via the duality relation. For this kernel we need to determine the corresponding operator of the martingale problem via the generator criterion for duality following \cite{DGP24}. But since we have the operator for the dual process, which is a Markov jump process, this is standard.
\end{proof}

\begin{proof}{[Proof of Lemma \ref{l.1786}]}
We use a standard $L^2$-argument. We return to \eqref{e1754}--\eqref{mfrndef}. The quantity
\begin{equation}
\label{e2051}
\left[\frac{1}{\delta} \int_0^\delta \dd s^\prime\, A_{s,s^\prime}(\eta,u)\right]
\mfh(\dd\eta) \, \mathfrak{n}(\dd u)
\end{equation}
is the \emph{different-type local time measure} at the time-space point $(s,\eta)$ averaged over the time interval $[s,s+\delta]$, which is an element of $\wt\CM(\Omega_N^\infty \times \K)$. This measure is applied to test functions $\varphi^2$ on $\Omega_N^\infty \times \K$ to obtain the volatility of the semi-martingale $\langle \wh X_s^\infty, \varphi \rangle$ at time $s$. In this way, for every choice of $\varphi$ from a countable dense subset, we obtain a random variable with values in $\R$, of which we must show convergence in probability as $\delta \downarrow 0$ to conclude convergence as $\delta \downarrow 0$ of the expression in \eqref{e2051} as a measure on $\Omega_N^\infty \times \K$. This runs as follows.

If in \eqref{e2051} we pass to the densities of the random measure, then we get an average of density functions on $\Omega_N^\infty \times \K$, namely, a space-time average. Inserting evaluations of the measure, we get averages of densities modulated with the test function $\varphi$. This is an average of an $(\Omega_N^\infty \times \K)$-indexed random field with values in $\R$ and with specific weights given by a random walk transition kernel. This setting is standard (see, in particular, Lemma ~\ref{lemmA4} in Appendix~\ref{appA:IPLT}, which provides a generalisation of the argument \cite{EF96}), so that the limit $\delta \downarrow 0$ can in fact be taken.
\end{proof}

\begin{proof}{[Proof of Lemma \ref{EffFV}]} 
We specify the operator coming from the Cannings resampling. The single-site resampling acts on level-$L$, while the block-resampling acts on all levels $k \geq -L+1$ simultaneously. The latter give rise to a jump process and plays no role for the diffusion part we deal with here. In the limit as $L\to\infty$, the terms on level $L$ give an additional contribution of the same form as the Fleming-Viot resampling, and the constant $d^\ast$ in front of the volatility has to be calculated for both jointly, and depends on $d$ and $|\Lambda_{-\infty}|$.

The single-site $\Lambda^L_{-\infty}$-resampling converges as $L \to \infty$ to a Fleming-Viot resampling, which follows by showing that the same holds for the dual process and that the resulting resampling rate is the sum of the two limiting resampling rates, i.e., $d^\ast=d+|\Lambda_{-\infty}|$, as we showed in the proof of Lemma~\ref{RWLprec}, Part (a). Therefore we must prove that $\wh X^\infty$ has a continuous part, which is the mass flow, and a diffusion part with increasing process as claimed in Lemma~\ref{EffFV}. We do this in Section~\ref{333*}. Consequently, the generator acts as a second-order differential operator, so that the quadratic variation is as claimed in the analysis of the limiting collision local times, which arise as the $L \to \infty$ limit of the collision local times on $\Omega_N^L$. Hence the compensators of the $(\langle \wh X_t^\infty,\varphi\rangle^2)_{t \geq 0}$, and therefore the evaluations $(\langle \wh X_t^\infty, \varphi\rangle)_{t \geq 0}$, satisfy the martingale problem conditions as a consequence of the finite-dimensional distribution convergence, and we get convergence of the dual processes $Y^L$ to $Y^\infty$ in path space.
\end{proof}

\begin{proof}{[Proof of \eqref{e1821}]}
The first claim in \eqref{e1821} is immediate from \eqref{e1754alt} because
\begin{equation}
\lim_{\delta \downarrow 0} \frac{1}{\delta} \int_0^\delta \dd s^\prime\, A_{s,s^\prime}(\eta,u)
= \int_{\K} Q^u_{\wh X^\infty_{s,\eta}}(\dd v).
\end{equation}
The second claim in \eqref{e1821} is \eqref{dinf}, which was settled via Corollary \ref{cor.1939}. An alternative route to the first claim in \eqref{e1821} is to show that \eqref{e1781} and \eqref{e1793} are the same. To that end, use that
\begin{equation}
\lim_{L \to \infty} \wh X^\infty_s[B_{-L}(\eta)] = \wh X_{s,\eta}^\infty
\end{equation}
and
\begin{equation}
\lim_{L \to \infty} N^L \int_{B_{-L}(\eta)} \mathfrak{h}(\dd\eta')\, \varphi^2(\eta',u) = \varphi^2(\eta,u).
\end{equation}
Substitution into \eqref{e1781} gives
\begin{equation}
\langle L^\ast_{\wh X^\infty_t}, \varphi^2 \rangle 
= \int_{\Omega_N^\infty} \mfh(\dd\eta) \int_0^t \dd s \int_{\K \times \K} Q_{\wh X_{s,\eta}^\infty}(\dd u,\dd v)\, \varphi^2(\eta,u),
\end{equation}
which is the same as saying that
\begin{equation}
L^\ast_{\wh X^\infty_t}(\dd\eta,\dd u) = \mfh(\dd\eta) \int_0^t \dd s \int_{\K} Q_{\wh X_{s,\eta}^\infty}(\dd u,\dd v),
\end{equation}
because $\varphi$ is an arbitrary test function. The latter is the same as \eqref{e1793}.
\end{proof}


\subsection{Combining different evolution mechanisms}
\label{333}

We have seen how to treat the jump part and the continuous part of $\wh X^\infty$ separately. Next, we deal with the problem how to \emph{combine} the martingale problems of the two parts into a single martingale problem. For that we decompose 
\begin{equation}
\label{e1118}
\wh X^\infty_t = \wh X_t^{\rm flow} + \wh X_t^{\rm jump} + \wh X_t^{\rm diff}, \qquad t \geq 0,
\end{equation}
with $\wh X_t^{\rm flow}$ the \emph{mass-flow part} due to migration (defined as the continuous finite-variation path), $\wh X_t^{\rm jump}$ the \emph{accumulated-jump part} (the sum of countably many jump-discontinuities in the finite-variation path), and $\wh X_t^{\rm diff}$ the \emph{diffusion part} (defined as the continuous infinite-variation path). The jump part $\wh X^{\rm jump}$ can be extracted from the c\`adl\`ag path $\wh X^\infty$ by forming the sum of the jumps up to time $t$ truncated at level-$L$ and afterwards taking the limit $L \to \infty$ (this is a Cauchy-sequence on $D([0,\infty),\Omega_N^\infty)$ by construction). The remaining \emph{continuous part} is decomposed uniquely into the finite-variation path $\wh X^{\rm flow}$ and the infinite-variation path $\wh X^{\rm diff}$, but with finite quadratic variation. 

For the path $\wh X^{\rm diff}+\wh X^{\rm flow}$ we use the quadratic variation and the compensating drift term for its characterisation, as in \eqref{e959}. For the jump part $\wh X_t^{\rm jump}$ we use the complete block-reshuffling-resampling operator to get the compensator for the accumulated jump process, as in \eqref{e1236}--\eqref{e1240}. Here we have to address the case of \emph{diverging} block-resampling rates, which requires that the process is locally mono-type, \emph{guaranteeing that the contribution is finite} and hence well defined. This requires constructing an \emph{entrance law} that for $t>0$ places the process in that subspace.

A final task is to piece the parts together. For that we use that they change on disjoint time sets, which are countable sets of jump times (whose complement has full Lebesgue measure) and which for finite $L$ is a union of intervals. Therefore the compensator of the pieces add up to the compensator of the whole path. This task is trivial for Markov pure-jump processes, but is more subtle when countably many jumps can occur in a finite time interval. However, we do have a limit of processes of the pure-jump type. Hence, for the truncated sum we have a well-defined expression and the contributions of the pieces form Cauchy sequences in $L$.


\subsection{Proof of the key proposition}
\label{333*}
 
\begin{proof}{[Proof of Proposition \ref{P.solwpmp}]}\\
The main difficulty sits in proving part (a) and identifying the local characteristics in part (b).

\medskip\noindent
(a) We use the results in \cite{DGP24} to show that, due to the duality equation, we can \emph{construct a Markov transition semigroup} $(\wh P^\infty_t)_{t \geq 0}$ for the process $\wh X^\infty$ that is \emph{Feller}, \emph{strong Markov}, and has a version with \emph{c\`adl\`ag paths}. Because the law of $(\wh X_t^\infty)_{t \geq 0}$ is uniquely determined by this semigroup, we have shown the \emph{strong} Markov property of $\wh X^\infty$, which we constructed as a limiting stochastic process via its finite-dimensional distributions by using duality. For that purpose we have to check the \emph{conditions} for the construction of the semigroup posed in \cite[Section 2.2]{DGP24} on the \emph{dual expectations}, which we next formulate and verify. 

Indeed, we have to \emph{adapt} the set-up of standard martingale problems and standard duality based on operator criteria. Here, the duality arises from the semigroup as a limit relation of dualities of \emph{standard} martingale problems, which are based on the \emph{generator criterion}. Namely, we are given the duality relations for the finite-dimensional distributions (and for the semigroups) via \eqref{ag99} and \eqref{e1547}, respectively. A key point is to verify a \emph{positivity condition} for the dual expectations, which allows us to obtain the Markov transition kernel from the duality equation \eqref{ag99} in connection with \eqref{e1431}--\eqref{e1474}. 

We define the transition kernel $(P_t)_{t \geq 0}$ on the state space $E$ by finding a transition kernel $(P_t)_{t \geq 0}$ and using the dual semigroup $(P^\ast_t)_{t \geq 0}$ with
\begin{equation}
\label{e2089}
\int_{E^\prime} P^\ast_t(e^{\,\prime},\dd \wt e^{\,\prime}) H(e,\wt e^{\,\prime}) 
= \int_E P_t(e,\dd \wt e\,) H (\wt e,e^{\,\prime}), \qquad e \in E,\,e^{\,\prime} \in E^\prime.
\end{equation}
According to \cite[Theorem 2.1]{DGP24}, we need to verify the \emph{existence} of the kernel above. The key to this is \cite[condition (2.8)]{DGP24}, which requires the \emph{positivity} of the functional $P_t^\ast (H(e,\cdot))$ on the space of $\{H(e,\cdot), e \in E\}$. We have to check this property by using moment techniques in combination with the fact that the dual expectations are the $L \to \infty$ limits of Markov processes and therefore \emph{do satisfy the positivity condition}. In \cite{DGP24} examples are discussed that include the Cannings model.

\medskip\noindent
(b) We have to verify the claimed \emph{local characteristics} of the semi-martingales given by linear evaluations. This means that we have to derive the local characteristics of the process $\wh X^\infty$ from the \emph{construction} of $\wh X^\infty$, by showing that the \emph{local characteristics of $\wh X^L$} converge to the claimed limiting expression, based on part (c) below. This is standard, except for the diffusion part. Nevertheless, some care is needed for the jump part, because we have an \emph{entrance law} rather than a Markov process with an initial state and without a transient part of the state space. However, from the literature we know the action of the martingale problem operators of $\wh X^L$ and of $\wh X^\infty$ via the duality relation for the limiting dual process with block-resampling as a jump process. 

For processes with finite quadratic variation, via Ito's formula we get the action of the martingale operators for the diffusion part. Only the changed form of the volatility coefficient has to be established, by calculating the quadratic variation of evaluations $\langle \wh X_t^\infty,\varphi\rangle$. In particular, the compensator of $\langle\wh X_t^\infty,\varphi\rangle^2$ has to be identified in terms of the \emph{different-type collision local time} by using the dual process based on the moment measures and carrying out \emph{moment calculations}. In other words, we have to establish \eqref{e959}--\eqref{e1750}. But these relations are immediate from duality.

\medskip\noindent
(c) We have to prove convergence of $\wh X^L$ to $\wh X^\infty$ by showing \emph{tightness}. Here we use the \emph{Joffe-Metivier tightness criteria} for semi-martingales (compare with \cite{D93}). Namely, it suffices to prove tightness of linear functions of the form $\langle u,f \rangle$ for a measure $u$ and for a convergence-determining set of test functions $f$. For the latter we can use the duality functions $\{H(\cdot,e^\prime),\, e^\prime \in E^\prime\}$ with $E,E^\prime$ the state space of the process, respectively, the dual. We use that these evaluations, when applied to $\wh X^\infty$, are \emph{semi-martingales}, so that the Joffe-Metivier criteria for tightness work. See \cite[Theorems 3.6.4 and 3.6.6 and Corollary 3.6.7]{D93} for detailed explanations and a list of sufficient conditions. These conditions are met in our setting as well.
\end{proof}


\subsection{Stochastic equation}
\label{334}


\paragraph*{Structure.}

Recall that, by the strong duality, $X_t^L$ for $t>0$ has at most countably many types from a random set $A^L$, where $A^L$ is non-decreasing in $L$ since we keep adding new positions and \emph{is contained in a countable set $A$} via truncation approximation. In particular, we can write $\wh X^L$ as a system of $\Omega_N^L \times A^L$-indexed \emph{stochastic differential equations}. We can label the random $A^L$ by $\N$, just as for the process on $\Omega_N^L$, as a system of stochastic differential equations for random probability measures on $\N$ given a countable random set $A^L$ of types contained in $\K$, where the driving white noise and Poisson field are independent of the choice of $A^L$. We can rephrase this, via the embedding of $\Omega_N^L$ in $\Omega_N^\infty$, as saying that we have a collection of measures on $\Omega_N^\infty$ with weights $\leq 1$ indexed by a countable set $A$ summing up to $\mfh$. Below we will write down a \emph{stochastic equation} for a version of the process $\wh X^\infty$ on $\Omega_N^\infty$ in the case of a \emph{stochastic} evolution. We will derive $A$ as a \emph{function of $\wh X^\infty_0$}, and we will look for an evolution given a realisation of the random set $A$.

We claim that, by the strong dual representation, $A^L \uparrow A$ as $L\to\infty$ when we pass to $X_t^\infty$, with $A$ a \emph{countable} set. This raises the question of how to view the evolution of the state as the evolution of a \emph{countable collection} of \emph{measures on $\Omega_N^\infty$} indexed by the types from $A$. The set $A$ depends on the events in $\mathop{\cap}_{\varepsilon>0} \mathop{\cap}_{t \in (0,\varepsilon)} \CB (\wh X_t^\infty)$. The point here is that the set $A$ is determined by the path observed for times in $(0,\delta]$, where $\delta$ can be chosen arbitrarily small. In particular, the evolution in $[\delta,\infty)$ is concerned only with properly defined subsets $A^\delta \subseteq A$ such that $A^\delta \uparrow A$ as $\delta \downarrow 0$. The same holds for the occupation measure $L_{\wh X}$. 

The \emph{form of the equation depends on $A$}, which in turn depends on the random path for which we write down stochastic evolution equations. These equations hold only for positive time and \emph{given $A$}. The form of the equation is independent of the realisation of $A$ and depends only on its cardinality, which equals that of $\N$.

We have seen in Theorem~\ref{T.char} that in the stochastic regime, $\wh X^\infty$ has \emph{densities} with respect to the Haar measure $\mathfrak{h}$ for the states of $\wh X_t^\infty$ projected on a \emph{specific type $a$} from the \emph{countable random set $A$} (the random set of types from $\K$ that survive up to some strictly positive time), and we get a \emph{collection of density functions} indexed by types $A \subseteq \K$, summing to $1$ in every site $\mfh$-a.e., that evolve in time. For strictly positive time we have a $\CP(\K)$-density $\wh X^{\infty,\dagger}$ on $\Omega_N^\infty$, which we can write as an \emph{atomic} probability measure on $\K$ with weight functions $\wh X^{\infty,\ddagger}(\cdot,\cdot)$ on $\Omega_N^\infty \times \K$ as
\begin{align}
\label{e982}
\wh X^{\infty,\dagger}_t(\eta,\cdot)
&= \sum_{a \in A} \wh X_t^{\infty,\ddagger}(\eta,a) \delta_a(\cdot) \in \CP (\K) \text{ for } \mathfrak{h}\text{-a.e.\ } \eta, t>0,
\end{align}
and we can represent $\wh X_t^\infty$ as
\begin{align}
\label{e1832}
\wh X^\infty_t(\dd\eta,\dd u) 
&= X^{\infty,\dagger}_t(\eta,\dd u)\,\mathfrak{h}(\dd\eta).
\end{align}

For the above $A$-indexed collection of processes of \emph{random fields} on $\Omega_N^\infty$, we obtain a description via a 
\begin{equation}
\label{e1907}
\text{\emph{system of stochastic equations} for $\{(\wh X_t^{\infty,\ddagger})(\eta,u)\colon\, 
\eta \in \Omega_N^\infty, \ u \in A,\; t > 0\}$.}
\end{equation}
If we code $A$ by an $\N$-indexed sequence, then we can write down a process of collections of random fields as
\begin{equation}\label{e2635}
\left(\big(\wh X_t^{\infty, \ddagger}(\cdot,m)\big)_{m \in \N}\right)_{t >0}
\end{equation}
and obtain a corresponding (SPDE+SIE)-stochastic equation (see Section~\ref{sss.sie}, in particular, equation~\eqref{a.1917}), since the \emph{evolution is independent} of the particular value in $\K$; and we can use a realisation of $A$ to obtain the paths of $\wh X^\infty$ as strong solutions. Since these functions turn out to be continuous in time and space, except for the countably many block-resampling times, the orbit in fact lies in a subset of the state space of the space-time process of $\wh X^{\infty,\ddag}$ and is entered \emph{instantaneously}. By introducing a \emph{stronger topology} on this subset, we can define a proper Markov process, namely, the \emph{restriction} of $\wh X^\infty$ to this subset.


\paragraph*{Ingredients.}

In summary, the solution of the martingale problem can be given by a \emph{stochastic equation} that is built and driven by three independent random elements:
\begin{itemize}
\item 
A \emph{random countable set $A \subseteq \K$} of all realised types for the population at strictly positive times whose law depends on $\wh X^\infty_0$. This is sampling locations from $\Omega_N^\infty$ by $\mfh$, and types by $\wh X_0^\infty(da \mid \eta)$.
\item
A \emph{white noise measure process} $W$ on $\Omega_N^\infty \times \K$, supported on the random set $A \times \Omega_N^\infty$, giving an occupation measure of noise on $[0,\infty) \times (\Omega_N^\infty \times \K)$ for the \emph{Fleming-Viot resampling} in the point-recurrent regime, chosen independently for \emph{any} realisation of the countable set $A$, to drive the stochastic equations for the \emph{diffusion part} of the path (i.e., the \emph{infinite-variation continuous} part of the path).
\item
A \emph{Poisson point process} $\Pi$ on $\Z \times [0,\infty) \times \Omega_N^\infty$ for the \emph{Cannings block-resampling}, giving the cadlag \emph{jump-process part} (with $\Z$ labelling the hierarchical levels of the blocks) and being independent of $A$. 
\item 
A \emph{deterministic} mass flow giving the \emph{continuous finite-variation} part of the path.
\end{itemize}
Below we spell out the details of how to obtain these ingredients. In Section \ref{sss.jpdp} we focus on the jump part and the diffusion part, in Section \ref{sss.sie} we focus on a stochastic integral equation, leading up to a key proposition (Proposition \ref{prop.1869}).


\subsubsection{Jump part and diffusion part}
\label{sss.jpdp}

(i) We start with the \emph{jump-process part}. Recall that $\Omega_N^\infty$ can be decomposed into $k$-balls for each $k \in \Z$, so that for every site there is exactly one $k$-ball where the site is involved in a $k$-block-resampling initiated by one of the sites, with rate $\lambda_k$. Therefore, for every $k \in \Z$ we need a Poisson point process on $[0,\infty)$ that places sites on the embedded $\Omega_N^k$-labelling of those $k$-blocks into which $\Omega_N^\infty$ decomposes, and characterises when and where a $k$-block-resampling happens. Thus,
\begin{equation}\label{e1913}
\text{$\Pi$ is a Poisson point process on } \bigcup_{k \in \Z} \Omega_N^k \times [0,\infty),
\end{equation}
where the intensity measure is specified via $(\mu_k)_{k \in \Z}$ such that $\mu_k\, (\CB_k(\eta))\,\mathfrak{n}_k(\cdot)$ is the intensity measure on level $k$. Each point of $\Pi$ gives rise to a $\Lambda_k$-block-resampling, and each of these transitions is independent of the \emph{set $A$} of types and the \emph{noise} that determines the set $A$.

\medskip\noindent
(ii) We continue with the \emph{diffusion-process part} (the infinite-variation part with finite quadratic variation). This uses the specific properties of $\wh X^\infty$, and uses that the $\K$-marginal of the state of the process is atomic and the $\Omega_N^\infty$-marginal is absolutely continuous with respect to $\mathfrak{h}$. 

Suppose that we have solved the martingale problem and observe the evolution beyond some positive time $\varepsilon >0$. On the interval $[\varepsilon,\infty)$ we order the atoms according to their size in a reference ball $B$ with radius $1$ around $0$ in $\Omega_N^\infty$, to obtain a unique representative of this vector in $\ell_1(\N)$, for an element of $[0,1]^\N$ for every $\eta \in \Omega_N^\infty$. The \emph{diffusion part} involves a measure state that has a density w.r.t.\ $h$ for fixed $u \in A$, which in \emph{``every site"} $\eta \in \Omega_N^\infty$ is driven by a family of Brownian motions for the different types. However, at a single site we have a speeding up of resampling and migration by $N^L$, which as $L \to \infty$ tends to $0$ or to $\infty$. Therefore we need to pass to the \emph{aggregated measure} to get a well-defined object. In particular, the densities have to be viewed as an aggregated object and not as a single-site value with a diffusive evolution. This means that, instead of sites $\eta \in \Omega_N^\infty$, we must consider test functions on $\Omega_N^\infty$, where we may choose a separable, countable set.

The type distributions in $A$ at $\eta \in \Omega_N^\infty$ of $X^{\infty,\ddagger}$ are probability measures, a property that is preserved over time, and the diffusion part of $\wh X^\infty$ is adapted to the noise. Furthermore, the random set $A$ depends on the evolution over times $t \in [0,\varepsilon]$ for arbitrary $\varepsilon>0$. The noise driving the evolution for larger times is, for given $A$, independent of the evolution for smaller times. On $\Omega_N^L$ the preservation property for probability measures can be achieved by using a Brownian motion $(W_1(t))_{t \geq 0}$ for an atom $a$ and a Brownian motion $W_2=(-W_1(t))_{t \geq 0}$ for all the other atoms in $A \setminus a$, in a proper combination. This needs to be formulated with care, because the set $A$ is random, while for every given time the noise is independent of the noise for later times. So the issue is to formulate an \emph{entrance law} and a suitable noise measure, i.e., an element of $\wt\CM(\Omega_N^\infty \times \K)$. 

In more detail, we associate with every site $\eta \in \Omega_N^L$ a Brownian motion $(W_\eta (t))_{t \geq 0}$, independent of $\eta \in \Omega_N^\infty$ and independent of the types in $\eta$ or elsewhere. This extends to $\Omega_N^{\infty, \rm disc}$, which gives us
\begin{equation}
\label{e2659}
\text{a \emph{space-time white-noise measure} $W(\dd s,\dd\eta,\dd u)$.}
\end{equation}
The white noise, respectively, $-1$ times the white noise have to be transformed into a noise measure on $\Omega_N^\infty \times \K$ by combining $W$ and $-W$ such that the effective noise over all types at a site add up to $0$ at every time point and every space point. For each $\eta\in\Omega_N^\infty$ we consider the set $A_{\eta,t}$ of types occurring at site $\eta$ at time $t$. We take an $a \in A_{\eta,t}$ that is maximal for a time interval $I$ that has positive length. As long as $a$ is maximal, we associate with $a$ the Brownian motion $W_\eta = (W_\eta(s))_{s \in I}$, and with all the other types we associate $-W_\eta$. In this way we obtain Brownian motions driving the various components, i.e., we drive our stochastic equation based on the white noise measure $W$. The extension of this construction from $\Omega_N^{\infty, \rm disc}$ to $\Omega_N^\infty$ is straightforward because $\Omega_N^{\infty, \rm disc}$ is a dense subset of $\Omega_N^\infty$. See \cite[Section 7]{D93} for more background.

For $u \in A$, set
\begin{equation}
\label{e2663}
W^u(\dd s,\dd\eta) = W(\dd s,\dd\eta,\{u\}).
\end{equation}
If we want to represent the process $\wh X^\infty$ as a solution to a stochastic equation, then we have to generate the Brownian motions from the path of the process, directly based on the increments of the continuous infinite-variation part of the process. We come to this in the proof of Proposition~\ref{prop.1869} below.


\subsubsection{Stochastic integral equation}
\label{sss.sie}

For $\wh X^L$ we can write down a system of \emph{stochastic differential} and \emph{integral equations} (indexed by $\Omega_N^L \times A$) based on the generator and Ito's formula (recall \eqref{e890}--\eqref{intTay}). The task is to obtain the analogue for $\wh X^\infty$ on $\Omega_N^\infty$. Namely, using the expression in \eqref{e1118} for $\wh X^\infty$, we want to write down a \emph{stochastic integral equation} with an integral expression for the changes of each of the three components of the path (continuous infinite-variation part, jump part and continuous finite-variation part), driven by $W$, a \emph{white noise measure process} on $\Omega_N^\infty \times \K$, and $\Pi$, a \emph{Poisson point process} that is \emph{independent} of the white noise measure $W$, and with a deterministic evolution part given by the \emph{mass flow} induced by a transition semigroup on $\Omega_N^\infty$, specified by the parameters $(c_k)_{k \in \Z}$ and $(\lambda_k)_{k \in \Z}$ that determine the motion of $\zeta^\infty$ and its transition semigroup on $\Omega_N^\infty \times \Omega_N^\infty$. We need to define these ingredients on a suitable probability space, \emph{together} with the solution to our stochastic equation. Note that, even though we have c\`adl\`ag functions as densities, we cannot write down an SPDE for $(\wh X_t^{\infty, \ddag}(\cdot,a))_{a \in A}$: we have to use the \emph{evaluations via test functions} integrated w.r.t.\ $\mfh$, since only these aggregated quantities evolve on a time scale of order $1$ (recall that the density is given via a sequence of local averages of ``single-site" states, different from the situation on $\Omega_N^L$). This means that we need a proper concept of a solution of such a stochastic equation, which we include in the following proposition. 

For a systematic account of \emph{stochastic equations on ultrametric spaces}, see \cite{K01}. For a detailed account of the concept of stochastic equations in general, see \cite{KS88}. Recall also \eqref{ar1749}--\eqref{e1750} in the definition of the martingale problem.

\begin{proposition}{\bf [Stochastic equation]}
\label{prop.1869}
Consider the process $\wh X^\infty$ constructed as limit of $\wh X^L$ as $L\to\infty$ on some probability space.
\begin{itemize}
\item[(a)] 
Let $A$ be the random set of types that survive up to some positive time in the random path $\wh X^\infty = (\wh X_t^\infty)_{t \geq 0}$. The latter is realised in law by running the dual backward process $\wh Y^{\ast,\infty}$ for some time horizon $t>0$, and observing the types surviving at some of these positive times. The set $A$ arises in law by sampling at each site in $\Omega_N^{\infty, \rm disc}$, using the restriction of $\wh Y^{\ast,\infty}$ to a countable set of types from $\wh X^\infty_0$ at these sites, and taking the union. 
\item[(b)] 
The evolution of $\wh X^\infty$ given $A$, encoded as an $A$-indexed set of functions on $\Omega_N^\infty$ with values in $\K$, is given on the probability space of $\wh X^\infty$, which supports $\Pi$ and $W$ (both independent of $A$) from \eqref{e1913} and \eqref{e2659}, which represents the unique weak solution of the martingale problem, and has the following properties:  
\begin{itemize}
\item
For $\wh X^\infty$ at every time $t>0$ there is for $\wh X_t^\infty$ a density with respect to $\mathfrak{h}$, in the sense that there exists a vector of density functions for the measures $\wh X^\infty_t (\cdot, \{a\})$ on $\Omega_N^\infty$ indexed by $a \in A$, written $\{\wh X_t^{\infty,\ddag}(\eta,a)\colon\,t>0,\,\eta \in \Omega_N^\infty\}_{a \in A}$, that are jointly continuous in $t$ and $\eta$, except at the countably many jump times in $t$ of block-resampling specified by $\Pi$ from~\eqref{e1913}. The sum over $a \in A$ of the densities must be equal to $1$ $\mfh$-a.e.\ and a.s.\ (to ensure that the $\K$-component of the solution is a probability measure on $A$).
\item
For every $a \in A$ there exists a white noise measure $W^a(\dd t,\dd\eta)$ on $(0,\infty) \times \Omega_N^\infty$, given via \eqref{e2659} and \eqref{e2663}, and based on \eqref{e1077} and \eqref{e1793}, defined on the same probability space as $\{\wh X_t^{\infty,\ddag}(\eta,a)\colon\,t >0,\,\eta \in \Omega_N^\infty\}$. 
\item
For every $\varphi \in \CD(\mathfrak{A}) \cap L^2(\mathfrak{h})$, given $A$ the following \emph{stochastic integro-differential equation} is satisfied:
\begin{equation}
\label{a.1917}\hspace{-1cm}
\begin{aligned}
&\langle \wh X_t^{\infty,\ddag},\varphi \rangle - \langle \wh X_0^{\infty,\ddag},\varphi \rangle\\ 
&\quad = \int_0^t \dd s\, \langle \wh X_s^{\infty,\ddag},\mathfrak{A} \varphi \rangle\\
&\qquad + \int_0^t \int_{\Omega_N^\infty} \sum_{u \in A} W^u(\dd s,\dd\eta)
\sqrt{\wh X_s^{\infty,\ddag} (\eta,u)(1-\wh X_s^{\infty,\ddag}(\eta,u))} \, \varphi(\eta,u)\\
&\qquad + \sum_{k \in \Z} \int_0^t \sum_{\eta^\prime \in \Omega_N^k} \Pi(k,\dd s,\eta^\prime)
\int_{(0,1]} r \Lambda_k(\dd r) \int_{B_k(\eta^\prime)} \mathfrak{h}(\dd\eta)\\
&\qquad \qquad \qquad \times \left(\sum_{a \in A} \wh X_s^{\infty,\ddag}(\eta^\prime,a)\, 
\varphi(\eta^\prime,a)- \sum_{u \in A \setminus \{a\}}
\wh X_s^{\infty,\ddag}(\eta,u)\, \varphi(\eta,u)\right)
\end{aligned}
\end{equation}
\item
$\Pi$ and $W$ are independent of $a \in \K$, and therefore the collection of densities $\{\wh X^{\infty,\dagger} (\cdot,a), a \in A\}$ with $\N$ can be relabelled, and \eqref{a.1917} can be rewritten for a process with values in $\CM_{\rm locfin} (\Omega_N^\infty \times \N)$ and unique solution in law on that space, which allows us to obtain the law of $\wh X^\infty$ as a functional by assigning a type $a \in \K$ to the labels in $\N$.
\end{itemize}
\end{itemize}
\end{proposition}

\begin{proof}{[Proof of Proposition \ref{prop.1869}]}
It suffices to show \emph{existence} of a solution, and to show that every solution \emph{satisfies the martingale problem} associated with the process $\wh X^\infty$ and hence is unique in law. 

The \emph{existence} of the process $(X_t^{\infty,\ddag})_{t \geq 0}$ follows in two ways, either from the existence and the properties of a solution of the martingale problem, or from the existence of the limit $\wh X^\infty$ of $\wh X^L$ as $L\to\infty$, provided we can construct the ingredients $\Pi$ and $W$ on the probability space for the solution of the $\wh X^\infty$-martingale problem, respectively, of the $L \to \infty$ limit. We focus on the latter possibility. 

Thus, to conclude the argument, we show that, starting with $\wh X^\infty$ and using it as a solution to the martingale problem, we get a solution to the stochastic equation in \eqref{a.1917}. The construction suggested in Proposition~\ref{prop.1869}, given the \emph{random} set $A$, has to be extended to cover the whole entrance law starting from any state, and gives the disintegration of this object w.r.t.\ the type variable $u \in \K$. For the path $\wh X^\infty$ we need to \emph{construct the ingredients $\Pi$ and $W$} on the probability space of the process $\wh X^\infty$, since having $\Pi$ and $W$ we can use \emph{Ito's formula} to go from the martingale problem to \eqref{a.1917}, using the formula for the operator of the martingale problem.

\medskip\noindent
(i) \emph{The point process $\Pi$}.
This is obtained via the process of resampling events (more specifically, the $\Lambda_k$-resampling events) in the form of a jump of the measure restricted to one of the $k$-balls of the decomposition of $\Omega_N^\infty$ into a collection of disjoint $k$-balls for $k \in \Z$, which identifies the jump times of $\wh X^\infty$,and thereby specifies the points of $\Pi$. We need to show that the jump points form a Poisson point process with intensity $\Lambda_k$ when the event concerns a $k$-ball of the decomposition of $\Omega_N^\infty$ into disjoint $k$-balls, which is read off from the path since the transition is given by a mapping of the states. The fact that this is a Poisson point process is standard and is read off from the form of the jump operator describing the jumps in the martingale problem. 

\medskip\noindent
(ii) \emph{The white noise measure $W$}.
The white noise $W$ has to be constructed based on the \emph{continuous} and \emph{infinite-variation part of the path of the martingale part}, to obtain the Brownian motions on the probability space on which $\wh X^\infty$ is defined, where we have an accumulated space-time volatility measure process (recall \eqref{ar1749}--\eqref{e1750})
\begin{equation}
\label{e2782}
\int_0^t \dd s\,L_{X^\infty_s} (\dd\eta, \dd u),
\end{equation}
representing the quadratic variation of the evaluations of test functions arising from the path $\wh X^\infty$. From here we have to obtain Brownian motions on the probability space of $\wh X^\infty$ via the evaluation processes time-transformed by the volatility measure.

The key object is a \emph{Gaussian field} related to $L_X(\dd\eta,\dd u)$, i.e., the \emph{different-type collision local time} describing the \emph{increasing process} of the evaluation process $\langle \wh X_t^\infty,\varphi \rangle_{t \geq 0}$ of $\wh X^\infty$. Via the duality relation the law of this field is related to the local time of the difference process of the effective migration and the associated \emph{difference local time functional} encountered earlier. This relation guarantees that in each time interval the quadratic variation increases when $\varphi$ is the indicator of some ball. We want to extend the latter property to the density process, so that the continuous part is a (random) time-transform of a Brownian motion. Phrased differently, we need to argue that we can use the expression for the density of the evaluation of the indicator of a ball to obtain the quadratic variation of the limit process, and that this process has positive quadratic variation over every positive time interval based on the collision local time of the dual path $\zeta$, so that we can construct the corresponding Brownian motion. 

The construction of a Brownian motion from the quadratic variation of a scalar continuous local martingale follows from the Dambis-Dubins-Schwarz theorem \cite[Theorem V.3.9]{RY99}. To see how this can be lifted to the measure-valued context, see \cite[Section 7.1]{D93}. For the particular setting of the hierarchical group, see \cite{Klenke1996} and \cite{EF96}.
\end{proof}


\section{Continuum-space scaling limit}
\label{s.results}

In this section we prepare for the proof of Theorem~\ref{T.CONVRES}, in which we claim that a special case of the continuum-space process in Example \ref{ex.1035} constructed in Sections \ref{sec:trlim}--\ref{sec:wpmp}, based on \emph{stable migration} on $\Omega_N^\infty$ for a single tagged individual, also appears as the \emph{scaling limit} of the Cannings process on $\Omega_N$. In particular, the scaling is indexed by a parameter $L$ such that $L\to\infty$, and the scaled processes on $\Omega_N$ can be viewed as a process on $\Omega_N^L$ after we appropriately scale space and time and scale the migration and resampling rates. This is done in such a way that the \emph{space-time scaled} process on $\Omega^L_N$ has the property that the underlying random walk converges to a L\'evy-process $\zeta^\infty$ on $\Omega_N^\infty$ and the mass in finite sets of $\Omega_N^\infty$ remains positive and finite as $L \to \infty$. This property singles out the \emph{stable} Cannings processes (i.e., stable under the scaling operation) as the canonical representatives of the \emph{universality classes} for the large space-time behaviour of the Fleming-Viot and Cannings processes with block-resampling on $\Omega_N$. Applying the analogous idea and scaling to the processes on $\Omega_N^\infty$, we obtain the same models based on stable migration on $\Omega_N^\infty$ as the \emph{representatives of the universality class} of the continuum dynamics, which are now \emph{fixed points} under the scaling operation.

Section \ref{sss.keypr} states a key proposition, Section~\ref{ss.resc} exhibits the scaling of space, time, migration rates and resampling rates, Section \ref{ss.rescon} proves convergence after scaling.


\subsection{Key proposition}
\label{sss.keypr}

We saw already in the previous sections that for our state spaces, to have convergence of the finite-dimensional distributions of the scaled processes, it suffices to show convergence of the \emph{mixed space-time moment measures} under the scaling to the ones of the limiting dynamics on $\Omega_N^\infty$. We conclude from the \emph{duality} for the scaled system and the analysis of the measure-valued versions of the Cannings processes constructed in Sections \ref{sec:trlim}--\ref{sec:wpmp} that, to get convergence of finite-dimensional distributions in the original measure-valued process $\wh X$ on $\Omega_N$ or $\wh X^\infty$ on $\Omega_N^\infty$, it suffices to show convergence of the space-time dual processes under the scaling to the dual process of the limiting process. 

Therefore we have to show that under the scaling the dual process $Y^L$ converges to $Y^\infty$ as $L\to\infty$. For that purpose we choose initial positions $\xi_1, \ldots,\xi_n \in \Omega_N^\infty$, $n\in\N$, and take the $L$-truncation of these as initial positions for the process $Y^L$ on $\Omega_N^L$, to get the dual process of $\wh X^L$ on $\Omega_N^L$, respectively, the dual of the space-time process. 

Furthermore, we need to consider the \emph{entrance law} of the coalescent with countably many individuals per site, since the duality relation on the continuum space requires the integrated version in space, i.e., both sides of the duality relation have to be formulated as \emph{measure-valued} processes and the types have to be assigned in the dual as well. The limiting dual process will again be denoted by $(\wh Y_t^{\infty,\ast})_{t \geq 0}$. Recall that, in order for this process on to be well defined, a key ingredient is the \emph{no-dust condition} for $\Lambda_{-\infty}$ assumed in \eqref{cond3}.

The key to the proof of Theorem \ref{T.CONVRES} in Section \ref{s.dual} is the following proposition.

\begin{proposition}{\bf [Convergence of the scaled dual]}
\label{P.convscaldual}
Subject to conditions \eqref{ag7a}--\eqref{ag7d} below,
\be{e1104}
\CL\big[(Y_t^L)_{t \geq 0}\big] \Lto \CL \big[(Y_t^\infty)_{t \geq 0} \big],
\ee
and
\be{e1104alt}
\begin{aligned}
\CL\big[(\wh Y_t^{\ast,L})_{t \geq 0}\big] &\Lto \CL\big[(\wh Y_t^{\ast,\infty})_{t \geq 0}\big],\\
\CL\big[(\wh Y_t^{\ast,L})_{t \geq 0}\big] &\Lto \CL\big[(\wh \xi_t^{\ast,\infty})_{t \geq 0}\big],
\end{aligned}
\ee
where in the last display the convergence holds for the recurrent regime (no dust), respectively, the transient regime (complete dust).
\end{proposition}

\noindent
We already know from the literature \cite{DGZpr} that the above claim holds for any collection of scaled independent L\'evy-processes. However, we have to do more to get a solution to the martingale problem in path space in the limit $L \to \infty$, namely, we need to prove \emph{convergence of the coalescence times} to be able to get convergence of finite-dimensional distributions, as well as \emph{tightness} in the space $D([0,\infty),\wt\CM(\Omega_N^\infty))$. 

The proof of Proposition~\ref{P.convscaldual} is organised as follows. In Section \ref{ss.resc} we introduce the scaling of the original process. In Section \ref{ss.rescon} we look at scaling and convergence of the dual process. In Section \ref{ss.proofprop} we use the latter to prove Proposition \ref{P.convscaldual}.


\subsection{Scaling of space, time, migration rates and resampling rates}
\label{ss.resc}

We focus on the case where both the migration and the migration induced by the block-resampling through the block-reshuffling are in the \emph{domain of attraction} of \emph{stable laws}. In Section \ref{s.dual} we briefly comment on possible generalisations to the domain of attraction of \emph{infinitely divisible} laws. 

In what follows we make the following \emph{regularity assumptions} on the migration and resampling rates
(recall \eqref{ag40b}):
\be{ag7a}
c_k = c^k \bar c_k, \qquad \lambda_k = \lambda^k \bar \lambda_k, 
\qquad \lambda^*_k = \lambda^k \bar{\lambda}^*_k, \qquad k \in \Z,
\ee
with
\be{ag7b}
\begin{aligned}
&c,\lambda \in (0,N), \quad \bar c_k,\bar \lambda_k, \bar{\lambda}^*_k \in (0,\infty), \quad k\in\Z,\\
&\lim_{k\to\infty} \frac{\bar c_{k+1}}{\bar c_k} = 1, 
\quad \lim_{k\to\infty} \frac{\bar \lambda_{k+1}}{\bar \lambda_k} = 1,
\quad \lim_{k\to\infty} \frac{\bar{\lambda}^*_{k+1}}{\bar{\lambda}^*_k} = 1,
\end{aligned}
\ee
and
\be{ag7c}
\frac{\Lambda_k}{|\Lambda_k|} \kto \Lambda,
\ee
with
\be{ag7d}
|\Lambda| = 1, \quad \Lambda(\{0\}) = 0, \quad \int_{(0,1]} r^{-1} \Lambda(\dd r) < \infty.
\ee

We next introduce the scaling of \emph{space, mass, time and migration/resampling rates}, leading for every $L \in \N$ to a Cannings model with block-resampling on $\Omega^L_N$. We choose the \emph{space-time scaling} such that the \emph{migration-reshuffling process} becomes tight and converges to the corresponding \emph{L\'evy-process} on $\Omega_N^\infty$. Namely,
\begin{itemize}
\item
\emph{Shrink space}: 
\begin{itemize}
\item[(i)] 
Replace migration jumps of size $k$ on $\Omega_N^L$ by
\be{ag2a}
\mbox{jumps of size $k-L$},
\ee
i.e., a jump in a block of size $k$ on $\Omega_N$ is turned into a jump in a block of size $k-L$ on $\Omega^L_N$.
\item[(ii)] 
Shift the resampling by using for the single-site resampling at level $L$ the measure $\Lambda_{-\infty}$ and by using the block-resampling measure $\Lambda_k$ at level $k$ on $\Omega_N$ for the block-resampling in $(k-L)$-blocks on $\Omega^L_N$.
\end{itemize} 
\item
\emph{Speed up time}: Replace time $t$ by
\be{ag2}
\begin{array}{lll}
&c \geq \lambda\colon &\mbox{time } \frac{N^L}{c_L} t,\\[0.2cm]
&c \leq \lambda\colon &\mbox{time } \frac{N^L}{\lambda_L} t.
\end{array}
\ee
\end{itemize}
Next, we define an additional scaling of the \emph{mass} and the \emph{resampling rates} to obtain tightness of the complete process.
\begin{itemize}
\item
\emph{Embedding and mass scaling:} Embed $\Omega^L_N$ into $\Omega^\infty_N$ and \emph{scale the total mass} by
\be{ag2b}
N^{-L}
\ee
to compensate for the increasing number of components in a fixed macroscopic component in $\Omega_N^\infty$ in the system on $\Omega^L_N$ as $L$ increases.
\item
\emph{Scale the rate of resampling}: Replace the rate $d$ of the Fleming-Viot resampling on $\Omega^L_N$ by rate $d c_L$ or rate $d \lambda_L$, depending on which of the two dominates. Do the same for the Cannings single-site resampling rate and measure at level $k = -L$.
\end{itemize}

In this way we obtain a process taking values in $\wt\CM(\Omega^{L}_N \times \K)$, which by the above embedding induces a new process for every index $L$, namely, the process denoted by 
\be{ag2c}
(\wh X_t^L)_{t \geq 0}
\ee
taking values in $\wt\CM(\Omega^\infty_N \times \K)$ (not to be confused with the object defined in Section \ref{(ii)}). Stated in a single formula, the complete scaling above reads
\be{ar1966}
\begin{array}{lll}
&c \geq \lambda\colon &\wh X_t^L(\dd\eta,\dd u) = N^{-L} \wh X_{(N^L/c_L)t}(\dd(\pi_L \eta), \dd u)),\\[0.2cm]
&c \leq \lambda\colon &\wh X_t^L(\dd\eta,\dd u) = N^{-L} \wh X_{(N^L/\lambda_L)t}(\dd(\pi_L \eta), \dd u)),
\end{array}
\ee
for $(\eta,u) \in \Omega_N^\infty \times \K$, with $\pi_L$ the shrinking operator in \eqref{shr}.

In view of \eqref{ag2a}--\eqref{ag2}, the rates of the scaled process on $\Omega_N^L$ become 
\be{e1975}
(c^L_k)_{k \geq -L}, \qquad (\lambda^L_k)_{k \geq -L},
\ee 
where
\be{ag10}
\begin{array}{lll}
&c \geq \lambda \colon &c^L_k = \frac{1}{c_L} c_{k+L},\\[0.2cm] 
&c \leq \lambda \colon &c^L_k = \frac{1}{\lambda_L} c_{k+L},\\[0.2cm]
&c \geq \lambda \colon &\lambda^L_k = \frac{1}{c_L} \lambda_{k+L},\\[0.2cm]
&c \leq \lambda \colon &\lambda^L_k = \frac{1}{\lambda_L} \lambda_{k+L}.
\end{array}
\ee
For the Fleming-Viot resampling we get
\be{ag12}
\begin{array}{lll}
&c \geq \lambda\colon &d^L_{-L} = d N^L c_L,\\[0.2cm]
&c \leq \lambda\colon &d^L_{-L} = d N^L\lambda_L,\\[0.2cm] 
&&d^L_k = 0, \,k> -L.
\end{array}
\ee
The above defines a sequence indexed by $L$ of Markov processes $\{\wh X^L, L \in \N\}$ on $\Omega_N^L$. 

We next observe that the parameters of this sequence of processes have limiting values as $L \to \infty$. Namely, from our assumptions \eqref{ag7a}--\eqref{ag7c} we get the following.

\begin{lemma}{\bf [Limiting parameter values]}
\label{l.2163}
\be{ag13}
\begin{aligned}
&\lim_{L\to\infty} c^L_k = \left\{\begin{array}{ll}
c^k, &c \geq \lambda,\\
0, &c< \lambda, 
\end{array}
\right.
\qquad k \in \Z,\\
&\lim_{L\to\infty} \lambda^L_k = \left\{\begin{array}{ll}
\lambda^k, &c \leq \lambda,\\
0, &c>\lambda, 
\end{array}
\right.
\qquad k \in \Z.
\end{aligned}
\ee
Thus, for $c < \lambda$ the migration is asymptotically negligible, for $c > \lambda$ the reshuffling and block-resampling are asymptotically negligible, while for $c=\lambda$ the two mechanisms are both present on all hierarchical levels.
\end{lemma}

The above observation suggests that the scaled processes $\wh X^L$ converges to the process $\wh X^\infty$ from Theorem \ref{T.CONV} with the choice of the parameters being the above values, as stated in Theorem \ref{T.CONVRES}. To prove this, we have to show convergence of the corresponding dual processes, as stated in Proposition \ref{P.convscaldual}. 


\subsection{Scaling and convergence of the finite spatial coalescents}
\label{ss.rescon}

We first recall the invariance principle for the random walks underlying our processes. After that we pass to the finite coalescents and their limit as $L \to \infty$, because the scaling must be such that the random walks converge to their scaling limit. Scaling of the volatility rates must lead to a stochastic limit where we see coalescence in finite coalescents.


\paragraph*{Invariance principle for the L\'evy-process.}

Let $\zeta^\infty=(\zeta^\infty_t)_{t \geq 0}$ be the path of a \emph{single individual} in the dual, which is a random walk with migration rates $(\bar c_k)_{k\in\Z}$ given by (recall \eqref{effmig})
\be{ag3a}
\bar c_k = c_k + \frac{\lambda^*_{k+1}}{N}, \qquad k \in \Z.
\ee
Let $\zeta^L=(\zeta^L_t)_{t \geq 0}$ be the corresponding scaled path on $\Omega^L_N$, obtained via \eqref{ag2a}--\eqref{ag2}, embedded into $\Omega^\infty_N$. From \cite[Section 2.1]{DGZpr} (which uses the Fourier transform) we know that, under the assumptions in \eqref{ag7a}--\eqref{ag7c}, 
\be{ag3}
\CL\big[(\zeta^L_t)_{t \geq 0}\big] \Lto \CL\big[(\zeta^\infty_t)_{t \geq 0}\big],
\ee
where $\zeta^\infty$ is the L\'evy-process on $\Omega^\infty_N$ with migration rates $(\wh c_k)_{k \in \Z}$ given by
\be{ag4}
\wh c_k = \left\{\begin{array}{ll}
c^k, &c > \lambda,\\
c^k \left(1 + \frac{\lambda}{N}\right), &c = \lambda,\\
\frac{1}{N} \lambda^{k+1}, &c < \lambda,
\end{array}
\right.
\qquad k \in \Z.
\ee
In Appendix \ref{appA:IPLT} we state the invariance principle in \eqref{ag3} (which is immediate consequence of the techniques developed in \cite[Section 2.1]{DGZpr} carried over to our setting). It is straightforward to extend this to an invariance principle for an $m$-tuple of L\'evy processes with $m \in \N \setminus \{1\}$.  


\paragraph*{Convergence for coalescing L\'evy-processes and entrance laws from countable many individuals per site.}

With duality in mind, the key question that remains to show is that the \emph{spatial coalescent} corresponding to the scaled system on $\Omega^L_N$, starting from a \emph{finite} number of partition elements, converges as $L\to\infty$ to a \emph{block-coalescing system of L\'evy-processes.} From pathwise convergence of the underlying random walks we obtain convergence of the spatial coalescents with the correct scaling of the resampling parameters. Since the migration converges after scaling, the question is whether the \emph{coalescent} $Y^{L,\ast}$ converges to $Y^{\infty,\ast}$ and the \emph{measure-valued process} $\wh Y^{L,\ast}$ converges to $\wh Y^{\infty,\ast}$ (with proper entrance laws). 

First, we settle the convergence to $Y^{\infty,\ast}$.

\begin{lemma}{\bf [Convergence of the finite dual under scaling]}
\label{lem:finscaldual}
As $L\to\infty$, $Y^{L,\ast}$ converges in path space to $Y^{\infty,\ast}$. The limiting coalescent $Y^{\infty,\ast}$ is a collection of interacting L\'evy-processes with: 
\begin{itemize}
\item[(i)]
migration at rates $(\wh c_k)_{k \in \Z}$ given by \eqref{ag4};
\item[(ii)] 
single-site coalescence at rate $d^\ast$ when $d \in (0,\infty)$ or $|\Lambda_{-\infty}| \neq 0$, and no single-site coalescence when $d=0$ and $|\Lambda_{-\infty}|=0$; 
\item[(iii)] 
block-coalescence at rates $(\lambda^k\Lambda)_{k \in \Z}$. 
\end{itemize}
\end{lemma}

\begin{proof}
We need to distinguish between the recurrent regime and the transient regime.

\medskip\noindent
\emph{(I) Recurrent $\zeta^\infty$.}
We scaled the random walks to get the L\'evy-process $\zeta^\infty$. We also scaled the coalescence so that the occupation time of two random walks shrinks to zero. In the Harris-recurrent regime the block-coalescence gives coalescence after a finite random time for each pair of individuals (because of the upscaling of the coalescence rate with $\lambda^L$ as $L\to\infty$). In the Harris-recurrent but not point-recurrent regime, we simply have block-coalescence as the limit. In the point-recurrent regime we have \emph{delayed} coalescence as the limit, based on the \emph{local time} upon hitting (see Sections~\ref{sss.first}--\ref{sss.second} and Section~\ref{s.third}). 

In Appendix \ref{appA:IPLT} we indicate that the local times of a single L\'evy-process on $\Omega_N^L$ \emph{exists} and \emph{converges}, in the limit as $L\to\infty$, to the local times of a single L\'evy-process on $\Omega_N^\infty$. The proof uses that the transition kernels of both processes can be explicitly identified (see Appendix \ref{appA:pot}). It is not hard to see that the same holds for the collision local times of all the pairs in an $m$-tuple with $m \in \N\setminus \{1\}$, at the expense of getting more involved formulas. We refrain from writing this out, because these formulas are not needed. See \cite{EF96} for further details (which has no block-reshuffling and no block-resampling).

\medskip\noindent
\emph{(II) Transient $\zeta^\infty$.}
Without block-resampling, the claim follows from the invariance principle describing the pathwise convergence of the random walks to the L\'evy-processes. With block-resampling, for jumps with a size above any small threshold the claim again follows from the same argument. The smaller jumps can, however, be included because the $L$-approximation is uniform in path space.  
\end{proof}

Next, we settle the convergence to $\wh Y^{\infty,\ast}$.

\begin{lemma}{\bf [Convergence of the finite measure-valued process under scaling]}
\label{lem:finscalmeas}
As $L\to\infty$, $\wh Y^{L,\ast}$ converges in path space to $\wh Y^{\infty,\ast}$ in the recurrent regime, respectively, to $\wh\zeta^\infty$ in the transient regime. The limiting measure-valued process $\wh Y^{\infty,\ast}$ is the stable Cannings process with block-resampling, with parameters $(c_k)_{k\in\Z}$ and $(\lambda_k)_{k\in\Z}$ given by \eqref{ag7a}--\eqref{ag7b} and \eqref{ag3a}.
\end{lemma}

\begin{proof}
We show that at positive times the entrance law turns into a locally finite coalescent, which subsequently evolves according to the limit coalescent that we get by combining the coalescences in (ii) and (iii). The argument follows the pattern in Sections~\ref{sec:dual}--\ref{s.third}.
\end{proof}

We have seen in Propositions \ref{prop.1010}--\ref{P.DUAL} and Corollary \ref{cor.1815} that convergence of the finite-dimensional distributions of $\wh X^L$ to those of $\wh X^\infty$ is implied by the same convergence for the dual, which we have settled in Lemma \ref{lem:finscalmeas}. From \eqref{ag3} and Proposition~\ref{prop.two-four-regimes} we get the following corollary.

\begin{corollary}{\bf [Deterministic versus stochastic]}
\label{cor.2333}
$\wh X^\infty$ is deterministic when $c \vee \lambda >1$ and stochastic when $c \vee \lambda \leq 1$.
\end{corollary}


\subsection{Proof of the key proposition}
\label{ss.proofprop}

The claim follows from Lemmas \ref{lem:finscaldual}--\ref{lem:finscalmeas}.


\section{Proofs of the main theorems}
\label{s.dual}

In this section we prove the 7 theorems stated in Section \ref{ss.outline} based on the key propositions and the corollaries stated Sections \ref{sec:dual}--\ref{s.results}. The proofs of the first 6 theorems is immediate. The proof of the last theorem requires some further considerations.


\subsection{Continuum-space Cannings process: Existence and well-posedness}

Theorems \ref{T.CONV}(a)--(c) follow from Propositions \ref{prop.1010} and \ref{P.solwpmp}, while Theorem \ref{T.CONV}(d) follows from Proposition \ref{prop.1869}.


\subsection{Coincidence of dichotomies}

Theorem \ref{T.det} follows from Propositions \ref{P.DUAL} and \ref{P.1430}, Corollary \ref{C.sdual} and Remark \ref{rem.threecases}. 


\subsection{Criterion for the dichotomy}

Theorems \ref{T.STOCHDET}(a)--(b) follow from Corollary \ref{cor.two-four-regimes}. 


\subsection{Properties of continuum-space Cannings process}

Theorem \ref{T.char}(a),(b1)--(b2) follow from Proposition \ref{P.DUAL} and Corollary \ref{C.sdual}. 


\subsection{Convergence to equilibrium}

Theorem \ref{T.conveq} follows from Proposition \ref{P.DUAL}. 


\subsection{Coincidence of all dichotomies}

Theorem \ref{T.coindichs} follows from Proposition \ref{P.DUAL} and Corollary \ref{C.sdual}.


\subsection{Convergence in the stable domain of attraction}

To prove Theorem \ref{T.CONVRES}, we argue that the scaled process on $\Omega^L_N$ embedded into $\Omega^\infty_N$ converges as $L \to \infty$ to the \emph{stable} continuum-space process on $\Omega^\infty_N$, of which we have determined the parameters. To the candidate for the limit process we can apply all the results obtained in Sections \ref{sec:trlim}--\ref{sec:wpmp}, where we constructed the continuum-space process. In particular, we have to argue that, by the duality and the \emph{scaling limit of the dual obtained in Proposition~\ref{P.convscaldual}}, we get a limit process $\wh X^\infty$ (in the sense of finite-dimensional distributions) that is precisely the process constructed in Sections \ref{sec:trlim}--\ref{sec:wpmp}, with a specific choice of the parameters, given explicitly in Lemma \ref{l.2163} as consequence of our conditions on the parameters.

We still have to argue that the scaled process also \emph{converges in the path space} $D([0,\infty),\wt\CM(\Omega_N^\infty \times \K))$. But we have already established \emph{tightness}, and since we have identified a unique candidate for the limit, whose path properties we know, this stronger convergence follows. We have to show that the dual system has the \emph{no-dust property}, which is needed to obtain a non-degenerate stochastic process $\wh X^\infty$ on $\Omega_N^\infty$. But we can settle this with the help of the general criteria in \cite{EMS13}. Namely, the strong form of the duality in Corollary \ref{C.sdual}, in combination with Proposition~\ref{P.convscaldual}, gives us criteria for dust versus no-dust formation in $\wh X^\infty$.

\begin{remark}{\bf [Convergence in the domain of attraction of infinite-divisibility]}
{\rm We may apply the same reasoning to the setting where the rates of migration and resampling depend on the scaling parameter $L$ in a more general manner. The question is whether we can obtain a scaling limit with an \emph{array} of rates
\be{ag59b}
(c_{k,L})_{k \in \N_0}, \quad (\lambda_{k,L})_{k \in \N_0}.
\ee
We can speed up time by
\be{ag59c}
\frac{N^L}{c_{L,L}}, \mbox{ respectively, } \frac{N^L}{\lambda_{L,L}}, 
\ee
depending on whether
\be{ag59d}
\lim_{L \to \infty} \frac{c_{L,L}}{\lambda_{L,L}} \geq 1 \quad \mbox{ or } \quad 
\lim_{L \to \infty} \frac{c_{L,L}}{\lambda_{L,L}} \leq 1.
\ee
If
\be{ag59e}
\lim_{L \to \infty} \frac{c_{L+k,L}}{c_{L,L}} = \bar c_k, \quad 
\lim_{L \to \infty} \frac{\lambda_{L+k,L}}{\lambda_{L,L}} = \bar \lambda_k, \quad 
\lim_{L \to \infty} \frac{\lambda^*_{L+k,L}}{\lambda_{L,L}} = \bar{\lambda}^*_k, \quad k \in \Z,
\ee
then we get convergence of the scaled process to a Cannings process on $\Omega_N^\infty$ with parameters $(\bar c_k)_{k\in\Z}$ and $(\bar \lambda_k)_{k\in\Z}$. We will not pursue this line of thought further.} \hfill$\spadesuit$
\end{remark} 

\begin{remark}{\bf [Euclidean setting]}
{\rm On $\Z^d$, a non-trivial and stochastic continuum-space limit for the Fleming-Viot process exists only for $d=1$. For $d \geq 2$ the limiting dual process does not allow for collision and hence for coalescence of the Brownian motions, so that the continuum-space limit becomes deterministic. The latter is different for the Cannings process with block-coalescence. The lack of collisions can be overcome through the block-coalescence, which only requires the paths to come sufficiently close to each other sufficiently often.}\hfill$\spadesuit$
\end{remark}


\begin{appendix}


\section{Basic facts about random walks on $\Omega_N$ and L\'evy-processes on $\Omega_N^\infty$}
\label{appA}

In Appendix \ref{appA:pot} we recall the \emph{potential theory} associated with the random walk on $\Omega_N$ with migration rates $(c_k)_{k\in\N_0}$ (which represents the migration of a dual individual) and see how this extends to the L\'evy-process on $\Omega_N^\infty$ with migration rates $(c_k)_{k\in\Z}$. In particular, formulas are derived for the transition kernel and the Green function on $\Omega_N^\infty$. In Appendix \ref{appA:IPLT} we derive an invariance principle for the migration on $\Omega_N^\infty$ and show that it has a jointly continuous local time. In Appendix \ref{appA:oneLP} we identify when the motion of a \emph{single} dual individual is recurrent or transient. In Appendix \ref{appA:twoLP} we study the coalescence of \emph{two} dual individuals. In Appendix \ref{app:stable} we zoom in on stable migration. In Appendix \ref{app:digression} we offer a digression that provides additional insight. 


\subsection{Potential theory}
\label{appA:pot}

We first look at the random walk on $\Omega_N$ with migration rates $(c_k)_{k\in\N_0}$ satisfying \eqref{ag35}. After that we look at the L\'evy-process on $\Omega_N^\infty$ with migration rates $(c_k)_{k\in\Z}$ satisfying \eqref{cond1}--\eqref{cond2}.


\paragraph*{Discrete hierarchical group.} 

The time-$t$ transition kernel on $\Omega_N$ can be computed in closed form (see \cite{DGW05} for related examples).

\begin{lemma}{\bf [Transition kernel on $\Omega_N$]}
For every $t \geq 0$ and $\xi,\eta \in \Omega_N$,
\be{e2553}
p_t(\xi,\eta) = \sum_{\ell \geq k(\xi,\eta)} \frac{1}{N^\ell} \left[\eee^{-s_\ell t} - \eee^{-s_{\ell-1}t}\right],
\ee
where $k(\xi,\eta) \in \N_0$ is the hierarchical distance between $\xi$ and $\eta$ (recall \eqref{ag28}), and $(s_k)_{k\in\N_0}$ is given by (recall \eqref{fdhkernel})  
\be{e2561}
s_k = \sum_{\ell>k} \frac{c_{\ell-1}}{N^{\ell-1}} = \sum_{\ell \geq k} \frac{c_\ell}{N^\ell} 
\ee  
with the convention that $s_{-1} = \infty$ (to kill the last term between square brackets in \eqref{e2553} when $k=\ell=0$). 
\end{lemma}

\begin{proof}
Simply note that $\eee^{-s_\ell t}$ is the probability that up to time $t$ no block-migration larger than size $\ell$ has taken place in $\zeta$, and $\eee^{-s_\ell t} - \eee^{-s_{\ell-1}t}$ is the probability that up to time $t$ the largest block-migration has size $\ell$. Conditional on the latter event, the distribution of $\zeta$ is uniform over the $\ell$-block around $\xi$, which has $N^\ell$ sites.
\end{proof}

Note that the terms under the sum in \eqref{e2553} are non-negative and that the kernel is properly normalised:
\be{}
\begin{aligned}
\sum_{\eta \in \Omega_N} p_t(\xi,\eta) 
&= \eee^{-s_0 t} + \sum_{\ell\in\N} \frac{1}{N^\ell} \left[\eee^{-s_\ell t} - \eee^{-s_{\ell-1}t}\right]
+ \sum_{k\in\N} (N^k-N^{k-1}) \sum_{\ell \geq k} \frac{1}{N^\ell} \left[\eee^{-s_\ell t} - \eee^{-s_{\ell-1}t}\right]\\
&= \eee^{-s_0 t} + \sum_{\ell\in\N} \frac{1}{N^\ell} \left[\eee^{-s_\ell t} - \eee^{-s_{\ell-1}t}\right]
+ \sum_{l\in\N} \frac{1}{N^\ell} \left[\eee^{-s_\ell t} - \eee^{-s_{\ell-1}t}\right] \sum_{k=1}^\ell (N^k-N^{k-1})\\
&= \eee^{-s_0 t} + \sum_{\ell\in\N} \left[\eee^{-s_\ell t} - \eee^{-s_{\ell-1}t}\right] = 1,
\end{aligned}
\ee
where we twice use telescoping and note that $\lim_{\ell\to\infty} s_\ell = 0$ because of the summability condition in \eqref{ag35}. 

A similar formula holds for the random walk on the \emph{$L$-refined hierarchical group} $\Omega_N^L$, $L \in \N$, defined in \eqref{ag29}, namely,
\be{e2553L}
p^L_t(\xi,\eta) = \sum_{\ell \geq k^L(\xi,\eta)} \frac{1}{N^\ell} \left[\eee^{-s_\ell t} - \eee^{-s_{\ell-1}t}\right],
\qquad t \geq 0,\,\,\xi,\eta \in \Omega_N^L,
\ee
where $k^L(\xi,\eta) \in \{-L,\ldots,-1\} \cup \N_0$ is the hierarchical distance between $\xi$ and $\eta$ (recall \eqref{ag31}), and we employ the convention that $s_{-L-1} = \infty$.


\paragraph*{Continuum hierarchical group.}

We can easily extend \eqref{e2553L} for the random walk $\zeta^L$ on $\Omega_N^L$ to the L\'evy-process $\zeta^\infty$ on $\Omega_N^\infty$. 

\begin{lemma}{\bf [Transition kernel density on $\Omega_N^\infty$]}
\label{lem:kerneldensity}
For every $t \geq 0$ and $\xi,\eta \in \Omega_N^\infty$,
\be{e2553infty}
p_t(\xi,\dd\eta) = p_t^*(\xi,\eta)\, \mathfrak{h}(\dd\eta), 
\qquad p_t^*(\xi,\eta) = \sum_{\ell \geq k^\infty(\xi,\eta)} \frac{1}{N^\ell} \left[\eee^{-s_\ell t} - \eee^{-s_{\ell-1}t}\right],
\ee
where $k^\infty(\xi,\eta) \in \Z$ is the hierarchical distance between $\xi$ and $\eta$ (recall \eqref{hierdist}), $(s_k)_{k\in\Z}$ is given by \eqref{e2561}, and $\mathfrak{h}$ is the normed Haar measure on $\Omega_N^\infty$. (The sum represents a density because $\Omega_N^\infty$ is uncountable. The earlier convention is dropped because there is no boundary term.)
\end{lemma}

\begin{proof}
First observe that $\zeta^L$ and $\zeta^\infty$ can be coupled such that the distance between their paths is bounded by $\eee^{-L}$ uniformly over time (recall \eqref{ag29}--\eqref{ag30}). Namely, $\zeta^\infty$ is a refinement of $\zeta^L$ when the same clocks are used for the jumps at block level $\ell$ with $\ell \geq - L$, i.e., both can be realised on the same probability space that was used to construct $\zeta^\infty$. As a consequence, the Green functions $G^L(\xi,\cdot)$ and $G^\infty(\xi,\cdot)$ of $\zeta^L$ and $\zeta^\infty$, viewed as measures on $\Omega_N^\infty$, differ in the Prokhorov metric by at most by $\eee^{-L}$. This means that the potential associated with the migration can be controlled uniformly in space and time as $L\to\infty$. 
\end{proof}

Later we will see that in the point-recurrent regime the singularity of $p^\ast_t (0,0)$ as $t \downarrow 0$ is integrable, i.e.,
\begin{equation}
\label{e2949}
\int_0^T p^\ast_t(0,0)\, \dd t < \infty,
\end{equation}
and that, under extra regularity conditions, $ p^\ast_t(0,0) \asymp t^{-(1+\gamma)}$ with $\gamma$ the \emph{degree} of the effective migration.

We can introduce the \emph{Fourier transform} and \emph{characteristic function} $\wh p_t(\cdot)$ of the probability measure $p_t(0,\cdot)$ on $\Omega_N^\infty$ as follows. Denote by $(\Omega_N^\infty)^\ast$ the set of elements of $\Omega_N^\infty$ that have the form $(\xi_i)_{i \in \Z}$ with $\xi_i \equiv 0$ for all $i \leq -\ell$ for some $\ell \in \N$. Define
\begin{equation}
\label{2949}
\langle \xi,\xi^\ast \rangle = \sum_{i \in \Z} \xi_i \xi_i^\ast, 
\qquad \xi^\ast \in (\Omega_N^\infty)^\ast,\,\xi \in \Omega_N^\infty,
\end{equation}
and, for a finite measure $\mu$ on $\Omega_N^\infty$, set
\begin{equation}
\label{e2953}
\wh \mu(\xi^\ast) = \int_{\Omega_N^\infty} \eee^{- \frac{2\Pi i}{N}\langle \xi,\xi^\ast \rangle} \mu(\dd\xi), 
\qquad \xi^\ast \in \Omega_N^\infty.
\end{equation}
The normed Haar measure appearing in \eqref{e2553infty}, replacing the counting measure in \eqref{e2553} and \eqref{e2553L}, is such that 
\be{}
\int_{\Omega_N^\infty} \mathfrak{h}(\dd\eta)\,1_{\{k^\infty(\xi,\eta) \leq k\}} = N^k, \qquad k \in \Z,\, \xi \in \Omega_N^\infty.
\ee


\subsection{Invariance principle and local times}
\label{appA:IPLT}

We consider the $(c_k)_{k \in \Z}$-L\'evy-process $\zeta^\infty$ defined in Section \ref{Levydef} as a path process and establish an invariance principle as well as existence of local times.

\begin{lemma}{\bf [Invariance principle]}
\label{lemmA3}
The $(c_k)_{k \in \Z}$-L\'evy-process $\zeta^\infty$, under the regularity assumptions and the scaling of space, time, mass and migration/reshuffling rates stated in Section \ref{ss.resc}, satisfies an invariance principle with limit the stable migration defined in Remark \ref{R.levy}. 
\end{lemma}

\begin{proof}
The claim follows from the techniques developed in \cite[Section 2.1]{DGZpr}.
\end{proof}

\begin{lemma}{\bf [Local times]}
\label{lemmA4}
The $(c_k)_{k \in \Z}$-L\'evy-process $\zeta^\infty$ has jointly continuous local times, satisfying the properties in \eqref{e1926}--\eqref{e1936}. 
\end{lemma}

\begin{proof}
We need four \emph{regularity properties} of the \emph{$\lambda$-resolvent} $P_t^\lambda (\cdot)$ of $(P_t)_{t \geq 0}$ defined by
\begin{equation}
\label{e2967}
P^\lambda_t(\xi) = \int_0^\infty \eee^{-\lambda t} p_t^\ast (0,\xi)\, \dd t, \qquad \lambda \geq 0,\,\xi \in \Omega_N^\infty.
\end{equation}
The arguments in the proofs of \cite{EF96} of the statement analogous to Corollary~\ref{cor.1939} (which are given in\cite[Proposition 13 and Corollary 14]{EF96}) are based on the explicit calculation of the Fourier transform of $P_t^\lambda(\xi)$. 

This function $P^\lambda_t(\cdot)$ is shown to have a \emph{series representation} that has four regularity properties. The terms are \emph{positive and continuous} in $\xi$ and $t$, and are \emph{indicators of balls}, and the series is \emph{uniformly continuous} on $\Omega_N^\infty$. These arguments given for these properties in the stable case can therefore be carried over to get the same regularity properties for our model as well. 

Namely, we can use \eqref{e2553infty} to analyse the function $P^\lambda_t(\cdot)$ on $\Omega_N^\infty$ for $\lambda \geq 0$. The representation in \eqref{ptasymp}, expressing $p_t^\ast(0,0)$ as a weighted sum of indicators of balls, is crucial to carry through the necessary computations. Indeed, by the ultrametric structure of $\Omega_N^\infty$, such indicators are non-negative, bounded and uniformly continuous functions on $\Omega_N^\infty$.  
\end{proof}


\subsection{Recurrence and transience for a single dual individual}
\label{appA:oneLP}

We next investigate the recurrence versus transience properties of a \emph{single} individual in the coalescent on $\Omega_N^\infty$. We first consider the case without block-reshuffling. Afterwards we indicate what changes when block-reshuffling is added. 


\subsubsection{Migration without block-reshuffling} 
 
Consider a L\'evy-process $\zeta^\infty=(\zeta^\infty_t)_{t \geq 0}$ on $\Omega_N^\infty$ with coefficients $(c_k)_{k\in\Z}$ satisfying \eqref{cond1}--\eqref{cond2}. 
Define
\be{e2494}
\sigma^+_m = \sum_{k \geq m} \frac{1}{c_k}, \qquad
\sigma^-_m = \sum_{k < m} \frac{1}{c_k}, \qquad m \in \Z. 
\ee
Note that $\sigma^+_m$ and $\sigma^-_m$ are either finite for all $m \in \Z$ or infinite for all $m \in \Z$. We distinguish between four regimes:

\begin{proposition}{\bf [Classification of recurrence versus transience of a single L\'evy-process]}
\label{prop.one-four-regimes}
$\mbox{}$
\begin{itemize}
\item[(a)]
$\sigma^+_0<\infty$: $\zeta^\infty$ is transient. 
\item[(b)]
$\sigma^+_0=\infty$, $\sigma^-_0=\infty$: $\zeta^\infty$ is Harris-recurrent but not point-recurrent.
\item[(c)]
$\sigma^+_0=\infty$, $\sigma^-_0<\infty$: $\zeta^\infty$ is point-recurrent.
\end{itemize}
\end{proposition}

\begin{proof}
We first show that $\zeta^\infty$ is \emph{transient} when $\sigma^+_0 < \infty$ and \emph{Harris-recurrent} when $\sigma^+_0 = \infty$. Fix $\xi \in \Omega_N^\infty$ and write $\E_\xi$ to denote expectation over $\zeta^\infty$ given that $\zeta^\infty_0=\xi$. Let 
\be{e3020}
\CL_m = \int_0^\infty \dd t \,\,1_{\{k^\infty(\zeta^\infty_t,0) = m\}}, \qquad m \in \Z,
\ee
denote the total local time of $\zeta^\infty$ at hierarchical distance $m$ from $0$. From Lemma~\ref{lem:kerneldensity} we get
\be{a1}
\E_\xi[\CL_m] = \int_0^\infty \dd t \,\,\PP_\xi\big(k^\infty(\zeta^\infty_t,0) = m\big)
= \int_{\Omega_N^\infty} \mathfrak{h}(\dd\eta)\,1_{\{k^\infty(\eta,0) = m\}}\,\,G^\infty(\xi,\eta), \qquad m \in \Z,
\ee
with 
\be{a2}
G^\infty(\xi,\eta) = \int_0^\infty \dd t\,p^*_t(\xi,\eta)
\ee
the Green function density on $\Omega_N^\infty$. It follows from \eqref{e2553infty} that
\be{a3}
G^\infty(\xi,\eta) = \sum_{\ell \geq k^\infty(\xi,\eta)} \frac{1}{N^\ell} \left[\frac{1}{s_\ell} - \frac{1}{s_{\ell-1}}\right].
\ee
Abbreviate $k_\xi = k^\infty(\xi,0)$, and let $\CL_{<0} = \sum_{m<k_\xi} \CL_m$ be the total local time of $\zeta^\infty$ in the \emph{interior} of the ball around $0$ on whose boundary $\xi$ lies. Then \eqref{a1}--\eqref{a3} give
\be{threesum}
\E_\xi[\CL_{<0}] = \sum_{\ell \in \Z} \frac{1}{N^\ell} \left[\frac{1}{s_\ell} - \frac{1}{s_{\ell-1}}\right] 
\int_{\Omega_N^\infty} \mathfrak{h}(\dd\eta)\,\sum_{k \leq \ell}
1_{\{k^\infty(\xi,\eta) = k\}}\,1_{\{k^\infty(\eta,0) < k_\xi\}}.
\ee
Because $k^\infty(\cdot,\cdot)$ is an ultrametric, we have
\be{}
k^\infty(\xi,\eta) \vee k^\infty(\eta,0) = k^\infty(\xi,0) \qquad \forall\,\eta \in \Omega_N^\infty.
\ee 
Hence, in the double sum in the right-hand side of \eqref{threesum} the restriction $k = k_\xi$ applies. Consequently, \eqref{threesum} reduces to
\be{twosum}
\E_\xi[\CL_{<0}] = \sum_{\ell \geq k_\xi} \frac{1}{N^\ell} \left[\frac{1}{s_\ell} - \frac{1}{s_{\ell-1}}\right] 
\int_{\Omega_N^\infty} \mathfrak{h}(\dd\eta)\,1_{\{k^\infty(\eta,0) < k_\xi\}}, 
\ee
where we note that $k^\infty(\xi,\eta) = k_\xi$ is implied by $k^\infty(\eta,0) < k_\xi$. The integral in the right-hand side of \eqref{twosum} equals $N^{k_\xi}-N^{k_\xi-1}$, and so we get
\be{}
\E_\xi[\CL_{<0}] = \frac{N-1}{N^2} \sum_{\ell \geq k_\xi} \left[\frac{1}{s_\ell} - \frac{1}{s_{\ell-1}}\right] N^{-(\ell-k_\xi)}
= \frac{N-1}{N^2} \sum_{\ell \geq k_\xi} \frac{1}{s_\ell s_{\ell-1}}\, \frac{c_{\ell-1}}{N^{\ell-1}}\, N^{-(\ell-k_\xi)},
\ee
where we use \eqref{e2561}. 

Under the first condition in \eqref{regcond}, we have
\be{sapprox}
s_\ell \asymp \frac{c_\ell}{N^\ell} \quad \text{ uniformly in } \ell \in\Z,
\ee 
so that
\be{A21}
\E_\xi[\CL_{<0}] \asymp N^{k_\xi} \sum_{l \geq k_\xi} \frac{1}{c_l} = N^{k_\xi} \sigma^+_{k_\xi}.
\ee
It follows that, irrespective of the choice of $\xi$, if $\sigma^+_0=\infty$, then $\E_\xi[\CL_{<0}] = \infty$, implying Harris-recurrence, while if $\sigma^+_0<\infty$, then $\E_\xi[\CL_{<0}] < \infty$, implying transience. 

We next show that $\zeta^\infty$ is \emph{point-recurrent} when $\sigma^-_0 < \infty$ and \emph{not} point-recurrent when $\sigma^-_0 = \infty$. The argument uses that point-recurrence is equivalent to $t \mapsto p^*_t(0,0)$ being \emph{integrable} near $t=0$ (just like Brownian motion in $\R^d$ with $d<2$). Write
\be{A22}
\begin{aligned}
p^*_t(0,0) &= \sum_{\ell \in \Z} \frac{1}{N^\ell} \left[\eee^{-s_\ell t} -\eee^{-s_{\ell-1} t}\right]\\
&= \sum_{\ell \in \Z} \left[\frac{1}{N^\ell} \eee^{-s_\ell t} - \frac{1}{N^{\ell-1}}\eee^{-s_{\ell-1} t}\right]
+ \frac{N-1}{N} \sum_{\ell \in \Z} \frac{1}{N^{\ell-1}}\eee^{-s_{\ell-1} t}
\end{aligned}
\ee
and use that the first sum in the last line telescopes to zero for every $t>0$ because  
\be{e3075}
\lim_{\ell \to \infty} \frac{1}{N^\ell} \eee^{-s_\ell t} = 0, \qquad \lim_{\ell \to -\infty} \frac{1}{N^\ell} \eee^{-s_\ell t} = 0.
\ee
The first claim is obvious, the second follows because $s_\ell \geq c_\ell/N^\ell$ and the first condition in \eqref{condextra} implies that there exists an $M \in (0,N)$ such that $c_\ell \geq M^\ell c_0$ for all $\ell < 0$. The contribution from $\ell \in \N_0$ in the second sum is bounded by $1$. Hence 
\be{ptasymp}
p^*_t(0,0) = O(1) + \frac{N-1}{N} \sum_{\ell \in -\N} \frac{1}{N^{\ell-1}}\eee^{-s_{\ell-1} t}
\ee
and 
\be{e4088}
\int_0^T \dd t\,p^*_t(0,0) = O(T) + \frac{N-1}{N} \sum_{\ell \in -\N} \frac{1}{N^{\ell-1}s_{\ell-1}}(1-\eee^{-s_{\ell-1} T}).
\ee
Since $\lim_{\ell \to -\infty} s_\ell = \infty$, for every $T>0$ we have
\be{e4092}
\sum_{\ell \in -\N} \frac{1}{N^{\ell-1}s_{\ell-1}}(1-\eee^{-s_{\ell-1} T}) \asymp \sum_{\ell \in -\N} \frac{1}{N^{\ell-1}s_{\ell-1}} 
\asymp \sum_{\ell \in -\N} \frac{1}{c_{\ell-1}} = \sigma^-_0,
\ee
where we again use \eqref{sapprox}. Hence we see that $t \mapsto p^*_t(0,0)$ is integrable near $t=0$ if and only of $\sigma^-_0<\infty$.
\end{proof}


\subsubsection{Migration with block-reshuffling}
 
The reshuffling adds to the migration rates, namely, 
\be{e4104}
c_{k,N} = c_k + \frac{\lambda^*_{k+1}}{N}, \qquad k \in \Z,  
\ee
with $\lambda^*_\ell = \int_{(0,1]} r^{-1}\Lambda_\ell(\dd r)$ (recall \eqref{effmig}). With this change the same conclusions as in Proposition \ref{prop.one-four-regimes} apply, with the sums in \eqref{e2494} replaced by 
\be{e2494alt}
\sigma_{m,N}^+ = \sum_{k \geq m} \frac{1}{c_{k,N}}, \qquad \sigma_{m,N}^- = \sum_{k < m} \frac{1}{c_{k,N}}, \qquad m\in\Z.
\ee


\subsection{Coalescence for two dual individuals}
\label{appA:twoLP}

We next show that a similar computation can be done for the distance between \emph{two} individuals in the coalescent. Section~\ref{sss.migr} treats migration and block-reshuffling, Section~\ref{sss.block} adds the block-coalescence.


\subsubsection{Migration with block-reshuffling and without block-coalescence}
\label{sss.migr}

Consider two L\'evy-processes without block-resampling, $\zeta^1=(\zeta^1_t)_{t \geq 0}$ and $\zeta^2=(\zeta^2_t)_{t \geq 0}$, where for convenience we suppress the upper index $\infty$. Then the migration rates for the \emph{difference process} $\zeta^1-\zeta^2 = (\zeta^1_t-\zeta^2_t)_{t \geq 0}$ are
\be{1901pre}
\bar c_{k,N} = 2 c_k + \frac{\lambda^*_{k+1}-\lambda_{k+1}}{N}, \qquad k \in \Z.  
\ee
The same conclusions as in Proposition \ref{prop.one-four-regimes} apply to $\zeta^1-\zeta^2$, with the sums in \eqref{e2494} replaced by
\be{1901}
\bar{\sigma}^+_{m,N} = \sum_{k \geq m} \frac{1}{\bar{c}_{k,N}}, \qquad
\bar{\sigma}^-_{m,N} = \sum_{k < m} \frac{1}{\bar{c}_{k,N}}, \qquad m \in \Z. 
\ee
The difference process itself is of no interest to us, but it serves as a tool to include the block-coalescence.
  

\subsubsection{Migration with block-reshuffling and with block-coalescence}
\label{sss.block}

If $\zeta^1$ and $\zeta^2$ are at distance $m$, then they block-coalescence at rate $\sum_{\ell \geq m} \lambda_\ell/N^\ell$, because $\lambda_\ell/N^\ell$ is the rate at which two lineages in an $\ell$-block are \emph{both marked} and coalesce through resampling. Hence the total rate at which block-coalescence occurs from distance $m$ is 
\be{4140}
\bar{\CL}_m \sum_{\ell \geq m} \frac{\lambda_\ell}{N^\ell},
\ee
where $\bar{\CL}_m$ is the local time at $m$ for the L\'evy-process with migration rates $(\bar{c}_{k,N})_{k\in\Z}$. Consequently, the total rate of block-coalescence equals
\be{e4144}
\CR = \sum_{m\in\Z} \bar{\CL}_m \sum_{\ell \geq m} \frac{\lambda_\ell}{N^\ell}.
\ee
By analysing when $\CR$ diverges and when not, we come to the following classification. Abbreviate 
\be{Sdef}
S = \sum_{k\in\N_0} \frac{1}{\bar{c}_{k,N}} \sum_{\ell=0}^k \lambda_\ell.
\ee

\begin{proposition}{\bf [Classification of locally mono-type versus locally multi-type]}
\label{prop.two-four-regimes}
$\mbox{}$
\begin{itemize}
\item[(a)]
$\bar{\sigma}^+_{0,N}<\infty$: $\widehat{X}^\infty$ is locally mono-type if and only if $S=\infty$.
\item[(b)]
$\bar{\sigma}^+_{0,N}=\infty$, $\bar{\sigma}^-_{0,N}<\infty$: $\widehat{X}^\infty$ is locally mono-type if and only if $d^\ast>0$.
\item[(c)]
$\bar{\sigma}^+_{0,N}=\infty$, $\bar{\sigma}^-_{0,N}=\infty$: $\widehat{X}^\infty$ is locally mono-type if and only if $\lambda_k \not\equiv 0$.
\end{itemize}
\end{proposition}

\begin{proof}
(a) By transience, two L\'evy-processes drift apart without block-coalescence, and coalesce if and only if $\PP_0(\CR^+=\infty)=1$, where
\be{e4167}
\CR^+ = \sum_{m\in\N_0} \bar{\CL}_m \sum_{\ell \geq m} \frac{\lambda_\ell}{N^\ell}
\ee
is the total rate of block-coalescence in $\N_0$. Compute 
\be{e4171}
\begin{aligned}
\E_0[\CR^+] &= \sum_{m\in\N_0} \E_0[\bar{\CL}_m] \sum_{\ell \geq m} \frac{\lambda_\ell}{N^\ell}
= \sum_{m\in\N_0} N^m \bar{\sigma}^+_{m,N} \sum_{\ell \geq m} \frac{\lambda_\ell}{N^\ell}\\
&= \sum_{m\in\N_0} N^m \sum_{k \geq m} \frac{1}{\bar{c}_{k,N}} \sum_{\ell \geq m} \frac{\lambda_\ell}{N^\ell}
= \sum_{k\in\N_0} \frac{1}{\bar{c}_{k,N}} \sum_{\ell \in \N_0} \frac{\lambda_\ell}{N^\ell}
\sum_{m = 0}^{k \wedge \ell} N^m 
\end{aligned}
\ee
and note that 
\be{e4181}
N^{-\ell} \sum_{m = 0}^{k \wedge \ell} N^m = N^{(k-\ell) \wedge 0}\,\frac{1-N^{-(k \wedge \ell)-1}}{1-N^{-1}}
\asymp N^{(k-\ell) \wedge 0}.
\ee
Hence $\E_0[\CR^+] \asymp S'$ with 
\be{e4186}
S' = \sum_{k\in\N_0} \frac{1}{\bar{c}_{k,N}} \sum_{\ell \in \N_0} \lambda_\ell
\left[1_{\{0 \leq \ell \leq k\}} + 1_{\{\ell > k\}} N^{-(\ell-k)}\right]. 
\ee
By the second condition in \eqref{cond2}, we have $S' \asymp S$, and so $\E_0[\CR^+] = \infty$ if and only if $S=\infty$. To obtain a \emph{zero-one law} for $\CR^+$ we need the following. 

\begin{lemma}{\bf [Second moment bound]}
\label{lem:zero-one}
$\E_0([\CR^+]^2) \leq 2 [\E_0(\CR^+)]^2$. 
\end{lemma}

\begin{proof}
Compute
\be{e4199}
\begin{aligned}
\E_0[(\bar{\CL}_m)^2] &= 2 \int_0^\infty \dd t \int_t^\infty \dd s\, \PP_0\big(k^\infty(\zeta^1_t,\zeta^2_t) 
= m, k^\infty(\zeta^1_s,\zeta^2_s) = m\big)\\
&= 2 \int_0^\infty \dd t\, \PP_0\big(k^\infty(\zeta^1_t,\zeta^2_t) = m\big) 
\int_0^\infty \dd u\, \PP_m\big(k^\infty(\zeta^1_u,\zeta^2_u) = m\big)
= 2\,\E_0[\bar{\CL}_m]\,\E_m[\bar{\CL}_m]. 
\end{aligned}
\ee
Since $\E_0[\bar{\CL}_m] = \PP(\tau_m < \infty)\,\E_m(\wh\CL_m)$, and $\PP(\tau_m < \infty) = 1$ for all $m\in \N_0$, it follows that $\E_0[(\bar{\CL}_m)^2] = 2\,(\E_0[\bar{\CL}_m])^2$ for all $m \in \N_0$. Now use Cauchy-Schwarz to estimate
\be{e4209}
\begin{aligned}
\E_0[(\CR^+)^2] &= \sum_{m \in \N_0} \sum_{m' \in \N_0} \E_0[\bar{\CL}_m\,\bar{\CL}_{m'})
\sum_{\ell \geq m} \frac{\lambda_\ell}{N^\ell} \sum_{\ell' \geq m'} \frac{\lambda_{\ell'}}{N^{\ell'}}\\ 
&\leq \sum_{m \in \N_0} \sum_{m' \in \N_0} \sqrt{\E_0[(\bar{\CL}_m)^2]\,\E_0[(\bar{\CL}_{m'})^2]}
\sum_{\ell \geq m} \frac{\lambda_\ell}{N^\ell} \sum_{\ell' \geq m'} \frac{\lambda_{\ell'}}{N^{\ell'}}\\ 
&= 2 \sum_{m \in \N_0} \sum_{m' \in \N_0} \E_0[\bar{\CL}_m]\,\E_0[\bar{\CL}_{m'}]
\sum_{\ell \geq m} \frac{\lambda_\ell}{N^\ell} \sum_{\ell' \geq m'} \frac{\lambda_{\ell'}}{N^{\ell'}}
= 2\,(\E_0[\CR^+])^2. 
\end{aligned}
\ee
\end{proof}

\noindent
Lemma \ref{lem:zero-one} and the Paley-Zygmund inequality imply that a zero-one law holds for $\CR^+$, i.e., $\PP_0(\CR^+=\infty)=1$ if and only if $\E_0[\CR^+]=\infty$ (see e.g. \cite[Section 5.2]{GHK18}). Hence the claim follows.

\medskip\noindent
(b) By point-recurrence, two L\'evy-processes meet in finite time without the block-coalescence. Hence, if $d^\ast = d + |\Lambda_{-\infty}|>0$, then they coalesce without block-coalescence and therefore also coalesce when block-coalescence is added. 

\medskip\noindent
(c) By Harris recurrence, $\PP_0(\CR=\infty)=1$ if and only if $\lambda_k \not\equiv 0$.
\end{proof}

The following lemma simplifies the criteria in Proposition \ref{prop.two-four-regimes}. 

\begin{lemma}{\bf [Equivalence of classifications]}
\label{lem:equivalt}
Suppose that 
\be{equivcondalt}
\sup_{k \in \Z} \frac{\lambda^*_{k+1}-\lambda_{k+1}}{c_k} < \infty. 
\ee 
Then
\be{}
\sigma_0^+ = \infty \quad \Longleftrightarrow \quad \bar\sigma_{0,N}^+ = \infty, 
\qquad \qquad
\sigma_0^- = \infty \quad \Longleftrightarrow \quad \bar\sigma_{0,N}^- = \infty. 
\ee
\end{lemma} 

\noindent
Note that, subject to \eqref{regcond}, \eqref{equivcondalt} is equivalent to \eqref{criterion}. Moreover, subject to \eqref{regcond}, if \eqref{equivcondalt} fails, then $\bar c_{k,N} \asymp \lambda_k$ and $S \asymp \sum_{k \in \N_0} \frac{1}{\lambda_k} \sum_{\ell = 0}^k \lambda_\ell = \infty$. 

\medskip
With the help of Lemma \ref{lem:equivalt} we can summarise our findings in Proposition \ref{prop.two-four-regimes} by the following statement (recall \eqref{ag6} and \eqref{ag6alt}):

\begin{corollary}{\bf [Classification of locally mono-type versus locally multi-type]}
\label{cor.two-four-regimes}
Under \eqref{cond1}--\eqref{regcond} the following hold:
\begin{itemize}
\item
Suppose that $\lambda_k \not\equiv 0$. If \eqref{criterion} holds, then
\be{}
\widehat{X}^\infty \text{ is locally mono-type if and only if } \sum_{k\in\N_0} \frac{1}{c_k} \sum_{\ell=0}^k \lambda_\ell = \infty. 
\ee
If \eqref{criterion} fails, then $\widehat{X}^\infty$ is locally mono-type. 
\item
Suppose that $\lambda_k \equiv 0$. If $d^\ast>0$, then
\be{}
\widehat{X}^\infty \text{ is locally mono-type if and only if } 
\sum_{k\in\N_0} \frac{1}{c_k} = \infty, \quad \sum_{k\in -\N} \frac{1}{c_k} < \infty.
\ee
If $d^\ast=0$, then $\widehat{X}^\infty$ is locally multi-type.
\end{itemize}
\end{corollary}

\noindent
The former corresponds to cases (a) and (c), the latter to case (b). Note that $\sum_{k\in\N_0} (1/c_k) = \infty$ implies $S=\infty$ when $\lambda_k \not\equiv 0$, i.e., when there is some $k \in \Z$ for which $\lambda_k>0$. Note that the case $\lambda_k \equiv 0$ and $d^\ast = 0$ is trivial because there is no resampling at all. 

In conclusion, we see that $\widehat{X}^\infty$ is locally mono-type when the migration is recurrent and also when the migration is transient and the block-coalescence is strong enough, while $\widehat{X}^\infty$ is locally multi-type when the migration is transient and the block-coalescence is weak enough.


\subsection{Stable migration}
\label{app:stable}

In this appendix we zoom in on the case of stable migration. In case the migration rates are \emph{regularly varying at} $-\infty$, i.e., satisfy \eqref{ag7a}--\eqref{ag7c}, we can derive the rough asymptotics of $p^*_t(0,0)$ for $t \downarrow 0$. 

\begin{figure}[htbp]
\vspace{3.5cm}
\begin{center}
\setlength{\unitlength}{.6cm}
\begin{picture}(0,-10)
{
\thicklines
\qbezier(-5,0)(0,0)(5,0)
\qbezier(0,0)(0,3)(0,6)
\qbezier(-4,5)(-2,1)(4,0.5)
\qbezier[30](-2.95,0)(-2.95,1)(-2.95,3.4)
\qbezier[30](-2.95,3.4)(-1.5,3.4)(0,3.4)
}
\put(5.5,-.2){$\ell$}
\put(-.35,6.5){$s(\ell)$}
\put(-.1,-.75){$0$}
\put(0.2,3.3){$\tfrac{1}{t}$}
\put(-3.7,-.75){$s^{-1}(\tfrac{1}{t})$}
\end{picture}
\end{center}
\vspace{0.3cm}
\caption{\small Plot of $\ell \mapsto s(\ell)$ extrapolated from $\Z$ to $\R$. Recall that $\lim_{\ell\to\infty} s(\ell) = 0$ by \eqref{cond2} and $\lim_{\ell\to -\infty} s(\ell) = \infty$ by \eqref{condextra}.}
\label{fig-sfunction}
\end{figure}
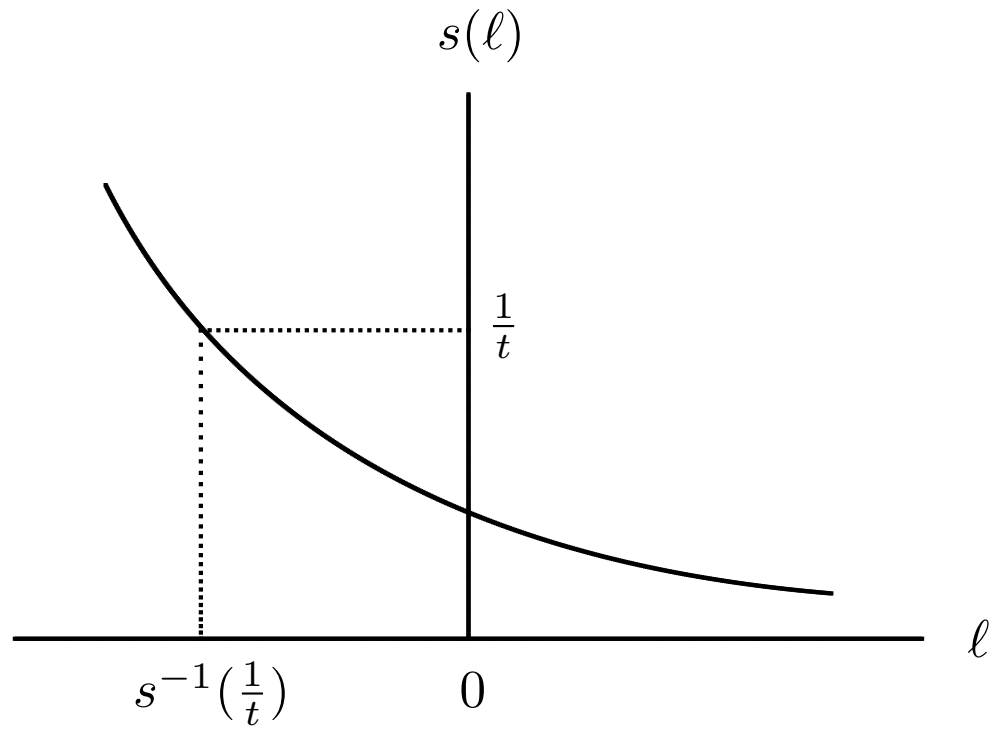

\begin{proposition}{\bf [Short-time asymptotics for the transition kernel density]}
Suppose that $(c_k)_{k\in\Z}$ is regularly varying at $-\infty$. Then
\be{pt*asymp}
p^*_t(0,0) \asymp \frac{1}{t} \left[\sum_{k \asymp s^{-1}(\frac{1}{t})} \frac{1}{c_k} + N^{-o\big(s^{-1}\big(\tfrac{1}{t}\big)\big)}\right], 
\qquad t \downarrow 0,
\ee
where $s^{-1}$ is the inverse of the function $s$ defined in \eqref{e2561} (see Fig.\ \ref{fig-sfunction}). 
\end{proposition}

\begin{proof}
Define $f(x) = x\,\eee^{-x}$, $x \in [0,\infty)$. From \eqref{sapprox} and \eqref{ptasymp} we have
\be{}
p^*_t(0,0) \asymp \sum_{\ell \in -\N} \frac{1}{N^\ell}\, \eee^{-s_\ell t} 
\asymp \sum_{\ell \in -\N} \frac{1}{c_\ell}\, \eee^{-s_\ell t} s_\ell
= \frac{1}{t} \sum_{\ell \in -\N} \frac{1}{c_\ell} f(s_\ell t).
\ee
Since $x \mapsto f(x)$ is bounded from above on $[0,\infty)$ and bounded from below on compact subsets on $[0,\infty)$, the sum splits into two parts: one over $\ell$ such that $s_\ell t \asymp 1$, i.e., $\ell \asymp s^{-1}(\tfrac{1}{t})$, and one over $\ell$ such that $s^{-1}(\tfrac{1}{t}) \ll \ell < 0$, which is to be viewed as an error term.
\end{proof}

\vspace{0.5cm}

For the case of exponentially varying migration coefficients
\be{expovar}
c_k = c^k, \quad \lambda_k = \lambda^k, \quad \lambda^*_k = C \lambda^k, \quad k \in \Z,
\qquad c \in (0,N), \lambda \in [0,N), C \in (1,\infty) 
\ee
the following dichotomy holds.  

\begin{lemma}{\bf [Dichotomy for exponentially varying migration coefficients]}
Subject to \eqref{expovar}, $S = \infty$ when $c \vee \lambda \leq 1$ and $c \leq \lambda < 1$, while $S<\infty$ when $c > 1 \geq \lambda$ and $c > \lambda > 1$.
\end{lemma}

\begin{proof}
The claim is immediate from \eqref{1901pre}--\eqref{1901} and \eqref{Sdef}. 
\end{proof}

Potential theory stipulates that $p^*_t(0,0) \asymp t^{-1-\gamma + o(1)}$, with $\gamma \in (-1,\infty)$ the \emph{degree} of the migration and $o(1)$ representing a sub-polynomial error term. Point-recurrent corresponds to $\gamma \in (-1,0)$, Harris-recurrent but not point-recurrent to $\gamma= 0$, transient to $\gamma \in (0,\infty)$. 

\begin{lemma}{\bf [Degree of the migration for exponentially varying migration coefficients]}
Subject to \eqref{expovar},
\be{}
p^*_t(0,0) \asymp t^{-(1+\gamma)+o(1)}, \quad t \downarrow 0, \qquad \gamma = \frac{\log c}{\log(N/c)}. 
\ee
\end{lemma}

\begin{proof}
Since $s_\ell \asymp (c/N)^\ell$, $\ell \to -\infty$, we have $s^{-1}(1/t) \sim -\log(1/t)/\log (N/c)$, and hence
\be{}
p^*_t(0,0) \asymp \frac{1}{t} \left[\sum_{\ell \sim -\log(1/t)/\log(N/c)} \frac{1}{c^\ell} + t^{o(1)}\right]
\asymp \frac{1}{t} \left[c^{\log(1/t)/\log(N/c)} + t^{o(1)}\right] = t^{-(1+\gamma)+o(1)}.
\ee 
\end{proof}

\noindent
Note that $\gamma$ changes sign at $c=1$, which is the critical threshold between point-recurrence and transience.  


\subsection{A brief digression}
\label{app:digression}

Since the formulas in Appendices \ref{appA:oneLP}--\ref{appA:twoLP} are technical, it is useful to provide some further insight into where they come from, which we do here. We again consider a single L\'evy-process, respectively, two L\'evy-processes. 


\subsubsection{Single individual}

Let $k^\infty(\zeta^\infty_t,0)$ denote the hierarchical distance between $\zeta^\infty_t$ and the origin. Then 
\begin{equation*}
K = \big(k^\infty(\zeta^\infty_t,0)\big)_{t \geq 0}
\end{equation*}
is the Markov chain on $\Z$ with transition rates $r(m,n)$ given by 
\be{}
\label{ratesmig}
r_N(m,n) = \frac{N-1}{N} \sum_{\ell \geq m \vee n} \frac{c_{\ell-1}}{N^{\ell-1}} \frac{1}{N^{\ell-n}},
\qquad m,n \in \Z,\,m \neq n. 
\ee
Indeed, $\frac{c_{\ell-1}}{N^{\ell-1}}$ is the rate at which the migration chooses the $\ell$-block, the value of $\ell$ must be at least $m$ otherwise the distance to the origin cannot change, it must also be at least $n$ otherwise the distance cannot reach $n$, while there are $N^{\ell-n}$ $n$-blocks inside each $\ell$-block, of which only one contains the origin. The factor $(N-1)/N$ comes from the fact that the distance must end up being exactly $n$ after the random relocation in the $n$-block around the origin is carried out. 

\begin{remark}{\bf [Harmonic function]}
{\rm Note that 
\be{}
N^m r_N(m,n) = N^n r_N(n,m), \qquad m,n \in \Z,\,m \neq n,
\ee
so that, after a time change where at site $m$ the rate is speeded up by $N^m$, we get a reversible Markov process with a constant invariant measure. The harmonic function is not affected by the time change.}\hfill$\spadesuit$
\end{remark}  

For large $N$, in the Markov process $K$ the transitions between \emph{nearest-neighbour distances} are dominant, and we have
\be{limrates}
N^m r_N(m,m-1) \to c_{m-1}, \qquad N^m r_N(m,m+1) \to c_m, \qquad N\to\infty. 
\ee 

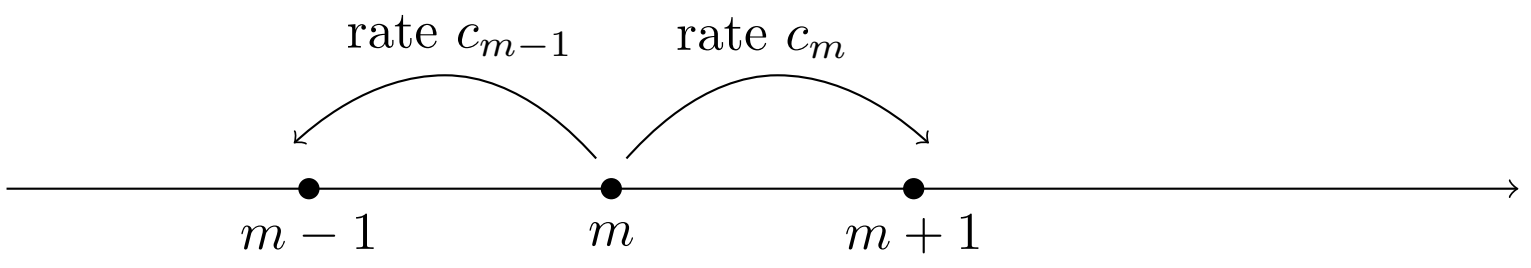
\begin{figure}[htbp]
\begin{center}
\begin{tikzpicture}[->]
\draw (0,0) -- (xyz cs:x=10);
\fill (2,0) circle (2pt);
\fill (4,0) circle (2pt);
\fill (6,0) circle (2pt);
\node (a) at (6,-0.3) {$m+1$};
\node (a) at (4,-0.3) {$m$};
\node (a) at (2,-0.3) {$m-1$};
\node (a) at (5,1) {rate $c_m$};
\node (a) at (3,1) {rate $c_{m-1}$};
\draw (4.1,0.2) parabola[parabola height=0.5cm] +(2,0.1);
\draw (3.9,0.2) parabola[parabola height=0.5cm] +(-2,0.1);
\end{tikzpicture}
\end{center}
\caption{Transition rates for the continuous-time random walk on $\Z$ obtained through \eqref{limrates}.}
\label{fig-apprtr}
\end{figure}

\noindent
The approximation made in \eqref{sapprox} amounts to \emph{pretending that these are the only transitions in} $K$ (see Fig.\ \ref{fig-apprtr}). As we will see below, this simplification does not qualitatively affect the recurrence versus transience properties of $K$. 

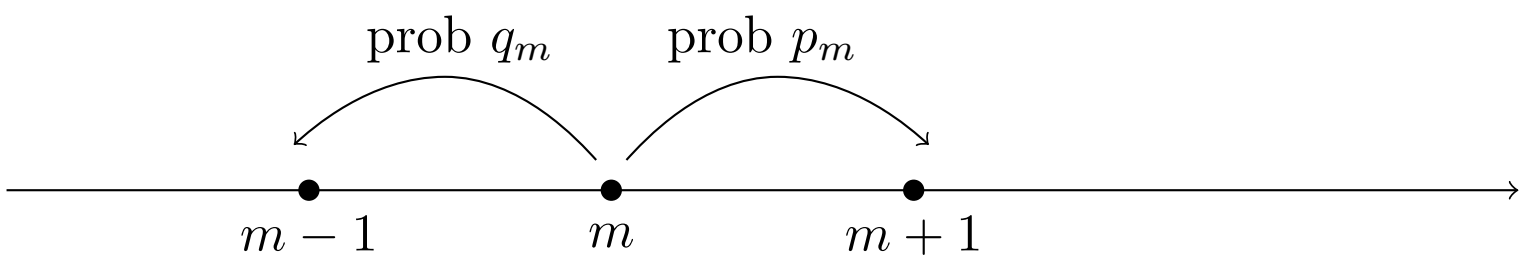
\begin{figure}[htbp]
\begin{center}
\begin{tikzpicture}[->]
\draw (0,0) -- (xyz cs:x=10);
\fill (2,0) circle (2pt);
\fill (4,0) circle (2pt);
\fill (6,0) circle (2pt);
\node (a) at (6,-0.3) {$m+1$};
\node (a) at (4,-0.3) {$m$};
\node (a) at (2,-0.3) {$m-1$};
\node (a) at (5,1) {prob $p_m$};
\node (a) at (3,1) {prob $q_m$};
\draw (4.1,0.2) parabola[parabola height=0.5cm] +(2,0.1);
\draw (3.9,0.2) parabola[parabola height=0.5cm] +(-2,0.1);
\end{tikzpicture}
\end{center}
\caption{Transition probabilities for the discrete-time embedded random walk on $\Z$ obtained from Fig.\ \ref{fig-apprtr}.}
\label{fig-trprob}
\end{figure}

Consider the discrete-time embedded random walk with transition probabilities (see Fig.~\ref{fig-trprob})
\be{e2486}
p_m = \frac{c_{m+1}}{c_m+c_{m+1}}, \qquad 
q_m = \frac{c_m}{c_m+c_{m+1}}
\ee

\noindent
Clearly,
\be{exploc}
\E_0[\CL_m]= \frac{N^m}{c_m+c_{m-1}}\,\, 
\frac{\PP_0(\tau_m<\infty)}{p_m\PP_{m+1}(\tau_m=\infty)+q_m\PP_{m-1}(\tau_m=\infty)}
\ee
where $\tau_m$ is the \emph{first hitting time} of $m$ for the discrete-time random walk. The three probabilities in \eqref{exploc} are easy to compute with the help of standard renewal theory (see e.g.\ \cite{KM57}). Namely,
\be{hits1}
\begin{aligned}
\PP_{m+1}(\tau_m=\infty) &= \frac{\sigma^+_m - \sigma^+_{m+1}}{\sigma^+_m}
= \frac{1}{c_m} \frac{1}{\sigma^+_m},\\
\PP_{m-1}(\tau_m=\infty) &= \frac{\sigma^-_m - \sigma^-_{m-1}}{\sigma^-_m}
= \frac{1}{c_{m-1}} \frac{1}{\sigma^-_m},
\end{aligned}
\ee
and
\be{hits2}
\PP_0(\tau_m<\infty)
= \left\{\begin{array}{ll}
\frac{\sigma^+_0}{\sigma^+_m}, &m<0,\\[0.3cm]
1, &m=0,\\[0.3cm]
\frac{\sigma^-_0}{\sigma^-_m}, &m>0.
\end{array}
\right.
\ee
Inserting \eqref{hits1}--\eqref{hits2} into \eqref{exploc}, we get 
\be{e2447}
N^{-m} \E_0[\CL_m] = \left\{\begin{array}{ll}
\frac{\sigma^+_0\sigma^-_m}{\sigma^+_m+\sigma^-_m}, &m<0,\\[0.3cm]
\frac{\sigma^+_0\sigma^-_0}{\sigma^+_0+\sigma^-_0}, &m=0,\\[0.3cm]
\frac{\sigma^+_m\sigma^-_0}{\sigma^+_m+\sigma^-_m}, &m>0.
\end{array}
\right.
\ee
In case $\sigma^+_m$ or $\sigma^-_m$ are infinite, the expressions in \eqref{hits2} and \eqref{e2447} simplify. Namely, if $\sigma^+_m=\infty$ and $\sigma^-_m<\infty$, then 
\be{int1}
\begin{aligned}
&\PP_0(\tau_m<\infty) \left\{\begin{array}{ll}
= 1, &m \leq 0,\\
< 1, &m>0,
\end{array} 
\right.\\
&N^{-m} \E_0[\CL_m] \left\{\begin{array}{ll}
= \sigma_m^-, &m \leq 0,\\
= \sigma_0^-, &m>0,
\end{array} 
\right.
\end{aligned}
\ee
while if $\sigma^+_m<\infty$ and $\sigma^-_m=\infty$, then
\be{int2}
\begin{aligned}
&\PP_0(\tau_m<\infty) \left\{\begin{array}{ll}
< 1, &m < 0,\\
= 1, &m \geq 0,
\end{array} 
\right.\\
&N^{-m} \E_0[\CL_m] \left\{\begin{array}{ll}
= \sigma_0^+, &m < 0,\\
= \sigma_m^+, &m \geq 0.
\end{array} 
\right.
\end{aligned}
\ee

\begin{proposition}{\bf [Classification for the distance to the origin of a single L\'evy-process]}
\label{prop.one-four-regimes-K}
$\mbox{}$
\begin{itemize}
\item[(a)]
$\sigma^+_0=\infty$, $\sigma^-_0<\infty$: $K$ is transient to the left and hits $-\infty$ in finite time with probability $1$.
\item[(b)]
$\sigma^+_0<\infty$, $\sigma^-_0=\infty$: $K$ is transient to the right and does not hit $\infty$ in finite mean time.
\item[(c)]
$\sigma^+_0=\infty$, $\sigma^-_0=\infty$: $K$ is recurrent.
\item[(d)]
$\sigma^+_0<\infty$, $\sigma^-_0<\infty$: $K$ is transient both to the left and to the right with positive probability. Conditional on being transient to the left, $K$ hits $-\infty$ in finite time with probability $1$. Conditional on being transient to the right, $K$ does not hit $\infty$ in finite mean time. 
\end{itemize}
\end{proposition}

\begin{proof}
(a) Use \eqref{int1}, and note that $\sum_{m \leq 0} \E_0[\CL_m] = \sum_{m \leq 0} \sigma_m^- N^{m} \leq \sigma_0^- \sum_{m \leq 0} N^{m} < \infty$.\\
(b) Use \eqref{int2}, and note that $\sum_{m \geq 0} \E_0[\CL_m] = \sum_{m \geq 0} \sigma_m^+ N^{m} \geq \sum_{m \geq 0} \frac{1}{c_m} N^{m} = \infty$ (by \eqref{ag35}).\\
(c) Use \eqref{hits1}.\\
(d) Use \eqref{hits1}--\eqref{hits2} and the computations in (a) and (b).
\end{proof}

\begin{remark}{\bf [Entrance law]}
{\rm After $K$ hits $-\infty$ (when $\sigma^+_0=\infty$ and $\sigma^-_0<\infty$), it returns to $\Z$ according to an entrance law that is easily deduced from \eqref{e2553infty}.}\hfill$\spadesuit$
\end{remark}


\subsubsection{Two individuals}

Let $k^\infty(\zeta^1_t,\zeta^2_t)$ denote the hierarchical distance between $\zeta^1_t$ and $\zeta^2_t$. Then 
\be{}
\bar{K} = \big(k^\infty(\zeta^1_t,\zeta^2_t)\big)_{t \geq 0}
\ee
is the Markov chain on $\Z$ with transition rates $\bar{r}_N(m,n)$ given by (compare with \eqref{ratesmig})
\be{fdh12}
\bar{r}_N(m,n) = 2 r_N(m,n) + \frac{N-1}{N} \sum_{\ell \geq m \vee n} 
\frac{\lambda^*_\ell-\lambda_\ell}{N^\ell}\,\frac{1}{N^{\ell-n}}, \qquad m,n \in \Z, m \neq n.
\ee
Indeed, $(\lambda^*_\ell-\lambda_\ell)/N^\ell$ is the rate at which two lineages in an $\ell$-block are \emph{not both marked} and make a random jump through resampling and reshuffling. For the distance process this corresponds to migration rates
\be{}
\bar{c}_{k,N} = 2 c_k + \frac{\lambda^*_{k+1}-\lambda_{k+1}}{N}, \qquad k \in \Z.
\ee

We again distinguish between four regimes.

\begin{proposition}{\bf [Classification for the distance between two L\'evy-processes]}
\label{prop.two-four-regimes-K}
$\mbox{}$
\begin{itemize}
\item[(a)]
$\bar{\sigma}^+_{0,N}=\infty$, $\bar{\sigma}^-_{0,N}<\infty$: $\bar{K}$ is transient to the left and hits $-\infty$ in finite time with probability $1$.
\item[(b)]
$\bar{\sigma}^+_{0,N}<\infty$, $\bar{\sigma}^-_{0,N}=\infty$: $\bar{K}$ is transient to the right and does not hit $\infty$ in finite mean time.
\item[(c)]
$\bar{\sigma}^+_{0,N}=\infty$, $\bar{\sigma}^-_{0,N}=\infty$: $\bar{K}$ is recurrent.
\item[(d)]
$\bar{\sigma}^+_{0,N}<\infty$, $\bar{\sigma}^-_{0,N}<\infty$: $\bar{K}$ is transient both to the left and to the right with positive probability. Conditional on being transient to the left, $\bar{K}$ hits $-\infty$ in finite time with probability $1$. Conditional on being transient to the right, $\bar{K}$ does not hit $\infty$ in finite mean time.
\end{itemize}
\end{proposition}

\begin{proof}
The claims follow from the same type of argument as in the proof of Proposition \ref{prop.one-four-regimes-K}. 
\end{proof}


\section{Alternative routes to the construction of the continuum-space process}
\label{appB}


\paragraph*{Formulation with graphical representation.}

An alternative way to obtain the continuum process is via a \emph{graphical construction}, representing the genealogy of the individuals alive at a given time. In \cite{GSW16} the genealogy of the continuum spatial Fleming-Viot process on $\R$ was described as a system of coalescing Brownian motions starting with an individual in every space-time point, called the \emph{Brownian web} (see \cite{A79}, \cite{FINR04}, \cite{MT95}). In \cite{gh2} we will show that the genealogy of the continuum spatial Cannings process can be described as a system of coalescing L\'evy-processes with block-coalescence, which we call the \emph{L\'evy-web}. We will see that the dichotomy of locally mono-type versus locally multi-type plays out in the structure of the L\'evy-web, in the sense that the latter either locally is a \emph{tree} (corresponding to the stochastic mode of evolution) or degenerates to \emph{dust} (corresponding to the deterministic mode of evolution).


\paragraph*{Formulation with atomic measure states.}

We may ask whether the martingale problem can be rephrased in a more classical form by coding the state as a certain \emph{atomic measure} on $\Omega_N^\infty \times \K$. 

\begin{itemize}
\item
For the case where the geographic space of the Cannings process is $\R$, this is indeed the case, as was first investigated by \cite{A79}. Such a dynamics can be specified by using that $\R$ is totally ordered. The dynamics is driven by what are called \emph{annihilating Brownian motions}. In the case where all the individuals initially have a different type, the mono-type regions are given by intervals whose boundaries evolve as annihilating Brownian motions. The intervals have a total length equal to the full Lebesgue measure on $\R$. This setting can also be described as an \emph{evolving atomic random measure} on $\R$, the atom being located at the middle of the interval and carrying a weight that is equal to the length of the interval. The two evolve as a pair of independent Brownian motions coupled by the fact that the left boundary of one interval is the right boundary of the other interval. However, a proper martingale problem description still must work with the \emph{space-time view}, since typically there are countably many intervals in any macroscopic interval close to potentially countably many points on the left and the right of an interval. This has not been carried out yet.
\item
For the Cannings model on $\Omega_N^\infty$ we have a different structure, yet the above description looks promising, the core being the collection of $\Z$-valued processes given by the objects from Lemma \ref{l.1562}(b), which provides the analogue of the annihilating Brownian motions. For our dynamics $\wh X^\infty$ on $\wt\CM(\Omega_N^\infty \times \K)$, we consider the state space
\begin{equation*}
\CM^{\ast \ast}_{\rm loc}(\Omega_N^\infty \times \K)
\end{equation*}
consisting of all measures $x \in \wt\CM(\Omega_N^\infty \times \K)$ that are locally mono-type, are regular with respect to $\mathfrak{h}$, and are locally uniform (compare with \eqref{e1664}). We can decompose $\Omega_N^\infty$ into maximal mono-type balls, which are locally uniform, and obtain an \emph{atomic measure} on $\Omega_N^\infty \times \K$, provided we identify all the space points in $L$-balls. i.e., take for an $L$-ball the element in $\Omega_N^L$ representing the center of that ball after the embedding in $\Omega_N^\infty$. This way the decomposition in maximal $L$-balls is characterised by an element in $(\Omega_N^\infty \times \Z)^\N$ specifying the pairs (center and radius) of the decomposition, which is in one-to-one correspondence with an atomic measure on $\Omega_N^\infty$ and on $\Omega_N^\infty \times \K$, provided we code the type as well. Considering
\begin{equation}
\label{e2126}
\text{the \emph{weak-atomic topology} on $\CM^{\ast \ast}(\Omega_N^\infty \times \K)$},
\end{equation}
we obtain a \emph{Polish space} as state space for atomic measure-valued functionals on this subset of states, which are locally mono-type and uniform and are given by $\wh X_t^\infty$, $t >0$. This (atomic measure)-valued process on $\Omega_N^\infty \times \K$ is the analogue of the annihilating Brownian motion picture for Fleming-Viot on $\R$, except that now a type can be assigned to a \emph{multitude} of balls. 
\item
We can ask whether the corresponding version of $(\wh X_t^{\ast,\infty})_{t >0}$ satisfies, in this \emph{new} coding with purely \emph{atomic} measures as states, a martingale problem of a more classical form, where \emph{atom locations} and \emph{atom sizes} fluctuate. For this dynamics we can ask whether we have one or several spatial balls for a type. Does the answer depend on the parameters, in that we expect the degree of the effective motion to be the relevant quantity? In fact, we will show that with positive probability we find more than one such clusters. This comes from the fact that we have jump processes as migration. Therefore we need to study the dynamics for the $\Z$-valued components and their locations, i.e., a set of pairs $(\eta,k) \in \Omega_N^\infty \times \K$, with elements in $(\Omega_N^\infty \times \Z)^\N$. In this process new balls can appear and they can become extinct, meaning that new atom locations are created or atom locations disappear. These are birth and death processes with values in $\Omega_N^\infty \times \K$, with a jump canceling a pair or a new pair appearing at a random time, where for a given time $t$ these births or deaths accumulate. Also the process of the ball radii can jump upwards or downwards in the set $\Z$ by a jump process on $\Z$ paired with a fixed location. This process moves eventually with the macroscopic evolution for $t\to\infty$ and for $t \downarrow 0$, forming an entrance law from $-\infty$.
\end{itemize}

Note that the above description will also have to work with the \emph{space-time view} to describe the fluctuations, possibly in a more intuitive way. This would give a richer picture of spatial patterns of the types for jump processes, similarly to what was done in \cite{GSW16} for the Fleming-Viot model on $\R$. Note that this \emph{cannot} replace the martingale problem on the full state space $\CM(\Omega_N^\infty \times \K)$, even though it is informative.


\end{appendix}


\bibliography{gh1}

\bibliographystyle{alpha}


\end{document}